\documentclass[10pt]{article}
\usepackage[utf8]{inputenc}

\usepackage[T1]{fontenc}
\usepackage[english]{babel}
\usepackage{lmodern}
\usepackage{amsmath}
\usepackage{amsthm}
\usepackage{amssymb}
\usepackage{yhmath}
\usepackage{stmaryrd}
\usepackage{mathrsfs}
\usepackage{fullpage}
\usepackage[all]{xy}
\usepackage{float}
\usepackage{enumitem}
\usepackage{comment}
\usepackage{graphicx}
\usepackage{verbatim}
\usepackage{latexsym}
\usepackage{hyperref}
\usepackage{overpic}
\usepackage{mathtools}
\usepackage{caption}
\usepackage{xcolor}
\usepackage{tikz-cd}
\usepackage[bottom]{footmisc}
\usepackage{accents}

\usetikzlibrary{calc}

\definecolor{green}{RGB}{34,139,34}   
\definecolor{red}{RGB}{220,20,60}     
\definecolor{blue}{RGB}{28,123,220} 
\definecolor{sectioncolor}{RGB}{155,135,185} 

\newcommand{\fonction}[5]{\begin{array}{lccl}
#1: & #2 & \longrightarrow & #3 \\
    & #4 & \longmapsto & #5 \end{array}}

\newcommand{\fonctionsansnom}[4]{\begin{array}{ccc}
#1 & \longrightarrow & #2 \\
#3 & \longmapsto & #4 \end{array}}

\newcommand{\limit}[4]{#1 \underset{#2 \rightarrow #3}{\longrightarrow} #4}

\newcommand{\liftingproblem}[9]{
\xymatrix{
#1 \ar[d]^-{#6} \ar[r]^-{#5} & #2 \ar[d]^-{#7} \\
#3 \ar[r]^-{#8} \ar@{-->}[ru]^-{#9} & #4
}
}

\newcommand{\squarediagram}[8]{
\xymatrix{
#1 \ar[d]_-{#6} \ar[r]^-{#5} & #2 \ar[d]^-{#7} \\
#3 \ar[r]^-{#8} & #4
}
}

\newcommand{\Z}{\mathbb{Z}}

\newcommand{\R}{\mathbb{R}}
\newcommand{\C}{\mathbb{C}}

\newcommand{\K}{\mathbb{K}}

\newcommand{\SM}{\mathbb{S}}

\newcommand{\Nl}{\mathcal{N}}
\newcommand{\Cl}{\mathcal{C}}
\newcommand{\Ml}{\mathcal{M}}

\newcommand{\Ql}{\mathcal{Q}}

\newcommand{\Pl}{\mathcal{P}}

\newcommand{\Sl}{\mathcal{S}}

\newcommand{\El}{\mathcal{E}}

\newcommand{\Fl}{\mathcal{F}}

\newcommand{\Hl}{\mathcal{H}}
\newcommand{\Dl}{\mathcal{D}}
\newcommand{\Ul}{\mathcal{U}}
\newcommand{\Vl}{\mathcal{V}}
\newcommand{\Wl}{\mathcal{W}}
\newcommand{\Kl}{\mathcal{K}}
\newcommand{\Cj}{\mathscr{C}}
\newcommand{\Dj}{\mathscr{D}}

\newcommand{\Lj}{\mathscr{L}}
\newcommand{\Gj}{\mathscr{G}}

\newcommand{\PR}{\mathfrak{P}}
\newcommand{\Cr}{\mathfrak{C}}

\newcommand{\Gd}{\mathsf{G}}
\newcommand{\Ld}{\mathsf{L}}

\newcommand{\sd}{\mathrm{sd}}
\newcommand{\Pre}{\mathrm{Pre}}
\newcommand{\id}{\mathrm{id}}
\newcommand{\ev}{\mathrm{ev}}

\newcommand{\Cat}{\mathrm{Cat}}
\newcommand{\Set}{\mathrm{Set}}
\newcommand{\sSet}{\mathrm{sSet}}

\newcommand{\Hom}{\mathrm{Hom}}

\newcommand{\Ch}{\mathrm{Ch}}
\newcommand{\op}{\mathrm{op}}

\newcommand{\Mod}{\mathrm{Mod}}

\newcommand{\Cyl}{\mathrm{Cyl}}

\newcommand{\Obj}{\mathrm{Obj}}
\newcommand{\Spine}{\mathrm{Spine}}
\newcommand{\Pos}{\mathrm{Pos}}

\newcommand{\Fun}{\mathrm{Fun}}

\newcommand{\colim}{\mathrm{colim}}

\newcommand{\Exit}{\mathrm{Exit}}
\newcommand{\Top}{\mathrm{Top}}

\newcommand{\Sing}{\mathrm{Sing}}

\newcommand{\Det}{\mathrm{Det}}

\newcommand{\Loc}{\mathrm{Loc}}

\newcommand{\ag}{\mathbf{a}}
\newcommand{\bg}{\mathbf{b}}

\DeclareRobustCommand{\Mi}{\accentset{\circ}{\mathcal{M}}}
\newcommand{\tu}{\underline{t}}
\newcommand{\qo}{\overline{q}}
\newcommand{\consel}{_{f(a_i),f(a_{i+1})}^{\ag}}
\newcommand{\Int}{\mathrm{Int}}
\newcommand{\Strat}{\mathrm{Strat}}
\newcommand{\Spaces}{\mathrm{Spaces}}
\newcommand{\CW}{\mathrm{CW}}
\newcommand{\hKan}{\mathrm{hKan}}
\newcommand{\hCW}{\mathrm{hCW}}
\newcommand{\inftyCat}{\infty \mathrm{Cat}}

\newtheoremstyle{normalstyle}{}{}{}{}{\bfseries}{.}{ }{}

\newtheorem{Thms}{Theorem}[subsection]
\newtheorem{Lemmas}[Thms]{Lemma}
\newtheorem{Proposition-Definition}[Thms]{Proposition-Definition}
\newtheorem{Proposition}[Thms]{Proposition}
\newtheorem{Thm}{Theorem}
\newtheorem*{Thm*}{Theorem}

\newtheorem{Corollary}[Thms]{Corollary}

\newtheorem*{Th}{Theorem}

\theoremstyle{normalstyle}
\newtheorem{Definition}[Thms]{Definition}
\newtheorem{Notation and Terminology}[Thms]{Notation and Terminology}
\newtheorem{Example}[Thms]{Example}
\newtheorem{Remark}[Thms]{Remark}
\newtheorem{Notation}[Thms]{Notation}

\newtheorem{Construction}[Thms]{Construction}
\newtheorem*{examplecont}{Example \ref{example: flow category of the other sphere} (continued)}

\usepackage{pgfplots}
\pgfplotsset{compat=1.18}
\usetikzlibrary{arrows.meta}
\usetikzlibrary{patterns, patterns.meta}
\usetikzlibrary{svg.path}
\usetikzlibrary{decorations.markings}

\tikzset{
  dotted vertical lines/.style={
    pattern={
      Lines[angle=90, distance=8pt, line width=0.4pt,
            dash pattern=on 1pt off 2pt]
    }
  }
}

\begin{document}

\begin{center}

\Large{Morse flow categories as exit path categories}

\normalsize

Colin Fourel \footnote{ colin.fourel@math.unistra.fr \\ Institut de Recherche Mathématique Avancée, \\
UMR 7501 Université de Strasbourg et CNRS \\ 7 rue René-Descartes, \\
67000 Strasbourg, France}

\end{center}

\vspace{1cm}

\textbf{Abstract.} We prove that the topological flow category $\Ml$ arising from a Morse-Smale pair $(f,\xi)$ on a smooth closed manifold $X$ is equivalent, as an $\infty$-category, to Lurie's $\infty$-category $\Sing_A(X)$ of exit paths in $X$ with respect to the stratification by the stable manifolds of $\xi$.

The objects of $\Ml$ are the critical points of $f$, and for every pair of critical points, the space of morphisms of $\Ml$ between these is the space of possibly broken trajectories of $\xi$ connecting them; it can be identified up to homotopy with the space of unbroken ones. The latter maps naturally to the space of exit paths connecting these critical points; we prove this map to be a weak homotopy equivalence. Then, we combine these ingredients with several others to construct a zigzag of equivalences between the homotopy coherent nerve of $\Ml$, denoted $\Nl(\Ml)$, and $\Sing_A(X)$. The $n$-simplices of $\Nl(\Ml)$ are homotopy coherent diagrams of $n$ composable morphisms of $\Ml$; we introduce the notion of unbroken diagram, yielding an $\infty$-subcategory of $\Nl(\Ml)$, which we refer to as the flow coherent nerve of $\Ml$. The simplices of the latter give rise to stratified maps out of a family of stratified cubes, into $X$. We organize this family into a functor from the category of finite ordered sequences of critical points, to the category of $A$-stratified topological spaces, and we prove a comparison result with the usual stratified geometric realization functor. We finally use a theorem of Tanaka that associates a functor to a semi-simplicial map between $\infty$-categories that satisfies certain conditions.

Our theorem has implications regarding constructible sheaves and the description of homotopy types in terms of flow categories.

\vspace{2cm}

\textbf{Reader's guide.} After recalling the definitions of the main objects and stating the main result in section \ref{section: main objects and main resut}, we present some motivations and consequences in section \ref{section: motivations and consequences}, we briefly present the idea of the proof in section \ref{section: idea of the proof and organization} and we give a detailed summary of the proof in section \ref{section: detailed summary of the proof}. This paper offers two perspectives on the same proof. The first perspective, which is the proof in itself, consists of sections \ref{section: infinity categories, stratified spaces and Morse-Smale pairs} through \ref{section: turning a semi simplicial map into functor}. On the other hand, our theorem is a consequence of a more general statement which does not refer to Morse flow categories, but instead to what we call \emph{stratified categories}. This second perspective is presented in section \ref{section: stratified categories as exit path categories}; this section is designed to be directly accessible to the reader with a good familiarity with the basics of $(\infty,1)$-category theory (as developed in \cite{HigherTopos}), stratified homotopy theory and Morse flow categories. The necessary background on these topics is presented in section \ref{section: infinity categories, stratified spaces and Morse-Smale pairs}.  

\newpage

\tableofcontents

\newpage

\section{Introduction}\label{section: Introduction}

\subsection{Main objects and main result}\label{section: main objects and main resut}

Consider a closed smooth manifold $X$ together with a Morse function $f : X \rightarrow \R$, and denote by $A$ the set of critical points of $f$. Let $\xi$ be a negative Morse-Smale pseudo-gradient vector field on $X$ adapted to $f$ \footnote{Although we will prove the main result of this paper under this assumption on $\xi$, one can think of $\xi$, as a first approximation in this introduction, as the negative gradient of $f$ for some Riemannian metric on $X$.}. Such a datum $(f,\xi)$ will be called a \textit{Morse-Smale pair} on $X$.

The central object of study of this paper is the \textit{flow category} of $(X,f,\xi)$, denoted $\Ml$. It is defined as follows:

\begin{itemize}
    \item[$\bullet$] The set of objects of $\Ml$ is $A$, the set of critical points of $f$.

    \item[$\bullet$] Given two critical points $a$ and $b$, $\Ml(a,b)$ is the set of possibly broken, and unparameterized, trajectories of $\xi$ linking $a$ to $b$.

    \item[$\bullet$] The composition law is defined as the concatenation of trajectories of $\xi$.
\end{itemize}

\begin{Example}\label{example: flow category of the other sphere}
    Let us illustrate this construction with figures \ref{figure: the other sphere} and \ref{figure: some trajectories on other sphere}.

\begin{figure}[H]
    \centering
    \includegraphics[width=0.5\linewidth]{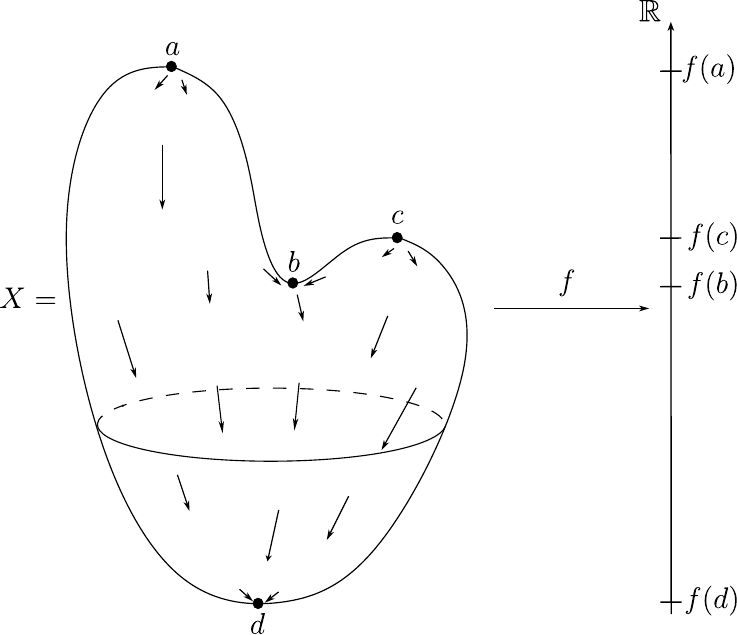}
    \caption{The manifold $X$ is the image of an embedding of $\SM^2$ in $\R^3$ and $f$ is the restriction to $X$ of the projection to the vertical axis. It has four critical points; one maximum $a$, one local maximum $c$, one saddle point $b$ and one minimum $d$. The vector field $\xi$ on $X$ is defined to be the negative gradient of $f$ with respect to the restriction to $X$ of the Euclidean metric on $\R^3$.}
    \label{figure: the other sphere}
\end{figure}

\begin{figure}[H]
    \centering
    \includegraphics[width=0.23\textwidth]{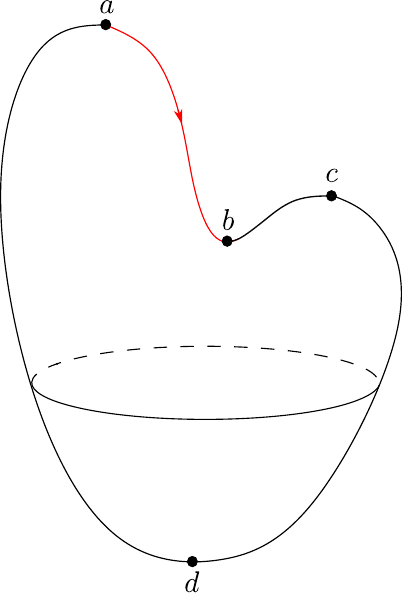}
    \, \, \,
    \includegraphics[width=0.23\textwidth]{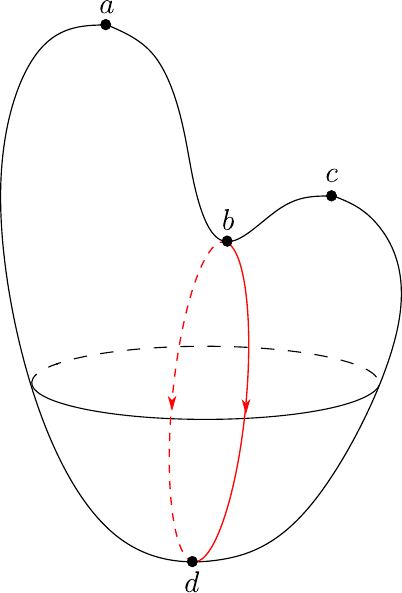}
    \, \, \,
    \includegraphics[width=0.23\textwidth]{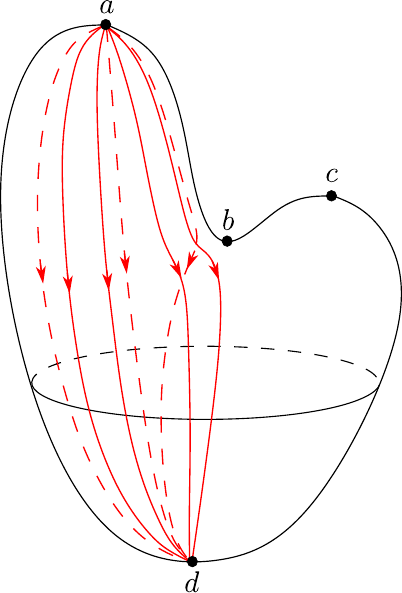}
    \caption{On the left: the only trajectory of $\xi$ from $a$ to $b$. In the middle: the two trajectories from $b$ to $d$. On the right: some unbroken trajectories from $a$ to $d$. In addition to these, there are two broken trajectories from $a$ to $d$, obtained by composing the one from $a$ to $b$ with the two from $b$ to $d$. We say that these two trajectories are broken at $b$. This allows to identify $\Ml(a,d)$ with the closed interval $[0,1]$, by identifying the two boundary components of $[0,1]$ with the two broken trajectories from $a$ to $d$, and by identifying the interior of $[0,1]$ with the one-parameter family of unbroken trajectories connecting the two broken ones. Similarly, $\Ml(c,b)$ consists of one morphism only, and $\Ml(c,d)$ can be identified with $[0,1]$. Finally, all the endomorphism sets of $\Ml$ are reduced to the constant trajectory at the corresponding critical points, and all the other morphism sets of $\Ml$ are empty. This gives a complete description of the category $\Ml$ in that case.}
    \label{figure: some trajectories on other sphere}
\end{figure}

\end{Example}

One can equip the morphism sets of $\Ml$ with several structures compatible with the composition operation in $\Ml$. In this paper, we will mainly view $\Ml$ as a \textit{topological category}. We will adopt the following definition of topological category.

\begin{Definition}\label{definition: intro topological category}

    A topological category is a category enriched over the category of (compactly generated and weakly Hausdorff) topological spaces.

\end{Definition}

We defer to later the precise description of the topology on the morphism sets of $\Ml$ (see section \ref{section: the flow category of a Morse-Smale pair}). For the moment we just describe it in the case of example \ref{example: flow category of the other sphere}.

\begin{examplecont}\label{example continued: topology on morphism sets of flow category}
    The topology on $\Ml(b,d)$ is the discrete topology, and the topology on $\Ml(a,d)$ is identified with the usual topology on the closed interval $[0,1]$. Similarly, the topology on $\Ml(c,d)$ is identified with the usual topology on $[0,1]$.
\end{examplecont}

Topological categories form one model for the theory of $(\infty,1)$-categories, and the goal of this paper is to understand $\Ml$ as an $\infty$-category \footnote{All the $\infty$-categories considered in this paper are $(\infty,1)$-categories. We will usually omit to write the $1$.}.

To this end we will consider a \textit{stratification} of $X$ determined by $(f,\xi)$. We will define later what we mean by a stratification (see section \ref{section: stratified topological spaces and exit paths}), for the moment one can just think of this stratification as a partition of $X$ indexed by the set of critical points of $f$ (which we denote by $A$)

$$
X = \bigsqcup_{a \in A} X_a.
$$

For every $x \in X$ and $t \in \R$, let us denote by $\xi^t(x)$ the image of $x$ by the flow of $\xi$ at time $t$. The subspace $X_a \subset X$ is called the $a-$\textit{stratum} and is defined to be the \textit{stable manifold of} $a$, defined by

$$
X_a=\{ x \in X \, | \limit{\xi^t(x)}{t}{+ \infty}{a} \}.
$$

Note that for every $x \in X$, the trajectory of the flow of $\xi$ starting at $x$ is defined for all time and converges to a critical point.

\begin{Example}\label{example: stratification on the other sphere}
In the case of example \ref{example: flow category of the other sphere}, this stratification is as in figure \ref{figure: stratification on other sphere}.

\begin{figure}[H]
    \centering
    \includegraphics[width=0.3\linewidth]{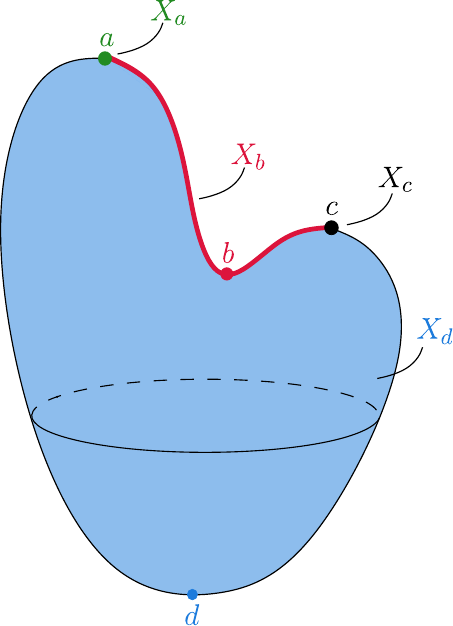}
    \caption{The strata associated with the two local maxima $a$ and $c$ are reduced to the respective critical point; indeed, the only way for a trajectory of $\xi$ to converge to a local maximum is to be constant at that point. The $b$-stratum is an embedded copy of $\R$ in $X$. The $d$-stratum consists of all the rest, and is an embedded copy of $\R^2$ in $X$.}
    \label{figure: stratification on other sphere}
\end{figure}

\end{Example}

Let $x \in X$ and consider the orbit of $x$ under the flow of $\xi$:

$$
\{\xi^t(x) \mid t \in \R\}.
$$

If $x$ is a critical point, this orbit reduces to the single point $x$. Otherwise, it is a submanifold diffeomorphic to $\R$ via the map
\[
t \longmapsto \xi^t(x).
\]
One can compactify this subspace by adding the critical points

$$
\underset{t \rightarrow - \infty}{\lim} \xi^t(x) \qquad \text{and} \qquad \underset{t \rightarrow + \infty}{\lim} \xi^t(x).
$$

In this way, one obtains a subspace homeomorphic to the closed interval $[0,1]$, and therefore a \emph{path} in $X$ (we will explain in detail below how to parameterize this subspace by $[0,1]$, see section \ref{section: stratification by stable manifolds}). The paths defined in this way belong to the following class of paths in $X$.

\begin{Definition}
    An \emph{exit path} is a path $\delta$ such that $\delta((0,1])$ is contained in a single stratum.
\end{Definition}

\begin{figure}[H]
    \centering
    \includegraphics[width=0.25\linewidth]{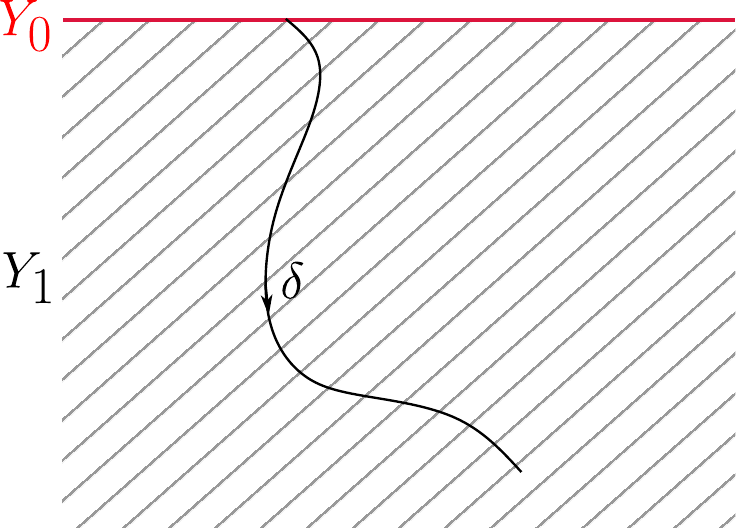}
    \caption{An exit path $\delta$ in a stratified space $Y$ with two strata $Y_0$ and $Y_1$.}
    \label{figure: exit path in intro}
\end{figure}

The class of exit paths in a given stratified topological space is not closed under concatenation in general, since the concatenation of two exit paths that each visit two strata, will visit three strata, which an exit path cannot.
However, for every stratified topological space satisfying a certain regularity condition — namely, being \emph{conically stratified} — Lurie (\cite[Appendix A]{HigherAlgebra}) constructed an $\infty$-category satisfying the following property

$$
\left\{
\begin{array}{ll}
   \text{Objects}  &  = \text{points of the stratified space} \\
    \Hom(x,y) & = \text{space of exit paths from } x \text{ to } y.
\end{array}
\right.
$$

This $\infty$-category is referred to as the $\infty$-\emph{category of exit paths} in the stratification, or the \textit{exit path $\infty$-category} associated with the stratification

By a theorem of Nicolaescu (theorem \ref{theorem: Smale if and only if Whitney}), the stratification on $X$ that we are interested in satisfies Whitney's regularity condition (b) and therefore belongs to the class of conically stratified spaces. The main theorem of this paper is the following:

\begin{Thm}\label{Main theorem}
    The $\infty$-category defined by the topological category $\Ml$ is equivalent to the exit path $\infty$-category associated with the stratification of $X$ by the stable manifolds of $\xi$.
\end{Thm}

The precise formulation of this theorem as well as additional properties of the equivalence that we construct will be given in section \ref{section: precise formulation of the main theorem}, after we have set up the necessary preliminaries in sections \ref{section: Left Kan extensions} through \ref{section: Morse-Smale pairs and flow categories}.

In the next section, we present some motivations for this theorem and some of its consequences.

\subsection{Motivations and consequences}\label{section: motivations and consequences}

\subsubsection{Morse homology and constructible sheaves}\label{section: Morse homology and constructible sheaves}

For every critical point $a$ of $f$ we denote by $|a|$ the Morse index of $a$. The \emph{Morse complex} with coefficients in $\Z / 2 \Z$ associated with $(f,\xi)$ is the graded abelian group

$$
\bigoplus_{a \in A} \Z / 2\Z \cdot a
$$

where the element $a$ is homogeneous of degree $|a|$, equipped with the differential defined for every $a \in A$ by

$$
\partial a = \sum_{\substack{b \in A \\ |b|=|a|-1}} \# \Ml(a,b) \cdot b.
$$

Note that the condition $|b| = |a|-1$ and the Smale transversality condition imply that $\Ml(a,b)$ is a finite set. We refer the reader to \cite[Section 3]{AudinDamianMorseFloer} for more details.

This construction can actually be extended to a larger class of coefficients. Consider the category $\pi_0(\Ml)$ obtained from the topological category $\Ml$ by the construction

$$
\pi_0(\Ml) = \left\{
\begin{array}{ll}
    \mathrm{Objects} & = A, \\
    \pi_0(\Ml)(a,b) & = \pi_0(\Ml(a,b)).
\end{array} \right. 
$$

For every functor $G : \pi_0(\Ml) \rightarrow \Mod_{\Z/2\Z}$, where $\Mod_{\Z / 2 \Z}$ denotes the category of $\Z / 2 \Z$-modules, we define the \emph{Morse complex with coefficients in} $G$ as the graded abelian group

$$
\bigoplus_{a \in A} G(a) \cdot a
$$

where the element $a$ is homogeneous of degree $|a|$, equipped with the differential defined for every $a \in A$ and every $g \in G(a)$ by

$$
\partial g = \sum_{|b|=|a|-1} \sum_{\gamma \in \Ml(a,b)} G(\gamma)(g) \cdot b.
$$

Note that, given a locally constant sheaf of $\Z / 2 \Z$-modules $\Fl$ on $X$, taking the monodromy of $\Fl$ along the trajectories of $\xi$ gives rise to a functor

$$
m_{\Fl} : \pi_0 (\Ml) \rightarrow \Mod_{\Z / 2\Z}
$$

and the homology of the resulting complex is known to be the homology of $X$ with local coefficients in $\Fl$. However, not every functor $G : \pi_0(\Ml) \rightarrow \Mod_{\Z / 2 \Z}$ is obtained by taking the monodromy of a local system; indeed, a necessary condition for this to hold is that for every $a,b \in A$ and $\gamma \in \Ml(a,b)$, the morphism $G(\gamma) : G(a) \rightarrow G(b)$ is an isomorphism. As long as there exists a non-constant trajectory of $\xi$ (i.e. as long as $X$ has dimension at least $1$), there exists a functor $G$ that does not satisfy this property. 

The above discussion was formulated in the setting of mod $2$ coefficients in order to keep it as simple as possible, but once one takes into account orientations on the morphism spaces of $\Ml$, it becomes valid for integer coefficients as well. Let us work with integer coefficients from now on. In view of the above discussion one can ask:

\begin{itemize}
    \item Denote by $\Loc_X$ the category of locally constant sheaves of $\Z$-modules on $X$. Does the functor
    
    $$
    \fonction{m}{\mathrm{Loc}_X}{\Fun(\pi_0(\Ml),\Mod_{\Z}),}{\Fl}{m_{\Fl}}
    $$
    
    lose information? What is its essential image?

    \item How can we describe the category $\Fun(\pi_0(\Ml), \Mod_{\Z})$?
\end{itemize}

This construction and these questions are $1$-categorical versions of analogous construction and questions formulated at the $\infty$-categorical level, as we now briefly explain. In \cite{BDHOMorse}, building on \cite{BarraudCorneaSerre}, Barraud, Damian, Humilière and Oancea introduced \emph{Morse homology with differential graded coefficients} as a generalization of classical Morse homology. Assume $X$ is connected and fix a basepoint in $X$. The category of coefficients of this theory is the category of differential graded (DG) modules over the DG algebra of cubical chains on the space of based Moore loops on $X$, denoted $C_*(\Omega X)$. This category can be described as the category of $\infty$-\emph{local systems} on $X$, i.e., the category of representations over $\Z$ of the $\infty$-groupoid associated with $X$ (this is discussed in \cite[Theorem 26]{HolsteinMorita} and \cite[Corollary 1.5]{PortaTeyssierExodromy}). In order to present the analogue of the above construction in this context, consider the differential graded category $C_{*}(\Ml)$ defined as

$$
C_{*}(\Ml) = \left\{
\begin{array}{ll}
    \text{Objects} & = \text{Objects of } \Ml, \\
    C_{*}(\Ml)(a,b) & = C_{*}(\Ml(a,b)),
\end{array}
\right.
$$

where $C_{*}(-)$ denotes the singular simplicial chain complex functor. The identity of $a$ in $C_{*}(\Ml)$ is the chain $\{c_a\}$ in $C_{0}(\Ml(a,a))$, where $c_a$ denotes the constant trajectory at $a$. The composition operation in $C_{*}(\Ml)$ is defined for every $a,b,c \in \Ml$ as:

$$
C_{*}(\Ml(a,b)) \otimes C_{*}(\Ml(b,c)) \overset{EZ}{\longrightarrow} C_{*}(\Ml(a,b)\times \Ml(b,c)) \overset{c_*}{\longrightarrow} C_{*}(\Ml(a,c)),
$$

where $EZ$ is the Eilenberg-Zilber map and $c_*$ is the map induced at the chain level by the composition operation at the level of morphism spaces of $\Ml$.

The construction of Morse homology with DG coeffiients generalizes naturally to the category of \emph{DG modules over} $C_*(\Ml)$. By a DG module over $C_{*}(\Ml)$ we mean a differential graded functor $F : C_{*}(\Ml) \rightarrow \Ch_{\Z}$, where $\Ch_{\Z}$ denotes the DG category of chain complexes of abelian groups. Equivalently, a DG module over $C_*(\Ml)$ is the datum of a chain complex $F(a)$ for every $a \in \Ml$, and of a morphism of chain complexes $F(a) \otimes C_{*}(\Ml(a,b)) \rightarrow F(b)$ for every $a,b \in \Ml$, such that the identity morphisms in $C_{*}(\Ml)$ induce identity morphisms in $\Ch_{\Z}$, and associativity holds (in the strict sense).

To summarize, once $(f,\xi)$ is fixed, the natural domain of coefficients for Morse homology with abelian coefficients seems to be the category $\Fun(\pi_0(\Ml),\Mod_{\Z})$ and, more generally, the natural domain of coefficients for Morse homology with abelian DG coefficients seems to be the category of DG modules over $C_*(\Ml)$. The following will be consequences of our main result combined with work of Ørsnes Jansen, Lurie and Porta-Teyssier on exit path categories. The proofs will be given in section \ref{section: end of the proof}.

\begin{Corollary}\label{corollary: 1-category of constructible sheaves}

The category $\Fun(\pi_0(\Ml),\Mod_{\Z})$ is equivalent to the full subcategory of the category of sheaves of abelian groups on $X$ whose objects are those sheaves that are constructible with respect to the stratification by the stable manifolds of $\xi$ (i.e., that are locally constant in restriction to each stratum).

\end{Corollary}

\begin{Corollary}\label{corollary: local systems as representations of flow category}

    The functor $m : \Loc_X \rightarrow \Fun(\pi_0(\Ml),\Mod_{\Z})$ induces an equivalence between $\Loc_X$ and the full subcategory of the right-hand side whose objects are those functors that send every morphism to an isomorphism.

\end{Corollary}

Corollary \ref{corollary: local systems as representations of flow category} is a special case of a more general statement. In order to explain this, consider the topological category $P_X$ defined as

$$
P_X = \left\{
\begin{array}{ll}
    \mathrm{\mathrm{Objects}} & = \text{points of } X, \\
    \Hom(x,y) & = \text{space of (Moore) paths from } x \text{ to } y.
\end{array} \right. 
$$

Associating to a critical point of $f$ the corresponding point of $X$, and to a (possibly broken) trajectory of $\xi$ between $x$ and $y$ the corresponding path in $X$ defines a functor of topological categories

$$\Ml \rightarrow P_X.$$

Consider further the category $\Pi_1(X) = \pi_0(P_X)$ defined from the topological category $P_X$ by taking $\pi_0$ of the morphism spaces. This is the category whose objects are the points of $X$, and whose sets of morphisms are the sets of homotopy classes of paths with fixed endpoints; it is called the \emph{fundamental groupoid} of $X$. The functor $\Ml \rightarrow P_X$ thus determines a functor $\pi_0(\Ml) \rightarrow \Pi_1(X)$.

Taking the monodromy along paths determines an equivalence of categories $\Loc_X \xrightarrow{\simeq} \Fun(\Pi_1(X),\Mod_{\Z})$ and the functor $m$ in corollary \ref{corollary: local systems as representations of flow category} identifies with the restriction

$$
\Fun(\Pi_1(X),\Mod_{\Z}) \rightarrow \Fun(\pi_0(\Ml), \Mod_{\Z})
$$

along the functor $\pi_0(\Ml) \rightarrow \Pi_1(X)$. The following is a generalization of corollary \ref{corollary: local systems as representations of flow category}.

\begin{Corollary}\label{corollary: fundamental groupoid as localization of flow category}

    The functor $\pi_0(\Ml) \rightarrow \Pi_1(X)$ exhibits $\Pi_1(X)$ as the localization of $\pi_0(\Ml)$ with respect to the set of all morphisms. In other words, for every category $\Cl$ the restriction functor
    
    $$
    \Fun(\Pi_1(X),\Cl) \rightarrow \Fun(\pi_0(\Ml),\Cl)
    $$
    
    induces an equivalence between the left-hand side and the full subcategory of the right-hand side whose objects are those functors that send every morphism to an isomorphism.

\end{Corollary}

At the $\infty$-categorical level we have the following results.

\begin{Corollary}\label{corollary: infinity category of constructible sheaves}

    The $\infty$-category of functors from the $\infty$-category associated with $\Ml$ to the derived $\infty$-category $D(\Z)$ is equivalent to the full $\infty$-subcategory of the $\infty$-category of sheaves on $X$ with values in $D(\Z)$ whose objects are those sheaves that are constructible with respect to the stratification by the stable manifolds of $\xi$. More generally, the statement holds for any compactly generated $\infty$-category in place of $D(\Z)$.

\end{Corollary}

In future work, we plan to provide a description of this $\infty$-category as the derived $\infty$-category of DG-modules over $C_*(\Ml)$. Given a DG-module over $C_*(\Ml)$, we plan to prove that the associated Morse complex with DG coefficients is a model for the homotopy colimit of the corresponding $D(\Z)$-valued constructible sheaf, regarded as a $D(\Z)$-valued functor from the $\infty$-category associated with $\Ml$.

\begin{Remark}\label{remark: Trygve and Alice}

    The formalism of \emph{Floer homotopy theory} developed in \cite{AbouzaidBlumbergI} allows one to define Morse homology with more general coefficients, such as the sphere spectrum. It was already suggested in \cite{BarraudCorneaHomotopical} that it is possible to define a Morse complex with coefficients in a spectral local system on $X$. It would be interesting to combine our main result with the results of \cite{TrygveAlice} to define Morse homology with coefficients in a constructible sheaf of spectra.

\end{Remark}

Finally, we have the following consequence of our main result.

\begin{Corollary}\label{corollary: localization of infinity category associated with M}

    The functor $\Nl(\Ml) \rightarrow \Sing(X)$, obtained by composing the equivalence $\Nl(\Ml) \rightarrow \Sing_A(X)$ provided by theorem \ref{Main theorem} with the inclusion $\Sing_A(X) \rightarrow \Sing(X)$, exhibits the fundamental $\infty$-groupoid of $X$ as the localization of $\Nl(\Ml)$ with respect to the class of all morphisms.

\end{Corollary}

Corollary \ref{corollary: localization of infinity category associated with M} is the $\infty$-categorical analogue of corollary \ref{corollary: fundamental groupoid as localization of flow category} but it is less precise in that it does not describe as explicitly the functor exhibiting the fundamental $\infty$-groupoid of $X$ as the localization of the $\infty$-category associated with $\Ml$ with respect to the class of all morphisms. We expect this functor to be the functor of $\infty$-categories induced by the functor of topological categories $\Ml \rightarrow P_X$ above.

\subsubsection{Flow categories and homotopy types}\label{section: flow categories and homotopy types}

In \cite{CohenJonesSegalMorse}, Cohen, Jones and Segal considered a Morse function $f$ on a smooth closed Riemannian manifold $X$. In this context, taking $\xi = -\nabla f$ gives rise to a flow category $\Ml$ in the same way as presented in section \ref{section: main objects and main resut}, even without assuming Smale transversality on $\xi$. They also considered the classifying space construction, which associates a topological space to every topological category. They claimed the following.

\begin{enumerate}[label=(\roman*)]
    \item The classifying space of $\Ml$, denoted $B \Ml$, is homotopy equivalent to $X$.
    \item If $\xi$ satisfies the Smale transversality condition then $B \Ml$ is homeomorphic to $X$.
\end{enumerate}

The proof of claim (i) contains a mistake. A counterexample to claim (i) as well as a corrected proof under some assumptions on $\xi$ are given in \cite{CalleLiuClassifyingSpace}. On the other hand, the proof of claim (ii) relies on the existence of compatible structures of smooth manifolds with corners on the morphism spaces of $\Ml$, which were proven later by Qin (\cite[Theorem 3.3]{QinModuli}) under the further assumption that the metric is locally trivial, in other words, under the same assumptions on $\xi$ as the ones that we make in this paper (see section \ref{section: Morse-Smale pairs and flow categories} for the precise formulation of these assumptions). The proof of (ii) also relies on an unproven lemma (\cite[Lemma 3.5]{CohenJonesSegalMorse}), which is a consequence of one of our theorems (theorem \ref{theorem: isomorphism between |-|'_A and C_A}, see also remark \ref{remark: relation to unproven lemma of CJS}).

Now suppose that $\xi$ satisfies the hypotheses of the main theorem of the present paper. Under this assumption, our results thus complete the proof of claim (ii), and therefore also that of claim (i). However, a different approach to claim (i) is possible, as we now explain. To every $\infty$-category is associated an $\infty$-groupoid (i.e., an $\infty$-category all of whose morphisms are invertible) called the \emph{groupoïd completion}, obtained by formally inverting all morphisms. The homotopy theory of $\infty$-groupoids is equivalent to the homotopy theory of CW-complexes via the fundamental $\infty$-groupoid construction and, according to corollary \ref{corollary: localization of infinity category associated with M}, the group completion of the $\infty$-category associated with $\Ml$ is the homotopy type of the manifold $X$.

Given a topological category $\Cj$, it is not always the case that its classifying space and the group completion of its associated $\infty$-category determine the same homotopy type. However, we will prove in \cite{ClassifyingSpaceVSGroupCompletion} that this is the case for $\Ml$, thus recovering (i) from our main theorem.

Property (ii) illustrates the fact that $\Ml$, as a topological category, contains more information than the homotopy type of $X$. In a similar spirit, by further taking into account compatible normal framings of the morphism spaces of $\Ml$, Cornea showed in \cite{CorneaHomotopicalII} that certain Poincaré duality spaces are not smoothable. The strategy is to derive a contradiction from the fact that any smooth realization would admit a Morse-Smale pair.

To finish this section, we would like to make a simple observation showing that the $\infty$-category associated with $\Ml$ also contains more information than the homotopy type of $X$. To this end, suppose we have another Morse-Smale pair $(f',\xi')$ on $X$ and let us denote by $\Ml'$ the associated flow category. We then have two $\infty$-categories whose groupoïd completions are the same.

$$
\xymatrix@C=0.5em{
(\infty,1)-\mathrm{Categories} \ar[d] & [\Ml] \ar@{|->}[rd] && [\Ml'] \ar@{|->}[ld] \\
\mathrm{Homotopy \, types} && [X].
}
$$

Here we denoted by $[\Ml]$, $[\Ml']$ the $\infty$-categories associated to $\Ml$ and $\Ml'$ respectively, by $[X]$ the homotopy type of $X$, and the left vertical arrow is the groupoïd completion functor.

However, it is not the case in general that $\Ml$ and $\Ml'$ are equivalent as $\infty$-categories, as we now discuss.

\begin{Example}\label{example: non equivalent flow categories on same manifold}

Let us consider the Morse-Smale pair on $\SM^2$ described in figure \ref{figure: height function on the round sphere}.

\begin{figure}[H]
    \centering
    \includegraphics[width=0.25\linewidth]{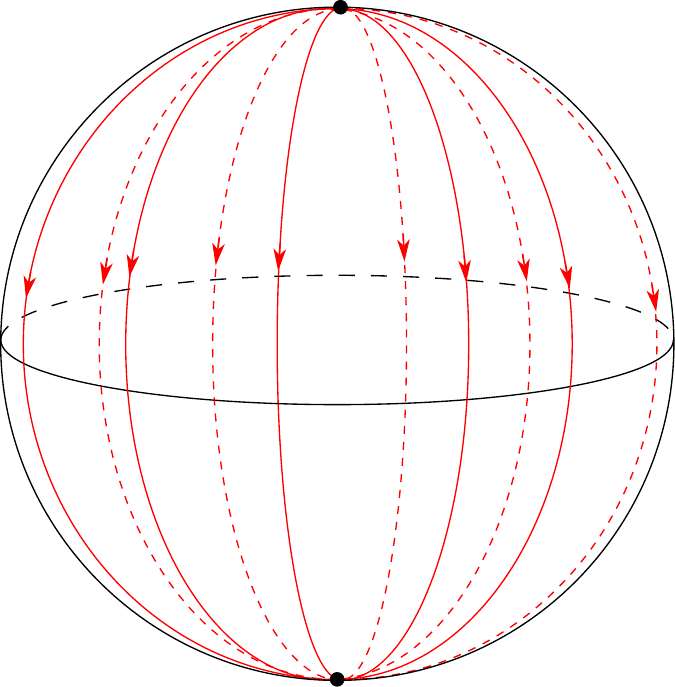}
    \caption{The Morse function is the restriction to the round $2$-sphere $\SM^2 \subset \R^3$ of the projection to the vertical axis. It has two critical points, one maximum and one minimum. The Morse-Smale pseudo-gradient vector field is the negative gradient of the function with respect to the restriction to $\SM^2$ of the Euclidean metric on $\R^3$. Mapping a point of the equator of $\SM^2$ to the trajectory that passes through it defines a homeomorphism between $\SM^1$ and the morphism space of the flow category between the maximum and the minimum.}
    \label{figure: height function on the round sphere}
\end{figure}

\end{Example}

The characterization of equivalences of ordinary categories as functors which are essentially surjective and fully faithful extends to the $(\infty,1)$-categorical setting. In the $(\infty,1)$-categorical world, the full faithfulness condition is the condition that the map induced between morphism \textit{spaces} is a \textit{weak homotopy equivalence}. The definition of essential surjectivity in the $\infty$-categorical setting will be given in section \ref{section: infinity categories}. The flow category from example \ref{example: non equivalent flow categories on same manifold} has a morphism space which is homeomorphic to $\SM^1$, while none of the morphism spaces of the flow category from example \ref{example: flow category of the other sphere} are weakly homotopy equivalent to $\SM^1$. There can therefore not exist a fully faithful functor between the $\infty$-categories defined by these two flow categories, and in particular these two $\infty$-categories are not equivalent.

It is therefore natural to ask whether $\Ml$, as an $\infty$-category, reflects some extra structure on $X$. Theorem \ref{Main theorem} can be regarded as an answer to this question.

\subsection{Idea of the proof and organization of the paper}\label{section: idea of the proof and organization}

Let us explain briefly the principle of our proof. We will work in the model of simplicial sets for $(\infty,1)$-categories, and we will only consider $\infty$-categories that are $(\infty,1)$-categories. For us, an $\infty$-\emph{category} will thus be a simplicial set that satisfies the inner horn filling property. The $\infty$-category defined by the topological category $\Ml$ is the \emph{homotopy coherent nerve} of $\Ml$, denoted $\Nl(\Ml)$. The exit path $\infty$-category associated with the stratification of $X$ by the stable manifolds of $\xi$ is the \emph{stratified singular simplicial set} of the stratified space $X$, denoted $\Sing_A(X)$. For small values of $n$, here is a description of these simplicial sets and how we compare them.

\underline{$n=0$}: The $0$-simplices of $\Sing_A(X)$ are the points of $X$. On the other hand, the $0$-simplices of $\Nl(\Ml)$ are the objects of $\Ml$, i.e., the critical points of $f$.  There is therefore an inclusion

$$
\Nl(\Ml)_0 \subseteq \Sing_A(X)_0.
$$

\underline{$n=1$}: The $1$-simplices of $\Sing_A(X)$ are the exit paths in $X$. On the other hand, a $1$-simplex $\sigma \in \Nl(\Ml)$ consists of the following datum:

$$
\sigma =
\left\{
\begin{array}{l}
    a_0, a_1 \in A  \\
    \gamma_0 \in \Ml(a_0,a_1). 
\end{array}
\right.
$$

To the trajectory $\gamma_0$ is associated a path in $X$, but the latter is an exit path only when $\gamma_0$ is \emph{not} broken. We therefore introduce the subset $(S_{\Ml})_1 \subset \Nl(\Ml)_1$ defined as

$$
(S_{\Ml})_1 = \{ \sigma \in \Nl(\Ml)_1 \mid \gamma_0 \text{ is unbroken} \},
$$

and we have an inclusion

$$
(S_{\Ml})_1 \subseteq \Sing_A(X)_1.
$$

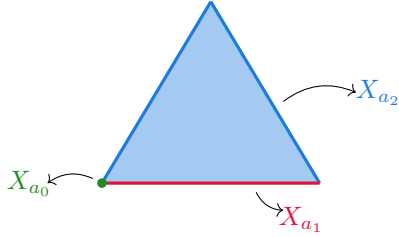
\begin{figure}[H]
    \centering
    \begin{tikzpicture}[scale=1.2]
        \draw[green] node at (-2,0) {$X_{a_0}$};
        \draw[red] node at (1,-0.4) {$X_{a_1}$};
        \draw[blue] node at (1.85,1) {$X_{a_2}$};
        
        \draw[->, bend right]  (-1.3,0.05) to (-1.8,0);
        \draw[->, bend right]  (0.5,-0.1) to (0.8,-0.3);
        \draw[->, bend left]  (0.8,0.9) to (1.6,1);
        
        \draw[fill, blue, opacity=0.4] (-1.2,0)--(1.2,0)--(0,2)--cycle;
        
        \draw[blue, line width=1.2pt] (0,2)--(1.2,0);
        \draw[blue, line width=1.2pt] (0,2)--(-1.2,0);
        
        \draw[red, line width=1.2pt] (-1.2,0)--(1.2,0);
        
        \fill[green] (-1.2,0) circle (1.5pt);
        
    \end{tikzpicture}
    \caption{Each of the three colored parts of $\Delta^2$ is mapped to the corresponding stratum of $X$.}
    \label{figure: homotopy coherent composition of two exit paths}
\end{figure}

\underline{$n=2$}: The $2$-simplices of $\Sing_A(X)$ are the \emph{homotopy coherent compositions of two exit paths}. These are the singular $2$-simplices of $X$ such that there exist three strata $X_{a_0}$, $X_{a_1}$, $X_{a_2}$ satisfying the condition described in figure \ref{figure: homotopy coherent composition of two exit paths}.

On the other hand, a $2$-simplex $\sigma \in \Nl(\Ml)_2$ consists of the following datum:

$$
\sigma =
\left\{
\begin{array}{l}
    a_0, a_1, a_2 \in A  \\
    \gamma_0 \in \Ml(a_0,a_1), \gamma_1 \in \Ml(a_1,a_2), \gamma_2 \in \Ml(a_0,a_2) \\
    h : [0,1] \rightarrow \Ml(a_0,a_2) \text{, } h(0)=\gamma_2 \text{ and } h(1) = \gamma_1 \circ \gamma_0 \text{ (i.e., a homotopy between } \gamma_2 \text{ and } \gamma_1 \circ \gamma_0).
\end{array}
\right.
$$

In order to explain how we compare $\Nl(\Ml)_2$ and $\Sing_A(X)_2$, let us first precise how the trajectories of $\xi$ are parameterized in order to yield paths in $X$. To every trajectory $\gamma \in \Ml(a_0,a_2)$ is associated a parameterization $\widetilde{\gamma} : [f(a_2),f(a_0)] \rightarrow X$ uniquely characterized by the condition that $f(\widetilde{\gamma}(s))=s$. Combining this with the map $h$ we get a map (see figure \ref{figure: stratified 2 simplex out of flow simplex})

$$
\fonction{f_{\sigma}}{[0,1] \times [f(a_2),f(a_0)]}{X}{(t,s)}{\widetilde{h(t)}(s).}
$$

When attempting to associate to $f_{\sigma}$ a $2$-simplex of $\Sing_A(X)$, we encounter two difficulties. Firstly, the presence of trajectories in $\Ml(a_0,a_2)$ that are possibly broken many time, has the effect that $f_{\sigma}$ may visit other strata than those of $a_0$, $a_1$ and $a_2$, which a $2$-simplex of $\Sing_A(X)$ cannot, by definition. To solve this, we introduce the subset $(S_{\Ml})_2 \subset \Nl(\Ml)_2$ defined as

$$
(S_{\Ml})_2 = \{ \sigma \in \Nl(\Ml)_2 \mid \gamma_0, \gamma_1 \text{ are unbroken and for every } 0 \le t < 1, h(t) \text{ is unbroken}\}.
$$

Secondly, $f_{\sigma}$ has a rectangular domain, not a simplicial one. To solve this, we observe that, by definition, $f_{\sigma}$ is constant equal to $a_2$ on $[0,1] \times \{f(a_2)\}$ and equal to $a_0$ on $[0,1] \times \{f(a_0)\}$. Hence, $f_{\sigma}$ factors through the space obtained from $[0,1] \times [f(a_2),f(a_0)]$ by collapsing $[0,1] \times \{f(a_2)\}$ and $[0,1] \times \{f(a_0)\}$ to distinct points. We denote by $\widetilde{f_{\sigma}}$ the map induced by $f_{\sigma}$ on this quotient space. Altogether, when $\sigma \in (S_{\Ml})_2$ and $a_0, a_1, a_2$ are distinct, the situation is described in figure \ref{figure: stratified 2 simplex out of flow simplex}.

\begin{figure}[H]
\centering
\begin{tikzcd}
\tikz[scale=1]{
        \draw node at (1.44,1.12) {$-$ \small{$f(a_0)$}};
        \draw node at (1.44,0.22) {$-$ \small{$f(a_1)$}};
        \draw node at (1.44,-1.28) {$-$ \small{$f(a_2)$}};
        \draw node at (-1,-1.28) {\small{|}};
        \draw node at (1,-1.28) {\small{|}};
        \draw node at (-0.2,-1.28) {\small{|}};
        \draw node at (-0.2,-1.65) {$t$};
        \draw node at (-1,-1.65) {$0$};
        \draw node at (1,-1.65) {$1$};

        \draw[fill, blue, opacity=0.5] (-1,1.2)--(1,1.2)--(1,-1.2)--(-1,-1.2)--cycle;
    
        \draw[blue] (1,0.3)--(1,-1.2)--(-1,-1.2)--(-1,1.2);
        \draw[blue] (-0.2,1.2)--(-0.2,-1.2);
        \draw[red, line width=1.4pt] (1,1.2)--(1,0.3);
        \draw[green, line width=1.4pt] (-1,1.2)--(1,1.2);
        
        
        \draw[red, ->, shift={(1,0.67)}] (0,0) -- (0,-0.0001);
        \draw node at (1.3,0.67) {\small{$\gamma_0$}};
        \draw[blue, ->, shift={(1,-0.53)}] (0,0) -- (0,-0.0001);
        \draw node at (1.3,-0.53) {\small{$\gamma_1$}};
        \draw[blue, ->, shift={(-1,0)}] (0,0) -- (0,-0.0001);
        \draw node at (-1.3,0) {\small{$\gamma_2$}};
        \draw[blue, ->, shift={(-0.2,0)}] (0,0) -- (0,-0.0001);
        \draw node at (0.2,0) {\small{$h(t)$}};
        
   } \ar[rd, swap, "{\displaystyle f_{\sigma}}"] \ar[rr, "/ \sim"] &&
   \tikz[scale=1]{

        \draw[fill, blue, opacity=0.5] (0,1.2)--(1,0.3)--(0,-1.2)--cycle;
        \draw[blue] (1,0.3)--(0,-1.2)--(0,1.2);
        \fill[blue] (0,-1.2) circle (1.3pt);
    
        \draw[red, line width=1.4pt] (0,1.2)--(1,0.3);
        
        \fill[green] (0,1.2) circle (1.3pt);
        
        \draw[->, bend left]  (0.1,1.3) to (0.4,1.4);
        \draw[->, bend left]  (0.6,0.75) to (0.9,0.85);
        \draw[->, bend right]  (0.6,-0.55) to (0.9,-0.65);
        
        
        \draw[green] node at (0.7,1.3) {$X_{a_0}$};
        \draw[red] node at (1.2,0.75) {$X_{a_1}$};
        \draw[blue] node at (1.2,-0.75) {$X_{a_2}$};
        
    } \ar[ld, "{\displaystyle \widetilde{f_{\sigma}}}"] \\
& X
\end{tikzcd}
\caption{}
\label{figure: stratified 2 simplex out of flow simplex}
\end{figure}

The quotiented rectangle is homeomorphic to the $2$-simplex, and the condition that $\sigma \in (S_{\Ml})_2$ is equivalent to the condition that $\widetilde{f_{\sigma}}$ belongs to $\Sing_A(X)_2$. We thus obtain a map

$$
\fonctionsansnom{(S_{\Ml})_2}{\Sing_A(X)_2}{\sigma}{\widetilde{f_{\sigma}}.}
$$

Extending these definitions to every nonnegative integer $n$, we will define a simplicial subset $S_{\Ml} \subseteq \Nl(\Ml)$, which we will call the \emph{flow coherent nerve} of $\Ml$, and which we will prove to be an $\infty$-category. Its $n$-simplices can be thought of as generalizing the notion of unbroken trajectory of $\xi$ to higher simplices of $\Nl(\Ml)$. We have (roughly) a zigzag of $\infty$-categories

$$
\xymatrix{
\Nl(\Ml) && \Sing_A(X). \\
& S_{\Ml} \ar[lu]^-{i_1} \ar[ru]_-{i_2}
}
$$

\begin{Remark}\label{remark: simplification}

    The assignment $\sigma \mapsto \widetilde{f_{\sigma}}$ only yields a morphism of \emph{semi-simplicial sets} from $S_{\Ml}$ to $\Sing_A(X)$. Following \cite{HiroFunctors}, we can think of this morphism as a functor of $\infty$-categories that isn't strictly unital, i.e., that preserves identity morphisms only up to homotopy. Fortunately, the results of \cite{HiroFunctors} guarantee the existence of a genuine functor that we denote here $i_2$.

\end{Remark}

We will prove that $i_1$ and $i_2$ are both equivalences of $\infty$-categories, by proving that these are essentially surjective and fully faithful. The idea behind the full faithfulness is the following. Given two critical points $a$ and $b$, we will give an explicit description of the space of morphisms of $S_{\Ml}$ between $a$ and $b$ as the space of unbroken trajectories between $a$ and $b$. Passing to the morphism spaces between $a$ and $b$, the above zigzag thus becomes:

$$
\xymatrix{
\parbox{4cm}{
\begin{center}
$\Nl(\Ml)(a,b)$ = Space of possibly broken trajectories from $a$ to $b$
\end{center}
}
&&
\parbox{4cm}{\begin{center} 
$\Sing_A(X)(a,b)$ =
Space of exit paths from $a$ to $b$.
\end{center}
}
\\
&
\parbox{4.8cm}{
\begin{center} 
$S_{\Ml}(a,b)$ =
Space of unbroken trajectories from $a$ to $b$
\end{center}
}
\ar[lu]^-{i_1(a,b)}
\ar[ru]_-{i_2(a,b)}
}
$$

We will show that the map $i_1(a,b)$ is the inclusion and the map $i_2(a,b)$ is the one that carries an unbroken trajectory to the corresponding exit path. These are both weak homotopy equivalences: this is an already known result in the case of $i_1(a,b)$, and we will prove it in the case of $i_2(a,b)$.

The paper is organized as follows. Section \ref{section: detailed summary of the proof} consist of a very detailed summary of our proof. It does not contain detailed constructions or proofs, but rather aims to explain the ideas underlying the overall structure of the proof, together with the constructions and results involved. It also illustrates these with numerous figures.

The remainder of the paper consists of carrying out the constructions and completing the arguments. Section \ref{section: infinity categories, stratified spaces and Morse-Smale pairs} is a brief introduction to $(\infty,1)$-category theory, stratified homotopy theory and Morse flow categories. We recommend to consult it depending on the reader's background. In section \ref{section: spaces of exit paths as spaces of pseudo gradient trajectories}, we review the theory of Whitney stratifications and use it to prove that the map $i_2(a,b)$ above is a weak homotopy equivalence. In section \ref{section: the homotopy coherent nerve}, we study the homotopy coherent nerve functor on the category of topological categories. As is already apparent in the sketch of proof given above, the simplices of the homotopy coherent nerve of $\Ml$ naturally give rise to \emph{higher rectangles} in $X$, not simplices on the nose. In section \ref{section: cubical stratified geometric realization}, we  construct from these higher rectangles an alternative to the usual stratified geometric realization functor, and give the definition of the flow coherent nerve of $\Ml$. In section \ref{section: comparison between realizations}, we prove a comparison result between the two geometric realizations. This can be understood as an instance of the comparison between the theory of simplicial sets and the theory of cubical sets. Although a mere isomorphism between these functors does not exist, our result is sufficient for our purposes, thanks to a theorem of Tanaka that we review in section \ref{section: end of the proof}. In section \ref{section: unbroken simplices}, we prove that $S_{\Ml}$ is an $\infty$-category and, roughly, that the two functors $i_1$ and $i_2$ above are equivalences.

In section \ref{section: stratified categories as exit path categories}, we give a general framework to which our methods extend. We prove a generalization of our main result where Morse flow categories are replaced by another class of objects called \emph{stratified categories} inspired from them, which we introduce.

\subsection{Detailed summary of the proof}\label{section: detailed summary of the proof}

The exit path $\infty$-category construction is defined by Lurie in the model of \emph{simplicial sets} for $\infty$-categories. It is modeled by a simplicial set called the \emph{stratified singular simplicial set} of the stratified space $X$, and denoted $\Sing_A(X)$ \footnote{Note that there is an abuse of notation here, since this simplicial depends on more than simply the topological space $X$ and the set $A$.}. Our strategy is to first bring $\Ml$ in the world of simplicial sets. This is done through the \emph{homotopy coherent nerve functor}

$$
\Nl : \Cat_{\Top} \longrightarrow \sSet
$$

from the category of topological categories to the category of simplicial sets, which identifies the theory of $\infty$-categories formulated in the language of topological categories, with the theory of $\infty$-categories formulated in the language of simplicial sets. In other words, the $\infty$-category defined by $\Ml$ is the same as that defined by $\Nl(\Ml)$. We are thus reduced to comparing the simplicial sets $\Nl(\Ml)$ and $\Sing_A(X)$. Our proof consists in constructing a zigzag of morphisms of simplicial sets that are equivalences of $\infty$-categories

\[\label{zigzag} \tag{Z}
\Nl(\Ml) \leftarrow S_{\Ml} \rightarrow \Sing_{C_A,r}(X) \rightarrow \Sing'_{A,r}(X) \rightarrow \Sing_A(X).
\]

In this section, we summarize the construction of these $\infty$-categories and equivalences.

To begin with, let us describe the simplicial set $\Sing_A(X)$. This is a simplicial subset of $\Sing(X)$, that is, the $n$-simplices of $\Sing_A(X)$ are some particular singular $n$-simplices in $X$. Let us describe them for $n$ going from $0$ to $3$ (complete definitions are given in section \ref{section: stratified spaces, exit paths and stratified simplicial sets}).

\underline{$n = 0$:} $\Sing_A(X)_0 = \Sing(X)_0$, in other words, the $0$-simplices of $\Sing_A(X)$ are the points of $X$.

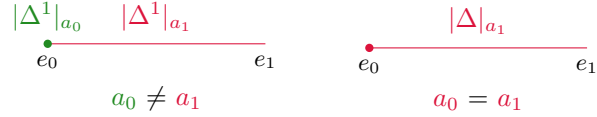
\begin{figure}[H]
    \centering
    \begin{tabular}{cccc}
    \begin{tikzpicture}[scale=1.2]
        \draw[green] node at (-1.2,0.3) {$|\Delta^1|_{a_0}$};
        \draw[red] node at (0,0.3) {$|\Delta^1|_{a_1}$};
        \draw node at (-1.2,-0.2) {\small{$e_0$}};
        \draw node at (1.2,-0.2) {\small{$e_1$}};
        \draw node at (0,-0.6) {\textcolor{green}{$a_0$} $\neq$ \textcolor{red}{$a_1$}};
        
        \draw[red] (-1.2,0)--(1.2,0);
        
        \fill[green] (-1.2,0) circle (1.3pt);
    \end{tikzpicture}

    & &

    \begin{tikzpicture}[scale=1.2]
        \draw node at (-1.2,-0.2) {\small{$e_0$}};
        \draw node at (1.2,-0.2) {\small{$e_1$}};
        \draw[red] node at (0,0.3) {$|\Delta|_{a_1}$};
        \draw node at (0,-0.6) {\textcolor{red}{$a_0$} = \textcolor{red}{$a_1$}};
    
        \draw[red] (-1.2,0)--(1.2,0);
        \fill[red] (-1.2,0) circle (1.3pt);
    \end{tikzpicture}
    \end{tabular}
    \caption{The possible stratifications of the standard $1$-simplex depending on whether $a_0$ and $a_1$ are different or not.}
    \label{figure: stratification of the standard 1-simplex}
\end{figure}

\underline{$n = 1$:} The $1$-simplices of $\Sing_A(X)$ are the \emph{exit paths} in $X$. In other words, a $1$-simplex of $\Sing_A(X)$ is the datum of a sequence of two (not necessarily different) critical points $a_0, a_1$ and a path $|\Delta^1| \rightarrow X$ that is compatible with the stratification of $|\Delta^1|$ by $\{a_0,a_1\}$ described in figure \ref{figure: stratification of the standard 1-simplex}.

Let us clarifiy that, by "compatible with the stratification", we mean that the $a_0$-stratum (resp. $a_1$-stratum) of $|\Delta^1|$ is mapped to the $a_0$-stratum (resp. $a_1$-stratum) of $X$. We indeed recover the definition of an exit path in this way.

Before moving on to higher $n$, we would like to point out that given two critical points $a_0$ and $a_1$, there does not always exist an exit path starting in the $a_0$-stratum and ending in the $a_1$-stratum. It is useful to endow $A$ with the relation $\le$ defined as $a_0 \le a_1$ if and only if such an exit path exists (equivalent characterizations of this condition are given in corollary \ref{corollary: equivalent characterizations of partial order on A}). As the notation suggests, this turns out to be a partial order.
 
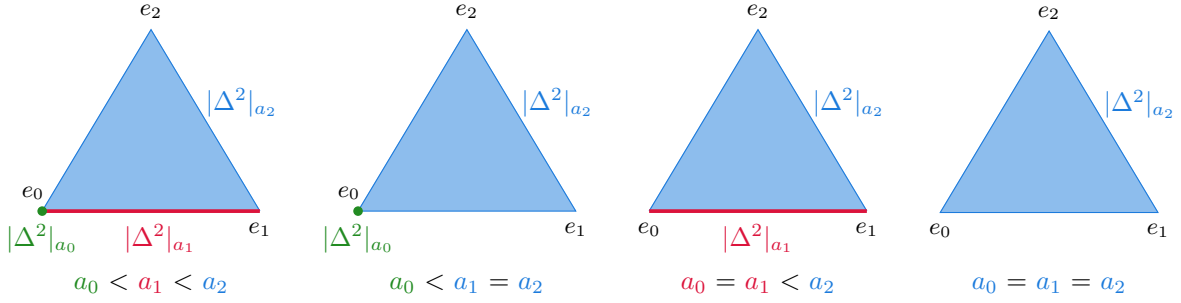
\begin{figure}[H]
    \centering
    $$\begin{array}{cccccc}
    \begin{tikzpicture}[scale=1.2]
        \draw[red] node at (0.1,-0.3) {$|\Delta^2|_{a_1}$};
        \draw[green] node at (-1.2,-0.3) {$|\Delta^2|_{a_0}$};
        \draw[blue] node at (1,1.2) {$|\Delta^2|_{a_2}$};
        \draw node at (-1.3,0.2) {\small{$e_0$}};
        \draw node at (1.2,-0.2) {\small{$e_1$}};
        \draw node at (0,2.2) {\small{$e_2$}};
        \draw node at (0,-0.8) {\textcolor{green}{$a_0$} < \textcolor{red}{$a_1$} < \textcolor{blue}{$a_2$}};
        
        \draw[blue] (-1.2,0)--(0,2)--(1.2,0);
        \draw[fill, blue, opacity=0.5] (-1.2,0)--(0,2)--(1.2,0);
        
        \draw[red, line width=1.4pt] (-1.2,0)--(1.2,0);
        
        \fill[green] (-1.2,0) circle (1.5pt);
        
    \end{tikzpicture}
    &
    \begin{tikzpicture}[scale=1.2]
        \draw[green] node at (-1.2,-0.3) {$|\Delta^2|_{a_0}$};
        \draw[blue] node at (1,1.2) {$|\Delta^2|_{a_2}$};
        \draw node at (-1.3,0.2) {\small{$e_0$}};
        \draw node at (1.2,-0.2) {\small{$e_1$}};
        \draw node at (0,2.2) {\small{$e_2$}};
        \draw node at (0,-0.8) {\textcolor{green}{$a_0$} < \textcolor{blue}{$a_1$} = \textcolor{blue}{$a_2$}};
        
        \draw[blue] (-1.2,0)--(0,2)--(1.2,0)--cycle;
        \draw[fill, blue, opacity=0.5] (-1.2,0)--(0,2)--(1.2,0);
        
        \fill[green] (-1.2,0) circle (1.5pt);
    \end{tikzpicture}
    &
    \begin{tikzpicture}[scale=1.2]
        \draw node at (-1.2,-0.2) {\small{$e_0$}};
        \draw node at (1.2,-0.2) {\small{$e_1$}};
        \draw node at (0,2.2) {\small{$e_2$}};
        \draw[red] node at (0,-0.3) {$|\Delta^2|_{a_1}$};
        \draw[blue] node at (1,1.2) {$|\Delta^2|_{a_2}$};
        \draw node at (0,-0.8) {\textcolor{red}{$a_0$} = \textcolor{red}{$a_1$} < \textcolor{blue}{$a_2$}};
        
        \draw[blue] (-1.2,0)--(0,2)--(1.2,0);
        \draw[fill, blue, opacity=0.5] (-1.2,0)--(0,2)--(1.2,0);
        
         \draw[red, line width=1.4pt] (-1.2,0)--(1.2,0);
    \end{tikzpicture}
     &
    \begin{tikzpicture}[scale=1.2]
        \draw node at (-1.2,-0.2) {\small{$e_0$}};
        \draw node at (1.2,-0.2) {\small{$e_1$}};
        \draw node at (0,2.2) {\small{$e_2$}};
        \draw[blue] node at (1,1.2) {$|\Delta^2|_{a_2}$};
        \draw node at (0,-0.8) {\textcolor{blue}{$a_0$} = \textcolor{blue}{$a_1$} = \textcolor{blue}{$a_2$}};
        
        \draw[white] node at (0,-0.3) {$|\Delta|_{a_0}$};
        
        \draw[blue] (-1.2,0)--(0,2)--(1.2,0)--cycle;
        \draw[fill, blue, opacity=0.5] (-1.2,0)--(0,2)--(1.2,0);
    \end{tikzpicture}
    \end{array}$$
    \caption{The possible stratifications of the standard $2$-simplex determined by a sequence of three critical points. On the left: the case when the critical points are distinct. When they are not distinct, the stratification is obtained by merging, for each group of identical points, the strata corresponding to these points.}
    \label{figure: stratification of the standard 2-simplex}
\end{figure}

\underline{$n=2$:} A $2$-simplex of $\Sing_A(X)$ is the datum of a sequence of three (not necessarily distinct) critical points $a_0, a_1, a_2$ and a singular $2$-simplex $|\Delta^2| \rightarrow X$ that is compatible with the stratification of $|\Delta^2|$ by $\{a_0,a_1,a_2\}$ described in figure \ref{figure: stratification of the standard 2-simplex}.

We note that if such a singular $2$-simplex in $X$ exists, then $a_0 \le a_1 \le a_2$.

\begin{figure}[H]
    \centering
    $$\begin{array}{ccccc}
    \begin{tikzpicture}[scale=1.2]
        
        \draw (-1.2,0)--(0,2)--(1.2,0);
        \draw[blue] (1.2,0)--(1.7,0.6);
        \draw (1.7,0.6)--(0,2);
        \draw[dashed, blue] (-1.2,0)--(1.7,0.6);
        
        \draw[fill, opacity = 0.1] 
            (-1.2,0)--(0,2)--(1.7,0.6)--(1.2,0)--cycle;
         \draw[fill, blue, opacity=0.4] 
            (-1.2,0)--(1.7,0.6)--(1.2,0)--cycle;

        \draw[red, line width=1.4pt] (-1.2,0)--(1.2,0);
        \fill[green] (-1.2,0) circle (1.5pt);

        \draw[green] node at (-1.2,-0.3) {$|\Delta^3|_{a_0}$};
        \draw[blue] node at (1.9,0.2) {$|\Delta^3|_{a_2}$};
        \draw node at (1.2,1.5) {$|\Delta^3|_{a_3}$};
        \draw[red] node at (0.1,-0.3) {$|\Delta^3|_{a_1}$};
        \draw node at (-1.3,0.2) {\small{$e_0$}};
        \draw node at (1.2,-0.2) {\small{$e_1$}};
        \draw node at (1.9,0.6) {\small{$e_2$}};
        \draw node at (0,2.2) {\small{$e_3$}};
        \draw node at (0.25,-0.8) {\textcolor{green}{$a_0$} < \textcolor{red}{$a_1$} < \textcolor{blue}{$a_2$} < $a_3$};
        
    \end{tikzpicture}
    &
    \begin{tikzpicture}[scale=1.2]
        
        \draw (-1.2,0)--(0,2);
        \draw[blue] (-1.2,0)--(1.2,0)--(1.7,0.6);
        \draw (1.7,0.6)--(0,2);
        \draw[dashed, blue] (-1.2,0)--(1.7,0.6);
        \draw (1.2,0)--(0,2);
        
        \draw[fill, opacity = 0.1] 
                (-1.2,0)--(0,2)--(1.7,0.6)--(1.2,0)--cycle;
         \draw[fill, blue, opacity=0.4] 
                (-1.2,0)--(1.7,0.6)--(1.2,0)--cycle;

        \fill[green] (-1.2,0) circle (1.5pt);

        \draw[green] node at (-1.2,-0.3) {$|\Delta^3|_{a_0}$};
        \draw[blue] node at (1.9,0.2) {$|\Delta^3|_{a_2}$};
        \draw node at (1.2,1.5) {$|\Delta^3|_{a_3}$};
        \draw node at (-1.3,0.2) {\small{$e_0$}};
        \draw node at (1.2,-0.2) {\small{$e_1$}};
        \draw node at (1.9,0.6) {\small{$e_2$}};
        \draw node at (0,2.2) {\small{$e_3$}};
        \draw node at (0.25,-0.8) {\textcolor{green}{$a_0$} < \textcolor{blue}{$a_1$} = \textcolor{blue}{$a_2$} < $a_3$};
        
    \end{tikzpicture}
    &
    \begin{tikzpicture}[scale=1.2]
        
        \draw (-1.2,0)--(1.2,0)--(1.7,0.6)--(0,2)--cycle;
        \draw[dashed] (-1.2,0)--(1.7,0.6);
        \draw (1.2,0)--(0,2);
        
        \draw[fill, opacity = 0.1] (-1.2,0)--(1.2,0)--(1.7,0.6)--(0,2)--cycle;

        \fill[green] (-1.2,0) circle (1.5pt);

        \draw[green] node at (-1.2,-0.3) {$|\Delta^3|_{a_0}$};
        \draw node at (1.2,1.5) {$|\Delta^3|_{a_3}$};
        \draw node at (-1.3,0.2) {\small{$e_0$}};
        \draw node at (1.2,-0.2) {\small{$e_1$}};
        \draw node at (1.9,0.6) {\small{$e_2$}};
        \draw node at (0,2.2) {\small{$e_3$}};
        \draw node at (0.25,-0.8) {\textcolor{green}{$a_0$} < $a_1$ = $a_2$ = $a_3$};
        
    \end{tikzpicture}
    \end{array}$$
    \caption{Some stratifications of the standard $3$-simplex determined by a sequence of four critical points $a_0,a_1,a_2,a_3$. On the left: the case when they are distinct. As in the previous cases, when they are not distinct, the stratification is obtained by merging, for each group of identical points, the strata corresponding to these points.}
    \label{figure: some stratifications of the 3 simplex}
\end{figure}
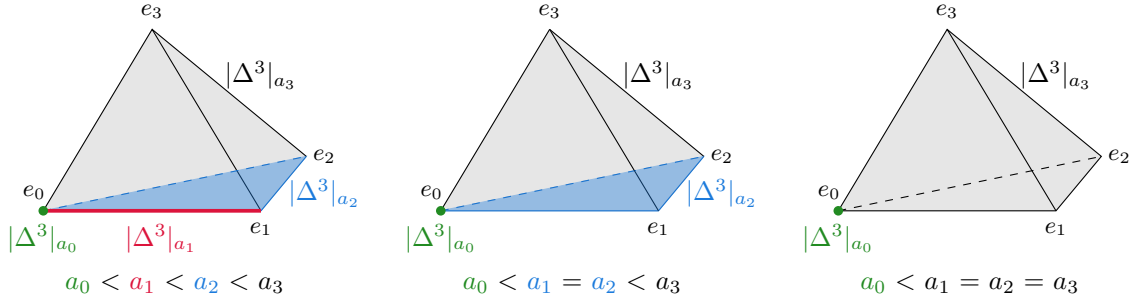

\underline{$n=3$:} A $3$-simplex of $\Sing_A(X)$ is the datum of a sequence of four (not necessarily distinct) critical points $a_0, a_1, a_2, a_3$ and a singular $3$-simplex $|\Delta^3| \rightarrow X$ that is compatible with the stratification of $|\Delta^3|$ determined by $\{a_0,a_1,a_2, a_3\}$. Some examples of such stratifications are presented in figure \ref{figure: some stratifications of the 3 simplex}.

We note again that if such a singular $3$-simplex in $X$ exists, then $a_0 \le a_1 \le a_2 \le a_3$.

These definitions extend to every nonnegative integer $n$, namely:

\underline{$n \ge 0$ arbitrary:} Every sequence of $n+1$ critical points $a_0,\hdots,a_n$ determines a stratification of the standard $n$-simplex $|\Delta^n|$ by $\{a_0,\hdots,a_n\}$, and an $n$-simplex of $\Sing_A(X)$ is a singular $n$-simplex $|\Delta^n| \rightarrow X$ that is compatible with the stratifications, for some (necessarily increasing \footnote{By saying "increasing", we allow repetitions. If no repetitions are allowed, we will say "strictly increasing".} with respect to $\le$) sequence of $n+1$ critical points.

Note that this simplicial set can be associated to any stratified space. Extra regularity conditions on the stratification are needed in order for this simplicial set to be an $\infty$-category (see in particular theorem \ref{theorem: stratified simplicial set of conically stratified space is infty cat} and proposition \ref{proposition: Sing_A infinity category in terms of starting point evaluation}).

In order to compare $\Sing_A(X)$ with $\Nl(\Ml)$, we need to understand what $\Nl(\Ml)$ is. The functor $\Nl$ admits a left adjoint, which we denote by $F$. By adjunction we have for every nonnegative integer $n$

$$
\begin{aligned}
    \Nl(\Ml)_n & = \Hom_{\sSet}(\Delta^n,\Nl(\Ml)) \\
    & = \Hom_{\Cat_{\Top}}(F(\Delta^n),\Ml).
\end{aligned}
$$

For every nonnegative integer $n$, describing explicitly what are the $n$-simplicies of the homotopy coherent nerve of a topological category therefore amounts to describing explicitly the topological category $F(\Delta^n)$. We define it and prove a number of properties in section \ref{section: the homotopy coherent nerve}. Its set of objects is $\{0, 1, \hdots ,n \}$. Denoting by $I$ the closed interval $[0,1]$, there are homeomorphisms

$$
F(\Delta^n)(i,j) \simeq \left\{
    \begin{array}{ll}
        * & \mbox{if } i=j \\
        I^{j-i-1} & \mbox{if } i < j.
    \end{array}
\right.
$$

Consider an $n$-simplex $\sigma : F(\Delta^n) \rightarrow \Ml$ of $\Nl(\Ml)$ and denote the critical point $\sigma(i)$ by $a_i$. This $n$-simplex comes with a map $\sigma(0,n) : I^{n-1} \rightarrow \Ml(a_0,a_n)$. Combining this with the map $\Ml(a_0,a_n) \rightarrow C^0([f(a_n),f(a_0)],X)$ that associates, to every trajectory from $a_0$ to $a_n$, its parametrization by the interval $[f(a_n),f(a_0)]$ (see \ref{construction: parametrization of all trajectories by values of the function} for the precise construction), we obtain a map

$$
f_{\sigma} : I^{n-1} \times [f(a_n),f(a_0)] \rightarrow X.
$$

\begin{figure}[H]
    \centering
    \includegraphics[width=0.23\linewidth]{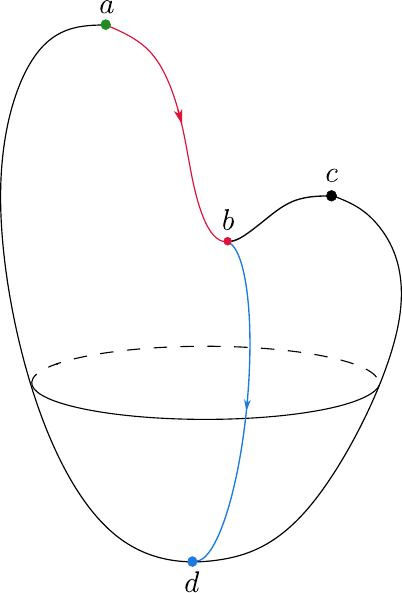}
    \caption{Coming back to example \ref{example: flow category of the other sphere}, if $n=1$ and $\sigma(0,1)(*)$ is the above trajectory from $a$ to $d$, broken at $b$, then $f_{\sigma}$ visits the strata of $a$ and $d$, but also that of $b$. It therefore does not yield an exit path.}
    \label{figure: example of broken simplex}
\end{figure}

One can see $f_{\sigma}$ as an $n$-dimensional singular cube of $X$ \footnote{This cube is degenerate if $a_0=a_n$.}. Recall that our goal is to compare $\Nl(\Ml)$ and $\Sing_A(X)$. We would therefore wish to associate to $f_{\sigma}$ an $n$-simplex of $\Sing_A(X)$, that is, a stratified singular $n$-simplex in $X$. Moreover, the $i^{\text{th}}$ vertex of this simplex should be the critical point $a_i$. But observe that, by definition of $\Sing_A(X)$, such an $n$-simplex can only visit the $a_i$-strata for $0 \le i \le n$. In comparison, the presence of \emph{broken trajectories} in the morphism spaces of $\Ml$ has the effect that one does not have such control on the strata which the map $f_{\sigma}$ can visit. Consider for example the $1$-simplex of $\Nl(\Ml)$ described in figure \ref{figure: example of broken simplex}.

More precisely, that map $f_{\sigma}$ has to visit \emph{at least} the $a_i$-strata for $0 \le i \le n$. Indeed, by functoriality of $\sigma : F(\Delta^n) \rightarrow \Ml$, the image by $\sigma(0,n)$ of the morphism of $F(\Delta^n)$ between $0$ and $n$ obtained as the composition

$$
* \simeq F(\Delta^n)(0,1) \times F(\Delta^n)(1,2) \times \hdots \times F(\Delta^n)(n-1,n) \rightarrow F(\Delta^n)(0,n)
$$

is broken at least at each of the $a_i$, and therefore the $a_i$'s belong to the image of $f_{\sigma}$.

The condition preventing $f_{\sigma}$ from visiting "too many" strata can be formulated as follows: we should restrict ourselves to those $\sigma$ that map morphisms of $F(\Delta^n)$ to trajectories that are not broken more times than imposed by the functoriality of $\sigma$.

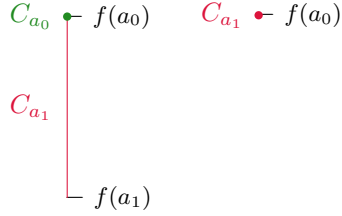
\begin{figure}[H]
    \centering
    \begin{tabular}{cc}
    \begin{tikzpicture}[scale=1.2]
        \draw node at (0.44,1) {$-$ \small{$f(a_0)$}};
        \draw[green] node at (-0.4,1) {$C_{a_0}$};
        \draw node at (0.45,-1) {$-$ \small{$f(a_1)$}};
        \draw[red] node at (-0.4,0) {$C_{a_1}$};
        \draw[red] (0,-1)--(0,1);
        \fill[green] (0,1) circle (1.3pt);
    \end{tikzpicture}

    &

    \begin{tikzpicture}[scale=1.2]
        \draw node at (0.44,1) {$-$ \small{$f(a_0)$}};
        \draw[red] node at (-0.4,1) {$C_{a_1}$};
        \fill[red] (0,1) circle (1.3pt);
        \draw[white] node at (0.48,-1) {$-$ $f(a_1)$};
    \end{tikzpicture}
    \end{tabular}
    \caption{The stratification of the cube $I^0 \times [f(a_1),f(a_0)]$ arising from an increasing sequence of two critical points $a_0 \le a_1$. To simplify the notation, we denote the corresponding cube by $C$ on each picture. In the case when the critical points are different (on the left), the condition that $f_{\sigma}$ preserves stratification is equivalent to the condition that the trajectory from $a_0$ to $a_1$ given by $\sigma$ is unbroken. On the right: the case when $a_0 = a_1$.}
    \label{figure: stratifications of the 1-cube}

\end{figure}

\begin{figure}[H]
    \centering
    $$\begin{array}{ccc}
         \begin{tikzpicture}[scale=1]
            \draw node at (2.64,-1) {$-$ \small{$f(a_3)$}};
            \draw node at (2.64,-0.5) {$-$ \small{$f(a_2)$}};
            \draw node at (2.64,0.9) {$-$ \small{$f(a_1)$}};
            \draw node at (2.64,2) {$-$ \small{$f(a_0)$}};
            \draw node at (-1.2,-1.5) {\small{$|$}};
            \draw node at (-1.2,-1.8) {\small{$(0,0)$}};
            \draw node at (-0.2,-1) {\small{$|$}};
            \draw node at (-0.2,-1.3) {\small{$(0,1)$}};
            \draw node at (1.2,-1.5) {\small{$|$}};
            \draw node at (1.2,-1.8) {\small{$(1,0)$}};
            \draw node at (2.2,-1) {\small{$|$}};
            \draw node at (2.2,-1.3) {\small{$(1,1)$}};
            \draw node at (-1.5,0) {$C_{a_3}$};
            \draw[green] node at (1,2.2) {$C_{a_0}$};
            \draw[red] node at (2.5,1.45) {$C_{a_1}$};
            \draw[blue] node at (2.5,0.2) {$C_{a_2}$};
            \draw node at (0.5,-2.3) {\textcolor{green}{$a_0$} < \textcolor{red}{$a_1$} < \textcolor{blue}{$a_2$} < $a_3$};
         
             \draw (-1.2,1.5)--(-1.2,-1.5)--(1.2,-1.5)--(1.2,0.4);
             \draw[dashed] (-1.2,-1.5)--(-0.2,-1);
             \draw[dashed] (-0.2,-1)--(2.2,-1);
             \draw[dashed] (-0.2,-1)--(-0.2,-0.5);
             \draw (2.2,-1)--(2.2,-0.5);
             \draw (1.2,-1.5)--(2.2,-1);

             \draw[blue] (2.2,-0.5)--(2.2,0.9);
             \draw[dashed, blue] (2.2,-0.5)--(-0.2,-0.5)--(-0.2,2);
             \draw[fill, blue, opacity=0.4] (-0.2,2)--(-0.2,-0.5)--(2.2,-0.5)--(2.2,2)--cycle;

             \draw[green] (-1.2,1.5)--(1.2,1.5)--(2.2,2)--(-0.2,2)--cycle;
             \draw[fill, green, opacity=0.4] (-1.2,1.5)--(1.2,1.5)--(2.2,2)--(-0.2,2)--cycle;

             \draw[red] (1.2,1.5)--(1.2,0.4)--(2.2,0.9)--(2.2,2);
             \draw[fill, red, opacity=0.4] (1.2,1.5)--(1.2,0.4)--(2.2,0.9)--(2.2,2)--cycle;

             \draw[fill, opacity = 0.1] (-1.2,-1.5)--(1.2,-1.5)--(2.2,-1)--(2.2,2)--(-0.2,2)--(-1.2,1.5)--cycle;
         \end{tikzpicture}
        &
         \begin{tikzpicture}[scale=1]
            \draw node at (2.64,-1) {$-$ \small{$f(a_3)$}};
            \draw node at (2.64,0.9) {$-$ \small{$f(a_1)$}};
            \draw node at (2.64,2) {$-$ \small{$f(a_0)$}};
            \draw node at (-1.2,-1.5) {\small{$|$}};
            \draw node at (-1.2,-1.8) {\small{$(0,0)$}};
            \draw node at (-0.2,-1) {\small{$|$}};
            \draw node at (-0.2,-1.3) {\small{$(0,1)$}};
            \draw node at (1.2,-1.5) {\small{$|$}};
            \draw node at (1.2,-1.8) {\small{$(1,0)$}};
            \draw node at (2.2,-1) {\small{$|$}};
            \draw node at (2.2,-1.3) {\small{$(1,1)$}};
            \draw node at (-1.5,0) {$C_{a_3}$};
            \draw[green] node at (1,2.2) {$C_{a_0}$};
            \draw[blue] node at (2.5,1.45) {$C_{a_2}$};
            \draw node at (0.5,-2.3) {\textcolor{green}{$a_0$} < \textcolor{blue}{$a_1$} = \textcolor{blue}{$a_2$} < $a_3$};
         
             \draw (-1.2,1.5)--(-1.2,-1.5)--(1.2,-1.5)--(1.2,0.4);
             \draw[dashed] (-1.2,-1.5)--(-0.2,-1);
             \draw[dashed] (-0.2,-1)--(2.2,-1);
             \draw[dashed] (-0.2,-1)--(-0.2,0.9);
             \draw (2.2,-1)--(2.2,0.9);
             \draw (1.2,-1.5)--(2.2,-1);

             \draw[green] (-1.2,1.5)--(1.2,1.5)--(2.2,2)--(-0.2,2)--cycle;
             \draw[fill, green, opacity=0.4] (-1.2,1.5)--(1.2,1.5)--(2.2,2)--(-0.2,2)--cycle;

             \draw[blue] (1.2,1.5)--(1.2,0.4)--(2.2,0.9)--(2.2,2);
             \draw[fill, blue, opacity=0.4] (1.2,1.5)--(1.2,0.4)--(2.2,0.9)--(2.2,2)--cycle;
             \draw[dashed, blue] (2.2,0.9)--(-0.2,0.9)--(-0.2,2);
             \draw[fill, blue, opacity=0.4] (-0.2,2)--(-0.2,0.9)--(2.2,0.9)--(2.2,2)--cycle;

             \draw[fill, opacity = 0.1] (-1.2,-1.5)--(1.2,-1.5)--(2.2,-1)--(2.2,2)--(-0.2,2)--(-1.2,1.5)--cycle;
         \end{tikzpicture}
        &
         \begin{tikzpicture}[scale=1]
            \draw node at (2.64,-1) {$-$ \small{$f(a_3)$}};
            \draw node at (2.64,2) {$-$ \small{$f(a_0)$}};
            \draw node at (-1.2,-1.5) {\small{$|$}};
            \draw node at (-1.2,-1.8) {\small{$(0,0)$}};
            \draw node at (-0.2,-1) {\small{$|$}};
            \draw node at (-0.2,-1.3) {\small{$(0,1)$}};
            \draw node at (1.2,-1.5) {\small{$|$}};
            \draw node at (1.2,-1.8) {\small{$(1,0)$}};
            \draw node at (2.2,-1) {\small{$|$}};
            \draw node at (2.2,-1.3) {\small{$(1,1)$}};
            \draw node at (-1.5,0) {$C_{a_3}$};
            \draw[green] node at (1,2.2) {$C_{a_0}$};
            \draw node at (0.5,-2.3) {\textcolor{green}{$a_0$} < $a_1$ = $a_2$ = $a_3$};
         
             \draw (-1.2,1.5)--(-1.2,-1.5)--(1.2,-1.5)--(1.2,1.5);
             \draw[dashed] (-1.2,-1.5)--(-0.2,-1);
             \draw[dashed] (-0.2,-1)--(2.2,-1);
             \draw[dashed] (-0.2,-1)--(-0.2,2);
             \draw (2.2,-1)--(2.2,2);
             \draw (1.2,-1.5)--(2.2,-1);

             \draw[green] (-1.2,1.5)--(1.2,1.5)--(2.2,2)--(-0.2,2)--cycle;
             \draw[fill, green, opacity=0.4] (-1.2,1.5)--(1.2,1.5)--(2.2,2)--(-0.2,2)--cycle;

             \draw[fill, opacity = 0.1] (-1.2,-1.5)--(1.2,-1.5)--(2.2,-1)--(2.2,2)--(-0.2,2)--(-1.2,1.5)--cycle;
         \end{tikzpicture}
    \end{array}$$
    \vspace{-0.5cm}
    \caption{Some stratifications of the cube $I^2 \times [f(a_3),f(a_0)]$ arising from an increasing sequence of four critical points $a_0 \le a_1 \le a_2 \le a_3$. To simplify the notation, we denote the corresponding cube by $C$ on each picture. On the left: the case when the critical points are distinct. For an arbitrary $\sigma \in \Nl(\Ml)_3$ such that the $a_i$'s are distinct, the trajectories in $\sigma(0,3)(\{1\} \times I)$ are broken at least at $a_1$, the trajectories in $\sigma(0,3)(I \times \{1\})$ are broken at least at $a_2$ and the trajectory $\sigma(0,3)((1,1))$ is broken at least at $a_1$ and $a_2$. The condition that $f_{\sigma}$ be compatible with stratifications is equivalent to the condition that these trajectories are not broken elsewhere, and all the other trajectories in the image of $\sigma(0,3)$ are unbroken. An analogous description follows from the two other pictures in the two cases $a_0 < a_1 = a_2 < a_3$ and $a_0 < a_1 = a_2 = a_3$.}
    \label{figure: some stratifications of 3-cubes}
\end{figure}
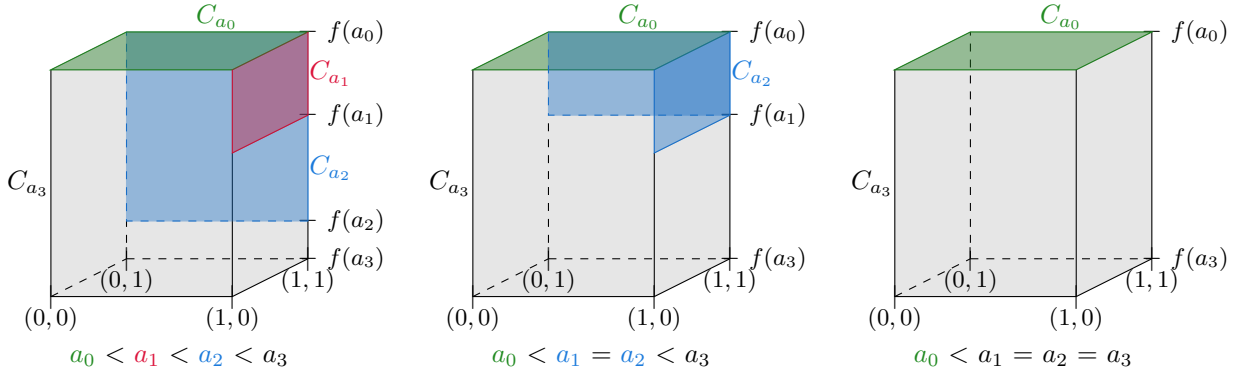

This condition is introduced in section \ref{section: cubical stratified geometric realization}. We formulate it here in a different, but equivalent way. For all finite increasing sequence of critical points $a_0 \le \hdots \le a_n \in A$, we define a \emph{stratification} of the cube $I^{n-1} \times [f(a_n),f(a_0)]$ by $\{a_0,\hdots,a_n\}$ \footnote{What we really need in order to define this stratification of $I^{n-1} \times [f(a_n),f(a_0)]$ is that the sequence satisfies $f(a_0) \ge f(a_1) \ge \hdots \ge f(a_n)$, which is in particular implied by the condition that it is increasing with respect to $\le$.} such that the restrictive condition we want on $\sigma \in \Nl(\Ml)_n$ is implemented by the requirement that $f_{\sigma}$ is \emph{compatible} with the respective stratifications of $I^{n-1} \times [f(a_n),f(a_0)]$ and $X$. By "compatible", we mean that the $a_i$-stratum of the cube is mapped to the $a_i$-stratum of $X$ for every $i$. Some of these stratifications of cubes for small values of $n$ are presented in figures \ref{figure: stratifications of the 1-cube}, \ref{figure: stratifications of the 2-cube} and \ref{figure: some stratifications of 3-cubes}.

\begin{figure}[H]
    \centering
    $$\begin{array}{c@{\!}c@{\!}c@{\!}c}
    \begin{tikzpicture}[scale=1.2]
        \draw node at (1.44,1.2) {$-$ \small{$f(a_0)$}};
        \draw node at (1.44,0.3) {$-$ \small{$f(a_1)$}};
        \draw node at (1.44,-1.2) {$-$ \small{$f(a_2)$}};
        \draw node at (-1,-1.2) {\small{|}};
        \draw node at (1,-1.2) {\small{|}};
        \draw node at (-1,-1.5) {$0$};
        \draw node at (1,-1.5) {$1$};
        \draw[green] node at (0,1.4) {$C_{a_0}$};
        \draw[red] node at (1.4,0.72) {$C_{a_1}$};
        \draw[blue] node at (-1.3,0) {$C_{a_2}$};
        \draw node at (0,-2) {\textcolor{green}{$a_0$} < \textcolor{red}{$a_1$} < \textcolor{blue}{$a_2$}};

        \draw[fill, blue, opacity=0.5] (-1,1.2)--(1,1.2)--(1,-1.2)--(-1,-1.2)--cycle;
    
        \draw[blue] (1,0.3)--(1,-1.2)--(-1,-1.2)--(-1,1.2);
        \draw[red, line width=1.4pt] (1,1.2)--(1,0.3);
        \draw[green, line width=1.4pt] (-1,1.2)--(1,1.2);
    \end{tikzpicture}
    &
    \begin{tikzpicture}[scale=1.2]
        \draw node at (1.44,1.2) {$-$ \small{$f(a_0)$}};
        \draw node at (1.44,0.3) {$-$ \small{$f(a_1)$}};
        \draw node at (-1,0.3) {\small{|}};
        \draw node at (1,0.3) {\small{|}};
        \draw node at (-1,0) {$0$};
        \draw node at (1,0) {$1$};
        \draw[white] node at (-1,-1.5) {$0$};
        \draw[green] node at (0,1.4) {$C_{a_0}$};
        \draw[blue] node at (1.4,0.72) {$C_{a_2}$};
        \draw node at (0,-2) {\textcolor{green}{$a_0$} < \textcolor{blue}{$a_1$} = \textcolor{blue}{$a_2$}};

        \draw[fill, blue, opacity=0.5] (-1,1.2)--(1,1.2)--(1,0.3)--(-1,0.3)--cycle;
    
        \draw[blue] (1,1.2)--(1,0.3)--(-1,0.3)--(-1,1.2);
        \draw[green, line width=1.4pt] (-1,1.2)--(1,1.2);
    \end{tikzpicture}
    &
    \begin{tikzpicture}[scale=1.2]
        \draw node at (1.44,1.2) {$-$ \small{$f(a_0)$}};
        \draw node at (1.44,-1.2) {$-$ \small{$f(a_2)$}};
        \draw node at (-1,-1.2) {\small{|}};
        \draw node at (1,-1.2) {\small{|}};
        \draw node at (-1,-1.5) {$0$};
        \draw node at (1,-1.5) {$1$};
        \draw[red] node at (0,1.4) {$C_{a_1}$};
        \draw[blue] node at (-1.3,0) {$C_{a_2}$};
        \draw node at (0,-2) {\textcolor{red}{$a_0$} = \textcolor{red}{$a_1$} < \textcolor{blue}{$a_2$}};

        \draw[fill, blue, opacity=0.5] (-1,1.2)--(1,1.2)--(1,-1.2)--(-1,-1.2)--cycle;
    
        \draw[blue] (1,1.2)--(1,-1.2)--(-1,-1.2)--(-1,1.2);
        \draw[red, line width=1.4pt] (-1,1.2)--(1,1.2);
    \end{tikzpicture}
     &
    \begin{tikzpicture}[scale=1.2]
        \draw node at (1.44,1.2) {$-$ \small{$f(a_0)$}};
        \draw node at (-1,1.2) {\small{|}};
        \draw node at (1,1.2) {\small{|}};
        \draw node at (-1,0.9) {$0$};
        \draw node at (1,0.9) {$1$};
        \draw[white] node at (-1,-1.5) {$0$};
        \draw[blue] node at (0,1.4) {$C_{a_2}$};
        \draw node at (0,-2) {\textcolor{blue}{$a_0$} = \textcolor{blue}{$a_1$} = \textcolor{blue}{$a_2$}};
        
        \draw[blue, line width=1.4pt] (-1,1.2)--(1,1.2);
    \end{tikzpicture}
    \end{array}$$
    \vspace{-0.5cm}
    \caption{The possible stratifications of the cube $I \times [f(a_2),f(a_0)]$ arising from an increasing sequence of three critical points $a_0 \le a_1 \le a_2$. To simplify the notation, we denote the corresponding cube by $C$ on each picture. On the left: the case when the critical points are distinct. In that case, by functoriality, the trajectory $\sigma(0,2)(1)$ is the composition of $\sigma(0,1)(*)$ and $\sigma(1,2)(*)$, and is therefore necessarily broken at $a_1$. The condition that the map $f_{\sigma}$ is compatible with the stratifications is equivalent to the conditions that $\sigma(0,1)(*)$ and $\sigma(1,2)(*)$ are unbroken, and $\sigma(0,2)(t)$ is unbroken for every $0 \le t < 1$. In the cases when $a_0 \neq a_1 = a_2$ and $a_0 = a_1 \neq a_2$, the condition that $f_{\sigma}$ is compatible with the stratifications is equivalent to the condition that $\sigma(0,2)(t)$ is unbroken for every $0 \le t \le 1$. In the case when $a_0 = a_1 = a_2$, $\sigma(0,2)$ maps every morphism to the constant trajectory at $a_0$, and requesting compatibility with stratifications for $f_{\sigma}$ does not impose any constraints on $\sigma$.}
    \label{figure: stratifications of the 2-cube}
\end{figure}
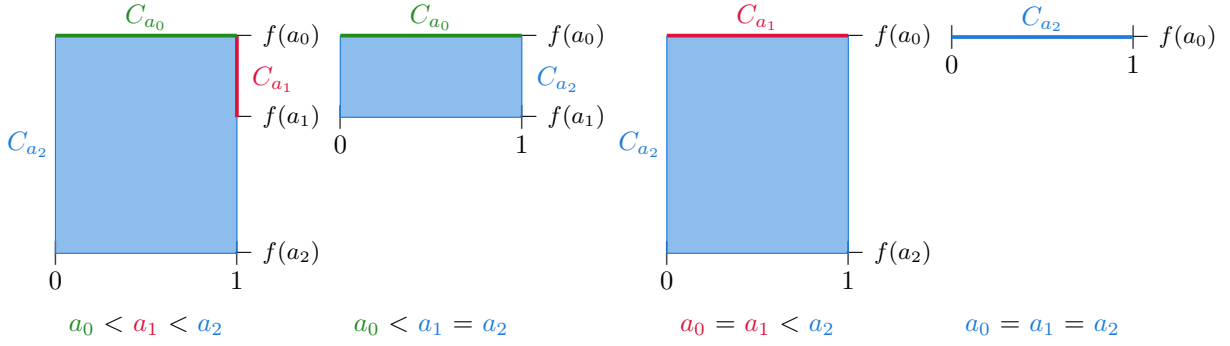

When $f_{\sigma}$ preserves the stratifications, we say that $\sigma$ is \emph{unbroken}. Unbroken simplices form a simplicial subset $S_{\Ml} \subset \Nl(\Ml)$, which we call the \emph{flow coherent nerve of} $\Ml$. We prove in section \ref{section: flow coherent nerve is an infinity category} that this is an $\infty$-category. This is done using structures of smooth manifolds with corners on the morphism spaces of $\Ml$, which are presented in section \ref{section: Morse-Smale pairs and flow categories}. We prove in section \ref{section: equivalence between flow and homotopy coherent nerve} that the functor $S_{\Ml} \rightarrow \Nl(\Ml)$ defined by the morphism of inclusion, is an equivalence of $\infty$-categories. For this, we use the characterization of equivalences of $\infty$-categories as those functors that are essentially surjective and fully faithful. In turns out that the morphism spaces of $S_{\Ml}$ can be explicitly described as the spaces of unbroken trajectories connecting critical points (see proposition \ref{proposition: right morphism spaces of flow coherent nerve}). Now, for every pair of critical points, the inclusion of the space of unbroken trajectories between them, into the space of all trajectories between them, is a homotopy equivalence (see corollary \ref{Corollary: space of broken trajectories equivalent to space of trajectories}).

Until now, we have imposed restrictions on the simplices of $\Nl(\Ml)$ that we consider, with the perspective to associate to them simplices of $\Sing_A(X)$, that is, stratified singular simplices of $X$. When the $n$-simplex $\sigma$ of $\Nl(\Ml)$ is unbroken, the associated map $f_{\sigma}$ visits the correct strata. Moreover, unbroken simplices form a simplicial subset $S_{\Ml} \subset \Nl(\Ml)$ and this inclusion is an equivalence of $\infty$-categories: we therefore aim to compare the simplicial sets $S_{\Ml}$ and $\Sing_A(X)$ from now on.

At this point, the following difficulty arises: when comparing the stratifications that we have introduced on cubes (figures \ref{figure: stratifications of the 1-cube}, \ref{figure: stratifications of the 2-cube} and \ref{figure: some stratifications of 3-cubes}) with the stratifications on the standard simplices (figures \ref{figure: stratification of the standard 1-simplex}, \ref{figure: stratification of the standard 2-simplex} and \ref{figure: some stratifications of the 3 simplex}), we see that, in most cases, these are not homeomorphic. Therefore, even with this restriction imposed, $f_{\sigma}$ does not yield a stratified singular $n$-simplex of $X$ on the nose.

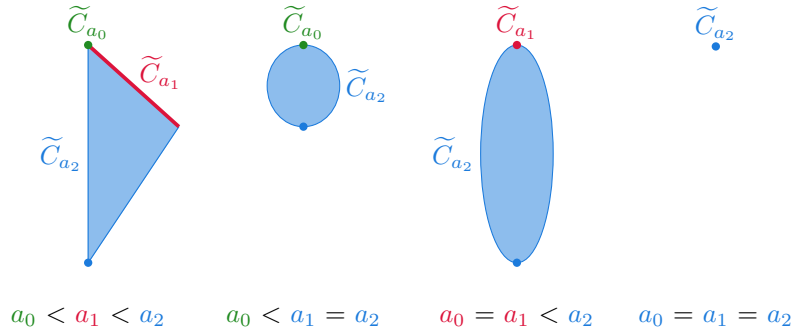
\begin{figure}[H]
    \centering
    $$\begin{array}{cccc}
    \begin{tikzpicture}[scale=1.2]
        \draw[green] node at (0,1.45) {$\widetilde{C}_{a_0}$};
        \draw[red] node at (0.8,0.9) {$\widetilde{C}_{a_1}$};
        \draw[blue] node at (-0.3,0) {$\widetilde{C}_{a_2}$};
        \draw node at (0,-1.8) {\textcolor{green}{$a_0$} < \textcolor{red}{$a_1$} < \textcolor{blue}{$a_2$}};

        \draw[fill, blue, opacity=0.5] (0,1.2)--(1,0.3)--(0,-1.2)--cycle;
        \draw[blue] (1,0.3)--(0,-1.2)--(0,1.2);
        \fill[blue] (0,-1.2) circle (1.3pt);
    
        \draw[red, line width=1.4pt] (0,1.2)--(1,0.3);
        
        \fill[green] (0,1.2) circle (1.3pt);
    \end{tikzpicture}
    &
    \begin{tikzpicture}[scale=1.2]
        \draw[green] node at (0,1.45) {$\widetilde{C}_{a_0}$};
        \draw[blue] node at (0.7,0.75) {$\widetilde{C}_{a_2}$};
        \draw node at (0,-1.8) {\textcolor{green}{$a_0$} < \textcolor{blue}{$a_1$} = \textcolor{blue}{$a_2$}};
    
        \fill[blue, opacity=0.5] (0,0.75) ellipse (0.4 and 0.45);
        \draw[blue] (0,0.75) ellipse (0.4 and 0.45);
        \fill[blue] (0,0.3) circle (1.3pt);
        
        \fill[green] (0,1.2) circle (1.3pt);
        
        \draw node at (0,-1.2) {$ $};
        \draw node at (-0.7,0) {$ $};
        
    \end{tikzpicture}
    &
    \begin{tikzpicture}[scale=1.2]
        \draw[red] node at (0,1.45) {$\widetilde{C}_{a_1}$};
        \draw[blue] node at (-0.7,0) {$\widetilde{C}_{a_2}$};
        \draw node at (0,-1.8) {\textcolor{red}{$a_0$} = \textcolor{red}{$a_1$} < \textcolor{blue}{$a_2$}};
        
        \draw[blue] (0,0) ellipse (0.4 and 1.2);
        \fill[blue, opacity=0.5] (0,0) ellipse (0.4 and 1.2);
        \fill[blue] (0,-1.2) circle (1.3pt);
        
        \fill[red] (0,1.2) circle (1.3pt);
    \end{tikzpicture}
    &
    \begin{tikzpicture}[scale=1.2]
        \draw[blue] node at (0,1.45) {$\widetilde{C}_{a_2}$};
        \draw node at (0,-1.8) {\textcolor{blue}{$a_0$} = \textcolor{blue}{$a_1$} = \textcolor{blue}{$a_2$}};
        
        \fill[blue] (0,1.2) circle (1.3pt);
        
        \draw node at (0,-1.2) {$ $};
        \draw node at (-0.3,1.2) {$ $};
        \draw node at (0.3,-1.2) {$ $};
    \end{tikzpicture}
    \end{array}$$
    \caption{The possible quotients of the cube arising from a sequence of three critical points. We denote the quotient by $\widetilde{C}$. Note that the stratification on the cube descends to the quotient.}
    \label{figure: quotients of stratified 2-cubes}
\end{figure}

To solve this problem, we further notice that, by construction, the singular cubes of $X$ arising from simplices of the homotopy coherent nerve of the flow category, are \emph{constant} on some subsets of the cube. As a basic illustration of this, we observe that all the trajectories in the image of $\sigma(0,n)$ start at $a_0$ and end at $a_n$. Consequently, $f_{\sigma}$ is constant equal to $a_0$ on $I^{n-1} \times \{f(a_0)\}$, and is constant equal to $a_n$ on $I^{n-1} \times \{f(a_n)\}$. The map $f_{\sigma}$ thus induces a map out of the cube $I^{n-1} \times [f(a_n),f(a_0)]$ with each of the two subsets $I^{n-1} \times \{f(a_n)\}$ and $I^{n-1} \times \{f(a_0)\}$ collapsed to a point. In the four cases of figure \ref{figure: stratifications of the 2-cube}, corresponding to the case $n=2$, these are described in figure \ref{figure: quotients of stratified 2-cubes}.

Now, when we compare figure \ref{figure: stratification of the standard 2-simplex} with figure \ref{figure: stratifications of the 2-cube}, it appears that we have partially resolved the problem of the difference between stratified $2$-cubes and stratified $2$-simplices. In the case of a sequence of three distinct critical points, the quotiented stratified $2$-cube is clearly stratified homeomorphic to the corresponding stratified $2$-simplex. This is however not the case for all sequences. For example, when it comes to constant sequences, the quotiented $2$-cube is a point, while the $2$-simplex is not. We will explain below how this problem is solved.

More generally, for every nonnegative integer $n$ and every increasing sequence of $n+1$ critical points $a_0 \le \hdots \le a_n$, we introduce in section \ref{section: cubical stratified geometric realization} an \emph{equivalence relation} on $I^{n-1} \times [f(a_n),f(a_0)]$ with which both $f_{\sigma}$ and the stratification are compatible. The compatibility with $f_{\sigma}$ is valid whether $\sigma$ is unbroken or not. In practice, section \ref{section: cubical stratified geometric realization} is organized as follows: we first define in section \ref{section: non stratified functor C A} the equivalence relations on cubes associated with finite increasing sequences of critical points, and then we introduce the stratifications in section \ref{section: lift to stratified spaces}.

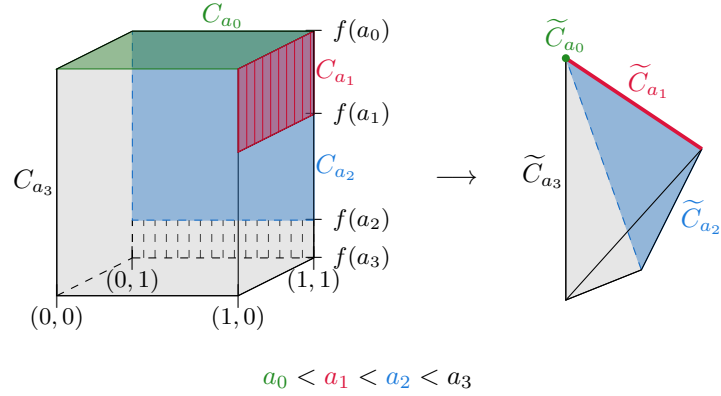
\begin{figure}[H]
    \centering
    $$\begin{array}{ccc}
        \begin{tikzpicture}[scale=1]
            \draw node at (2.64,-1) {$-$ \small{$f(a_3)$}};
            \draw node at (2.64,-0.5) {$-$ \small{$f(a_2)$}};
            \draw node at (2.64,0.9) {$-$ \small{$f(a_1)$}};
            \draw node at (2.64,2) {$-$ \small{$f(a_0)$}};
            \draw node at (-1.2,-1.5) {\small{$|$}};
            \draw node at (-1.2,-1.8) {\small{$(0,0)$}};
            \draw node at (-0.2,-1) {\small{$|$}};
            \draw node at (-0.2,-1.3) {\small{$(0,1)$}};
            \draw node at (1.2,-1.5) {\small{$|$}};
            \draw node at (1.2,-1.8) {\small{$(1,0)$}};
            \draw node at (2.2,-1) {\small{$|$}};
            \draw node at (2.2,-1.3) {\small{$(1,1)$}};
            \draw node at (-1.5,0) {$C_{a_3}$};
            \draw[green] node at (1,2.2) {$C_{a_0}$};
            \draw[red] node at (2.5,1.45) {$C_{a_1}$};
            \draw[blue] node at (2.5,0.2) {$C_{a_2}$};
            
             \draw (-1.2,1.5)--(-1.2,-1.5)--(1.2,-1.5)--(1.2,0.4);
             \draw[dashed] (-1.2,-1.5)--(-0.2,-1);
             \draw[dashed] (-0.2,-1)--(2.2,-1);
             \draw[dashed] (-0.2,-1)--(-0.2,-0.5);
             \draw (2.2,-1)--(2.2,-0.5);
             \draw (1.2,-1.5)--(2.2,-1);
             \foreach \x in {-0.2,-0.05,...,2.2} {
            \draw[dashed, opacity = 0.7] (\x,-1) -- (\x,-0.5);
            }
            
             \draw[blue] (2.2,-0.5)--(2.2,0.9);
             \draw[dashed, blue] (2.2,-0.5)--(-0.2,-0.5)--(-0.2,2);
             \draw[fill, blue, opacity=0.4] (-0.2,2)--(-0.2,-0.5)--(2.2,-0.5)--(2.2,2)--cycle;

             \draw[green] (-1.2,1.5)--(1.2,1.5)--(2.2,2)--(-0.2,2)--cycle;
            \fill[green, fill opacity=0.4]
            (-1.2,1.5)--(1.2,1.5)--(2.2,2)--(-0.2,2)--cycle;

             \draw[red] (1.2,1.5)--(1.2,0.4)--(2.2,0.9)--(2.2,2);
             \draw[pattern=vertical lines, pattern color=red] (1.2,1.5)--(1.2,0.4)--(2.2,0.9)--(2.2,2)--cycle;
             \draw[fill, red, fill opacity=0.4] (1.2,1.5)--(1.2,0.4)--(2.2,0.9)--(2.2,2)--cycle;

             \draw[fill, fill opacity = 0.1] (-1.2,-1.5)--(1.2,-1.5)--(2.2,-1)--(2.2,2)--(-0.2,2)--(-1.2,1.5)--cycle;
    \end{tikzpicture}
    & 
    \begin{tikzpicture}
            \draw node at (0,0.2) {$\longrightarrow$};
            \draw node at (0,2.2) {$ $};
            \draw node at (0,-1.8) {$ $};
    \end{tikzpicture}
    &
    \begin{tikzpicture}[scale=1]
            \draw[green] node at (0,1.9) {$\widetilde{C}_{a_0}$};
            \draw[red] node at (1.1,1.3) {$\widetilde{C}_{a_1}$};
            \draw[blue] node at (1.8,-0.5) {$\widetilde{C}_{a_2}$};
            \draw node at (-0.3,0.1) {$\widetilde{C}_{a_3}$};
            
            \draw node at (0,2) {$ $};
            \draw node at (0,-2) {$ $};
        
             \draw[blue] (1,-1.2)--(1.8,0.4);
             \draw[dashed, blue] (1,-1.2)--(0,1.6);
             \draw[fill, blue, opacity=0.4] (1,-1.2)--(0,1.6)--(1.8,0.4)--cycle;

             \draw[fill, fill opacity = 0.1] (0,-1.6)--(1,-1.2)--(1.8,0.4)--(0,1.6)--cycle;
              
             \draw (0,1.6)--(0,-1.6)--(1,-1.2);
             \draw (0,-1.6)--(1.8,0.4);

             \draw[red, line width=1.4pt] (1.8,0.4)--(0,1.6);
             
            \fill[green] (0,1.6) circle (1.5pt);
    \end{tikzpicture}
    \end{array}$$
    
    \vspace{-0.7cm}
    
    $$\textcolor{green}{a_0} < \textcolor{red}{a_1} < \textcolor{blue}{a_2} < a_3$$
    
    \caption{The quotient of the cube $I^2 \times [f(a_3),f(a_0)]$ in the case of a strictly increasing sequence of four critical points. We denote it by $\widetilde{C}$. The map $f_{\sigma}$ is constant equal to $a_0$ on $I^2 \times \{f(a_0)\}$ and equal to $a_3$ on $I^2 \times \{f(a_3)\}$. In addition, on the vertical segments $(\{(t,1)\} \times [f(a_3),f(a_2)])_{0 \le t \le 1}$ depicted in the figure, $f_{\sigma}$ corresponds to the evaluation of the same trajectory for every $t$. The same holds on the vertical segments $(\{(1,t)\} \times [f(a_1),f(a_0)])_{0 \le t \le 1}$ also depicted in the figure. Consequently, the map $f_{\sigma}$ is compatible with the equivalence relation that collapses to a point each of the two subsets $I^2 \times \{f(a_3)\}$ and $I^2 \times \{f(a_0)\}$, as well as each of the horizontal segments $I \times \{1\} \times \{s\}$ for $f(a_3) \le s \le f(a_2)$, and $\{1\} \times I \times \{s\}$ for $f(a_1) \le s \le f(a_0)$. By comparing with figure \ref{figure: some stratifications of the 3 simplex}, we see that the stratified quotient is stratified homeomorphic to the stratified $3$-simplex corresponding to the given sequence of critical points.}
    \label{figure: quotient of cube for a0<a1<a2<a3}
\end{figure}

\begin{figure}[H]
    \centering
    $$\begin{array}{ccc}
         \begin{tikzpicture}[scale=1]
            \draw node at (2.64,-1) {$-$ \small{$f(a_3)$}};
            \draw node at (2.64,2) {$-$ \small{$f(a_0)$}};
            \draw node at (-1.2,-1.5) {\small{$|$}};
            \draw node at (-1.2,-1.8) {\small{$(0,0)$}};
            \draw node at (-0.2,-1) {\small{$|$}};
            \draw node at (-0.2,-1.3) {\small{$(0,1)$}};
            \draw node at (1.2,-1.5) {\small{$|$}};
            \draw node at (1.2,-1.8) {\small{$(1,0)$}};
            \draw node at (2.2,-1) {\small{$|$}};
            \draw node at (2.2,-1.3) {\small{$(1,1)$}};
            \draw node at (-1.5,0) {$C_{a_3}$};
            \draw[green] node at (1,2.2) {$C_{a_0}$};
         
             \draw (-1.2,1.5)--(-1.2,-1.5)--(1.2,-1.5)--(1.2,1.5);
             \draw[dashed] (-1.2,-1.5)--(-0.2,-1);
             \draw[dashed] (-0.2,-1)--(2.2,-1);
             \draw[dashed] (-0.2,-1)--(-0.2,2);
             \draw (2.2,-1)--(2.2,2);
             \draw (1.2,-1.5)--(2.2,-1);
             \foreach \x in {-0.2,-0.05,...,2.2} {
            \draw[dashed, opacity = 0.7] (\x,-1) -- (\x,2);
            }

             \draw[green] (-1.2,1.5)--(1.2,1.5)--(2.2,2)--(-0.2,2)--cycle;
             \draw[fill, green, opacity=0.4] (-1.2,1.5)--(1.2,1.5)--(2.2,2)--(-0.2,2)--cycle;

             \draw[fill, opacity = 0.1] (-1.2,-1.5)--(1.2,-1.5)--(2.2,-1)--(2.2,2)--(-0.2,2)--(-1.2,1.5)--cycle;
         \end{tikzpicture}
         & 
        \begin{tikzpicture}
                \draw node at (0,0.2) {$\longrightarrow$};
                \draw node at (0,2.2) {$ $};
                \draw node at (0,-1.8) {$ $};
        \end{tikzpicture}
        &
        \begin{tikzpicture}
                \draw[green] node at (0,1.9) {$\widetilde{C}_{a_0}$};
                \draw node at (-0.9,0.1) {$\widetilde{C}_{a_3}$};
                
                \draw node at (0,2) {$ $};
                \draw node at (0,-2) {$ $};
                
                \draw (0,0) ellipse (0.6 and 1.6);
                \draw[fill, opacity=0.1] (0,0) ellipse (0.6 and 1.6);
                \draw[dashed] (0,1.6) arc [start angle=90, end angle=-90, x radius=0.3, y radius=1.6];
                \draw[dashed] (-0.6,0) arc [start angle=180, end angle=0, x radius=0.6, y radius=0.1];
                \draw (-0.6,0) arc [start angle=-180, end angle=0, x radius=0.6, y radius=0.1];
                
                \fill[green] (0,1.6) circle (1.5pt);
        \end{tikzpicture}
    \end{array}$$
    
    \vspace{-0.7cm}
    
    $$\textcolor{green}{a_0} < a_1 = a_2 = a_3$$
    
    \caption{The quotient of the cube $I^2 \times [f(a_3),f(a_0)]$ in the case of a sequence of four critical points satisfying $a_0 < a_1 = a_2 = a_3$. We denote this quotient by $\widetilde{C}$. The map $f_{\sigma}$ is constant equal to $a_0$ on $I^2 \times \{f(a_0)\}$ and equal to $a_3$ on $I^2 \times \{f(a_3)\}$. In addition, on the vertical segments $(\{(t,1)\} \times [f(a_3),f(a_0)])_{0 \le t \le 1}$ depicted in the figure, $f_{\sigma}$ corresponds to the evaluation of the same trajectory for every $t$. Consequently, the map $f_{\sigma}$ is compatible with the equivalence relation that collapses to a point each of the two subsets $I^2 \times \{f(a_3)\}$ and $I^2 \times \{f(a_0)\}$, as well as each of the horizontal segments $I \times \{1\} \times \{s\}$ for $f(a_3) \le s \le f(a_0)$. We will explain below how to compare the resulting stratified quotiented cube with the corresponding stratified $3$-simplex.}
    \label{figure: quotient of cube for a0<a1=a2=a3}
\end{figure}
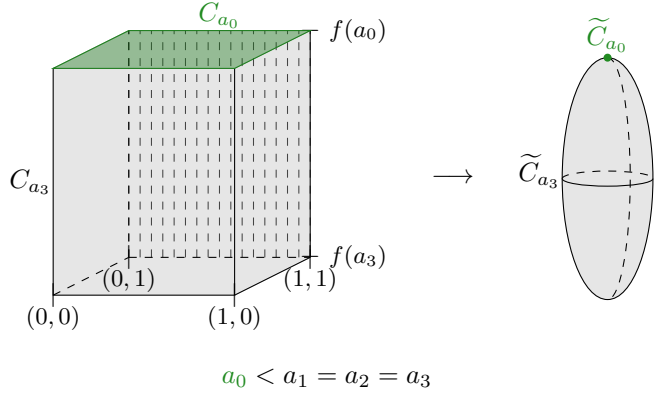

In figures \ref{figure: quotient of cube for a0<a1<a2<a3}, \ref{figure: quotient of cube for a0<a1=a2<a3} and \ref{figure: quotient of cube for a0<a1=a2=a3} we describe this equivalence relation in the three cases from figure \ref{figure: some stratifications of 3-cubes}, as well as the resulting quotients.

\begin{figure}[H]
    \centering
    $$\begin{array}{ccc}
         \begin{tikzpicture}[scale=1]
            \draw node at (2.64,-1) {$-$ \small{$f(a_3)$}};
            \draw node at (2.64,0.9) {$-$ \small{$f(a_1)$}};
            \draw node at (2.64,2) {$-$ \small{$f(a_0)$}};
            \draw node at (-1.2,-1.5) {\small{$|$}};
            \draw node at (-1.2,-1.8) {\small{$(0,0)$}};
            \draw node at (-0.2,-1) {\small{$|$}};
            \draw node at (-0.2,-1.3) {\small{$(0,1)$}};
            \draw node at (1.2,-1.5) {\small{$|$}};
            \draw node at (1.2,-1.8) {\small{$(1,0)$}};
            \draw node at (2.2,-1) {\small{$|$}};
            \draw node at (2.2,-1.3) {\small{$(1,1)$}};
            \draw node at (-1.5,0) {$C_{a_3}$};
            \draw[green] node at (1,2.2) {$C_{a_0}$};
            \draw[blue] node at (2.5,1.45) {$C_{a_2}$};
         
             \draw (-1.2,1.5)--(-1.2,-1.5)--(1.2,-1.5)--(1.2,0.4);
             \draw[dashed] (-1.2,-1.5)--(-0.2,-1);
             \draw[dashed] (-0.2,-1)--(2.2,-1);
             \draw[dashed] (-0.2,-1)--(-0.2,0.9);
             \draw (2.2,-1)--(2.2,0.9);
             \draw (1.2,-1.5)--(2.2,-1);
             \foreach \x in {-0.2,-0.05,...,2.2} {
            \draw[dashed, opacity = 0.7] (\x,-1) -- (\x,0.9);
            }

             \draw[green] (-1.2,1.5)--(1.2,1.5)--(2.2,2)--(-0.2,2)--cycle;
             \draw[fill, green, opacity=0.4] (-1.2,1.5)--(1.2,1.5)--(2.2,2)--(-0.2,2)--cycle;

             \draw[dashed, blue] (2.2,0.9)--(-0.2,0.9)--(-0.2,2);
             \draw[fill, blue, opacity=0.4] (-0.2,2)--(-0.2,0.9)--(2.2,0.9)--(2.2,2)--cycle;
             \draw[blue] (1.2,1.5)--(1.2,0.4)--(2.2,0.9)--(2.2,2);
             \draw[pattern=vertical lines, pattern color=blue] (1.2,1.5)--(1.2,0.4)--(2.2,0.9)--(2.2,2)--cycle;
             \draw[fill, blue, fill opacity=0.4] (1.2,1.5)--(1.2,0.4)--(2.2,0.9)--(2.2,2)--cycle;

             \draw[fill, opacity = 0.1] (-1.2,-1.5)--(1.2,-1.5)--(2.2,-1)--(2.2,2)--(-0.2,2)--(-1.2,1.5)--cycle;
         \end{tikzpicture}
         & 
        \begin{tikzpicture}
                \draw node at (0,0.2) {$\longrightarrow$};
                \draw node at (0,2.2) {$ $};
                \draw node at (0,-1.8) {$ $};
        \end{tikzpicture}
        &
        \begin{tikzpicture}
                \draw[green] node at (0,1.9) {$\widetilde{C}_{a_0}$};
                \draw[blue] node at (1.1,1.3) {$\widetilde{C}_{a_2}$};
                \draw node at (-0.3,0.1) {$\widetilde{C}_{a_3}$};
                
                \draw node at (0,2) {$ $};
                \draw node at (0,-2) {$ $};
                
                \pgfmathsetmacro{\L}{veclen(1.8,-1.2)};
                \pgfmathsetmacro{\theta}{atan2(-1.2,1.8)};
                \draw[blue, line width=1.4pt] (0,1.6)--(1.8,0.4);
                \draw[shift={(0,1.6)}, rotate=\theta, dashed, blue] (0,0) arc [start angle=-180, end angle=0, x radius=\L/2, y radius=0.3];
                \draw[fill, blue, opacity=0.4] (1.8,0.4)--(0,1.6)--++(0,0) [shift={(0,1.6)}, rotate=\theta] (0,0) arc [start angle=-180, end angle=0, x radius=\L/2, y radius=0.3];
                
                \draw (1.8,0.4)--(0,-1.6)--(0,1.6);
                \pgfmathsetmacro{\K}{veclen(1.8,2)};
                \pgfmathsetmacro{\nu}{atan2(2,1.8)};
                \draw[shift={(0,-1.6)}, rotate=\nu] (0,0) arc [start angle=-180, end angle=0, x radius=\K/2, y radius=0.3];
                \draw[fill, opacity=0.1] (1.8,0.4)--(0,1.6)--(0,-1.6)--++(0,0) [shift={(0,-1.6)}, rotate=\nu] (0,0) arc [start angle=-180, end angle=0, x radius=\K/2, y radius=0.3];
                
                \fill[green] (0,1.6) circle (1.5pt);
        \end{tikzpicture}
    \end{array}$$
    
    \vspace{-0.7cm}
    
    $$\textcolor{green}{a_0} < \textcolor{blue}{a_1} = \textcolor{blue}{a_2} < a_3$$
    
    \caption{The quotient of the cube $I^2 \times [f(a_3),f(a_0)]$ in the case of a sequence of four critical points satisfying $a_0 < a_1 = a_2 < a_3$. We denote this quotient by $\widetilde{C}$. The map $f_{\sigma}$ is constant equal to $a_0$ on $I^2 \times \{f(a_0)\}$ and equal to $a_3$ on $I^2 \times \{f(a_3)\}$. In addition, on the vertical segments $(\{(t,1)\} \times [f(a_3),f(a_1)])_{0 \le t \le 1}$ depicted in the figure, $f_{\sigma}$ corresponds to the evaluation of the same trajectory for every $t$. The same holds on the vertical segments $(\{(1,t)\} \times [f(a_1),f(a_0)])_{0 \le t \le 1}$ also depicted in the figure. Consequently, the map $f_{\sigma}$ is compatible with the equivalence relation that collapses to a point each of the two subsets $I^2 \times \{f(a_3)\}$ and $I^2 \times \{f(a_0)\}$, as well as each of the horizontal segments $I \times \{1\} \times \{s\}$ for $f(a_3) \le s \le f(a_1)$, and $\{1\} \times I \times \{s\}$ for $f(a_1) \le s \le f(a_0)$. We will explain below how to compare the resulting stratified quotiented cube with the corresponding stratified $3$-simplex.}
    \label{figure: quotient of cube for a0<a1=a2<a3}
\end{figure}
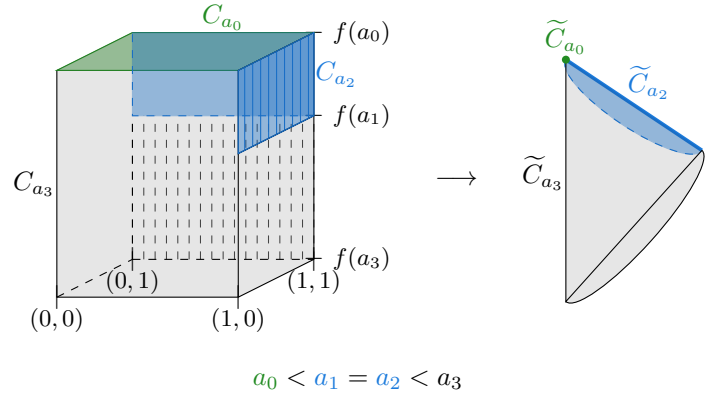

We would now like to point out that we do not merely wish to compare stratified simplices and stratified cubes \emph{individually}. We want to compare the simplicial sets $S_{\Ml}$ and $\Sing_A(X)$, and consequently, in addition to comparing the sets of $n$-simplices for every $n$, we also need some compatibility with the simplicial structure.

The reason why $\Sing_A(X)$ is a \emph{simplicial subset} of $\Sing(X)$ is that the face and degeneracy maps between standard simplices are \emph{compatible} with the stratifications on these. This can be precisely formulated as follows: the usual geometric realization functor $|-| : \Delta \rightarrow \Top$ lifts to a functor

$$
|-|_A : \Delta_A \rightarrow \Top_A.
$$

Here, $\Top_A$ is the category of $A$-stratified topological spaces, and $\Delta_A$ is the category of finite increasing sequences of critical points, with arrows generated by the face and degeneracy maps. More precise definitions of these notions are given in section \ref{section: stratified spaces, exit paths and stratified simplicial sets}. In the present section, we will just use the following notation: given a finite increasing sequence of critical points $\ag = a_0 \le a_1 \le \hdots \le a_n$, we denote by $\Delta^{\ag}$ the corresponding object of $\Delta_A$.

In section \ref{section: cubical stratified geometric realization}, we organize stratified quotiented cubes into a \emph{functor}

$$C_A : \Delta_A \rightarrow \Top_A.$$

This is another reason why we work with these \emph{quotients} of cubes, as these quotients are essential to define the functor $C_A$. One cannot realize this functor as a quotient of another functor whose objectwise values would be the cubes prior to quotienting.

We can define, thanks to $C_A$, a simplicial set $\Sing_{C_A}(X)$ in the same way as $\Sing_A(X)$ is defined from $|-|_A$ and the $A$-stratified space $X$, in other words:

$$
\begin{array}{ll}
\Sing_{C_A}(X)_n = & \{C_A(\Delta^{\ag}) \rightarrow X \text{ continuous and compatible with stratifications, for some increasing} \\
& \text{ sequence } \ag \text{ of } n+1 \text{ critical points }\},
\end{array}
$$

with the simplicial structure induced by the functorial structure on $C_A$. The functor $C_A$ is constructed in such a way that the natural assignment $S_{\Ml} \rightarrow \Sing_{C_A}(X)$ defined as $\sigma \mapsto f_{\sigma}$, is a \emph{morphism} of simplicial sets.

At this point of our discussion towards a comparison between $\Nl(\Ml)$ and $\Sing_A(X)$, the following questions arise:

\begin{itemize}
    \item Is the simplicial set $\Sing_{C_A}(X)$ an $\infty$-category? If yes, is the functor determined by $S_{\Ml} \rightarrow \Sing_{C_A}(X)$ an equivalence?
    
    \item To what extent can we compare the simplicial sets $\Sing_{C_A}(X)$ and $\Sing_A(X)$?
\end{itemize}

Note that the first question is related to the second, because we know that $\Sing_A(X)$ is an $\infty$-category.

Our strategy to address these questions can be briefly summarized as follows.

\begin{itemize}
    \item We prove a comparison result between the functors $|-|_A$ and $C_A$, which in turn implies a comparison result between the simplicial sets  $\Sing_{C_A}(X)$ and $\Sing_A(X)$. This result is not a mere isomorphism, but is sufficient for our purposes.
    
    \item We prove that the morphism $S_{\Ml} \rightarrow \Sing_{C_A}(X)$ is an equivalence of $\infty$-categories after corestriction to a suitable subcategory of the target.
\end{itemize}

We now explain this in more detail. We first focus our attention on the functors $|-|_A$ and $C_A$.

Let us summarize the similarities and differences between $|-|_A$ and $C_A$ that we have noted so far. At this point, we have seen on the above pictures that, for $n=1, 2$ and $3$, and given a sequence of $n$ \emph{distinct} critical points, the corresponding stratified $n$-simplex and stratified quotiented $n$-cube are stratified homeomorphic. On the other hand, when the sequence is \emph{constant}, the quotiented stratified cube is a point, so this statement does not hold.

Let us take a closer look at the stratification of the standard $n$-simplex $|\Delta^n|$ defined by the \emph{constant} sequence of length $n$, at some critical point $a_0$. It consists of the single $a_0$-stratum (see figures \ref{figure: stratification of the standard 1-simplex} and \ref{figure: stratification of the standard 2-simplex} for the cases $n=1$ and $n=2$). An $n$-simplex of $\Sing_A(X)$ associated with the constant sequence at $a_0$ is therefore the same as a singular $n$-simplex of the $a_0$-stratum of $X$.

Roughly speaking, such simplices of $\Sing_A(X)$ serve to capture the \emph{homotopy types of the strata}. But in our case, the strata are contractible, so we don't really need those kind of simplices. We can thus quotient the functor $|-|_A$ by its restriction to constant sequences, and the resulting functor, denoted $|-|_A'$, is just as good as $|-|_A$ to define the exit path $\infty$-category of this stratification of $X$.

\begin{Remark}\label{remark: contractible strata is necessary}
    Note that the endomorphism spaces of $\Ml$ are reduced to points, while the endomorphism space of a point $x \in X$ as an object of $\Sing_A(X)$, is the space of loops of the stratum of $x$, based at $x$. In particular, the condition that the strata are (weakly) contractible is necessary for an equivalence of $\infty$-categories between $\Ml$ and $\Sing_A(X)$ to exist.
\end{Remark}

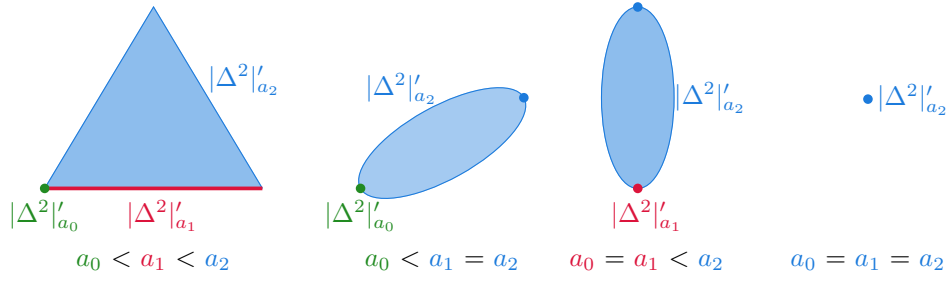
\begin{figure}[H]
    \centering
    $$\begin{array}{cccccc}
    \begin{tikzpicture}[scale=1.2]
        \draw node at (0,-0.8) {\textcolor{green}{$a_0$} < \textcolor{red}{$a_1$} < \textcolor{blue}{$a_2$}};
    
        \draw[blue] (-1.2,0)--(0,2)--(1.2,0);
        \draw[fill, blue, opacity=0.5] (-1.2,0)--(0,2)--(1.2,0);
        \draw[blue] node at (1,1.2) {$|\Delta^2|'_{a_2}$};
        \draw[red] node at (0.1,-0.3) {$|\Delta^2|'_{a_1}$};
        \draw[red, line width=1.4pt] (-1.2,0)--(1.2,0);
        \draw[green] node at (-1.2,-0.3) {$|\Delta^2|'_{a_0}$};
        \fill[green] (-1.2,0) circle (1.5pt);
    \end{tikzpicture}
    &
    \begin{tikzpicture}[scale=1.2]
        \draw node at (-0.3,-0.8) {\textcolor{green}{$a_0$} < \textcolor{blue}{$a_1$} = \textcolor{blue}{$a_2$}};
    
        \pgfmathsetmacro{\L}{veclen(1.8,1)};
        \pgfmathsetmacro{\theta}{atan2(1,1.8)};
        \draw[shift={(-0.3,0.5)}, rotate=\theta, blue] (0,0) ellipse (\L/2 and 0.4);
        \draw[fill, blue, opacity=0.4] [shift={(-0.3,0.5)}, rotate=\theta, blue] (0,0) ellipse (\L/2 and 0.4);
         \fill[blue] (0.6,1) circle (1.5pt);
        \draw[blue] node at (-0.75,1.1) {$|\Delta^2|'_{a_2}$};
        \draw[green] node at (-1.2,-0.3) {$|\Delta^2|'_{a_0}$};
        \fill[green] (-1.2,0) circle (1.5pt);
    \end{tikzpicture}
    &
    \begin{tikzpicture}[scale=1.2]
        \draw node at (0.1,-0.8) {\textcolor{red}{$a_0$} = \textcolor{red}{$a_1$} < \textcolor{blue}{$a_2$}};
    
        \draw[blue] (0,1) ellipse (0.4 and 1);
        \draw[fill, blue, opacity=0.5] (0,1) ellipse (0.4 and 1);
        \fill[blue] (0,2) circle (1.5pt);
        \draw[blue] node at (0.8,1) {$|\Delta^2|'_{a_2}$};
        \draw[red] node at (0.1,-0.3) {$|\Delta^2|'_{a_1}$};
        \fill[red] (0,0) circle (1.5pt);
    \end{tikzpicture}
     &
    \begin{tikzpicture}[scale=1.2]
        \draw node at (0,-0.8) {\textcolor{blue}{$a_0$} = \textcolor{blue}{$a_1$} = \textcolor{blue}{$a_2$}};
        \draw[blue] node at (0.1,-0.3) {$ $};
        \fill[blue] (0,1) circle (1.5pt);
        \draw[blue] node at (0.5,1) {$|\Delta^2|'_{a_2}$};
    \end{tikzpicture}
    \end{array}$$
    \caption{The possible quotients of the standard $2$-simplex determined by an increasing sequence of three critical points $a_0 \le a_1 \le a_2$. When the critical points are distinct, the constant subsequences of $a_0,a_1,a_2$ are each reduced to one point, so the quotiented simplex is the same as the original simplex. In the case when $a_0 < a_1 = a_2$, the segment joining $e_1$ to $e_2$, which corresponds to the constant subsequence $a_1 = a_2$, is collapsed to a point. In the case when $a_0 = a_1 < a_2$, the segment joining $e_0$ to $e_1$, which corresponds to the constant subsequence $a_0 = a_1$, is collapsed to a point. In the case when $a_0 = a_1 = a_2$, the whole $2$-simplex, which corresponds to the constant sequence $a_0 = a_1 = a_2$, is collapsed to a point.}
    \label{figure: quotiented stratified 2-simplices}
\end{figure}

\begin{figure}[H]
    \centering
    $$\begin{array}{ccccc}
    \begin{tikzpicture}[scale=1.2]
        
        \draw (-1.2,0)--(0,2);
        \draw[blue] (1.2,0)--(1.7,0.6);
        \draw (1.7,0.6)--(0,2);
        \draw[dashed, blue] (-1.2,0)--(1.7,0.6);
        \draw (1.2,0)--(0,2);
        
        \draw[fill, opacity = 0.1] 
                (-1.2,0)--(0,2)--(1.7,0.6)--(1.2,0)--cycle;
         \draw[fill, blue, opacity=0.4] 
                (-1.2,0)--(1.7,0.6)--(1.2,0)--cycle;

        \draw[red, line width=1.4pt] (-1.2,0)--(1.2,0);
        \fill[green] (-1.2,0) circle (1.5pt);

        \draw[green] node at (-1.2,-0.3) {$|\Delta^3|'_{a_0}$};
        \draw[blue] node at (1.9,0.2) {$|\Delta^3|'_{a_2}$};
        \draw node at (1.2,1.5) {$|\Delta^3|'_{a_3}$};
        \draw[red] node at (0.1,-0.3) {$|\Delta^3|'_{a_1}$};
        \draw node at (0.25,-0.8) {\textcolor{green}{$a_0$} < \textcolor{red}{$a_1$} < \textcolor{blue}{$a_2$} < $a_3$};
        
    \end{tikzpicture}
    &
    \begin{tikzpicture}[scale=1.2]
        

        \draw[blue] (-1.2,0) arc [start angle=-180, end angle=0, x radius=1.3, y radius=0.2];
        \draw[blue, dashed] (-1.2,0) arc [start angle=180, end angle=0, x radius=1.3, y radius=0.2];
        \draw[fill, blue, opacity=0.4] (0.1,0) ellipse (1.3 and 0.2);
        
        \pgfmathsetmacro{\theta}{atan2(-2,1.4)};
        \pgfmathsetmacro{\L}{veclen(1.4,-2)};
        \draw (-1.2,0)--(0,2);
        \draw [shift={(0.7,1)}, rotate=\theta] (0,0) ellipse (\L/2 and 0.2);
        \draw[fill, opacity=0.1] (1.4,0) 
        ([shift={(1.4,0)}] 0,0) 
        arc[start angle=0, end angle=-180, x radius=1.3, y radius=0.2]--(0,2)--++(0,0) [shift={(0,2)}, rotate=\theta] (0,0) arc [start angle=180, end angle=0, x radius=\L/2, y radius=0.2];

        \fill[green] (-1.2,0) circle (1.5pt);

        \draw[green] node at (-1.2,-0.3) {$|\Delta^3|'_{a_0}$};
        \draw[blue] node at (1.9,0.2) {$|\Delta^3|'_{a_2}$};
        \draw node at (1.2,1.5) {$|\Delta^3|'_{a_3}$};
        \draw node at (0.15,-0.8) {\textcolor{green}{$a_0$} < \textcolor{blue}{$a_1$} = \textcolor{blue}{$a_2$} < $a_3$};
        
    \end{tikzpicture}
    &
    \begin{tikzpicture}[scale=1.2]
        
        \draw[green] node at (-1.2,-0.3) {$|\Delta^3|'_{a_0}$};
        \draw node at (1.2,1.6) {$|\Delta^3|'_{a_3}$};
        \draw node at (0.15,-0.8) {\textcolor{green}{$a_0$} < $a_1$ = $a_2$ = $a_3$};
        
        \pgfmathsetmacro{\L}{veclen(2.65,1.3)};
        \pgfmathsetmacro{\theta}{atan2(1.3,2.65)};
        \draw [shift={(0.125,0.65)}, rotate=\theta] (0,0) ellipse (\L/2 and 0.45);
        \draw[fill, opacity=0.1] [shift={(0.125,0.65)}, rotate=\theta] (0,0) ellipse (\L/2 and 0.45);
        \draw [shift={(-1.2,0)}, rotate=\theta, dashed] (0,0) arc [start angle=-180, end angle=0, x radius=\L/2, y radius=0.25];

        \fill[green] (-1.2,0) circle (1.5pt);
        
    \end{tikzpicture}
    \end{array}$$
    \caption{Some possible quotients of the standard $3$-simplex determined by a sequence of three critical points $a_0 \le a_1 \le a_2 \le a_3$. When the sequence is strictly increasing, the constant subsequences of $a_0,a_1,a_2, a_3$ are each reduced to one point, so the quotiented simplex is the same as the original simplex. In the case when $a_0 < a_1 = a_2 < a_3$, the segment joining $e_1$ to $e_2$, which corresponds to the constant subsequence $a_1 = a_2$, is collapsed to a point. In the case when $a_0 < a_1 = a_2 = a_3$, the face generated by the vertices $e_1, e_2, e_3$, which corresponds to the constant subsequence $a_1 = a_2 = a_3$, is collapsed to a point.}
    \label{figure: quotiented stratified 3-simplices}
\end{figure}
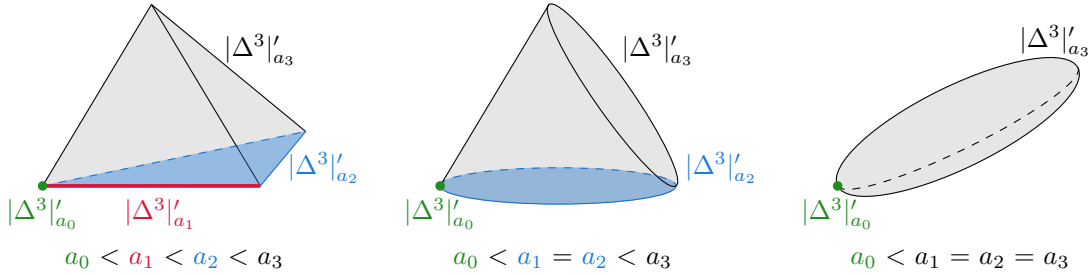

In practice, for every $\Delta^{\ag} \in \Delta_A$, $|\Delta^{\ag}|_A'$ is obtained from $|\Delta^{\ag}|_A$ by collapsing to a point each of the subsimplices of $|\Delta^{\ag}|_A$ corresponding to \emph{constant} subsequences of $\ag$ (the detailed construction is given as construction \ref{construction: functor |-|_A'}). In the cases presented above in figures \ref{figure: stratification of the standard 2-simplex} and \ref{figure: some stratifications of the 3 simplex}, the results are described in figures \ref{figure: quotiented stratified 2-simplices} and \ref{figure: quotiented stratified 3-simplices}.

We can now explain all the terms of the zigzag \eqref{zigzag} from the beginning of this section. Let us recall it for convenience

\[
\Nl(\Ml) \leftarrow S_{\Ml} \rightarrow \Sing_{C_A,r}(X) \rightarrow \Sing'_{A,r}(X) \rightarrow \Sing_A(X).
\]

Again, we define from $|-|_A'$ and the stratification on $X$, a simplicial set $\Sing'_A(X)$ as

$$
\begin{array}{ll}
\Sing'_A(X)_n = & \{|\Delta^{\ag}|'_A \rightarrow X \text{ continuous and compatible with stratifications,} \\
& \text{ for some sequence } \ag \text{ of } n+1 \text{ critical points }\}.
\end{array}
$$

Note that the quotient natural transformation $|-|_A \rightarrow |-|'_A$ induces an inclusion of simplicial sets $\Sing'_A(X) \rightarrow \Sing_A(X)$.

Let us explain the meaning of the subscript "$_r$" in the zigzag. For this, observe that the simplices of $\Sing'_A(X)$ and $\Sing_{C_A}(X)$ that are entirely contained in a given stratum, are constant. Said differently, the restriction of these simplicial sets to a given stratum is always discrete. In contrast, this is not the case for $\Sing_A(X)$ (as long as the stratum is of dimension at least 1). To remedy this, we introduce the \emph{restricted versions} of these two simplicial sets, denoted with a subscript "$_r$", and defined to be the respective simplicial subsets formed by those simplices whose vertices are critical points. This modification is harmless: using the fact that the strata of $X$ are contractible, we can show that the inclusion $\Sing'_{A,r}(X) \rightarrow \Sing_A(X)$ is an equivalence of $\infty$-categories (proposition \ref{proposition: restricted version into unrestricted is equivalence}).

Next, by construction, the morphism of simplicial sets $S_{\Ml} \rightarrow \Sing_{C_A}(X)$ corestricts to a morphism of simplicial sets $S_{\Ml} \rightarrow \Sing_{C_A,r}(X)$. We prove in section \ref{section: comparison between S M and Sing C A (X)} that this is an equivalence of $\infty$-categories (proposition \ref{proposition: S_M -> Sing_C_A,r is an equivalence}). Again, this is done by proving that it is essentially surjective and fully faithful. The geometric reason behind the full faithfulness of this functor is the result that the inclusion of the space of unbroken trajectories between two given critical points, into the space of exit paths between these, is a weak homotopy equivalence. The proof of this result is the main purpose of section \ref{section: spaces of exit paths as spaces of pseudo gradient trajectories}. Our proof combines a description of the spaces of exit path in cylindrically stratified spaces (which is presented in section \ref{section: geometric description of spaces of exit paths}), with cylindrical properties of Whitney stratifications (which are presented in section \ref{section: Whitney conditions and cylindrical structures}). It uses as an essential ingredient a theorem of Nicolaescu (\cite{NicolaescuInvitation}, theorems 4.32 and 4.33) asserting that the Smale transversality condition for the pseudo-gradient $\xi$ is equivalent to the condition that the partition of $X$ by the stable manifolds of $\xi$ is a Whitney stratification.

Lastly, it remains to compare the functors $C_A$ and $|-|_A'$. If these functors were isomorphic, the two simplicial sets $\Sing_{C_A}(X)$ and $\Sing_A'(X)$ would be isomorphic, and so would be the simplicial sets $\Sing_{C_A,r}(X)$ and $\Sing_{A,r}'(X)$, which would complete the construction of a zigzag of equivalences between $\Nl(\Ml)$ and $\Sing_A(X)$. However, as long as $\dim (X) > 0$, these two functors are \emph{not} isomorphic (proposition \ref{proposition: functors not isomorphic}). In other words, it is not possible to find a family of stratified homeomorphisms

$$
C_A(\Delta^{\ag}) \simeq |\Delta^{\ag}|'_A
$$

for all finite increasing sequences of critical points, which are compatible with both face and degeneracy maps. Still, a weaker, but sufficient result, does hold. Namely, we prove that such a family of stratified homeomorphisms exist if one requires compatibility with \emph{faces} only (theorem \ref{theorem: isomorphism between |-|'_A and C_A}). In other words, the restrictions of these two functors to the subcategory of $\Delta_A$ generated by faces, are isomorphic.

This has two essential consequences for us. Firstly, this result suffices to imply that $\Sing_{C_A}(X)$ is an $\infty$-category (proposition \ref{proposition: Sing_C_A is an infinity category}), and therefore also $\Sing_{C_A,r}(X)$. Secondly, it implies that $\Sing_{C_A}(X)$ and $\Sing_A'(X)$ are isomorphic as \emph{semi-simplicial sets}, and likewise for $\Sing_{C_A,r}(X)$ and $\Sing_{A,r}'(X)$. Tanaka shows in \cite{HiroFunctors} how to naturally associate a functor of $\infty$-categories to a morphism of semi-simplicial sets between $\infty$-categories, provided that the latter preserves identities up to homotopy. We discuss this result in section \ref{section: end of the proof}. It applies in our case and yields an equivalence of $\infty$-categories between $\Sing_{C_A,r}(X)$ and $\Sing_{A,r}'(X)$.  With this, the proof of our main theorem is complete.

To finish this section, we quickly outline our proof that $C_A$ and $|-|_A'$ are isomorphic in restriction to the subcategory generated by face morphisms. For every finite increasing sequence of critical points $\ag$ of length $n$, $C_A(\Delta^\ag)$ (resp. $|\Delta^\ag|_A'$) is constructed as some quotient of an $n$-cube (resp. an $n$-simplex). Both the cube and the simplex are homeomorphic to the closed Euclidean ball $D^n$. The key point of our proof is to show that this remains true after passing to the respective quotients (as long as $\ag$ is not constant), and that any homeomorphism $C_A(\Delta^\ag) \simeq D^n$ (resp. $|\Delta^\ag|_A' \simeq D^n$) identifies the union of the faces of $C_A(\Delta^\ag)$ (resp. $|\Delta^\ag|_A'$) with $\partial D^n = \SM^{n-1}$. This result is formulated as theorem \ref{theorem: stratified simplices and cubes as balls}; although it might be intuitively clear that such a statement holds, it takes quite some work to prove. The proof is the purpose of section \ref{section: proof of the ball theorem for realizations}. The argument analyzes the geometry of polytopes and shows that certain quotients of polytopes are themselves homeomorphic to polytopes. With this theorem at hand, we can then argue by induction: once a family of homeomorphisms compatible with faces is constructed for sequences of length at most $n$, this family can be extended to sequences of length at most $n+1$ as follows. If $\ag$ is of length $n+1$ and not constant, combining two homeomorphisms $C_A(\Delta^\ag) \simeq D^n$ and $|\Delta^\ag|_A' \simeq D^n$, with the family of homeomorphisms between the faces of $C_A(\Delta^\ag)$ and $|\Delta^\ag|_A'$, available to us by the induction hypothesis, we get a homeomorphism of the sphere $\SM^{n-1}$. A homeomorphism between $C_A(\Delta^\ag)$ and $|\Delta^\ag|_A'$, which is both stratified and compatible with the faces, is then the same as an extension of this homeomorphism of $\SM^{n-1}$ to a homeomorphism of $D^n$. The conclusion follows because such an extension exists.

\phantomsection
\addcontentsline{toc}{subsection}{Acknowledgements}
\subsection*{Acknowledgements}

I would first like to warmly thank my PhD advisor Alexandru Oancea for constantly supporting me throughout this project and for carefully reading this text. I also owe many thanks to Sylvain Douteau, Mauro Porta, Robin Riegel and Lukas Waas for precious discussions closely related to this paper, and to Jonathan Block for first telling me about exit path categories. I am also very grateful to the members of my PhD committee for their feedbacks on this work, and especially to Mohammed Abouzaid and Hiro Lee Tanaka for pointing out several minor mistakes and possible clarifications. Finally I would like to thank Mihai Damian, Grégory Ginot, Jean Gutt and Manuel Rivera for their encouraging interest in this work. I benefited from the ANR grant COSY 21-CE40-0002.

\newpage

\section{\texorpdfstring{$(\infty,1)$-categories, stratified topological spaces and Morse-Smale pairs}{infinity,1-categories, stratified spaces and Morse-Smale pairs}}\label{section: infinity categories, stratified spaces and Morse-Smale pairs}

\subsection{Categories of presheaves, left Kan extensions and simplicial sets}\label{section: Left Kan extensions}

In this paper, we make a very important use of the category of simplicial sets, as well as the more general notions of categories of presheaves, and left Kan extensions to these. This section is intended to serve as a reference about these concepts. We start by reviewing the basic properties of categories of presheaves, and left Kan extensions to categories of presheaves. We then focus on the particular case of simplicial sets. For a more comprehensive review, see for example \cite[Section 2.2]{DuggerSheaves} or \cite{RiehlCategory}.

We denote by $\Set$ the category of sets.

Let $\Cl$ be a small category.

\begin{Definition}\label{definition: category of presheaves}
    The \textit{category of presheaves} on $\Cl$ is the functor category $\Fun(\Cl^{\op},\Set)$. It is denoted $\Pre (\Cl)$.
\end{Definition}

\begin{Remark}\label{remark: presheaf categories admit all limits and colimits}
    The category $\Pre(\Cl)$ admits all (small) limits and colimits. For every diagram $D : J \rightarrow \Pre(\Cl)$, the limit and colimits of $D$ are given explicitly by the formulas

    $$
    (\underset{D}{\lim} \, J)(C) = (\underset{D}{\lim} \, J(C)) \qquad (\underset{D}{\colim} \, J)(C) = (\underset{D}{\colim} \, J(C)).
    $$
\end{Remark}

\begin{Definition}\label{definition: Yoneda functor}
    The \textit{Yoneda functor} is the functor

    $$
    \fonction{r}{\Cl}{\Pre (\Cl)}{C}{\Hom_{\Cl}(-,C).}
    $$

    The presheaf $r(C)$ is said to be \textit{represented} by $C$. Objects in the essential image of $r$ are called \textit{representables}.
\end{Definition}

The following lemma is known as the Yoneda lemma:

\begin{Lemmas}\label{Yoneda lemma}
    For every $C \in \Cl$ and $F \in \Pre (\Cl)$, the correspondence

    $$
    \fonctionsansnom{\Hom_{\Pre (\Cl)}(r(C),F)}{F(C)}{\varphi}{\varphi(C)(\id_C)}
    $$

    is a bijection. In particular, $r$ is fully faithful. \qed
\end{Lemmas}

An important example for us is the special case of the category of simplicial sets, which we now introduce.

\begin{Notation}\label{notation: posets}
    We adopt the following notation about partially ordered sets (posets).
    
    \begin{itemize}
        \item  We denote by $\Pos$ the category of posets. Recall that this is the category whose objects are posets, and whose morphisms are order-preserving maps between posets.
        \item Given $n \ge 0$, we denote by $[n]$ the ordered set $0 < 1 < \hdots < n$, which we regard as an object of $\Pos$.
    \end{itemize}
   
\end{Notation}

\begin{Definition}\label{definition: simplex category}
    The \textit{simplex category} $\Delta$ is the full subcategory of $\Pos$ whose objects are the posets $[n]$, for $n$ ranging through nonnegative integers. Given an integer $n \ge 0$, we denote by $\Delta^n$ the object of $\Delta$ corresponding to the poset $[n]$.
\end{Definition}

\begin{Definition}\label{definition: simplicial set}
    A \textit{simplicial set} is an object of the category of presheaves $\Pre(\Delta)$. This category of presheaves is called the \textit{category of simplicial sets} and denoted $\sSet$.
\end{Definition}

Let us now introduce some notation and terminology about simplicial sets.

\begin{Notation and Terminology}\label{notation and terminology about simplicial sets}

Let $n$ be a nonnegative integer and $0 \le i \le n$.

    \begin{itemize}
        \item For $n \ge 1$, the only injective arrow in $\Hom_{\Delta}(\Delta^{n-1},\Delta^n)$ whose image does not contain $i$ is denoted $\delta_n^i$. Arrows of this form are called \textit{face morphisms}, or simply \emph{faces}.

        \item The only surjective arrow in $\Hom_{\Delta}(\Delta^{n+1},\Delta^n)$ which maps both $i$ and $i+1$ to $i$ is denoted $\sigma_n^i$. Arrows of this form are called \textit{degeneracy morphisms}, or simply \emph{degeneracies}.

        \item The object of $\sSet$ represented by $\Delta^n$ is denoted $\Delta^n$ as well, and is called the \textit{standard} $n$-\textit{simplex}.
        
        \item Given $0 \le i \le n$, we refer to the corresponding morphism of simplicial sets $\Delta^0 \rightarrow \Delta^n$ as the \emph{vertex} $i$ of $\Delta^n$.

        \item We refer to the morphisms of simplicial sets $\delta_n^i : \Delta^{n-1} \rightarrow \Delta^n$ for $0 \le i \le n$ as the \textit{faces} of $\Delta^n$. We refer to the morphism $\delta_n^i : \Delta^{n-1} \rightarrow \Delta^n$ as the \textit{face opposite to the vertex} $i$.
        
    \end{itemize}
        
     Let $T$ be a simplicial set.

    \begin{itemize}
        \item The value of $T$ at the object $\Delta^n$ is denoted $T_n$.

        \item The elements of $T_n$ are called $n$-\textit{simplices} of $T$. By lemma \ref{Yoneda lemma}, they correspond to elements of $\Hom_{\sSet}(\Delta^n,T)$. We will often abuse notation and identify these two sets.

        \item The $0$-simplices of $T$ are called the \emph{vertices} of $T$. Given an $n$-simplex $\sigma : \Delta^n \rightarrow T$, the \textit{vertex} $i$ of $\sigma$ is the $0$-simplex of $T$ obtained by restriction of $\sigma$ along the morphism $\Delta^0 \rightarrow \Delta^n$ induced by $i$. 

        \item The map $T(\delta_n^i) : T_n \rightarrow T_{n-1}$ is denoted $d_i^n$. Maps of this form are called \textit{face operators} of $T$.

        \item The map $T(\sigma_n^i) : T_{n} \rightarrow T_{n+1}$ is denoted $s_i^n$. Maps of this form are called \textit{degeneracy operators} of $T$.
        
        \item A simplex of $T$ which is not in the image of a degeneracy operator is called \emph{nondegenerate}.
    \end{itemize}
    
    We will in particular use the following simplicial sets.
    
    \begin{itemize}

        \item The \textit{boundary} of $\Delta^n$ is, informally, the simplicial set obtained from $\Delta^n$ by removing the interior. It is denoted by $\partial \Delta^n$. Formally, it is the simplicial subset of $\Delta^n$ defined as

        $$
        (\partial \Delta^n)_m = \{ \varphi : \Delta^m \rightarrow \Delta^n \mid \varphi \text{ is not surjective} \}.
        $$

        We will later give another description of $\partial \Delta^n$ as the gluing of the faces of $\Delta^n$ along their common faces (proposition \ref{proposition: boundary of standard simplex as colimit of its faces}).

        \item The $i^{\mathrm{th}}$ \textit{horn in} $\Delta^n$ is, informally, the simplicial set obtained from $\Delta^n$ by removing the interior together with the face opposite to the vertex $i$. It is denoted by $\Lambda_i^n$. Formally, it is the simplicial subset of $\Delta^n$ defined as

        $$
        (\Lambda_i^n)_m = \{ \varphi : \Delta^m \rightarrow \Delta^n \mid \varphi(\{0,1,\hdots,m\}) \cup \{i\} \neq \{0,1,\hdots,n\}\}.
        $$

        We will later give another description of $\Lambda_i^n$ as the gluing of the faces of $\Delta^n$ different from the one opposite to the vertex $i$, along their common faces (proposition \ref{proposition: horn as colimit of faces}).

        \item The horn $\Lambda_i^n$ is called an \textit{inner} horn when  $0 < i < n$.
        
        \item The \emph{spine} of $\Delta^n$ is, informally, the simplicial subset of $\Delta^n$ comprised of all vertices of $\Delta^n$ together with those edges that join successive vertices. Formally, it is the simplicial subset of $\Delta^n$ whose $k$-simplices are morphisms $\sigma : \Delta^k \rightarrow \Delta^n$ satisfying $\sigma(k) \le \sigma(0) + 1$. It is denoted $\Spine(\Delta^n)$.
    
        \item Given two simplicial sets $K$ and $L$, one defines the simplicial set $\underline{\Hom}_{\sSet}(K,L)$ by the formula

         $$
            \underline{\Hom}_{\sSet}(K,L)_n = \Hom_{\sSet}(\Delta^n \times K,L),
        $$
    
        with face and degeneracy operators determined by face and degeneracy maps between standard simplices. Note that the set of $0$-simplices of the latter is the set of morphisms of simplicial sets between $K$ and $L$.
         
         \item Let $g : K \rightarrow L$ be a morphism of simplicial sets and $x$ a $0$-simplex of $L$, which we regard as a morphism of simplicial sets $\Delta^0 \rightarrow L$. The \emph{fiber} of $g$ over $x$ is the simplicial set defined as the pullback of the diagram
    
         $$
        \xymatrix{
        & K \ar[d]^-g \\
        \Delta^0 \ar[r]_-x & L.
        }
        $$
    \end{itemize}
   
\end{Notation and Terminology}

We will also use the notion of \textit{semi-simplicial set}, which we now introduce. 

\begin{Notation}\label{notation: Delta inj}
    Let $\Delta^+$ denote the subcategory of $\Delta$ with the same objects, and morphisms corresponding to injective maps only.
\end{Notation}

\begin{Remark}\label{remark: face morphisms generate delta +}
    The category $\Delta^+$ can also be described as the subcategory of $\Delta$ with the same objects and morphisms generated under composition by the face morphisms.
\end{Remark}

\begin{Definition}\label{definition: semi-simplicial set}
    A \textit{semi-simplicial set} is an object of the category of presheaves $\Pre(\Delta^+)$. This category is called the \emph{category of semi-simplicial sets} and denoted $\sSet^+$.
\end{Definition}

\begin{Remark}\label{remark: a simplicial set defines a semi-simplicial set}
    To every simplicial set $T$ is associated an underlying semi-simplicial set by restriction along the functor $(\Delta^+)^{\op} \rightarrow \Delta^{\op}$. This semi-simplicial set sends $\Delta^n$ to $T_n$ for every $n$, and has structural maps given by the faces operators of $T$ and their compositions.
\end{Remark}

We now come back to the general setting of a (small) category $\Cl$. We will see that every object of $\Pre (\Cl)$ can be naturally written as a colimit of representables.

\begin{Definition}\label{definition: overcategory}
    Given a category $\Dl$ and an object $X \in \Dl$, the \textit{overcategory} of $X$ is the category $\Dl_{/X}$ whose objects are morphisms of $\Dl$ with target $X$, and whose morphisms between two objects $\phi : Y \rightarrow X$ and $\psi : Z \rightarrow X$ are the commutative diagrams of the form

    $$
    \xymatrix{
    Y \ar[rr] \ar[rd]_{\phi} && Z \ar[ld]^{\psi} \\
    & X.
    }
    $$
\end{Definition}

\begin{Definition}\label{definition: diagram of presheaves induced by an object}
    Let $F \in \Pre(\Cl)$. The \textit{category of elements of} $F$, denoted $D_F$, is the subcategory of the overcategory $\Pre(\Cl)_{/F}$ defined as the pullback
    
    $$
    \squarediagram{D_F}{\Pre(\Cl)_{/F}}{\Cl}{\Pre(\Cl)}{}{}{}{r}
    $$
    
    where the right vertical arrow is the forgetful functor. Equivalently, by lemma \ref{Yoneda lemma}, this is the category described as

    $$
    D_F = \left\{
    \begin{array}{ll}
        \Obj & = (C,x) \text{ such that } x \in F(C) \\
        \Hom((C,x),(D,y)) & = \{f \in \Hom_{\Cl}(C,D) \text{ such that } F(f)(y) = x \}.  
    \end{array}
    \right.
    $$

    There is an associated functor

    $$
    \fonction{c_F}{D_F}{\Pre(\Cl),}{(C,x)}{r(C).}
    $$
\end{Definition}

Let $G \in \Pre(\Cl)$ and $c_F \rightarrow G$ be a natural transformation. By lemma \ref{Yoneda lemma}, this is the same as the datum of an object of $G(C)$ for every $C \in \Cl$ and every element of $F(C)$, compatible with  the map $F(f) : F(D) \rightarrow F(C)$ for every $C,D \in \Cl$ and $f \in \Hom_{\Cl}(C,D)$. In other words, this is the same as a morphism $F \rightarrow G$ in $\Pre(\Cl)$. Hence we have the following result:

\begin{Proposition}\label{proposition: a presheaf is a colimit of representables}
    $F$ is the colimit of $c_F$. \qed
\end{Proposition}

Now let $\Dl$ be an ordinary category which admits all (small) colimits and $H : \Cl \rightarrow \Dl$ a functor. We would like to extend $H$ as a colimit-preserving functor out of $\Pre(\Cl)$. This leads to the following construction.

\begin{Construction}\label{construction: left Kan extension}
    The \textit{left Kan extension} of $H$ to $\Pre(\Cl)$ is the functor given by the formula

    $$
    \fonctionsansnom{\Pre(\Cl)}{\Dl,}{F}{\underset{D_F}{\colim} \, (H \circ c_F).}
    $$
\end{Construction}

\begin{Notation}\label{notation: left Kan extension}
    Depending on the functor $H$ that we consider, we will denote the left Kan extension of $H$ to $\Pre(\Cl)$ by $H$ as well, or by $|-|_H$.
\end{Notation}

\begin{Remark}\label{Remark: concrete description of left Kan extension}
    Denote the left Kan extension of $H$ to $\Pre(\Cl)$ by $H$ as well. More concretely, given $F \in \Pre(\Cl)$, the object $H(F) \in \Dl$ can be explicitly described as obtained from the disjoint union
    
    $$
    \bigsqcup_{c \in \Cl} F(c) \times H(c)
    $$
    
    by making identifications along the morphisms of $\Cl$. Namely, for every morphism $\varphi : c \rightarrow d$ of $\Cl$, there are identifications determined by coequalizing the two arrows
    
    $$
    F(d) \times H(c) \underset{\id \times H(\varphi)}{\overset{F(\varphi) \times \id}{\rightrightarrows}} F(c) \times H(c) \sqcup F(d) \times H(d).
    $$
\end{Remark}

To prove that $H$ commutes with all colimits, we observe that it admits a right adjoint, and it is a general fact that functors that admit a right adjoint commute with colimits. Suppose such a right adjoint exists, and denote it $\Sing_H(-)$. We must have for every $D \in \Dl$ and $C \in \Cl$

$$
\begin{aligned}
    \Sing_H(D)(C) & = \Hom_{\Pre(\Cl)}(r(C),\Sing_H(D)) \text{ by lemma \ref{Yoneda lemma}} \\
    & = \Hom_{\Dl}(H(C),D) \text{ by adjunction.}
\end{aligned}
$$

\begin{Construction}\label{construction: right adjoint of left Kan extension}
    The functor $\Sing_H$ is given by the formula

    $$
    \Sing_H(D)(C) = \Hom_{\Dl}(H(C),D).
    $$
\end{Construction}

This definition is made so that $H$ and $\Sing_H$ are adjoint at the level of representable objects of $\Pre(\Cl)$. But since every object of $\Pre(\Cl)$ is a colimit of representables, these are in fact adjoint functors. In particular, the left Kan extension of $H$ to $\Pre(\Cl)$ commutes with colimits.

To illustrate these constructions, let us consider the following functor in the case $\Cl = \Delta$ (definition \ref{definition: simplex category}).

\begin{Definition}\label{definition: geometric realization at the level of simplex category}
    The \textit{geometric realization functor} is the functor, denoted $|-|$, from the category of simplicial sets to the category of topological spaces, defined as the left Kan extension of the following functor

    $$\fonction{|-|}{\Delta}{\Top}{
    \left\{
    \begin{array}{l}
        \Delta^n \\
        \Delta^n \overset{\varphi}{\longrightarrow} \Delta^m
    \end{array}
    \right.
    }{
    \left\{
    \begin{array}{l}
        \{(t_0,\hdots,t_n) \in \R^{n+1} \mid t_i \ge 0 \text{ for all } 0 \le i \le n \text{ and } \sum_i t_i = 1\} \\
        (t_0,\hdots,t_n) \mapsto (s_0, \hdots, s_m) \text{ with } s_i = \sum_{\varphi(j) = i} t_j.
    \end{array}
    \right.
    }$$
    
    Let us explain how remark \ref{Remark: concrete description of left Kan extension} specializes to this case. For every simplicial set $S$, the geometric realization of $S$ can be concretely described as the quotient of the disjoint union
    
    $$
    \bigsqcup_{n \ge 0} S_n \times |\Delta^n|
    $$
    
    by the equivalence relation generated by the following relations: 
    
    \begin{itemize}
    
    \item $(d_i^n(\tau),x) \sim (\tau, |\delta_n^i|(x))$ for every $n \ge 1$, $\tau \in S_n$, $x \in |\Delta^{n-1}|$ and $0 \le i \le n$.
    
    \item $(s_i^n(\tau),y) \sim (\tau, |\sigma_n^i|(y))$ for every $n \ge 0$, $\tau \in S_n$, $y \in |\Delta^n|$ and $0 \le i \le n$.
    
    \end{itemize}

    The right adjoint of the geometric realization functor is called the \textit{singular simplicial set functor} and is denoted $\Sing(-)$. By construction \ref{construction: right adjoint of left Kan extension}, for every topological space $Y$, the set of $n$-simplices of $\Sing(Y)$ can be described as

    $$
    \Sing(Y)_n = \text{set of continuous maps from } |\Delta^n| \text{ to } Y.
    $$
\end{Definition}

\subsection{\texorpdfstring{$(\infty,1)$-categories}{infinity,1-categories}}\label{section: infinity categories}

This section is devoted to the theory of $(\infty,1)$-categories and follows \cite{HigherTopos} and \cite{kerodon}. The theory of $(\infty,1)$-categories is part of a more general theory of $(\infty,n)$-categories, but $(\infty,1)$-categories are often referred to as $\infty$-categories in the literature. We will do the same from now on.

The theory of $\infty$-categories can be regarded as a common generalization of ordinary category theory and the classical homotopy theory of CW-complexes. The first two sections are dedicated to these, respectively. The definition of an $\infty$-category is given in section \ref{section: infinity categories: definitions}, as well as the $\infty$-categorical counterparts of several notions from classical category theory, such as compositions, functors and equivalences. Finally in section \ref{section: topological categories and infinity categories} we explain how topological categories provide a model for the theory of $\infty$-categories.

\subsubsection{Ordinary categories as simplicial sets}\label{section: ordinary categories as simplicial sets}

To every (small) category $\Cl$ is associated a simplicial set $N(\Cl)$ called the \textit{nerve} of $\Cl$. The $n$-simplices of $N(\Cl)$ are the commutative diagrams in $\Cl$ supported by the vertices and edges of the standard $n$-simplex. From $n = 0$ to $n = 3$, these look as follows (the face and degeneracy operators of this simplicial set will be described below).

\begin{figure}[H]
    \centering
    \begin{tikzpicture}
        \begin{scope}[scale=1.1, xshift=-7cm]
        \draw node at (0,0) {$X_0$};
        \draw node at (0,-1) {$n=0$};
        \end{scope}
        
        \begin{scope}[scale=1.1, xshift=-3cm]
        \draw node at (0,-0.3) {$f_{0,1}$};
        \draw node at (-1.2,0) {$X_0$};
        \draw node at (1.2,0) {$X_1$};
        
        \draw node at (0,-1) {$n=1$};

        \draw[->] (-0.9,0) to (0.9,0);
        
        \end{scope}
        
        \begin{scope}[scale=1.1, xshift=1cm]
        \draw node at (0,-0.3) {$f_{0,1}$};
        \draw node at (0.9,1.2) {$f_{1,2}$};
        \draw node at (-0.9,1.2) {$f_{0,2}$};
        \draw node at (-1.2,0) {$X_0$};
        \draw node at (1.2,0) {$X_1$};
        \draw node at (0,2) {$X_2$};
        
        \draw node at (0,-1) {$n=2$};
        
        \draw[->] (-0.9,0) to (0.9,0);
        \draw[->] (1,0.2) to (0.2,1.8);
        \draw[->] (-1,0.2) to (-0.2,1.8);
        
        \end{scope}
        
        \begin{scope}[scale=1.1, xshift=5cm]
        
        \draw node at (0,-0.3) {$f_{0,1}$};
        \draw node at (0.25,1.1) {$f_{1,3}$};
        \draw node at (-0.9,1.2) {$f_{0,3}$};
        \draw node at (1.2,1.5) {$f_{2,3}$};
        \draw node at (-0.2, 0.5) {$f_{0,2}$};
        \draw node at (1.9,0.2) {$f_{1,2}$};
        \draw node at (-1.2,0) {$X_0$};
        \draw node at (1.2,0) {$X_1$};
        \draw node at (0,2) {$X_3$};
        \draw node at (1.9,0.6) {$X_2$};
        
        \draw node at (0,-1) {$n=3.$};
        
        \draw[->] (-0.9,0) to (0.9,0);
        \draw[->] (-1,0.1) to (1.7, 0.6);
        \fill[white] (0.876,0.448) circle (1.5pt);
        \draw[->] (-1,0.2) to (-0.2,1.85);
        \draw[->] (1.35,0.15) to (1.8,0.45);
        \draw[->] (1,0.2) to (0.2,1.8);
        \draw[->] (1.7,0.7) to (0.3, 1.9);
        \end{scope}
    \end{tikzpicture}
    \label{figure: simplices of the nerve of a category}
\end{figure}

For $n = 2$, the datum of such a diagram is actually equivalent to the datum of the two arrows $f_{0,1}$ and $f_{1,2}$, since the arrow $f_{0,2}$ is equal to the composition $f_{1,2} \circ f_{0,1}$, because of the commutativity condition. More generally, the datum of an $n$-simplex of $N(\Cl)$, with arrows $f_{i,j} : X_i \rightarrow X_j$ for $i \le j$, is equivalent to the datum of the arrows $f_{i,i+1}$ from $i=0$ to $i = n-1$. Indeed, because of the commutativity condition, the other arrows are obtained by the formula

$$
f_{i,j} = f_{j-1,j} \circ \hdots \circ f_{i,i+1}.
$$

The datum of an $n$-simplex of $N(\Cl)$ is therefore equivalent to the datum of a \textit{string of} $n$ \textit{composable morphisms} of $\Cl$

$$
X_0 \overset{f_1}{\longrightarrow} X_1 \overset{f_2}{\longrightarrow} \hdots \overset{f_n}{\longrightarrow} X_n.
$$

The face and degeneracy operators of $N(\Cl)$ are described as follows. Consider two integers $n \ge 0$ and $0 \le i \le n$. If $0 < i < n$, the face operator $d_i^n$ corresponds to composing $f_{i-1,1}$ with $f_{i,i+1}$ and leaving the other morphisms unchanged. If $i=0$ or $i=n$, it corresponds to deleting $f_i$ and leaving the other morphisms unchanged. The degeneracy operator $s_i^n$ corresponds to inserting the identity morphism of $X_i$ in position $i$.

The nerve functor, defined from the category of small categories to the category of simplicial sets, is fully faithful (\cite[Proposition 1.3.3.1]{kerodon}). In other words, it allows to see ordinary categories as special kinds of simplicial sets.

Given an arbitrary simplicial set $S$, we can think of its $0$-simplices as objects, its $1$-simplices as morphisms, its $2$-simplices as witnessing compositions of two morphisms and more generally, its $n$-simplices as some kind of coherent compositions of $n$ composable morphisms. Compositions in $S$ may then not exist, and not be unique in general. A simplicial set is isomorphic to the nerve of an ordinary category precisely when compositions exist and are unique, in the sense that for every $n \ge 0$, the map

$$
\Hom_{\sSet}(\Delta^n,S) \rightarrow \Hom_{\sSet}(\Spine(\Delta^n),S)
$$

is a bijection (\cite[Corollary 1.5.7.9]{kerodon}). This condition is also equivalent to the condition that for every $n \ge 0$ and every $0 < k < n$, the map

$$
\Hom_{\sSet}(\Delta^n,S) \rightarrow \Hom_{\sSet}(\Lambda_k^n,S)
$$

is a bijection (\cite[Proposition 1.3.4.1]{kerodon}).

An $\infty$-\textit{category} is a simplicial set in which compositions exist and are unique \textit{up to a contractible space of choices}. We will give the precise definition in section \ref{section: infinity categories: definitions}. The next section is devoted to the classical homotopy theory of spaces, both in terms of topological spaces and in terms of simplicial sets. It will in particular allow to define precisely what we mean by a “space of choices.”

\subsubsection{The homotopy category of spaces}\label{section: homotopy category of spaces}

In this section, we discuss the classical homotopy category of spaces, first in terms of topological spaces and then in terms of simplicial sets. We start with the following definition.

\begin{Definition}\label{definition: weak homotopy equivalence}

    Let $X, Y$ be two topological spaces and $f : X \rightarrow Y$ a continuous map. The map $f$ is called a \emph{weak homotopy equivalence} when the following conditions are satisfied.
    
    \begin{itemize}
        \item The map $\pi_0(X) \rightarrow \pi_0(Y)$ induced by $f$ is a bijection.
        
        \item For every point $x \in X$ and every integer $n \ge 1$, the morphism $\pi_n(X,x) \rightarrow \pi_n(Y,f(x))$ induced by $f$ is an isomorphism.
    \end{itemize}

\end{Definition}

This condition is in particular verified when $f$ is a \emph{homotopy equivalence}, i.e., when there exists a continuous map $g : Y \rightarrow X$ such that $f \circ g$ and $g \circ f$ are homotopic to the identity. The condition that a map is a weak homotopy equivalence is easier to check in practice.

\begin{Definition}

    The homotopy category of spaces is the localization of the category of topological spaces at the class of weak homotopy equivalences. It is denoted $\mathrm{Spaces}$.

\end{Definition}

We refer to \cite[Definition 1.2.1]{HoveyModelCategories} for the definition of the localization of a category with respect to a class of morphisms.

Denote by $\Top$ the category of topological spaces. The category $\mathrm{Spaces}$ satisfies the following universal property: there exists a functor $\Top \rightarrow \Spaces$ that sends weak homotopy equivalences to isomorphisms and such that for every category $\Cl$, every functor $F : \Top \rightarrow \Cl$ that sends weak homotopy equivalences to isomorphisms factors uniquely through $\Top \rightarrow \Spaces$.

The category $\Spaces$ can be described in concrete terms thanks to the following theorem. This theorem can be deduced from the existence of the classical model structure on topological spaces first introduced by Quillen in \cite[Chapter 2, Section 3]{QuillenHomotopicalAlebra}.

\begin{Thms}\label{theorem: homotopy theory of spaces in terms of CW complexes}

    Denote by $\CW$ the full subcategory of $\Top$ whose objects are CW-complexes. Denote by $\hCW$ the category described as follows.
    
    \begin{itemize}
        \item The objects of $\hCW$ are the CW-complexes.
        
        \item For every pair of CW-complexes $(X,Y)$, the set of morphisms of $\hCW$ between $X$ and $Y$ is the set of homotopy classes of continuous maps between $X$ and $Y$.
        
        \item The composition operation in $\hCW$ is induced from the composition operation in $\CW$ by passing to homotopy classes.
    \end{itemize}
    
    Consider the composite of functors $\CW \rightarrow \Top \rightarrow \Spaces$, where the first one is the inclusion and the second one is the localization. For every $X, Y \in \CW$, the induced map $\Hom_{\CW}(X,Y) \rightarrow \Hom_{\Spaces}(X,Y)$ is compatible with the homotopy relation. The induced functor $\hCW \rightarrow \Spaces$ is an equivalence of categories.

\end{Thms}

This generalizes a classical theorem of Whitehead asserting that a continuous map between CW-complexes is a homotopy equivalence if and only if it is a weak homotopy equivalence.

It will be useful, while reading this paper, to be familiar with the description of the category $\Spaces$ in terms of simplicial sets. We now briefly present this description. We start with the notion of homotopy between two morphisms of simplicial sets.

\begin{Definition}\label{definition: homotopy between morphisms of simplicial sets}

Let $f,g : S \rightarrow$ be two morphisms of simplicial sets. A \emph{homotopy} between $f$ and $g$ is a morphism of simplicial sets $h : S \times \Delta^1 \rightarrow T$ such that $h_{S \times \{0\}} = f$ and $h_{S \times \{1\}} = g$. If such a homotopy exists, we say that $f$ is \emph{homotopic} to $g$.

\end{Definition}

\begin{Remark}\label{remark: homotopic morphism give homotopic maps}

Suppose given a homotopy $h : S \times \Delta^1 \rightarrow T$ from $f$ to $g$. Passing to the geometric realization, we get a morphism of topological spaces $|h| : |S \times \Delta^1| \rightarrow |T|$. Moreover there is a natural map of topological spaces $| S \times \Delta^1| \rightarrow |S| \times |\Delta^1|$ that is a homeomorphism (\cite[Corollary 3.6.2.2]{kerodon}). The two maps $|f|$ and $|g|$ are therefore homotopic.

\end{Remark}

The homotopy relation between morphisms of simplicial sets is not reflexive in general. For example, there exists a homotopy between the morphisms $\Delta^0 \xrightarrow{0} \Delta^1$ and $\Delta^0 \xrightarrow{1} \Delta^1$ given by the identity morphism $\Delta^1 \rightarrow \Delta^1$, but there does not exist a homotopy in the other direction. This relation is not transitive either in general. For example, consider the three morphisms of simplicial sets $\Delta^0 \xrightarrow{0,1,2} \Lambda_1^2$. The morphisms $0$ and $1$ are homotopic, as well as the morphisms $1$ and $2$, but the morphisms $0$ and $2$ are not homotopic.

The failure of the homotopy relation between morphisms of simplicial sets to be an equivalence relation disappears when the target is assumed to belong to some class of simplicial sets called \emph{Kan complexes} which we now introduce.

\begin{Definition}\label{definition: Kan complexes}

    A \emph{Kan complex} is a simplicial set $K$ satisfying the following property: for any two integers $n$ and $k$ such that $n \ge 0$ and $0 \le k \le n$, every morphism of the form $\Lambda_k^n \rightarrow K$ can be extended to a morphism $\Delta^n \rightarrow K$. Here, the simplicial set $\Lambda_k^n$ is the $k^{\mathrm{th}}$ horn in $\Delta^n$, introduced in \ref{notation and terminology about simplicial sets}.

\end{Definition}

\begin{Lemmas}\label{lemma: homotopy relation equivalence relation when target Kan complex}

Let $S$ be a simplicial set and $K$ a Kan complex. The homotopy relation on the set $\Hom_{\sSet}(S,K)$ is an equivalence relation.

\end{Lemmas}

\begin{proof}

    See \cite[Proposition 3.1.5.4]{kerodon}. \qedhere

\end{proof}

\begin{Definition}\label{definition: homotopy category of Kan complexes}
    The \emph{homotopy category of Kan complexes}, denoted $\hKan$, is the category defined as follows.
    
    \begin{itemize}
        \item The objects of $\hKan$ are Kan complexes.
        
        \item Given two Kan complexes $K$ and $L$, the set of morphisms of $\hKan$ between $K$ and $L$ is the set of homotopy classes of morphisms of simplicial sets between $K$ and $L$. Given $f \in \Hom_{\sSet}(K,L)$, we denote by $[f]$ its homotopy class.
        
        \item Given two composable morphisms in $\hKan$ represented by two morphisms of simplicial sets $f$ and $g$, their composition is defined by $[f] \circ [g] = [f \circ g]$.
    \end{itemize}
\end{Definition}

We know explain how the category $\hKan$ is related to the category $\Spaces$. The comparison is done via the pair of adjoint functors

$$
|-| : \sSet \rightleftarrows \Top : \Sing
$$

introduced in section \ref{section: Left Kan extensions}.

The geometric realization of every simplicial set has a CW structure (\cite[Remark 1.2.3.12]{kerodon}). Moreover, by remark \ref{remark: homotopic morphism give homotopic maps}, if two morphisms of Kan complexes are homotopic then their geometric realization are homotopic as maps of topological spaces. In particular, the geometric realization functor induces a functor

$$
\mathrm{h}|-| : \hKan \rightarrow \hCW \simeq \Spaces.
$$

Consider on the other hand the singular simplicial set functor $\Sing : \Top \rightarrow \sSet$. For every topological space $X$ the simplicial set $\Sing(X)$ is a Kan complex (\cite[Proposition 1.2.5.8]{kerodon}). Suppose we are given two topological spaces $X$ and $Y$ and a homotopy $h : X \times [0,1] \rightarrow Y$ from a map $f$ to a map $g$. We get a map of Kan complexes $\Sing(X \times [0,1]) \rightarrow \Sing(Y)$. We want to deduce from this the existence of a homotopy from the map $\Sing(f)$ to the map $\Sing(g)$. To this end note that since $\Sing$ is a right adjoint it commutes with products, and therefore the natural map $\Sing(X \times [0,1]) \rightarrow \Sing(X) \times \Sing([0,1])$ is an isomorphism of simplicial sets. Furthermore by definition of the geometric realization functor one has

$$
|\Delta^1| = \{(t_0,t_1) \in \R^2 \mid t_0,t_1 \ge 0 \text{ and } t_0 + t_1 = 1\}
$$

and consequently one has a homeomorphism defined by

$$
\fonctionsansnom{|\Delta^1|}{[0,1],}{(t_0,t_1)}{t_1.}
$$

By adjunction this corresponds to a morphism of simplicial set $\Delta^1 \rightarrow \Sing([0,1])$. Now the composite

$$
\Sing(X) \times \Delta^1 \rightarrow \Sing(X) \times \Sing([0,1]) \xrightarrow{\Sing(h)} \Sing(Y)
$$

is a homotopy from $\Sing(f)$ to $\Sing(g)$. In particular the functor $\Sing$ induces a functor

$$
\mathrm{h} \Sing : \hCW \rightarrow \hKan.
$$

We have the following lemma.

\begin{Lemmas}\label{lemma: singular realization adjunction passes to homotopy cat}

    The adjunction between the functors $|-|$ and $\Sing$ induces an adjunction between the functors $\mathrm{h} |-|$ and $\mathrm{h} \Sing$.

\end{Lemmas}

\begin{proof}

    By adjunction, for every topological space $X$ and every simplicial set $S$ there is a natural identification
    
    $$\Hom_{\Top}(|S|,X) = \Hom_{\sSet}(S,\Sing(X)).$$
    
    We want to show that when $S$ is a Kan complex and $X$ is a CW-complex, this identification is compatible with the homotopy relation on both sides. We prove that this actually holds without these restrictions on $S$ and $X$. Indeed, identifying $[0,1]$ and $|\Delta^1|$ as above we have
    
    $$
    \begin{aligned}
    \Hom_{\Top}(|S| \times [0,1],X) & = \Hom_{\Top}(|S| \times |\Delta^1|,X) \\
    & =\Hom_{\Top}(|S \times \Delta^1|,X) \\
    & = \Hom_{\sSet}(S \times \Delta^1,\Sing(X))
    \end{aligned}
    $$
    
    and therefore the homotopy relation on the left-hand side coincides with the homotopy relation on the right-hand side. \qedhere

\end{proof}

\begin{Thms}\label{theorem: equivalence between topological spaces and simplicial sets}

    The functor $\mathrm{h} |-|$ is an equivalence of categories.

\end{Thms}

\begin{proof}

A proof is given in \cite[Section 3.6]{kerodon}. \qedhere

\end{proof}

\begin{Remark}

    It follows from theorem \ref{theorem: equivalence between topological spaces and simplicial sets} that the counit of the adjunction between h$|-|$ and h$\Sing$ is an equivalence, in other words for every CW-complex $X$ the map of topological spaces $|\Sing(X)| \rightarrow X$ is a homotopy equivalence. Actually, for every topological space $X$ the map $|\Sing(X)| \rightarrow X$ is a weak homotopy equivalence. This is due to Milnor in \cite{MilnorGeometricRealization} and a proof is also given in \cite[Theorem 3.6.4.1]{kerodon}. This formally implies that for every topological spaces $X,Y$ a map $f : X \rightarrow Y$ is a weak homotopy equivalence if and only if the morphism of Kan complexes $\Sing(f)$ is a homotopy equivalence. Indeed, there is a commutative diagram of topological spaces
    
    $$
    \xymatrix@C=4em{
    |\Sing(X)| \ar[r]^-{|\Sing(f)|} \ar[d] & |\Sing(Y)| \ar[d] \\
    X \ar[r]_-f & Y
    }
    $$
    
    and since the vertical maps are weak homotopy equivalences, the condition that $f$ is a weak homotopy equivalence is equivalent to the condition that $|\Sing(f)|$ is a weak homotopy equivalence. Since $|\Sing(X)|$ and $|\Sing(Y)|$ admit CW structures, this is equivalent to the condition that $|\Sing(f)|$ is a homotopy equivalence and therefore, by theorem \ref{theorem: equivalence between topological spaces and simplicial sets}, to the condition that $\Sing(f)$ is a homotopy equivalence.

\end{Remark}

\subsubsection{\texorpdfstring{$\infty$-categories: definitions}{infinity,1-categories: definitions}}\label{section: infinity categories: definitions}

The goal of this section is to define $\infty$-categories and explain how several concepts from ordinary category theory generalize to this setting.

We begin with two definitions.

\begin{Definition}\label{definition: lifting problem}
    Let $\Cl$ be a category. A \textit{lifting problem} in $\Cl$ is a commutative solid arrow diagram of the form

    $$
    \xymatrix{
    A \ar[r]^u \ar[d]^f & X \ar[d]^g \\
    B \ar[r]^v \ar@{-->}[ru]^h & Y,
    }
    $$

   i.e., the data of morphisms $u,g,f,v$ such that $g \circ u = v \circ f$. A \textit{solution to the lifting problem} is a dotted arrow $h$ making the resulting diagram commute, i.e., such that $h \circ f = u$ and $g \circ h = v$.
\end{Definition}

\begin{Definition}\label{trivial Kan fibration}
    A morphism of simplicial sets $K \rightarrow L$ is a \emph{trivial Kan fibration} if every lifting problem in $\sSet$ of the form

    $$
    \xymatrix{
    \partial \Delta^n \ar[r] \ar[d] & K \ar[d] \\
    \Delta^n \ar@{-->}[ru] \ar[r] & L
    }
    $$

    has a solution.
\end{Definition}

One can define $\infty$-categories in several equivalent ways.

\begin{Definition}\label{definition: infinity category}
    An $\infty$\textit{-category} is a simplicial set $\Cl$ that satisfies the following equivalent conditions.
    
    \begin{enumerate}[label=(\roman*)]
        \item For any two integers $n$ and $k$ such that $n \ge 0$ and $0 < k < n$, every lifting problem of the form
        
        $$
         \xymatrix{
         \Lambda_k^n \ar[d] \ar[r] & \Cl \\
        \Delta^n \ar@{-->}[ru]
        }
        $$
        
        has a solution.
        
        \item The morphism of simplicial sets
        
        $$
        \underline{\Hom}_{\sSet}(\Delta^2,\Cl) \rightarrow \underline{\Hom}_{\sSet}(\Lambda_1^2,\Cl)
        $$

        induced by restriction along the inclusion $\Lambda_1^2 \rightarrow \Delta^2$ is a trivial Kan fibration (see notation and terminology \ref{notation and terminology about simplicial sets} for the definition of $\underline{\Hom}_{\sSet}$).
        
        \item For every integer $n \ge 0$, the morphism of simplicial sets
    
        $$
        \underline{\Hom}_{\sSet}(\Delta^n,\Cl) \rightarrow \underline{\Hom}_{\sSet}(\Spine(\Delta^n),\Cl)
        $$
    
        induced by restriction along the inclusion $\Spine(\Delta^n) \rightarrow \Delta^n$ is a trivial Kan fibration.
        
    \end{enumerate}

\end{Definition}

The equivalence between conditions (i) and (ii) is due to Joyal and a proof can be found in \cite[Theorem 1.5.6.1]{kerodon}. Note that condition (iii) reduces to condition (ii) in the case $n=2$. Conversely a proof of the implication (i) $\Rightarrow$ (iii) is given in \cite[Remark 1.5.7.7]{kerodon}.

\begin{Example}\label{example: nerve of category is infinity category}

    By the discussion of section \ref{section: ordinary categories as simplicial sets}, the nerve of every (small) category is an $\infty$-category.

\end{Example}

\begin{Example}\label{example: Kan complex is infinity category}

    Every Kan complex is an $\infty$-category, since it satisfies condition (i) for any two integers $n \ge 0$ and $0 \le k \le n$.

\end{Example}

Condition (i) is called the \emph{inner horn filling property} and is in practice the easiest to check. On the other hand, condition (ii) justifies the idea that the compositions in an $\infty$-category are well-defined up to a contractible space of choices, as we now explain. Let $\Cl$ be an $\infty$-category.

\begin{Definition}\label{definition: objects and morphisms of infinity category}

    The \emph{objects} of $\Cl$ are the $0$-simplices of $\Cl$. The \emph{morphisms} of $\Cl$ are the $1$-simplices of $\Cl$. Given a morphism $\gamma$ of $\Cl$, the \emph{source} of $\gamma$ is the object $x = d_1^1(\gamma)$ and the \emph{target} of $\gamma$ is the object $y = d_0^1(\gamma)$. The morphism $\gamma$ is called a \emph{morphism from} $x$ \emph{to} $y$.

\end{Definition}

Note that when $\Cl$ is the nerve of an ordinary category, this definition coincides with the usual one.

Compositions in an $\infty$-category are defined as follows.

\begin{Definition}\label{definition: composition in infinity category}

Suppose given two morphisms $f,g$ of $\Cl$ such that the target of $f$ is equal to the source of $g$. A morphism $h$ is called a \emph{composition} of $f$ and $g$ if there exists a $2$-simplex $\sigma$ of $\Cl$ such that $d_0^2(\sigma)=g$, $d_1^2(\sigma)=h$ and $d_2^2(\sigma)=f$, as in the following picture.

\begin{figure}[H]
\centering
\begin{tikzpicture}[scale=2]
    \coordinate (A) at (-1,0);
    \coordinate (B) at (1,0);
    \coordinate (C) at (0,1.7);
    \draw (A)--(B) node[midway, below] {$f$};
    \draw (B)--(C) node[midway, above, right] {$g$};
    \draw (C)--(A) node[midway, above, left] {$h$};
    \node at (0,0.7) {$\sigma$};
\end{tikzpicture}
\end{figure}

\end{Definition}

Condition (i) in definition \ref{definition: infinity category} in the case $n=2$ and $k=1$ implies that there exists a composition of $f$ and $g$. Something stronger actually holds: a $2$-simplex $\sigma$ as in definition \ref{definition: composition in infinity category} is the same as a $0$-simplex of the fiber (in the sense of notation and terminology \ref{notation and terminology about simplicial sets}) of the morphism

$$
\underline{\Hom}_{\sSet}(\Delta^2,\Cl) \rightarrow \underline{\Hom}_{\sSet}(\Lambda_1^2,\Cl)
$$

over the horn defined by $f$ and $g$. This fiber can be regarded as parametrizing the compositions between $f$ and $g$ and by condition (ii) in definition \ref{definition: infinity category}, this is a contractible Kan complex. In other words, the compositions of $f$ and $g$ are parameterized by a contractible space.

\begin{Remark}\label{remark: contractible space of choices}

    Condition (ii) in definition \ref{definition: infinity category} is a priori stronger than the condition that the fibers of this map are contractible. This is because one does not only want the spaces parametrizing the compositions of any two fixed morphisms to be contractible. One also want coherence as the morphisms vary.

\end{Remark}

\begin{Definition}\label{definition: identity morphisms in infinity categories}

    For an object $X \in \Cl$, the morphism $s_0^0(X)$ is called the \textit{identity morphism of} $X$.

\end{Definition}

\begin{Definition}\label{definition: isomorphisms in infinity categories}

    A morphism $f$ from $X$ to $Y$ is called an \textit{isomorphism} if there exists a morphism $g$ from $Y$ to $X$ such that the identity morphism of $X$ is a composition of $f$ and $g$ and the identity morphism of $Y$ is composition of $g$ and $f$. If there exists an isomorphism between $X$ and $Y$, we say that $X$ and $Y$ are \textit{isomorphic}.

\end{Definition}

The notion of functor of categories generalizes to $\infty$-categories as follows.

\begin{Definition}\label{definition: functor of infinity categories}

    A functor between $\infty$-categories is a morphism between the corresponding simplicial sets.

\end{Definition}

\begin{Proposition}\label{proposition: morphism simplicial set is infinity category}

    If $K$ is a simplicial set and $\Cl$ an $\infty$-category, the simplicial set $\underline{\Hom}(K,\Cl)$ is an $\infty$-category.
    
\end{Proposition}

\begin{proof}

    See \cite[Proposition 1.2.7.3]{HigherTopos}. \qedhere

\end{proof}

The notion of equivalence of categories generalizes to $\infty$-categories as follows.

\begin{Definition}\label{definition: equivalence of infinity categories}

    A functor $F$ between two $\infty$-categories $\Cl$ and $\Dl$ is called an \textit{equivalence} if there exists a functor $G$ between $\Dl$ and $\Cl$ such that $F \circ G$ and $G \circ F$ are isomorphic to the identity as objects of the $\infty$-categories $\underline{\Hom}(\Dl,\Dl)$ and $\underline{\Hom}(\Cl,\Cl)$ respectively.

\end{Definition}

It is often convenient to prove that a functor between ordinary categories is an equivalence by verifying that it is fully faithful and essentially surjective. This characterization of equivalences extends to $\infty$-categories as we now discuss.

\begin{Definition}\label{definition: essential surjectivity for infinity categories}

A functor of $\infty$-categories $F : \Cl \rightarrow \Dl$ is called \textit{essentially surjective} if for every object $y$ of $\Dl$, there exists an object $x$ of $\Cl$ such that $F(x)$ and $y$ are isomorphic.

\end{Definition}

To define full faithfulness we need the following definition.

\begin{Definition}\label{definition: morphism space}
    Let $\Cl$ be a simplicial set and $x,y$ two $0$-simplices of $\Cl$. We define the simplicial set $\Hom_{\Cl}(x,y)$ to be the fiber over $(x,y)$ of the morphism

    $$
    \underline{\Hom}_{\sSet}(\Delta^1,\Cl) \xrightarrow{d_1^1 \times d_0^1} \Cl \times \Cl.
    $$
    
    When $\Cl$ is an $\infty$-category, we refer to this simplicial set as the \emph{space of morphisms of} $\Cl$ \emph{between} $x$ \emph{and} $y$.
\end{Definition}

\begin{Remark}\label{remark: space of morphisms is Kan complex}

This terminology is justified by the fact that $\Hom_{\Cl}(x,y)$ is a Kan complex when $\Cl$ is an $\infty$-category (\cite[Proposition 4.6.1.10]{kerodon}).

\end{Remark}

\begin{Definition}\label{definition: full faithfulness for infinity categories}

A functor of $\infty$-categories $F : \Cl \rightarrow \Dl$ is called \textit{fully faithful} if for every pair of objects $(x, y)$ of $\Cl$, the morphism of Kan complexes induced by $F$

$$
\Hom_{\Cl}(x,y) \rightarrow \Hom_{\Dl}(F(x),F(y))
$$

is a homotopy equivalence. 

\end{Definition}

We then have the following theorem.

\begin{Thms}\label{theorem: full faithfulness and essential surjectivity}
    A functor between $\infty$-categories is an equivalence if and only if it is fully faithful and essentially surjective.
\end{Thms}

\begin{proof}
    See \cite[Theorem 4.6.2.21]{kerodon}.
\end{proof}

\begin{Remark}\label{remark: homotopy category}

    To every $\infty$-category $\Cl$ is associated an ordinary category $\mathrm{h}\Cl$ called the \emph{homotopy category} of $\Cl$. Its objects are the objects of $\Cl$ and given two objects $x,y$ of $\Cl$, the set of morphism of $\mathrm{h} \Cl$ between these is the set $\pi_0(\Hom_{\Dl}(x,y))$, which denotes the set of homotopy classes of maps $\Delta^0 \rightarrow \Hom_{\Dl}(x,y)$ (this set is well-defined since $\Hom_{\Dl}(x,y)$ is a Kan complex). The notion of composition in $\Cl$ from definition \ref{definition: composition in infinity category} induces a well-defined composition operation in $\mathrm{h} \Cl$. The isomorphisms of $\Cl$ are precisely those that yield isomorphisms in the homotopy category. For a discussion of the homotopy category construction, we refer to \cite[Section 1.4.5]{kerodon}.

\end{Remark}

Until now we have defined in this section the notion of $\infty$-category, functor between $\infty$-categories and equivalence of $\infty$-categories. The $\infty$-categories are actually the objects of a category in which equivalences become isomorphisms.

\begin{Definition}\label{definition: category of infinity categories}

    The \emph{category of} $\infty$-\emph{categories}, denoted $\inftyCat$, is the category described as follows.
    
    \begin{itemize}
        \item The objects of $\inftyCat$ are the $\infty$-categories.
        
        \item Given two $\infty$-categories $\Cl$ and $\Dl$, the set of morphisms of $\inftyCat$ between $\Cl$ and $\Dl$ is the set of functors of $\infty$-categories between $\Cl$ and $\Dl$ modded out by the relation of being \emph{naturally equivalent}, meaning isomorphic as objects of the $\infty$-category $\underline{\Hom} (\Cl,\Dl)$ (this is an equivalence relation according to remark \ref{remark: homotopy category}).
    \end{itemize}

\end{Definition}

\begin{Remark}\label{remark: category of infinity categories is homotopy category}

    Denote by $\sSet_{\inftyCat}$ the full subcategory of $\sSet$ whose objects are the $\infty$-categories. Note that the natural functor $\sSet_{\inftyCat} \rightarrow \inftyCat$ sends equivalences of $\infty$-categories to isomorphisms. By considering the Joyal model structure on $\sSet$, one can prove that this functor exhibits $\inftyCat$ as the localization of $\sSet_{\inftyCat}$ with respect to the class of equivalences of $\infty$-categories. For a discussion of the Joyal model structure, we refer the reader to \cite[Section 2.2.5]{HigherTopos}.

\end{Remark}

\begin{Remark}

    In definition \ref{definition: infinity category} we defined an $\infty$-category to be a simplicial set satisfying certain conditions. However, as definition \ref{definition: category of infinity categories} suggests, two $\infty$-categories should be considered "the same" not only when they are isomorphic as simplicial sets, but more generally when they are equivalent, i.e., isomorphic as objects of $\inftyCat$. Strictly speaking, simplicial sets satisfying the equivalent conditions of definition~\ref{definition: infinity category} should therefore be regarded as models for $\infty$-categories. Nevertheless, throughout this paper we will refer to them simply as "$\infty$-categories".

\end{Remark}

\subsubsection{\texorpdfstring{Topological categories and $\infty$-categories}{Topological categories and infinity categories}}\label{section: topological categories and infinity categories}

We saw in the previous section that for every two objects $x,y$ of an $\infty$-category, there is a Kan complex of morphisms between them. On the other hand, we presented in section \ref{section: homotopy category of spaces} a description of the homotopy category of spaces in terms of Kan complexes. This leads to a description of the theory of $\infty$-categories in terms of \emph{topological categories} instead of simplicial sets, which we briefly review in this section.

\begin{Definition}\label{definition: topological category}
    A \textit{topological category} is a category enriched over the category of compactly generated weakly Hausdorff topological spaces. We denote the category of topological categories by $\Cat_{\Top}$.
\end{Definition}

A topological category is thus a category such that the set of morphisms between any two given objects has a structure of compactly generated weakly Hausdorff topological space, and for every triple of objects $X, Y, Z$, the composition map $\Hom(X,Y) \times \Hom(Y,Z) \rightarrow \Hom(X,Z)$ is continuous. A functor between topological categories is a functor such that the associated maps between morphism topological spaces are continuous.

\begin{Remark}\label{remark: product in the category of compactly generated spaces}
    The product $\Hom(X,Y) \times \Hom(Y,Z)$ above has to be taken in the category of compactly generated topological spaces. It has the same underlying set as the usual product in the category of topological spaces, but not the same topology in general. The reason why we work with compactly generated topological spaces is because we need the geometric realization functor 
    
    $$
    |-| : \sSet \rightarrow \text{Topological spaces}
    $$
    
    to commute with products of simplicial sets and of topological spaces, in order to define the geometric realization of a simplicial category (in section \ref{section: definition of the homotopy coherent nerve}). Apart from this technical point, this subtlety can be ignored.
\end{Remark}

In order to explain how topological categories are related to $\infty$-categories, we will need an appropriate notion of isomorphism in a topological category, which is weaker than the usual one. In order to introduce it, we first need the following definition.

\begin{Definition}\label{definition: homotopy in topological category}
    Let $\Cj$ be a topological category, $X$, $Y$ two objects of $\Cj$ and $f,g \in \Hom_{\Cj}(X,Y)$. A \emph{homotopy} from $f$ to $g$ is a continuous map $h : [0,1] \rightarrow \Hom_{\Cj}(X,Y)$ such that $h(0) = f$ and $h(1) = g$. If such a homotopy exists, we say that $f$ and $g$ are \emph{homotopic}.
\end{Definition}

\begin{Definition}\label{definition: isomorphism in topological category}

    Let $\Cj$ be a topological category and $f : X \rightarrow Y$ a morphism in $\Cl$. The morphism $f$ is called an \emph{isomorphism} when there exists a morphism $g : Y \rightarrow X$ such that $f \circ g$ is homotopic to the identity of $Y$ and $g \circ f$ is homotopic to the identity of $X$.

\end{Definition}

In order to compare the theory of topological categories and the theory of $\infty$-categories, the appropriate notion of equivalence between topological categories is the following (compare with theorem \ref{theorem: full faithfulness and essential surjectivity}).

\begin{Definition}\label{definition: equivalence of topological categories}

    A functor between topological categories $F : \Cj \rightarrow \Dj$ is called an \emph{equivalence} when it satisfies the following two conditions.
    
    \begin{itemize}
        \item Essential surjectivity: for every object $Y \in \Dj$, there exists an object $X \in \Cj$ such that $F(X)$ and $Y$ are isomorphic (in the sense of definition \ref{definition: isomorphism in topological category}).
        
        \item Full faithfulness: for every pair of objects $(X,Y)$ of $\Cj$, the map of topological spaces $F(X,Y) : \Hom_{\Cj}(X,Y) \rightarrow \Hom_{\Dj}(F(X),F(Y))$ is a weak homotopy equivalence.
    \end{itemize}

\end{Definition}

Contrary to the notion of equivalence of $\infty$-categories introduced in definition \ref{definition: equivalence of infinity categories}, an equivalence of topological categories does not always admit an inverse, even up to homotopy. Two topological categories should therefore be considered equivalent when they are connected by a zigzag of equivalences. More generally, one can consider the localization of the category of topological categories at the class of equivalences. It turns out that the resulting category is equivalent to $\inftyCat$, as we now explain.

The category of topological categories is related to the category of simplicial sets by means of a functor called the \textit{homotopy coherent nerve functor}

$$
\Nl : \Cat_{\Top} \rightarrow \sSet.
$$

We will study this functor in detail in section \ref{section: the homotopy coherent nerve}. We have the following theorem.

\begin{Thms}\label{theorem: equivalence between infinity categories topological categories}
    The functor $\Nl$ sends every topological category to an $\infty$-category, and sends equivalences between topological categories to equivalences between $\infty$-categories. The induced functor $\Cat_{\Top} \rightarrow \inftyCat$ exhibits $\inftyCat$ as the localization of $\Cat_{\Top}$ with respect to the class of equivalences.
\end{Thms}

\begin{proof}
    A proof is given in \cite[Theorem 1.1.5.13]{HigherTopos}. Another proof is given in \cite{JoyalEquivalence}. \qedhere
\end{proof}

Let $\Cj$ be a topological category. On the one hand, for every pair of objects $(X,Y)$ of $\Cj$ we have a topological space of morphisms $\Hom_{\Cj}(X,Y)$ and therefore a Kan complex $\Sing(\Hom_{\Cj}(X,Y))$. On the other hand, according to theorem \ref{theorem: equivalence between infinity categories topological categories}, there is an $\infty$-category $\Nl(\Cj)$ associated to $\Cl$. We thus have, for every pair of objects $(x,y)$ of $\Nl(\Cj)$, the Kan complex $\Hom_{\Nl(\Cj)}(x,y)$ defined in \ref{definition: morphism space}. We can wonder if these two notions of morphism spaces coincide. This is indeed the case by virtue of the following proposition.

\begin{Proposition}\label{proposition: morphism spaces of a topological category and its coherent nerve}
There is a natural bijection between the sets of objects of the topological category $\Cj$ and the $\infty$-category $\Nl(\Cj)$. Moreover, for every pair of objects $(X,Y)$ of $\Cl$, the Kan complexes $\Sing(\Hom_{\Cj}(X,Y))$ and $\Hom_{\Nl(\Cj)}(X,Y)$ are naturally homotopy equivalent.
\end{Proposition}

\begin{proof}
    The natural bijection between the sets of objects follows from the construction of the functor $\Nl$ (see section \ref{section: the homotopy coherent nerve}). The natural homotopy equivalence of morphism spaces is discussed in \cite[Section 2.2]{HigherTopos}. We will come back to it in section \ref{section: equivalence between flow and homotopy coherent nerve} (see in particular remark \ref{remark: homotopy equivalence of morphism spaces btw cat and nerve}). \qedhere
\end{proof}

\begin{Remark}\label{remark: homotopy category of topological category versus quasicategory}

    Given a topological category $\Cj$, there are two associated constructions of a "homotopy category". One can either consider the category $\pi_0(\Cl)$ with the same set of objects as $\Cj$ and $\Hom_{\pi_0(\Cj)}(x,y) = \pi_0(\Hom_{\Cj}(x,y))$ for every pair of objects $x,y$, or one can consider the homotopy category of the $\infty$-category $\Nl(\Cj)$. By proposition \ref{proposition: morphism spaces of a topological category and its coherent nerve} these two categories are naturally isomorphic. A proof of this can also be found in \cite[Proposition 2.4.6.9]{kerodon}.

\end{Remark}

\subsection{Stratified topological spaces, exit paths and stratified simplicial sets}\label{section: stratified spaces, exit paths and stratified simplicial sets}

This section is devoted to the notion of \emph{stratification} and its homotopy-theoretic study. In section \ref{section: stratified topological spaces and exit paths} we introduce the notions of stratified topological space and exit path. In section \ref{section: exit path infty category of conically stratified space} we define conically stratified spaces, we introduce the exit path $\infty$-category associated with a conically stratified space, and we give a description of the morphism spaces of this $\infty$-category in terms of spaces of exit paths. In section \ref{section: stratified simplicial sets} we introduce the analogues of stratified topological spaces in the category of simplicial sets, called stratified simplicial sets. Finally, in section \ref{section: geometric description of spaces of exit paths}, we introduce cylindrically stratified topological spaces and we give for these an alternative description of the spaces of exit paths, up to homotopy.

\subsubsection{Stratified topological spaces and exit paths}\label{section: stratified topological spaces and exit paths}

\begin{Remark}\label{remark on the definition of stratified space}
    The notion of stratified topological space that we review in this section is that of a \textit{poset-stratified space}. Let us mention that the particular stratification that interests us in this paper falls into the more restricted realm of \textit{Whitney stratified spaces}, which are defined by extra geometric conditions. This will be the subject of sections \ref{section: Whitney conditions and cylindrical structures} and \ref{section: Whitney stratifications} (see in particular \ref{definition: Whitney stratified space} for the definition). On the other hand, the framework of poset-stratified spaces is the most general one and has the advantage of being well-suited for the study of stratified spaces from a homotopy theoretic perspective, as it permits the development of a parallel theory of \textit{stratified simplicial sets}, which we will review in section \ref{section: stratified simplicial sets}. \\[0pt]
    For an overview of the various notions of stratified topological spaces and relations between them, see \cite{WaasWoolfYokuraStratifications}.
\end{Remark}

Let $A$ be a poset.

\begin{Definition}\label{definition: sub posets of greater or lower elements}
    Given $a \in A$, we let $A_{\ge a}$ (resp. $A_{>a}$, $A_{\le a}$, $A_{<a}$) denote the sub-poset of $A$ consisting of those elements $b \in A$ such that $b \ge a$ (resp. $b >a$, $b \le a$, $b < a$).
\end{Definition}

\begin{Definition}\label{definition: Alexandrov topology}
    The \emph{Alexandrov topology} on $A$ is the topology generated by the open subsets $A_{\ge a}$.
\end{Definition}

\begin{Remark}\label{remark: open subsets Alexandrov topology}

A subset $U \subseteq A$ is open for the Alexandrov topology if and only if it is \emph{closed upwards}, meaning $a \in U$ and $b \ge a$ implies $b \in U$.

\end{Remark}

\begin{Definition}\label{definition: stratified space}
    Let $Y$ be a topological space. An $A$-\emph{stratification} on $Y$ is a continuous map $\pi : Y \rightarrow A$, where $A$ is endowed with the Alexandrov topology. For every $a \in A$, we denote $Y_a = \pi^{-1}(\{a\})$ and refer to it as the $a$-\emph{stratum} of $Y$. We will often abuse notation and denote the $A$-stratified space $\pi : Y \rightarrow A$ by $(Y,A)$ without specifying the stratifying map.
\end{Definition}

\begin{Remark}\label{remark: reformulation of definition of stratification}
    Said differently, a stratification of $Y$ is a decomposition of $Y$ as a disjoint union of strata
    
    $$
    Y = \bigsqcup_{a \in A} Y_a
    $$
    
    such that for every $a \in A$, the union of the strata $Y_b$ for $b \ge a$ is open in $Y$.
\end{Remark}

\begin{Example}\label{example: 0 < 1 stratified space}

    Every $(0 < 1)$-stratification of $Y$ determines a closed subset $Y_0 \subseteq Y$ and conversely, every closed subset $F \subseteq Y$ determines a $(0 < 1)$-stratification of $Y$ by the conditions that $Y_0 = F$ and $Y_1 = Y \backslash F$. This defines a bijection between the set of $(0 < 1)$-stratifications of $Y$, and the set of closed subsets of $Y$.

\end{Example}

\begin{Notation}\label{notation: category of A-stratified topological spaces}
    We denote by $\Top_A$ the overcategory $\Top_{/A}$.
\end{Notation}

The category $\Top_A$ is thus the category whose objects are the $A$-stratified topological spaces and whose morphisms are the stratum-preserving continuous maps between these.

\begin{Proposition}\label{proposition: the category of A-stratified spaces has limits and colimits}
    The category $\Top_A$ admits all small limits and small colimits. Moreover, the forgetful functor $\Top_A \rightarrow \Top$ commutes with colimits.
\end{Proposition}

\begin{proof}
    To prove that $\Top_A$ admits all small limits, we prove that it admits all equalizers and all small products. To prove that it admits all equalizers, consider two arrows in $\Top_A$

    $$
    \xymatrix{
    X \ar@<0.5ex>[rr]^{f} \ar@<-0.5ex>[rr]_{g} \ar[rd] && Z \ar[ld] \\
    & A.
    }
    $$

    The two arrows $f$ and $g$ admit an equalizer in $\Top$, given by the subset $\{x \in X \mid f(x)=g(x)\} \subseteq X$ endowed with the topology induced by the topology of $X$. This equalizer inherits an $A$-stratification from the one of $X$. The resulting $A$-stratified space is an equalizer of the above two arrows in $\Top_A$.

    To prove that $\Top_A$ admits all small products, consider a family of $A$-stratified spaces $(X_j,A)_{j \in J}$ for some set $J$. The fiber product in $\Top$ of the family of topological spaces $(X_j)_{j \in J}$ over $A$ is naturally endowed with an $A$-stratification. The resulting $A$-stratified space is the product of the family $(X_j,A)_{j \in J}$ in $\Top_A$.

    To prove that $\Top_A$ admits all small colimits, let us consider a small category $D$ and a functor $G : D \rightarrow \Top_A$. It induces a functor $D \rightarrow \Top$ by forgetting the stratifications. This functor admits a colimit $X \in \Top$, which is endowed with a natural map $X \rightarrow A$. The latter is an $A$-stratification of $X$ which is a colimit of $G$. This also shows that the forgetful functor $\Top_A \rightarrow \Top$ commutes with colimits.
\end{proof}

Let us fix an $A$-stratification of $Y$. We denote by $I$ the interval $[0,1]$.

\begin{Definition}\label{definition: exit path}

An \textit{exit path} in $Y$ is a path $\gamma : I \rightarrow Y$ such that there exist $a, b \in A$ with $\pi(\gamma(0)) = a$ and $\pi(\gamma(t)) = b$ for all $t > 0$ (note that this implies $a \le b)$.

\end{Definition}

See figure \ref{figure: exit path in intro} for an illustration of this definition.

\begin{Remark}\label{remark: an path in a stratum is an exit path}
    Note that in the statement of definition \ref{definition: exit path}, one can have $a=b$. In other words, a path in $Y$ that stays in the same stratum at all times, is an exit path.
\end{Remark}

The purpose of the next section is to introduce a notion of composition between exit paths.

\subsubsection{\texorpdfstring{The exit path $\infty$-category of a conically stratified space}{The exit path infty-category of a conically stratified space}}\label{section: exit path infty category of conically stratified space}
    
The naive way to compose two paths by concatenating them, does not yield a notion of composition of exit paths. Indeed, suppose given $a,b,c \in A$ and two exit paths $\gamma, \delta$ from $Y_a$ to $Y_b$ and from $Y_b$ to $Y_c$ respectively, such that $\gamma(1) = \delta(0)$. The concatenation of $\gamma$ and $\delta$ looks as follows.

\begin{center}

\begin{tikzpicture}[scale=1.2]
        \draw node at (0,-0.3) {$\gamma$};
        \draw[green] node at (-2,0) {$Y_a$};
        \draw node at (0.3,1.1) {$\delta$};
        \draw[red] node at (1,-0.4) {$Y_b$};
        \draw[blue] node at (1.8,1) {$Y_c$};
        
        \draw[->, bend right]  (-1.3,0.05) to (-1.8,0);
        \draw[->, bend right]  (0.5,-0.1) to (0.8,-0.3);
        \draw[->, bend left]  (0.8,0.9) to (1.6,1);
        
        \draw[blue, line width=1.2pt] (0,2)--(1.2,0);
        
        \draw[red, line width=1.2pt] (-1.2,0)--(1.2,0);
        
        \fill[green] (-1.2,0) circle (1.5pt);
        
    \end{tikzpicture}

\end{center}

As soon as $a,b,c$ are distinct, this concatenation visits three distinct strata and is therefore not an exit path.

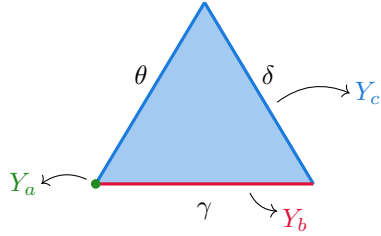
\begin{figure}[H]
    \centering
    \begin{tikzpicture}[scale=1.2]
        \draw node at (0,-0.3) {$\gamma$};
        \draw[green] node at (-2,0) {$Y_a$};
        \draw node at (0.7,1.2) {$\delta$};
        \draw[red] node at (1,-0.4) {$Y_b$};
        \draw node at (-0.7,1.2) {$\theta$};
        \draw[blue] node at (1.8,1) {$Y_c$};
        
        \draw[->, bend right]  (-1.3,0.05) to (-1.8,0);
        \draw[->, bend right]  (0.5,-0.1) to (0.8,-0.3);
        \draw[->, bend left]  (0.8,0.9) to (1.6,1);
        
        \draw[fill, blue, opacity=0.4] (-1.2,0)--(1.2,0)--(0,2)--cycle;
        
        \draw[blue, line width=1.2pt] (0,2)--(1.2,0);
        \draw[blue, line width=1.2pt] (0,2)--(-1.2,0);
        
        \draw[red, line width=1.2pt] (-1.2,0)--(1.2,0);
        
        \fill[green] (-1.2,0) circle (1.5pt);
        
    \end{tikzpicture}
    \caption{The composition law for exit paths.}
    \label{figure: composition of exit paths}
\end{figure}

The appropriate notion of composition of two exit paths is the following: we define a \emph{composition} of $\gamma$ and $\delta$ to be an exit path $\theta$ such that there exists a $2$-simplex $\sigma : \Delta^2 \rightarrow Y$ satisfying $d_0^2(\sigma) = \delta$, $d_1^2(\sigma) = \theta$, $d_2^2(\sigma)=\gamma$, and the restriction of $\sigma$ to the interior of $\Delta^2$ has image contained in $Y_c$. This is illustrated by figure \ref{figure: composition of exit paths}.

One can think of this as a notion of composition up to homotopy, with the homotopy being compatible with the stratification.

The question then arises as to what extent compositions exist and are unique. As explained in section \ref{section: infinity categories}, $\infty$-categories provide a framework to adress this question. This leads to Lurie's notion of exit path $\infty$-category, which we now introduce.

In order to give the definition of the exit path $\infty$-category, we first need to extend the definition of an exit path to higher simplices.

Recall that we denote by $\Pos$ the category of posets, and that we have the simplex category $\Delta$, which is a full subcategory of $\Pos$.

Let us fix a poset $A$.

\begin{Definition}\label{definition: finite increasing sequence}
    A \emph{finite increasing sequence} of elements of $A$ is a morphism of posets $[n] \rightarrow A$ for some integer $n \ge 0$.
\end{Definition}

\begin{Notation}\label{notation: finite increasing sequence}
    We will denote finite increasing sequences of elements of $A$ in lowercase bold letters. For every finite increasing sequence $\ag : [n] \rightarrow A$, we will denote $\ag(i) = a_i$ and also use the notation $[a_0 \le \hdots \le a_n]$ to denote the sequence $\ag$.
\end{Notation}

\begin{Definition}\label{definition: simplex category of a poset}
    The \emph{simplex category of} $A$, denoted $\Delta_A$, is the full subcategory of the overcategory $\Pos_{/ A}$ whose objects are the finite increasing sequences of elements of $A$. In other words, this is the category obtained as the pullback
    
    $$
    \squarediagram{\Delta_A}{\Pos_{/A}}{\Delta}{\Pos}{}{}{}{}
    $$
    
    where the right vertical arrow is the forgetful functor and the bottom horizontal one is the inclusion. Put another way, this is the category whose objects are the morphisms of posets $[n]\rightarrow A$ for $n$ ranging through nonnegative integers, and whose morphisms are the commutative diagrams of posets of the form

    $$
    \xymatrix{
    [n] \ar[rr] \ar[rd] && [m] \ar[ld] \\
    & A.
    }
    $$
\end{Definition}

\begin{Notation}\label{notation: standard stratified simplices}
    Given a finite increasing sequence $\ag$ of elements of $A$, we denote by $\Delta^{\ag}$ the corresponding object of $\Delta_A$.
\end{Notation}

\begin{Definition}\label{definition: length of a sequence}
    Given an object $\Delta^{\ag} \in \Delta_A$, it is mapped to some object $\Delta^n$ by the forgetful functor $\Delta_A \rightarrow \Delta$. We refer to $n$ as the \emph{length} of $\ag$ and denote it by $l(\ag)$.
\end{Definition}

\begin{Definition}\label{definition: stratified geometric realization}

An object $\Delta^{\ag} \in \Delta_A$ defines the following $A$-stratification of the geometric realization $|\Delta^n|$. For every $a \in A$, the $a$-stratum is defined to be

$$
|\Delta^n|_a = \bigcup_{0 \le i \le n, \, a_i = a} \{ (t_0,\hdots,t_n) \mid t_i \neq 0 \text{ and } t_j = 0 \text{ for all } j>i \}.
$$

This $A$-statified space is denoted $|\Delta^{\ag}|_A$.
\end{Definition}

\begin{Example}\label{example: stratifications of standard simplices}

    Examples of stratifications of standard simplices are given by figure \ref{figure: stratification of the standard 1-simplex} for $n=1$, by figure \ref{figure: stratification of the standard 2-simplex} for $n=2$, and by figure \ref{figure: some stratifications of the 3 simplex} for $n=3$. 

\end{Example}

\begin{Proposition}\label{proposition: stratified geometric realization is a functor}
    For every arrow $\Phi : \Delta^{\ag} \rightarrow \Delta^{\bg}$ in $\Delta_A$, denoting by $\varphi : \Delta^n \rightarrow \Delta^m$ the underlying arrow in $\Delta$, the induced map $|\varphi| : |\Delta^{\ag}|_A \longrightarrow |\Delta^{\bg}|_A$ is compatible with the stratifications. \qed
\end{Proposition}

\begin{Definition}\label{definition: stratified geometric realization at the level of simplex category}
    The functor

    $$
    \fonction{|-|_A}{\Delta_A}{\Top_A}{
    \left\{
    \begin{array}{l}
         \Delta^{\ag} \\
         \Delta^{\ag} \overset{\Phi}{\longrightarrow} \Delta^{\bg} 
    \end{array}
    \right.
    }{
     \left\{
    \begin{array}{l}
         |\Delta^{\ag}|_A \\
         |\Delta^{\ag}|_A \overset{|\varphi|}{\longrightarrow} |\Delta^{\bg}|_A 
    \end{array}
    \right.
    }
    $$

    is called the \textit{stratified geometric realization}. Here, again, $\varphi$ denotes the arrow in $\Delta$ underlying the arrow $\Phi$ in $\Delta_A$.
\end{Definition}

\begin{Definition}\label{definition: stratified singular simplicial set}

The singular simplicial set of $(Y,A)$ is the simplicial subset of $\Sing(Y)$ defined as

$$
\Sing_A(Y)_n = \{\sigma : |\Delta^{\ag}|_A \rightarrow Y \text{ compatible with stratifications for some } \Delta^{\ag} \in \Delta_A \}.
$$

\end{Definition}

\begin{Remark}\label{1-simplices of stratified singular simplicial set are exit paths}
    The $1$-simplices of $\Sing_A(Y)$ are the exit paths of $Y$.
\end{Remark}

\begin{Remark}\label{remark: what it means for stratified singular simplicial set to be an infinity category}
    Assume that $\Sing_A(Y)$ is an $\infty$-category. By the discussion of section \ref{section: infinity categories: definitions}, this means that the notion of composition of exit paths introduced at the beginning of this section, is well-defined up to a contractible space of choices.
\end{Remark}

It is the case that $\Sing_A(Y)$ is an $\infty$-category for a class of stratified topological spaces, called \emph{conically stratified spaces}, which we now introduce, following \cite[Section A.5]{HigherAlgebra}. We first need the following definition.

\begin{Definition}\label{definition: open cone on a stratified space}
    Let $A$ be a poset and let $A^{\triangleleft}$ be the poset obtained by adjoining a new smallest element $- \infty$ to $A$. Let $\pi : l \rightarrow A$ be an $A$-stratified space. We define the \emph{open cone on} $l$, denoted $C(l)$, to be the following $A^{\triangleleft}$-stratified space:
    
    \begin{enumerate}
        \item As a set, $C(l)$ is given by the union $* \cup (l \times \R_{>0})$.
        
        \item A subset $U \subset C(l)$ is open if and only if $U \cap (l \times \R_{>0})$ is open and, if $* \in U$, then $l \times ]0,\varepsilon[ \subset U$ for some positive real number $\varepsilon$.
        
        \item The $A^{\triangleleft}$-stratification on $C(l)$ is determined by the map $C(l) \rightarrow A^{\triangleleft}$ sending $*$ to $-\infty$ and $(x,t)$ to $\pi(x)$ for every $(x,t) \in l \times \R_{>0}$. 
    \end{enumerate}
\end{Definition}

\begin{Remark}\label{remark: topology on the open cone}
    The topology on the open cone on $l$ from definition \ref{definition: open cone on a stratified space} is called the \emph{teardrop topology}. It is in general coarser than the topology on the pushout $(l \times \R_{\ge 0})\bigsqcup_{l \times \{0\}}\{*\}$, but the two topologies coincide as long as $l$ is compact and Hausdorff. This will be the case in all the examples that we will consider.
\end{Remark}

The definition of a conically stratified space is the following.

\begin{Definition}\label{definition: conically stratified space}
    Let $A$ be a poset, let $Y$ be an $A$-stratified space and let $y$ be a point of $Y$. We let $a$ be the point of $A$ such that $y \in Y_a$. We say that $Y$ is \emph{conically stratified} at the point $y$ if there exists a $A_{>a}$-stratified space $l$, a topological space $V$, and an open embedding $C(l) \times V \hookrightarrow Y$ of $A$-stratified spaces whose image contains $y$. Here, we regard $C(l)$ as endowed with an $A_{\ge a}$ stratification by identifying $A_{>a}^{\triangleleft}$ with $A_{\ge a}$, and we regard $C(l) \times V$ as endowed with the $A_{\ge a}$-stratification defined as the composition of the $A_{\ge a}$-stratification on $C(l)$ with the projection $C(l) \times V \rightarrow C(l)$. Such an open embedding is called a \emph{local cone-like description} of the stratification near $y$. \\[0pt]
    We say that $Y$ is \emph{conically stratified} if it is conically stratified at every point $y \in Y$.
\end{Definition}

\begin{Example}\label{example: example of Morse stratification is conical}
    Let us illustrate definition \ref{definition: conically stratified space} with the example of the stratification by stable manifolds from the introduction.
    
    \begin{figure}[H]
        \centering
        \includegraphics[width=0.25\linewidth]{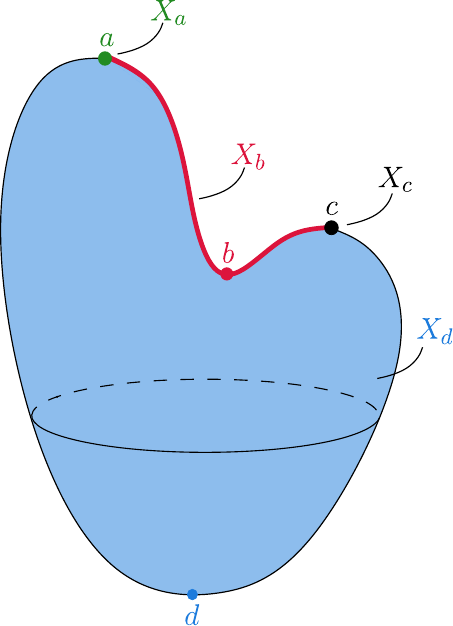}
    \end{figure}
    
    This is an example of a conical stratification. To see the conical structure at $a$, we can consider the following spaces.
    
    \begin{figure}[H]
        \centering
            \begin{tikzpicture}
            \begin{scope}[xshift=-4cm]
            \draw node at (-0.4,0) {$l=$};
            
            
            \draw[blue] (0,0) arc [start angle=-180, end angle=0, x radius=1.2, y radius=0.3];
            \draw[blue] (0,0) arc [start angle=180, end angle=0, x radius=1.2, y radius=0.3];
            
            
            \fill[red] (2.4,0) circle (1.5pt);
           \end{scope}
           
           \begin{scope}[xshift=0]
            \draw node at (0,0) {$V = *$};
            \end{scope}
            
            \begin{scope}[xshift=4.6cm]
            \draw node at (-1.6,0) {$C(l) \times V= C(l)=$};
            
            
            \draw[fill, blue, opacity=0.4] (0,0)--(1.2,1.8)-- (2.4,0)-- (2.4,0) arc [start angle=0, end angle=-180, x radius=1.2, y radius=0.3]--cycle; 
            \draw[blue, dashed, opacity=0.3] (0,0) arc [start angle=180, end angle=0, x radius=1.2, y radius=0.3];
            \draw[blue] (0,0)--(1.2,1.8);
            
            
            \draw[red, line width=1.4] (2.4,0)--(1.2,1.8);
            
            
            \fill[green] (1.2,1.8) circle (1.7pt);
            \end{scope}
            \end{tikzpicture}
    \end{figure}
    
Similarly, to see the conical structure at $c$, we can consider the following spaces.

    \begin{figure}[H]
        \centering
            \begin{tikzpicture}
            \begin{scope}[xshift=-4cm]
            \draw node at (-0.4,0) {$l=$};
            
            
            \draw[blue] (0,0) arc [start angle=-180, end angle=0, x radius=1.2, y radius=0.3];
            \draw[blue] (0,0) arc [start angle=180, end angle=0, x radius=1.2, y radius=0.3];
            
            
            \fill[red] (0,0) circle (1.5pt);
           \end{scope}
           
           \begin{scope}[xshift=0]
            \draw node at (0,0) {$V = *$};
            \end{scope}
            
            \begin{scope}[xshift=4.6cm]
            \draw node at (-1.6,0) {$C(l) \times V=C(l)=$};
            
            
            \draw[fill, blue, opacity=0.4] (0,0)--(1.2,1.8)-- (2.4,0)-- (2.4,0) arc [start angle=0, end angle=-180, x radius=1.2, y radius=0.3]--cycle; 
            \draw[blue, dashed, opacity=0.3] (0,0) arc [start angle=180, end angle=0, x radius=1.2, y radius=0.3];
            \draw[blue] (2.4,0)--(1.2,1.8);
            
            
            \draw[red, line width=1.4] (0,0)--(1.2,1.8);
            
            
            \fill (1.2,1.8) circle (1.7pt);
            \end{scope}
            \end{tikzpicture}
    \end{figure}
    
    To see the conical structure at any point of the $b$-stratum, we can consider the following spaces.
    
    \begin{figure}[H]
        \centering
            \begin{tikzpicture}
            \begin{scope}[xshift=-6cm]
            \draw node at (-0.5,0.05) {$l=$};
            
            
            \fill[blue] (0,0) circle (1.5pt);
            \fill[blue] (1,0) circle (1.5pt);
    
           \end{scope}
           
           \begin{scope}[xshift=-2.8cm]
                \node at (-0.7,0.05) {$C(l)=$};
                
                \draw[blue] (0,0)--(0.5,0.8);
                \draw[blue] (1,0)--(0.5,0.8);
                \fill[red] (0.5,0.8) circle (1.5pt);
           \end{scope}
           
           \begin{scope}[xshift=0]
            \draw node at (0,0.05) {$V = ]0,1[$};
            \end{scope}
            
            \begin{scope}[xshift=4cm]
            \draw node at (-1.3,0) {$C(l) \times V =$};
            
            \pgfmathsetmacro{\L}{9/80};
            \draw[dashed, blue, opacity=0.4] (0.6-\L,0.2)--(1.7-\L,0.2)--(2,0.7);
            \draw[fill, blue, opacity=0.4] (0.3,0.7)--(0,0.2)--(0.6-\L,0.2)--(0.6,-0.1)--(2.3,-0.1)--(2,0.7)--cycle;
             
            \draw[red, line width=1.4] (0.3,0.7)--(2,0.7);
            
            \end{scope}
            \end{tikzpicture}
    \end{figure}
    
    Finally, to see the conical structure at any point of the $d$-stratum, we can consider the following spaces.
    
    \begin{figure}[H]
        \centering
            \begin{tikzpicture}
            \begin{scope}[xshift=-5.5cm]
            \draw node at (-0.5,0.05) {$l=\varnothing$};
    
           \end{scope}
           
           \begin{scope}[xshift=-2.8cm]
                \node at (-0.8,0) {$C(l)=$};
                
                \fill[blue] (0,0) circle (1.5pt);
           \end{scope}
           
            \begin{scope}[xshift=0]
                \draw node at (-1.3,0) {$V=$};
                
                \draw[dashed, opacity=0.3] (0,0) circle (0.8cm);
                \draw[fill, opacity=0.1] (0,0) circle (0.8cm);
                
            \end{scope}
            
            \begin{scope}[xshift=4.5cm]
            \draw node at (-1.8,0) {$C(l) \times V =$};
            
            \draw[dashed, blue, opacity=0.7] (0,0) circle (0.8cm);
            \draw[fill, blue, opacity=0.4] (0,0) circle (0.8cm);
             
            \end{scope}
            \end{tikzpicture}
    \end{figure}
    
\end{Example}

\begin{Remark}\label{remark: we will see that every Morse stratification is conical}
    In section \ref{section: spaces of exit paths as spaces of pseudo gradient trajectories}, we will generalize this example by providing, for every point in a Whitney stratified space, a stratified space $l$ and a space $V$ giving a local cone-like description of the stratification near that point (see remark \ref{remark: cone like description of Whitney stratification}).
\end{Remark}

\begin{Thms}[Lurie]\label{theorem: stratified simplicial set of conically stratified space is infty cat}
    If $(Y,A)$ is conically stratified then $\Sing_A(Y)$ is an $\infty$-category.
\end{Thms}

\begin{proof}
    See \cite[Theorem A.6.4]{HigherAlgebra}.
\end{proof}

\begin{Definition}\label{definition: exit path infinity category}
    When $(Y,A)$ is conically stratified, the $\infty$-category $\Sing_A(Y)$ is called the \emph{$\infty$-category of exit paths in $Y$ with respect to the stratification $Y \rightarrow A$}, or the \emph{exit path} $\infty$-\textit{category of} $(Y,A)$.
\end{Definition}

\subsubsection{Stratified simplicial sets}\label{section: stratified simplicial sets}

Suppose we have an $A$-stratified space $Y$ and an $n$-simplex $\sigma$ of $\Sing_A(Y)$. This comes with a continuous map $|\Delta^n| \rightarrow Y$. Recall that $|\Delta^n|$ is the convex hull of the canonical basis $(e_0,\hdots,e_n)$ of $\R^{n+1}$. For every integer $0 \le i \le n$, let us define the element $a_i \in A$ by the condition that $\sigma(e_i) \in Y_{a_i}$. For every integer $0 \le i < n$, there is a corresponding morphism $\Delta^1 \rightarrow \Delta^n$ carrying $0$ to $i$ and $1$ to $i+1$. This allows to associate to $\sigma$ and $i$ a $1$-simplex of $\Sing_A(Y)$, and the latter is an exit path from $\sigma(e_i)$ to $\sigma(e_{i+1})$; we deduce that $a_i \le a_{i+1}$. We have thus a finite increasing sequence $[a_0 \le \hdots \le a_n]$ of elements of $A$. This observation leads to an analogue of the notion of stratified topological space, in the category of simplicial sets, which we now introduce.

\begin{Definition}\label{definition: nerve of a poset}
    The \emph{nerve} of $A$, denoted $N(A)$, is the nerve of the following category:
    
    $$
    \left\{
    \begin{array}{l}
        \mathrm{Objects} = \text{Elements of } A,  \\
        \Hom(a,b) = * \text{ if } a \le b \text{ and } \Hom(a,b) = \varnothing \text{ otherwise.} 
    \end{array}
    \right.
    $$
\end{Definition}

\begin{Remark}\label{remark: finite sequences as simplices of the nerve}
    We can alternatively regard the finite increasing sequences of elements of $A$ of length $n$ as the $n$-simplices of the nerve of $A$.
\end{Remark}

Assigning to a simplex of $\Sing_A(Y)$ the corresponding finite increasing sequence of elements of $A$ (as explained at the beginning of this section) determines a morphism of simplicial sets from $\Sing_A(Y)$ to the nerve of $A$. This is a particular case of the following definition.

\begin{Definition}
    The category of \textit{simplicial sets stratified by} $A$ is the overcategory $\sSet_{/N(A)}$. It is denoted $\sSet_A$.
\end{Definition}

\begin{Proposition}\label{proposition: stratified simplicial sets admit small limits and colimits}
    The category $\sSet_A$ admits all small limits and small colimits. Moreover, the forgetful functor $\sSet_A \rightarrow \sSet$ commutes with colimits.
\end{Proposition}

\begin{proof}
    The proof is similar to that of proposition \ref{proposition: the category of A-stratified spaces has limits and colimits}, replacing $\Top$ by $\sSet$ and $\Top_A$ by $\sSet_A$.
\end{proof}

\begin{Definition}\label{definition: stratum of a stratified simplicial set}
    Suppose we are given $T \in \sSet_A$ and $a \in A$. The $a$-\textit{stratum} of $T$, denoted $T_a$, is the simplicial set obtained as the fiber over $a$ of the morphism $T \rightarrow N(A)$. In other words, it is the pullback in $\sSet$

    $$
    \squarediagram{T_a}{T}{\Delta^0}{N(A).}{}{}{}{a}
    $$
\end{Definition}

\begin{Example}\label{example: stratification of Delta^n associated to an element of Delta_A}
    Let $n$ be a nonnegative integer. By definition, the finite increasing sequences of elements of $A$ of length $n$ are exactly the $n$-simplices of $N(A)$, or equivalently, the morphisms of simplicial sets from $\Delta^n$ to $N(A)$. The finite increasing sequences of elements of $A$ of length $n$ are therefore in one-to-one correspondence with the $A$-stratifications of the standard $n$-simplex $\Delta^n \in \sSet$.
\end{Example}

\begin{Remark}\label{remark: stratification determined by vertices}
    Example \ref{example: stratification of Delta^n associated to an element of Delta_A} implies that for every nonnegative integer $n$, an $A$-stratification of $\Delta^n$ is entirely determined by the images of the vertices of $\Delta^n$ in $A$. Consequently, the same statement holds replacing $\Delta^n$ by any simplicial set.
\end{Remark}

\begin{Remark}\label{remark: stratified geometric realization as Kan extension}
    Carrying a finite increasing sequence $\ag$ of elements of $A$ of length $n$, to the corresponding $A$-stratification of $\Delta^n$ (according to example \ref{example: stratifications of standard simplices}) determines a functor

    $$
    \Delta_A \rightarrow \sSet_A.
    $$

    Since $\sSet_A$ admits all small colimits (proposition \ref{proposition: stratified simplicial sets admit small limits and colimits}), this functor induces, by the constructions \ref{construction: left Kan extension} and \ref{construction: right adjoint of left Kan extension}, a pair of adjoint functors

    $$
    \Pre(\Delta_A) \rightleftarrows \sSet_A.
    $$

    These functors can be explicitly described as follows. Every presheaf $F \in \Pre(\Delta_A)$ is sent to the $A$-stratified simplicial set whose underlying simplicial set is defined as
    
    $$
    \Delta^n \mapsto \bigcup_{\ag \text{ of length }n} F(\Delta^{\ag}),
    $$
    
    and whose structural morphism to $N(A)$ carries an $n$-simplex $\sigma \in F(\Delta^{\ag})$ to the sequence $\ag$, regarded as an $n$-simplex of $N(A)$.
    
    Conversely, every $A$-stratified simplicial set $T$ is sent to the presheaf on $\Delta_A$ defined as
    
    $$
    \begin{array}{ll}
    \Delta^{\ag} \text{ of length } n \mapsto & 
        \text{Set of those morphisms of simplicial sets } \sigma : \Delta^n \rightarrow T \text{ such that } \\
        & \text{the induced stratification on } \Delta^n \text{ is that determined by } \ag. 
    \end{array}
    $$
\end{Remark}

\begin{Proposition}\label{proposition: stratified simplicial set equivalent to presheaves}
    These two functors are inverse to one another.
\end{Proposition}

\begin{proof}
    See \cite[Proposition 3.1.5]{DouteauThese}.
\end{proof}

These two functors thus identify the category of simplicial sets stratified by $A$ with the category of presheaves over $\Delta_A$. We will identify these two categories from now on.

\begin{Notation}\label{notation: stratifications of simplex, boundary and horns}
    Let $\ag$ be a finite increasing sequence of elements of $A$ of length $n$.
    
    \begin{itemize}
        \item Regarding $\sSet_A$ as the category of presheaves $\Pre(\Delta_A)$, the object $\Delta^{\ag} \in \Delta_A$ represents an object of $\sSet_A$, by the Yoneda embedding. We denote the latter by $\Delta^{\ag}$ as well. Regarding $\sSet_A$ as the overcategory $\sSet_{/ N(A)}$, the object $\Delta^{\ag} \in \sSet_A$ corresponds to the stratification of the standard $n$-simplex $\Delta^n$ determined by $\ag$.

        \item We let $\partial \Delta^{\ag}$ denote the simplicial set $\partial \Delta^n$ endowed with the restriction of the $A$-stratification $\Delta^{\ag}$ on $\Delta^n$.

        \item For every $0 \le k \le n$ we let $\Lambda_k^{\ag}$ denote the simplicial set $\Lambda_k^n$ endowed with the restriction of the $A$-stratification $\Delta^{\ag}$ on $\Delta^n$.
    \end{itemize}
\end{Notation}

\begin{Definition}\label{definition: stratified geometric realization at the level of stratified simplicial sets}
    Since the category $\Top_A$ admits all small colimits (proposition \ref{proposition: the category of A-stratified spaces has limits and colimits}), we can consider the left Kan extension of the functor $|-|_A : \Delta_A \rightarrow \Top_A$ to $\sSet_A$. It is called the \textit{stratified geometric realization functor} and still denoted $|-|_A$. This functor associates an $A$-stratified topological space to every $A$-stratified simplicial set.
\end{Definition}

\begin{Remark}\label{remark: explicit description of stratified realization}
    Recall from construction \ref{construction: left Kan extension} that the functor $|-|_A : \sSet_A \rightarrow \Top_A$ is determined (up to isomorphism) by its restriction to $\Delta_A$ and by the condition that it commutes with colimits. Let us also describe more explicitly the stratified geometric realization of a given $A$-stratified simplicial set $S$. For every $n \ge 0$ we denote by $S_n$ the set of $n$-simplices of $S$. For every $\Delta^{\ag} \in \Delta_A$, denoting by $n$ the length of $\ag$, we denote by $S_{\ag}$ the subset of $S_n$ described in the following two equivalent ways.
    
    \begin{itemize}
    
    \item Regarding $\sSet_A$ as the overcategory $\sSet_{/ A}$, $S_{\ag}$ is the set of those $n$-simplices $\Delta^n \rightarrow S$ such that the induced stratification on $\Delta^n$ is that determined by $\ag$.
    
    \item Regarding $\sSet_A$ as the category of presheaves $\Pre(\Delta_{A})$, $S_{\ag}$ is the image of $\Delta^{\ag}$ by $S$.
    
    \end{itemize}
    
    The underlying topological space of the stratified geometric realization of $S$ is the geometric realization of the underlying simplicial set of $S$. Recall from section \ref{section: Left Kan extensions} that the latter is obtained as a quotient of the disjoint union
    
    $$
    \bigsqcup_{n \ge 0} S_n \times |\Delta^n| = \bigsqcup_{\Delta^{\ag} \in \Delta_A} S_{\ag} \times |\Delta^{l(\ag)}|.
    $$
    
    The functor $|-|_A : \Delta_A \rightarrow \Top_A$ endows each of the $|\Delta^{l(\ag)}|$ with an $A$-stratification. These stratifications determine an $A$-stratification of the disjoint union. Unwinding the definitions and constructions, one sees that this stratification passes to the quotient, and that the stratified geometric realization of $S$ is obtained in this way.
\end{Remark}

\begin{Remark}\label{remark: stratified sing is right adjoint}

By construction \ref{construction: right adjoint of left Kan extension} we know that the functor $|-|_A : \sSet_A \rightarrow \Top_A$ admits a right adjoint. Unwinding the definitions and constructions, one sees that this right adjoint is the functor $\Sing_A(-)$ of definition \ref{definition: stratified singular simplicial set}.

\end{Remark}

Let $Y$ be an $A$-stratified topological space and $x, y$ two points of $Y$. We finish this section by studying the simplicial set $\Hom_{\Sing_A(Y)}(x,y)$ (definition \ref{definition: morphism space}). Before stating the result, we introduce the following definition.

\begin{Definition}\label{definition: spaces of exit paths}

\begin{itemize}
    \item The \emph{space of exit paths in $(Y,A)$ between} $x$ \emph{and} $y$ is the set of exit paths in $(Y,A)$ starting at $x$ and ending at $y$, endowed with the compact-open topology. It is denoted $\Exit(x,y)$.
    
    \item Given two elements $a,b \in A$, the \emph{space of exit paths in $(Y,A)$ between} $Y_a$ \emph{and} $Y_b$ is the set of exit paths in $(Y,A)$ starting in $Y_a$ and ending in $Y_b$, endowed with the compact-open topology. It is denoted $\Exit(Y_a,Y_b)$.
\end{itemize}

\end{Definition}

\begin{Proposition}\label{proposition: space of morphism is space of exit paths}
    There is a natural identification of simplicial sets
    
    $$
    \Hom_{\Sing_A(Y)}(x,y) = \Sing(\Exit(x,y)).
    $$
\end{Proposition}

In particular, when the stratification of $Y$ by $A$ is conical, the exit path $\infty$-category of $(Y,A)$ is an $\infty$-category whose points are the points of $Y$, and whose spaces of morphisms are the spaces of exit paths between the corresponding points.

\begin{proof}[Proof of proposition \ref{proposition: space of morphism is space of exit paths}]
Denote by $a, b$ the respective images of $x,y$ in $A$. Assume that $\Hom_{\Sing_A(Y)}(x,y)$ is nonempty. Then $x$ and $y$ are connected by an exit path, so $a \le b$. Let us endow $\Delta^1$ and $|\Delta^1|$ with the corresponding stratifications by $A$. For every integer $n \ge 0$ we also endow $\Delta^1 \times \Delta^n$ and $|\Delta^1| \times |\Delta^n|$ with the stratifications induced by the first projection. Consider an $n$-simplex of $\Hom_{\Sing_A(Y)}(x,y)$

$$
\sigma : \Delta^1 \times \Delta^n \rightarrow \Sing_A(Y).
$$

The map $\sigma$ is compatible with the arrows to $N(A)$ on both sides. Indeed, by remark \ref{remark: stratification determined by vertices} it suffices to verify this compatibility on vertices, where it follows from the fact that $\sigma$ is constant equal to $x$ in restriction to $0 \times \Delta^n$ and equal to $y$ in restriction to $1 \times \Delta^n$.

By adjunction we have

$$
\begin{aligned}
    \Hom_{\sSet_A}(\Delta^1 \times \Delta^n,\Sing_A(Y)) & = \Hom_{\Top_A}(|\Delta^1 \times \Delta^n|,Y) \\
    & = \Hom_{\Top_A}((|\Delta^1| \times |\Delta^n|,Y) \\
    & = \Hom_{\Top}(|\Delta^n|,\underline{\Hom}_{\Top_A}(|\Delta|^1,Y)),
\end{aligned}
$$

where $\underline{\Hom}_{\Top_A}$ denotes the subset of stratified maps endowed with the topology induced by the compact-open topology.

The above identification restricts to an identification of sets

$$
\Hom_{\Sing_A(Y)}(x,y)_n = \Hom_{\Top}(|\Delta^n|,\Exit(x,y)).
$$

These identifications are functorial in $n$, so that we finally get an identification of simplicial sets

\begin{equation*}
\Hom_{\Sing_A(Y)}(x,y) = \Sing(\Exit(x,y)).\qedhere
\end{equation*}

\end{proof}

\begin{Remark}\label{remark: homotopy category of exit path infinity category}

    By theorem \ref{theorem: stratified simplicial set of conically stratified space is infty cat} one can associate to every conically stratified topological space $(Y,A)$ an ordinary category defined to be the homotopy category of the $\infty$-category $\Sing_A(Y)$. The objects of this category are the points of $X$ and, by proposition \ref{proposition: space of morphism is space of exit paths}, the set of morphisms between every pair of points $(x,y)$ is the set $\pi_0(\Exit(x,y))$. The composition law is defined as in figure \ref{figure: composition of exit paths}; it is in particular a consequence of theorem \ref{theorem: stratified simplicial set of conically stratified space is infty cat} that this composition law is well-defined.

\end{Remark}

We close this section with the following proposition, which will be useful to us later.

\begin{Proposition}\label{proposition: Sing_A infinity category in terms of starting point evaluation}
    Let $Y$ be a metrizable $A$-stratified space. Then the following two conditions are equivalent:
    
    \begin{enumerate}[label=(\roman*)]
        \item the simplicial set $\Sing_A(Y)$ is an $\infty$-category,
        
        \item for every pair of elements $a < b$ of $A$, the starting point evaluation map $\Exit(Y_a,Y_b) \rightarrow Y_a$ is a Serre fibration.
    \end{enumerate}
\end{Proposition}

\begin{proof}

This a particular case of \cite[Proposition 5.1.4]{WaasPresenting}. Note that, in that paper, "categorically fibrant" is equivalent to condition $(i)$ (see \cite[Definition 5.1.2]{WaasPresenting}). \qedhere

\end{proof}

\subsubsection{A geometric description of the spaces of exit paths in cylindrical stratifications}\label{section: geometric description of spaces of exit paths}

Since the exit path $\infty$-category of a conically stratified space has spaces of exit paths as spaces of morphisms (proposition \ref{proposition: space of morphism is space of exit paths}), we are interested in spaces of exit paths up to weak homotopy equivalence. The main result of this section is proposition \ref{proposition: space of exit paths as boundary of regular neighborhood}, which provides an alternative description of the space of exit paths between two strata, up to weak homotopy equivalence, in what is called a \emph{cylindrical stratification}. This is based on \cite[Construction 2.41]{MaderWaasSamples}, and we give details for arguments that are only sketched there.

Except at the very end, this section is only concerned with topological spaces stratified by the poset $0 < 1$. Recall from example \ref{example: 0 < 1 stratified space} that such a datum is equivalent to the datum of a topological space and a closed subset.

We will work with the following definition of mapping cylinder.

\begin{Definition}\label{definition: mapping cylinder teardrop topology}

Let $r : L \rightarrow Y$ be a continuous map between topological spaces. The \emph{mapping cylinder} of $r$, denoted $\Cyl(r)$, is the set $Y \, \cup \, L \times (0,1]$ endowed with the minimal topology such that the following two conditions hold:

\begin{enumerate}[label=(\roman*)]
    \item the inclusion $L \times (0,1] \rightarrow \Cyl(r)$ is an open embedding,
    
    \item the map $R : \Cyl(r) \rightarrow Y \times [0,1]$ obtained as the union of the two maps
    
    $$
    \fonctionsansnom{Y}{Y \times [0,1]}{y}{(y,0)}
    $$
    
    and
    
    $$
    \fonctionsansnom{L \times (0,1]}{Y \times [0,1]}{(l,t)}{(r(l),t),}
    $$
    
    is continuous.

\end{enumerate}

This topology can be explicitly described as follows: a subset $U \subset \Cyl(r)$ is open if and only if $U \, \cap \, L \times (0,1]$ is open and, for every $y \in Y \cap U$, $U$ contains a set of the form $V \cup r^{-1}(V) \times (0,\varepsilon)$ for some neighborhood $V$ of $y$ in $Y$ and $0 < \varepsilon \le 1$.

\end{Definition}

\begin{Remark}\label{remark: teardrop topology generalizes open cone}
    Compare with definition \ref{definition: open cone on a stratified space}, which is the same except that $Y = *$ and $(0,1]$ is replaced by $\R_{>0}$.
\end{Remark}

\begin{Definition}\label{definition: stratified mapping cylinder}

    The \emph{stratified mapping cylinder} of $r$ is the space $\Cyl(r)$ endowed with the stratification by $0 < 1$ defined as $\Cyl(r)_0 = Y$. 

\end{Definition}

\begin{Remark}\label{remark: teardrop topology}
    The topology on the mapping cylinder from definition \ref{definition: mapping cylinder teardrop topology} is called the \emph{teardrop topology}. It is contained in the topology on the pushout in $\Top$ of the diagram
    
    $$
    \xymatrix{
    L \times \{0\} \ar[r]^r \ar[d] & Y, \\ 
    L \times [0,1]
    }
    $$
    
    but has in general fewer open sets. The two topologies coincide as long as $r$ is proper and $L, Y$ are locally compact and Hausdorff. The teardrop topology on $\Cyl(r)$ has the following essential feature: for any subspace $L' \subseteq L$, the restriction of the teardrop topology on $\Cyl(r)$ to the subset $\Cyl(r_{\mid L'})$ coincides with the teardrop topology. Thus, for example, taking $L=[0,1]$, $L' = (0,1)$ and $Y=*$, the topology on $\Cyl(r)$ coincides with the pushout $[0,1] \times [0,1] \sqcup_{[0,1] \times \{0\}} *$, but the induced topology on the subset $\Cyl(r')$ is the teardrop topology, and one can check that it is different from the pushout $(0,1) \times [0,1] \sqcup_{(0,1) \times \{0\}} *$.\\[-0.1em]
    We will use the teardrop topology in an essential way in the proof of lemma \ref{lemma: exit paths in mapping cylinder as pullback} below, where we will describe the space of exit paths $\Exit(\Cyl(r)_0,\Cyl(r)_1)$ as a pullback. This description leads to an elegant verification of the continuity of many maps constructed in the sequel.
\end{Remark}

Our next goal is to better understand the topological space $\Exit(\Cyl(r)_0,\Cyl(r)_1)$. For this, we consider the following commutative square of topological spaces
    
$$
    \xymatrix{
    \Exit(\Cyl(r)_0,\Cyl(r)_1) \ar[r]^-{\nu_1} \ar[d]_-{\nu_2} & \left( L \times (0,1] \right)^{(0,1]} \ar[d]^-{\nu_3} \\
    \Exit(Y \times \{0\},Y \times (0,1]) \ar[r]_-{\nu_4} & \left( Y \times (0,1] \right)^{(0,1]}.
    }
$$
    
In this commutative square, the space $Y \times [0,1]$ is endowed with the stratification by $0 < 1$ defined as $(Y \times [0,1])_0 = Y \times \{0\}$. The maps $\nu_1$ and $\nu_4$ are obtained by restriction to $(0,1]$, the map $\nu_2$ is obtained by postcomposition with the map $R : \Cyl(r) \rightarrow Y \times [0,1]$ from definition \ref{definition: mapping cylinder teardrop topology}, and the map $\nu_3$ is obtained by postcomposition with $r \times \id_{(0,1]}$. The following lemma provides a useful description of $\Exit(\Cyl(r)_0,\Cyl(r)_1)$.

\begin{Lemmas}\label{lemma: exit paths in mapping cylinder as pullback}
    This commutative square is a pullback square in the category of topological spaces.
\end{Lemmas}

\begin{proof}

    To simplify the notation, we denote $\Cyl(r)$ by $C$. We denote by $P$ the pullback of the diagram, explicitly described as
    
    $$
    P = \Exit(Y \times \{0\}, Y \times (0,1]) \times_{(Y \times (0,1])^{(0,1]}} (L \times (0,1])^{(0,1]}
    $$
    
    with the topology induced by the product topology on

    $$
    \Exit(Y \times \{0\}, Y \times (0,1]) \times (L \times (0,1])^{(0,1]}.
    $$
    
    Given a compact subset $K$ of $[0,1]$ (resp. $[0,1]$, $(0,1]$), and an open subset $V$ of $C$ (resp. $Y \times [0,1]$, $L \times (0,1]$), we will denote by $\Ul(K,V)$ (resp. $\Vl(K,V)$, $\Wl(K,V)$) the open subset of $\Exit(C_0,C_1)$ (resp. $\Exit(Y \times \{0\},Y \times (0,1])$, $(L \times (0,1])^{(0,1]}$) consisting of those $\gamma$ such that $\gamma(K) \subset V$. Finally, in our notation, we will not distinguish between $\Vl(K,V) \times \Wl(K',V')$ and its intersection with $P$.
    
    By the universal property of $P$, there is a continuous map $\varphi : \Exit(C_0,C_1) \rightarrow P$; our goal is to show that this is a homeomorphism. To this end, we will prove that the following expression defines the inverse of $\varphi$
    
    $$
    \fonction{\psi}{P}{\Exit(C_0,C_1)}{(\gamma_1,\gamma_2)}{
    t \mapsto \left\{
    \begin{array}{ll}
        \gamma_1(t) & \text{ if } t=0,  \\
        \gamma_2(t) & \text{ if } t > 0. 
    \end{array}
    \right.
    }
    $$
    
    We first prove that $\psi(\gamma_1,\gamma_2)$ is an exit path for every $(\gamma_1,\gamma_2) \in P$. Observe that a map $\gamma : [0,1] \rightarrow C$ is an exit path if and only if
    
    \begin{itemize}
    
        \item it is continuous, or equivalently, by definition of the topology on $C$ (definition \ref{definition: mapping cylinder teardrop topology}), the restriction $\gamma_{|(0,1]} : (0,1] \rightarrow L \times (0,1]$ and the composition $R \circ \gamma$ are continuous,
    
        \item $\gamma(0) \in Y$ and the restriction $\gamma_{|(0,1]}$ has image contained in $L \times (0,1]$.
    \end{itemize}
    
    These two conditions are indeed fulfilled by $\psi(\gamma_1,\gamma_2)$.
    
    Note that by construction, $\psi$ is inverse to $\varphi$ at the level of underlying sets. We conclude the proof by showing that $\psi$ itself is continuous. Let $K \subset [0,1]$ compact and $V \subset C$ open. We want to show that $\psi^{-1}(\Ul(K,V))$ is an open subset of $P$. We consider two cases.
    
    First assume that $0 \notin K$. Then $\Ul(K,V) = \Ul(K,V \cap L \times (0,1])$ and therefore $\psi^{-1}(\Ul(K,V)) = \Exit(Y \times \{0\},Y \times (0,1]) \times \Wl(K,V \cap L \times (0,1])$ is an open subset of $P$, as desired.
    
    Assume now that $0 \in K$ and $\psi^{-1}(\Ul(K,V))$ is not empty. Let $\gamma \in \psi^{-1}(\Ul(K,V))$. Then $\gamma(0) \in V$. By definition of the topology on $C$, there exists $\varepsilon > 0$ and a neighborhood $U$ of $\gamma(0)$ in $Y$ such that $U \cup r^{-1}(U) \times (0,\varepsilon) \subset V$. Moreover, by continuity of $\gamma$, there exists $\eta > 0$ such that for every $0 \le t \le \eta$, $\gamma(t) \in U \cup r^{-1}(U) \times (0,\varepsilon)$. Then the set $\Vl(K \cap [0,\eta],U \times [0,\varepsilon)) \times \Wl(K \cap [\eta,1],V \cap L \times (0,1]))$ is an open neighborhood of $\gamma$ in $\psi^{-1}(\Ul(K,V))$. Hence $\psi^{-1}(\Ul(K,V))$ is an open subset of $P$, as desired. \qedhere
    
\end{proof}

\begin{Notation}\label{notation: exit paths in mapping cylinder}

We consider the map

    $$
    \fonction{e_r}{L}{\Exit(\Cyl(r)_0,\Cyl(r)_1)}{l}{t \mapsto 
    \left\{
    \begin{array}{ll}
    r(l) \in Y & \text{if } t = 0, \\
    (l,t) \in L \times (0,1] & \text{if } t > 0.
    \end{array}
    \right.
    }
    $$
    
\end{Notation}

\begin{Lemmas}\label{lemma: link to exit paths map is continuous}
    The map $e_r$ is continuous.
\end{Lemmas}

\begin{proof}

    We use the description of $\Exit(\Cyl(r)_0),\Cyl(r)_1)$ provided by lemma \ref{lemma: exit paths in mapping cylinder as pullback}. By universal property of $\Exit(\Cyl(r)_0,\Cyl(r)_1)$, the map $e_r$ coincides with the continuous map induced by the commutative square of topological spaces
    
    $$
    \xymatrix{
    L \ar[r]^-{f_1} \ar[d]_-{f_2} & \left( L \times (0,1] \right)^{(0,1]} \ar[d]^-{\nu_3} \\
    \Exit(Y \times \{0\},Y \times (0,1]) \ar[r]_-{\nu_4} & \left( Y \times (0,1] \right)^{(0,1]}
    }
    $$
    
    where $f_1$ is the continuous map defined by
    
    $$
    l \mapsto (t \mapsto (l,t))
    $$
    
    and $f_2$ is the continuous map defined by
    
    \begin{equation*}
    l \mapsto (t \mapsto (r(l),t)). \qedhere
    \end{equation*}
    
\end{proof}

Now let $Z$ be a topological space stratified by $0 < 1$. We say that $Z$ is \textit{cylindrically stratified} if $Z_0$ has a neighborhood in $Z$ that is homeomorphic to a stratified mapping cylinder. More precisely:

\begin{Definition}\label{definition: cylindrically stratified space}
    Let $Z$ be a topological space stratified by $0<1$. We say that $Z$ is \textit{cylindrically stratified} if there exists a neighborhood $N$ of $Z_0$ in $Z$, a topological space $L$ and a continuous map $L \overset{r}{\rightarrow} Z_0$ such that there exists a homeomorphism between $N$ and $\Cyl(r)$ which is the identity on $Z_0$ and preserves the stratifications.
\end{Definition}

\begin{figure}[H]
    \centering
    \begin{tikzpicture}
    
        \coordinate (G) at (-3,0);
        \coordinate (D) at (3,0);
        
        
        
        \coordinate (A1) at (-3.5,2.5);
        
        \path (G) -- (A1) coordinate[pos=0.75] (B1);
        \draw (G)--(B1) coordinate[pos=0.5] (M1);
        \draw[dashed, opacity=0.3] (B1)--(A1);
        
        
        \coordinate (A2) at (-2.5,2);
        
        \path (G) -- (A2) coordinate[pos=0.75] (B2);
        \draw (G)--(B2) coordinate[pos=0.5] (M2);
        \draw[dashed, opacity=0.3] (B2)--(A2);

        
        \coordinate (A3) at (2.5,2.5);
        
        \path (D) -- (A3) coordinate[pos=0.75] (B3);
        \draw[densely dashed] (D)--(B3) coordinate[pos=0.5] (M3);
        \draw[dashed, opacity=0.3] (B3)--(A3);
        
        
        \coordinate (A4) at (3.5,2);
        
        \path (D) -- (A4) coordinate[pos=0.75] (B4);
        \draw (D)--(B4) coordinate[pos=0.5] (M4);
        \draw[dashed, opacity=0.3] (B4)--(A4);
        
        \draw[opacity=0.5] node at (3,2.5) {$Z_1$};
        
        
        \draw[green] (M1)--(M2)--(M4);
        \draw[green, dashed] (M1)--(M3)--(M4);
        \draw[fill, green, opacity=0.6] (M1)--(M2)--(M4)--(M3)--cycle;
        \draw[green] node at (0,1.2) {$L$};
        
        
        \fill[pattern=dots, opacity=0.5] (G)--(M1)--(M2)--(M4)--(D)--cycle;
        \fill[white] (0,0.4) circle (6pt);
        \draw node at (0,0.4) {$N$};
        
        
        \draw[red, line width=1.2] (G)--(D);
        \draw[red] node at (3.3,0) {$Z_0$};
        
        
        \draw[fill, opacity=0.1] (A1)--(G)--(D)--(A4)--(A3)--cycle;
        \draw[dashed, opacity=0.3] (A1)--(A2)--(A4)--(A3)--cycle;
        
        
        \path (G)--(D) coordinate[pos=0.7] (L1);
        \path (M2)--(M4) coordinate[pos=0.7] (L2);
        \path (M1)--(M3) coordinate[pos=0.7] (L3);
        \path (L2)--(L3) coordinate[pos=0.5] (L4);

        \node[inner sep=0pt] at (L4) {\tiny{$\times$}};
        \node[right] at (L4) {$l$};
        \node[below] at (L1) {$r(l)$};

        \draw[densely dashed] (L1) -- (L4);

    \end{tikzpicture}
    \caption{An example of cylindrical stratification. To see the cylindrical structure, one can take $r$ to be the restriction to $L$ of the orthogonal projection to $Z_0$. A homeomorphism between $N$ and $\Cyl(r)$ fulfilling the conditions of definition \ref{definition: cylindrically stratified space} is obtained by identifying the segment $\{l\} \times [0,1] \subset L \times [0,1]$ with the segment in $Z_1$ connecting $r(l)$ to $l$.}
    \label{figure: cylindrical stratification}
\end{figure}
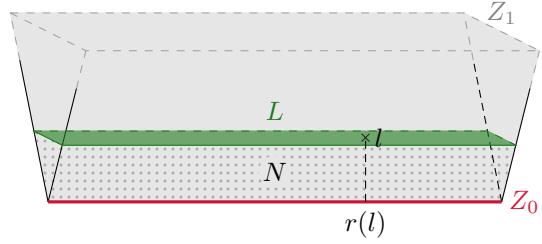

\begin{Notation}\label{notation: map from link to exit paths}

    Assume $Z$ is cylindrically stratified and fix a homeomorphism $H : N \xrightarrow{\simeq} \Cyl(r)$ as in definition \ref{definition: cylindrically stratified space}. Combining $H$ with the inclusion $N \subset Z$, we get a map 
    $$
    i_H : \Exit(\Cyl(r)_0,\Cyl(r)_1) \rightarrow \Exit(Z_0,Z_1).
    $$
    
    Composing $i_H$ with the map $e_r : L \rightarrow \Exit(\Cyl(r)_0,\Cyl(r)_1)$ (notation \ref{notation: exit paths in mapping cylinder}), we obtain a map 
    $$
    \alpha_H : L \rightarrow \Exit(Z_0,Z_1).
    $$
    
\end{Notation}

The goal of this section is to provide a proof of the following result.

\begin{Proposition}\label{proposition: space of exit paths as boundary of regular neighborhood}
If $Z$ is metrizable, then the map $\alpha_H$ (notation \ref{notation: map from link to exit paths}) is a weak homotopy equivalence.
\end{Proposition}

The proof will be based on the following two lemmas. The first one is essentially a reformulation of \cite[Lemma 2.4.(1)]{QuinnHomotopically}.

\begin{Lemmas}\label{lemma: exit paths in a neighborhood}
    For every neighborhood $U$ of $Z_0$ in $Z$, endow $U$ with the stratification induced by that of $Z$. The map
    
    $$
    \Exit(Z_0,U \cap Z_1) \rightarrow \Exit(Z_0,Z_1)
    $$
    
    induced by the inclusion $U \cap Z_1 \rightarrow Z_1$, is a homotopy equivalence.
\end{Lemmas}

\begin{proof}

Let $d$ be a distance defining the topology of $Z$. For every $z \in Z$ and $r \ge 0$, we denote by $B(z,r)$ the \emph{closed} ball of radius $r$ centered at $z$. Let $\delta : Z_0 \rightarrow (0,+\infty)$ be a continuous function such that for every $z \in Z_0$, $B(z,\delta(z)) \subset U$ (for example, for every $z_0 \in Z_0$ one can define $\delta(z_0)$ to be half of the distance between $Z_0$ and the complement of $U$ in $Z$). Denote by $\Exit^{\delta}(Z_0,Z_1)$ the subspace of $\Exit(Z_0,Z_1)$ consisting of those $\gamma$ such that the image of $\gamma$ is contained in $B(\gamma(0),\delta(\gamma(0)))$. We have a commutative diagram

$$
\xymatrix{
& \Exit(Z_0,U \cap Z_1) \ar[dd] \\
\Exit^{\delta}(Z_0,Z_1) \ar[ru]^{i_1} \ar[rd]_{i_2} \\
& \Exit(Z_0,Z_1).
}
$$

By \cite[Lemma 2.4.(1)]{QuinnHomotopically}, the maps $i_1$ and $i_2$ are homotopy equivalences, which allows to conclude that the vertical map is a homotopy equivalence. 
\end{proof}

\begin{Lemmas}\label{lemma: exit paths in mapping cylinder up to homotopy}
    For every map of topological spaces $r : L \rightarrow Y$, the map $e_r : L \rightarrow \Exit(\Cyl(r)_0,\Cyl(r)_1)$ (notation \ref{notation: exit paths in mapping cylinder}), is a homotopy equivalence.
\end{Lemmas}

We will give the proof of this lemma after the proof of proposition \ref{proposition: space of exit paths as boundary of regular neighborhood}.

\begin{proof}[Proof of proposition \ref{proposition: space of exit paths as boundary of regular neighborhood}]
 The map $\alpha_H$ is the composite
 
 $$
 L \xrightarrow{e_r} \Exit(\Cyl(r)_0,\Cyl(r)_1) \xrightarrow{i_H} \Exit(Z_0,Z_1).
 $$
 
 The first is a homotopy equivalence by lemma \ref{lemma: exit paths in mapping cylinder up to homotopy} and the second is a homotopy equivalence by lemma \ref{lemma: exit paths in a neighborhood}. Hence $\alpha_H$ is a homotopy equivalence.
\end{proof}

\begin{proof}[Proof of lemma \ref{lemma: exit paths in mapping cylinder up to homotopy}]
    To simplify the notation, we denote $\Cyl(r)$ by $C$. Recall that our goal is to show that the map
    
    $$
    \fonction{e_r}{L}{\Exit(C_0,C_1)}{l}{t \mapsto 
    \left\{
    \begin{array}{ll}
    r(l) \in Y & \text{if } t = 0, \\
    (l,t) \in L \times (0,1] & \text{if } t > 0
    \end{array}
    \right.
    }
    $$
    
    is a homotopy equivalence. Let us introduce the map
    
    $$
    \fonction{p_L}{L \times (0,1] \subset C}{L,}{(l,s)}{l,}
    $$
    
    as well as the maps $\rho : C \rightarrow [0,1]$ and $p_Y : C \rightarrow Y$ obtained as the compositions of the map $R : C \rightarrow Y \times [0,1]$ (from definition \ref{definition: mapping cylinder teardrop topology}) with the projections $Y \times [0,1] \rightarrow [0,1]$ and $Y \times [0,1] \rightarrow Y$ respectively.
    
    We will show that the continuous map
    
    $$
    \fonction{\beta}{\Exit(C_0,C_1)}{L}{\gamma}{p_L(\gamma(1))}
    $$
    
    is a homotopy inverse to $e_r$. First note that $\beta \circ e_r = \id$. We then aim to show that the following expression defines a homotopy between $\id$ and $e_r \circ \beta$:
    
    $$
    \fonction{H}{\Exit(C_0,C_1) \times [0,1]}{\Exit(C_0,C_1)}{(\gamma,s)}{t \mapsto 
    \left\{
    \begin{array}{ll}
    p_Y (\gamma(s)) \in Y & \text{if } t = 0, \\
    \left( p_L(\gamma(s+(1-s)t)),st+(1-s)\rho(\gamma(t)) \right) \in L \times (0,1] & \text{if } t > 0.
    \end{array}
    \right.}
    $$
    
    The challenge is to prove that this expression defines a continuous map. For this, we use the description of $\Exit(C_0,C_1)$ provided by lemma \ref{lemma: exit paths in mapping cylinder as pullback}. The continuity of $H$ follows from the fact that $H$ coincides with the map induced by the universal property of $\Exit(C_0,C_1)$ from the following commutative square of topological spaces

    $$
    \xymatrix{
    \Exit(C_0,C_1) \times [0,1] \ar[r]^-{f_1} \ar[d]_-{f_2} & \left( L \times (0,1] \right)^{(0,1]} \ar[d]^-{\nu_3} \\
    \Exit(Y \times \{0\},Y \times (0,1]) \ar[r]_-{\nu_4} & \left( Y \times (0,1] \right)^{(0,1]},
    }
    $$
    
    where $f_1$ is the continuous map defined by
    
    $$
    (\gamma,s) \mapsto (t \mapsto (p_L(\gamma(s+(1-s)t)),st+(1-s)\rho(\gamma(t)))
    $$
    
    and $f_2$ is the continuous map defined by
    
    \begin{equation*}
    (\gamma,s) \mapsto (t \mapsto (p_Y \circ R(\gamma(s+(1-s)t),st+(1-s)\rho(\gamma(t))). \qedhere
    \end{equation*}
    
\end{proof}

For later use, we record the following lemma about the map $\alpha_H$, in the case of a cylindrical stratification whose $0$-stratum is a point.

\begin{Lemmas}\label{lemma: link to exit paths map does not depend on homeomorphism}
    Let $Z$ be a topological space and $z_0 \in Z$ a closed point. Suppose that the $(0<1)$-stratification of $Z$ defined by $Z_0 = \{z_0\}$ is cylindrical. Let $N$ be a neighborhood of $z_0$ in $Z$, $L$ a topological space together with a homeomorphism $H : C(L) \xrightarrow{\simeq} N$ such that $H(*) = z_0$ (where $C(L)$ denotes the cone of $L$, that is, the cylinder of the unique map $r : L \rightarrow *$). Assume that we have another such homeomorphism $G : C(L) \xrightarrow{\simeq} N$ for which $H_{|L \times \{1\}} = G_{|L \times \{1\}}$. Then the two maps $\alpha_H, \alpha_G : L \rightarrow \Exit(Z_0,Z_1)$ are homotopic.
\end{Lemmas}

\begin{proof}

    We denote by $\varphi$ the homeomorphism $H^{-1} \circ G : C(L) \xrightarrow{\simeq} C(L)$. It has the property of being equal to the identity in restriction to $L \times \{1\}$. Consider the map $A_{\varphi}$ defined by
    
    $$
    \fonction{A_{\varphi}}{L}{\Exit(C(L)_0,C(L)_1),}{l}{t \mapsto
    \left\{
    \begin{array}{ll}
         r(l) & \text{if } t=0,  \\
         \varphi(l,t) & \text{if } t > 0. 
    \end{array}
    \right.
    }
    $$
    
    The continuity of this map follows again from an application of lemma \ref{lemma: exit paths in mapping cylinder as pullback}. Recall the map $e_r$ from notation \ref{notation: exit paths in mapping cylinder} and the map $i_H : \Exit(C(L)_0,C(L)_1) \rightarrow \Exit(Z_0,Z_1)$ from notation \ref{notation: map from link to exit paths}. Note that we have $i_H \circ A_{\varphi} = \alpha_G$. Since $i_H \circ e_r = \alpha_H$, it suffices to prove that $A_{\varphi}$ is homotopic to $e_r$, which we now aim to do.
    
    We further argue that it suffices to prove that there exists a homotopy $h$ between $\varphi$ and the identity satisfying $h_s(*) = *$ and $h_s(L \times (0,1]) \subset L \times (0,1]$ for all $0 \le s \le 1$. Indeed, the following expression then defines a homotopy between $A_{\varphi}$ and $e_r$:
    
    $$
    \fonctionsansnom{L \times [0,1]}{\Exit(C(L)_0,C(L)_1),}{(l,s)}{t \mapsto
    \left\{
    \begin{array}{ll}
        r(l) & \text{if } t=0,  \\
        h_s(l,t) & \text{if } t > 0. 
    \end{array}
    \right.
    }
    $$
    
    The continuity of this expression follows again from an application of lemma \ref{lemma: exit paths in mapping cylinder as pullback}.
    
    We finish the proof by constructing a map $h$ satisfying the desired properties. We denote by $p_L$ and $p_{(0,1]}$ the respective projections out of $L \times (0,1]$. The map $h$ is defined as
    
    $$
    \fonction{h}{C(L) \times [0,1]}{C(L),}{
    \left\{
    \begin{array}{l}
        (*,s) \in * \times [0,1] \\
        (l,t,s) \in L \times (0,1] \times [0,1]
    \end{array}
    \right.
    }{
    \left\{
    \begin{array}{l}
        * \\
        (p_L(\varphi(l,(1-s)t+s)),(1-s)p_{(0,1]}(\varphi(l,t))+st) \in L \times (0,1]. 
    \end{array}
    \right.
    }
    $$
    
    To conclude, we need to prove that $h$ is continuous; the other desired properties on $h$ are satisfied by construction. By definition of the topology on $C(L)$, $h$ is continuous if and only if the following two conditions are satisfied: the preimage by $h$ of any open subset of $L \times (0,1]$ is open, and the composition $R \circ h : C(L) \times [0,1] \rightarrow Z_0 \times [0,1]$ is continuous. To see that the first condition is satisfied, notice that the corestriction of $h$ to $L \times (0,1] \subset C(L)$ is the map $L \times (0,1] \times [0,1] \rightarrow L \times (0,1]$ given by the expression of the second line in the definition of $h$, and this map is continuous. To see that the second condition is satisfied, observe that $Z_0 \times [0,1] = [0,1]$ and the composition $R \circ h$ is given by the expression
    
    $$
    \fonctionsansnom{C(L) \times [0,1]}{[0,1]}{(c,s)}{(1-s)R(\varphi(c))+sR(c),}
    $$
    
    and this expression is continuous. This concludes the proof. \qedhere

\end{proof}

Finally, we also introduce a notion of cylindrical stratification with possibly more than two strata.

\begin{Definition}\label{definition: cylindrical stratification with many strata}

Suppose given a poset $A$ and a topological space $Z$ stratified by $A$. Given $a \in A$, we say that $Z$ is \emph{cylindrically stratified along the} $a$-\emph{stratum} if there exists a neighborhood $N$ of $Z_a$ in $Z$, a $A_{>a}$-stratified topological space $L$ and a continuous map $L \rightarrow Z_a$ such that there exists a homeomorphism between $N$ and $\Cyl(r)$ which is the identity on $Z_a$ and is compatible with the stratifications. Here, $\Cyl(r)$ is endowed with the $A_{\ge a}$-stratification whose restriction to $L \times (0,1]$ is induced by the first projection and whose $a$-stratum is $Z_a$.

\end{Definition}

\begin{Remark}\label{lemma: restricition of cylindrical stratification}

This more general notion will be useful to us below because it will allow us to prove the existence of a cylindrical structure for certain stratifications with two strata, from the existence of a cylindrical structure for certain stratifications with more than two strata. Indeed, suppose given a $A_{>a}$-stratified topological space $L$ and a continuous map $r : L \rightarrow Z_a$. Suppose given $b>a$ and consider the subspace $\Cyl(r)_a \cup \Cyl(r)_b \subseteq \Cyl(r)$ endowed with the $(0<1)$-stratification whose $0$-stratum is $\Cyl(r)_a$. According to remark \ref{remark: teardrop topology}, the evident set-theoretic identification $\Cyl(r_{\mid L_b}) \simeq \Cyl(r)_a \cup \Cyl(r)_b$ is a homeomorphism and therefore this $(0<1)$-stratified space is cylindrically stratified. Consequently, if $Z$ is cylindrically stratified along the $a$-stratum then for every $b > a$ the $(0<1)$-stratified space $Z_a \cup Z_b \subseteq Z$ is cylindrically stratified.

\end{Remark}

\subsection{Morse-Smale pairs and flow categories}\label{section: Morse-Smale pairs and flow categories}

The goal of this section is to present the definitions and results of Morse theory that we are going to use. We introduce the notions of Morse function, Morse-Smale pair, and prove that the stable manifolds of a Morse-Smale pair form a stratification in the sense of definition \ref{definition: stratified space} (proposition \ref{proposition: stable manifolds form a stratification}). We then define the flow category of a Morse-Smale pair and discuss smoothness of its morphism spaces. As an important consequence, we obtain that the space of morphisms between two critical points is homotopy equivalent to the space of unbroken pseudo-gradient trajectories connecting them (corollary \ref{Corollary: space of broken trajectories equivalent to space of trajectories}).

A general reference for the unproven results of this section is the first part of the book \cite{AudinDamianMorseFloer}.

We let $X$ be a smooth manifold of dimension $n$.

\subsubsection{Morse-Smale pairs and unbroken trajectories}

\begin{Definition}\label{definition: Morse function}
    A \emph{Morse function} on $X$ is a smooth function $f : X \rightarrow \R$ whose critical points are nondegenerate, meaning that for any critical point $a$ of $f$, the Hessian of $f$ at $a$ is a nondegenerate symmetric bilinear $2-$form on the tangent space of $X$ at $a$.
\end{Definition}

Morse functions always exist on $X$, and they form an open dense subset of $C^{\infty}(X,\R)$ for the $C^{\infty}$ topology (\cite[Theorem 1.2.5]{AudinDamianMorseFloer}).

The following lemma is known as the \textit{Morse lemma} (\cite[Theorem 1.3.1]{AudinDamianMorseFloer}).

\begin{Lemmas}\label{lemma: Morse lemma}
    Let $a \in X$ be a nondegenerate critical point of a real smooth function. There exist local coordinates $(x_1,\hdots,x_n)$ at $a$ in which this function is written as follows

    $$
    f(x) = f(a) -\sum_{i=1}^{k}x_i^2 + \sum_{i=k+1}^n x_i^2.
    $$
    \qed
\end{Lemmas}

\begin{Definition}\label{definition: Morse chart}
    Let $a \in X$ be a critical point of a real smooth function. A local chart at $a$ with coordinates as in lemma \ref{lemma: Morse lemma} is called a \emph{Morse chart} for this function at $a$.
\end{Definition}

\begin{Definition}\label{definition: Morse index}

The integer $0 \le k \le n$ from the expression in lemma \ref{lemma: Morse lemma} is the maximal dimension of a negative definite subspace for the Hessian of $f$ at $a$, and therefore does not depend on the Morse chart at $a$. It is called the \textit{Morse index}, or simply the \emph{index}, of the critical point $a$ with respect to the function, and is denoted by $|a|$.

\end{Definition}

\begin{Definition}\label{definition: pseudo-gradient}
    Let $f : X \rightarrow \R$ be a Morse function. A \emph{negative pseudo-gradient adapted to} $f$ is a vector field $\xi$ on $X$ satisfying the following conditions:

    \begin{itemize}
        \item $\xi(x)=0$ if and only if $x$ is a critical point of $f$, and if not, $\emph{d} f_x(\xi(x)) < 0$.

        \item For every critical point $a$, there exists a Morse chart for $f$ at $a$ in which $\xi$ coincides with the Euclidean negative gradient of $f$. In other words, $\xi$ is written in this chart as follows

        $$
        \sum_{i=1}^{k} 2x_i \frac{\partial}{\partial x_i} - \sum_{i=k+1}^n 2x_i \frac{\partial}{\partial x_i}.
        $$

    \end{itemize}
\end{Definition}

A negative pseudo-gradient adapted to $f$ always exists (\cite[Theorem 2.2.5]{AudinDamianMorseFloer}).

Assume from now on that $X$ is compact without boundary, and fix a negative pseudo-gradient $\xi$ adapted to $f$. The flow of $\xi$ is then defined for all times $t \in \R$.

\begin{Notation}\label{notation: flow of the pseudo-gradient}

We denote by $\phi^t(x)$ the flow of $\xi$ at time $t$, starting from the point $x \in X$.

\end{Notation}

\begin{Definition}\label{definition: flow line}
    A \emph{flow line of} $\xi$ is a map of the form

    $$
    \fonctionsansnom{(p,q)}{X}{t}{\phi^t(x)}
    $$

    for some $p, q \in \R \cup \{-\infty, +\infty\}$ such that $p < q$, and some $x \in X$.
\end{Definition}

\begin{Proposition}\label{proposition: limits of flow line are critical points}
    The limits at $- \infty$ and $+ \infty$ of any flow line of $\xi$ defined on $(-\infty,+\infty)$ exist and are critical points of $f$.
\end{Proposition}

\begin{proof}
    See \cite[Theorem 2.1.6]{AudinDamianMorseFloer}.
\end{proof}

\begin{Definition}\label{definition: stable and unstable manifolds}

For every critical point $a$ of $f$ we define the \textit{stable  manifold} of $a$ as

$$
W^s(a) = \{ x \in X \mid \limit{\phi^t(x)}{t}{+ \infty}{a} \},
$$

and the \textit{unstable manifold} of $a$ as

$$
W^u(a) = \{ x \in X \mid \limit{\phi^t(x)}{t}{- \infty}{a} \}.
$$

\end{Definition}

\begin{Example}\label{example: stable manifolds}
    Figure \ref{figure: stratification on other sphere} provides examples of stable manifolds.
\end{Example}

\begin{Proposition}\label{proposition: stable and unstable manifolds are R^k}

The stable and unstable manifolds of every critical point are smooth submanifolds of $X$, diffeomorphic to $\R^{n-|a|}$ and $\R^{|a|}$, respectively.

\end{Proposition}

\begin{proof}

See \cite[Section 2.1.d]{AudinDamianMorseFloer}. \qedhere

\end{proof}

\begin{Remark}\label{remark: points of stable or unstable manifolds as trajectories}
    Given a critical point $a$ of $f$, one can alternatively think of the points of $W^s(a)$ (resp. $W^u(a)$) as (unbroken) trajectories of $\xi$ \emph{ending at} $a$ (resp. \emph{starting} from $a$).
\end{Remark}

\begin{Definition}\label{definition: space of unbroken trajectories}

For any two critical points $a$ and $b$, the quotient space of $W^u(a) \cap W^s(b)$ by the action of $\R$ induced by the flow of $\xi$ is called the \emph{space of unbroken trajectories} of $\xi$ connecting $a$ to $b$, and is denoted $\Mi(a,b)$ \footnote{This notation is justified by the fact that we will consider the compactified version of this space below (proposition \ref{proposition: structure of smooth manifold with corners on space of broken traj}).}. An element of this quotient is called an \emph{unbroken trajectory} of $\xi$ connecting $a$ to $b$.

\end{Definition}

Suppose that there exists an unbroken trajectory $\gamma \in \Mi (a,b)$, and take a representing element $x \in W^u(a) \cap W^s(b)$ for it. The orbit of $x$ under the flow of $\xi$ is the image of $\R$ by the map $t \mapsto \phi^t(x)$. This map tends to $a$ as $t$ tends to $-\infty$, and tends to $b$ as $t$ tends to $+ \infty$.

\begin{Definition}\label{definition: image of a trajectory}
    The \emph{image} of $\gamma$ is the union of $\{a,b\}$ with the orbit of $x$ under the flow of $\xi$.
\end{Definition}

Note that this definition does not depend on the choice of the representing element $x$.

\begin{Remark}\label{remark: a critical point is connected ot itself by a trajectory}
    Every critical point is connected to itself by the constant trajectory of $\xi$ at that point.
\end{Remark}

\begin{Definition}\label{definition: Smale's transversality condition}
    We say that $\xi$ satisfies \textit{Smale's transversality condition} if for every two distinct critical points $a$ and $b$, the submanifolds $W^u(a)$ and $W^s(b)$ intersect transversely.
\end{Definition}

\begin{Remark}\label{remark: smoothness and dimension of spaces of unbroken trajectories}
     When $\xi$ satisfies Smale's transversality condition, the intersection $W^u(a) \cap W^s(b)$ is a submanifold of $X$ of dimension $|a| - |b|$. Moreover, when $f(a) > f(b)$, the space $\Mi(a,b)$ (definition \ref{definition: space of unbroken trajectories}) is naturally identified with the transverse intersection of $W^u(a) \cap W^s(b)$ with any regular level set of $f$ of level strictly comprised between $f(a)$ and $f(b)$. It is endowed in this way with the structure of a smooth manifold of dimension $|a| - |b| - 1$.

\end{Remark}

\begin{Remark}\label{remark: index increases along trajectories}
    By remark \ref{remark: smoothness and dimension of spaces of unbroken trajectories}, if $\xi$ satisfies Smale's transversality condition and $a, b$ are two critical points of $f$ such that there exists a non constant trajectory of $\xi$ starting at $a$ and ending at $b$, then $|a| > |b|$.
\end{Remark}

\begin{Definition}\label{definition: Morse-Smale pair}

The datum $(f,\xi)$ of a Morse function on $X$ and of a negative pseudo-gradient adapted to $f$ satisfying Smale's transversality condition is called a \textit{Morse-Smale pair} on $X$.

\end{Definition}

Let us fix a Morse-Smale pair $(f,\xi)$ on $X$ and denote by $A$ the set of critical points of $f$.

\begin{Proposition}\label{proposition: Smale ordering}
    The relation $\le$ on $A$ defined as: $a \le b$ if and only if $\Mi(a,b)$ is not empty, is a partial order.
\end{Proposition}

\begin{proof}
    This relation is reflexive by remark \ref{remark: a critical point is connected ot itself by a trajectory}. It is antisymmetric because, if $a$ and $b$ are two critical point such that $a$ is connected to $b$ by an unbroken trajectory of $\xi$ and $b$ is connected to $a$ by an unbroken trajectory of $\xi$, then, since the values of $f$ are strictly decreasing along non-constant flow lines of $\xi$, we must have that these two trajectories are constant, and therefore that $a = b$. Transitivity is more delicate and follows from proposition \ref{proposition: structure of smooth manifold with corners on space of broken traj} below.
\end{proof}

\begin{Remark}\label{remark: smale order oppositely compatible with order on R}

Because of our choice to consider a \emph{negative} pseudo-gradient of $f$, the map $f : A \rightarrow \R$ is \emph{decreasing.}

\end{Remark}

\begin{Definition}\label{definition: partial order on set of critical points}

The partial order $\le$ on $A$ from proposition \ref{proposition: Smale ordering} is called the \emph{Smale ordering}.

\end{Definition}

\begin{Remark}\label{remark: equivalence characterizations of partial order}
    We will see below equivalent characterizations of this partial order relation (corollary \ref{corollary: equivalent characterizations of partial order on A}).
\end{Remark}

\begin{Construction}\label{construction: filtration on usual realization}
    Let $\ag = [a_0 \le \hdots \le a_n]$ be a finite increasing sequence of elements of $A$ of length $n$ and denote by $(e_0,\hdots,e_n)$ the canonical basis of $\R^{n+1}$. The sequence $\ag$ determines a map
    
    $$
    \fonctionsansnom{|\Delta^n|}{[f(a_n),f(a_0)] \subset \R}{\sum_{i=0}^n t_ie_i}{\sum_{i=0}^n t_if(a_i).}
    $$
    
    Given another sequence $\bg$ of length $m$ and a morphism $\Phi : \Delta^{\ag} \rightarrow \Delta^{\bg}$ in $\Delta_A$, denoting by $\varphi$ the underlying morphism $\Delta^n \rightarrow \Delta^m$ in $\sSet$, the following diagram commutes
    
    $$
    \xymatrix{
    |\Delta^n| \ar[rd] \ar[rr]^-{|\varphi|} && |\Delta^m| \ar[ld] \\
    & \R.
    }
    $$
    
    Denote by $\underline{\R}$ the \emph{constant functor} at $\R$ on $\Delta_A$, defined as
   
   $$
   \fonction{\underline{\R}}{\Delta_A}{\Top}{
    \left\{\begin{array}{c}
         \Delta^{\ag}  \\
         \Phi : \Delta^{\ag} \rightarrow \Delta^{\bg} 
    \end{array}
    \right.
   }{
   \left\{\begin{array}{c}
         \R  \\
         \id : \R \rightarrow \R. 
    \end{array}
    \right.
   }
   $$
   
    We have thus constructed a natural transformation $|-|_A \rightarrow \underline{\R}$ of functors from $\Delta_A$ to $\Top$. Furthermore, the latter extends to a natural transformation between the left Kan extensions of these functors to $\sSet_A$. Note that the left Kan extension of $\underline{\R}$ to $\sSet_A$ is (canonically isomorphic to) the constant functor at $\R$ on $\sSet_A$.
\end{Construction}

\begin{Definition}\label{definition: filtration on usual realization}
    We refer to the natural transformation $|-|_A \rightarrow \underline{\R}$ of functors from $\sSet_A$ to $\Top$, constructed in \ref{construction: filtration on usual realization}, as the \emph{filtration} on $|-|_A$ induced by $f$, and we denote it by $f_A$.
\end{Definition}

\begin{Remark}\label{remark: explicit description of filtration}
    Given an $A$-stratified simplicial set $S$, let us describe explicitly the map
    $f_A(S) : |S|_A \rightarrow \R$. Recall from remark \ref{remark: explicit description of stratified realization} that the $A$-stratified topological space $|S|_A$ is obtained as a quotient of the following $A$-stratified topological space
    
    $$
    \bigsqcup_{\Delta^{\ag} \in \Delta_A} S_{\ag} \times |\Delta^{l(\ag)}|.
    $$
    
    The maps $|\Delta^{l(\ag)}| \rightarrow \R$ from construction \ref{construction: filtration on usual realization} determine a map from this disjoint union to $\R$, and the latter factors through $|S|_A$, thus defining the map $f_A(S)$.
\end{Remark}

\subsubsection{The stratification by the stable manifolds of a Morse-Smale pair}\label{section: stratification by stable manifolds}

Assume $X$ is closed and let $(f,\xi)$ be a Morse-Smale pair on $X$. The goal of this section is to introduce the stratification associated with this data, that we will study throughout the paper.

\begin{Proposition}\label{proposition: stable manifolds form a stratification}
    Recall the notation $\phi$ for the flow of $\xi$ (notation \ref{notation: flow of the pseudo-gradient}). The map

    $$\fonction{\pi}{X}{A}{x}{\underset{t\rightarrow + \infty}{\lim} \phi^t (x),}$$

    which is well-defined by proposition \ref{proposition: limits of flow line are critical points}, is a stratification of $X$ by $A$ (in the sense of definition \ref{definition: stratified space}), satisfying $X_a = W^s(a)$ for all $a \in A$. We refer to this stratification as the \emph{stratification by the stable manifolds of} $\xi$.
\end{Proposition}

\begin{Example}\label{example: stratification by stable manifolds}

    Example \ref{example: stratification on the other sphere} illustrates this stratification for a specific example of $(X,f,\xi)$.
    
\end{Example}

\begin{proof}
Recall that we denote by $A_{\ge a}$ the open subset of $A$ formed by those $b \in A$ such that $b \ge a$. We have to prove that for every $a \in A$, $\pi^{-1}(A_{\ge a})$ is an open subset of $X$, that is, for every $x \in X$ such that $\pi(x) \ge a$ there exists a neighborhood $U$ of $x$ in $X$ such that for every $y \in U$, $\pi(y) \ge a$. 

We first argue that it suffices to treat the case where $x$ is a critical point. Assume indeed that the property holds for the critical point $\pi(x)$. It means that there exists a neighborhood $V$ of $\pi (x)$ such that for every $y \in V$, $\pi(y) \ge \pi(x)$. Since the flow of $\xi$ is continuous, there exists a neighborhood $U$ of $x$ and a time $t \in \R$ such that $\phi^t(U) \subset V$. Then for every $y \in U$ we have $\pi(y) \ge \pi(x) \ge a$.

We assume from now on that $x$ is the critical point $a$. Our proof is by induction on $|a|$. For $|a| = 0$, by definition \ref{definition: pseudo-gradient}, there exists a neighborhood $V$ of $a$ in $X$ and coordinates $(y_1,\hdots,y_n)$ on $V$ in which $\xi$ is written

$$-\sum_{i=1}^{n} 2x_i \frac{\partial}{\partial x_i}.$$

Hence, for every $y \in V$ with coordinates $(y_1, \hdots, y_n)$, the flow line starting from $y$ is given in these coordinates by the formula

$$\phi^t(y) = (y_0e^{-2t},\hdots,y_ne^{-2t}).$$

This expression converges to $0$ as $t$ goes to $+ \infty$, in other words we have $\pi(y)=a$. This completes the first step of the induction proof.

Now let $ 1 \le k \le n$ and assume that for every critical point $b$ satisfying $|b| < k$, there exists a neighborhood $\Ul_b$ of $b$ in $X$ such that for every $y \in \Ul_b$, $\pi(y) \ge b$. Let $a$ be a critical point of index $k$, we wish to prove that there exists such a neighborhood of $a$. By definition \ref{definition: pseudo-gradient}, there exists a neighborhood $V$ of $a$ in $X$ and coordinates $(y_1,\hdots,y_n)$ on $V$ in which $\xi$ is written

$$
\sum_{i=1}^{k} 2x_i \frac{\partial}{\partial x_i} - \sum_{i=k+1}^n 2x_i \frac{\partial}{\partial x_i}.
$$

The flow line starting from $y=(y_1,\hdots,y_n) \in V$ is given in $V$ by the formula

$$\phi^t(y) = (y_1e^{2t},\hdots,y_ke^{2t},y_{k+1}e^{-2t},\hdots,y_ne^{-2t}).$$

Denote $H_- = \R^k \times 0$ and $H_+ = 0 \times \R^{n-k}$. 

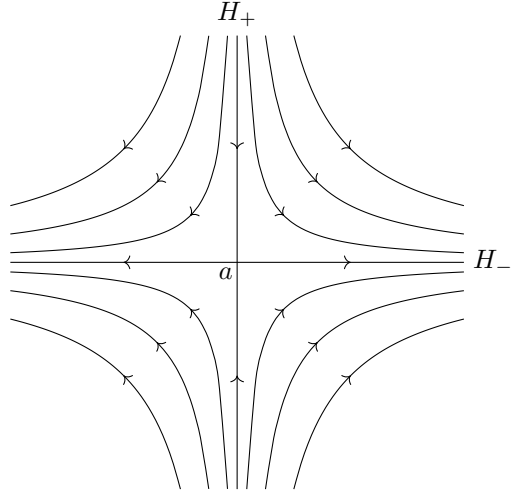
\begin{figure}[H]
    \centering
\begin{tikzpicture}[scale=1.5]

\draw (-2,0) -- (2,0) node[right] {$H_-$};
\draw[->, shift={(1,0)}] (0,0) -- (0.0001,0);
\draw[->, shift={(-1,0)}] (0,0) -- (-0.0001,0);
\draw (0,-2) -- (0,2) node[above] {$H_+$};
\draw[->, shift={(0,1)}] (0,0) -- (0,-0.0001);
\draw[->, shift={(0,-1)}] (0,0) -- (0,0.0001);

\draw node at (-0.1,-0.1) {$a$};


    
    \draw[domain=0.085:2,smooth,variable=\x] 
        plot ({\x},{0.17/\x});
    \draw[->, shift={({sqrt(0.17)},{sqrt(0.17)})}, rotate=-45] (0,0) -- (0.0001,0);
    
    \draw[domain=-2:-0.085,smooth,variable=\x] 
        plot ({\x},{0.17/\x});
    \draw[->, shift={(-{sqrt(0.17)},-{sqrt(0.17)})}, rotate=135] (0,0) -- (0.0001,0);
    
    \draw[domain=-2:-0.085,smooth,variable=\x] 
        plot ({\x},{-0.17/\x});
    \draw[->, shift={(-{sqrt(0.17)},{sqrt(0.17)})}, rotate=-135] (0,0) -- (0.0001,0);
        
    \draw[domain=0.085:2,smooth,variable=\x] 
        plot ({\x},{-0.17/\x});
    \draw[->, shift={({sqrt(0.17)},-{sqrt(0.17)})}, rotate=45] (0,0) -- (0.0001,0);

    
    \draw[domain=0.25:2,smooth,variable=\x] 
        plot ({\x},{0.5/\x});
    \draw[->, shift={({sqrt(0.5)},{sqrt(0.5)})}, rotate=-45] (0,0) -- (0.0001,0);
        
    \draw[domain=-2:-0.25,smooth,variable=\x] 
        plot ({\x},{0.5/\x});
    \draw[->, shift={(-{sqrt(0.5)},-{sqrt(0.5)})}, rotate=135] (0,0) -- (0.0001,0);
        
    \draw[domain=-2:-0.25,smooth,variable=\x] 
        plot ({\x},{-0.5/\x});
    \draw[->, shift={(-{sqrt(0.5)},{sqrt(0.5)})}, rotate=-135] (0,0) -- (0.0001,0);
        
    \draw[domain=0.25:2,smooth,variable=\x] 
        plot ({\x},{-0.5/\x});
    \draw[->, shift={({sqrt(0.5)},-{sqrt(0.5)})}, rotate=45] (0,0) -- (0.0001,0);
        
    
    \draw[domain=0.5:2,smooth,variable=\x] 
        plot ({\x},{1/\x});
    \draw[->, shift={(1,1)}, rotate=-45] (0,0) -- (0.0001,0);
        
    \draw[domain=-2:-0.5,smooth,variable=\x] 
        plot ({\x},{1/\x});
    \draw[->, shift={(-1,-1)}, rotate=135] (0,0) -- (0.0001,0);
        
    \draw[domain=-2:-0.5,smooth,variable=\x] 
        plot ({\x},{-1/\x});
    \draw[->, shift={(-1,1)}, rotate=-135] (0,0) -- (0.0001,0);
        
    \draw[domain=0.5:2,smooth,variable=\x] 
        plot ({\x},{-1/\x});
    \draw[->, shift={(1,-1)}, rotate=45] (0,0) -- (0.0001,0);

\end{tikzpicture}
    \caption{The flow lines of $\xi$ in $V$ for $n=2$ and $k=1$.}
    \label{figure: flow lines in Morse chart}
\end{figure}

Observe that for $z \in H_+$ one has $\pi(z)=a$ and for $z \in H_-$ one has $\underset{t \rightarrow - \infty}{\lim} \phi^t(z) = a$, hence $\pi(z) \ge a$, by definition. In particular, by remark \ref{remark: index increases along trajectories}, we have $|\pi(z)| < |a|$ as long as $z \neq a$, and thus the induction hypothesis applies to $\pi(z)$ in that case. Our strategy is to show that there exists a neighborhood $\Ul$ of $a$ such that if $y \notin H_- \cup H_+$ and $y \in \Ul$, then the flow line through $y$ passes sufficiently close to $H_-$ so that there exists a $z \in H_-$ such that the flow line starting from $y$ passes through $\Ul_{\pi(z)}$. It follows from this that $\pi(y) \ge \pi(z) \ge a$, as desired.

There exists $\varepsilon > 0$ such that $V$ contains $D^k(\varepsilon) \times D^{n-k}(\varepsilon)$ (where $D^k(\varepsilon)$ denotes the closed unit ball of radius $\varepsilon$). To simplify, we will assume in the sequel that $\varepsilon = 1$ and denote this product of disks $D^k \times D^{n-k}$. We will use the following lemma.

\begin{Lemmas}\label{structure du flot au voisinage d'un point critique}
    For every $0 < r \le 1$, there exists a neighborhood of $a$ such that for every $y$ in that neighborhood, either $y \in H_+$ or the flow line starting at $y$ passes through $\SM^{k-1} \times D^{n-k}(r)$.
\end{Lemmas}

\begin{proof}
We prove this using the explicit formula for the flow in $V$ given above.

For any $y \in \R^n$ one has $(y_1e^{2t},\hdots,y_ke^{2t},y_{k+1}e^{-2t},\hdots,y_ne^{-2t}) \in \SM^{k-1} \times \R^{n-k} $ if and only if $t = t_y$ where $t_y$ satisfies

$$
e^{-2t_y} = ||(y_1,\hdots,y_k)||_2.
$$

Let $y \in D^k \times D^{n-k}$. For $0 \le t \le t_y$ one has

$$
||(y_1 e^{2t},\hdots,y_k e^{2t})||_2 \le ||(y_1 e^{2t_y};\hdots,y_k e^{2t_y})||_2 = 1
$$

and

$$
||(y_{k+1}e^{-2t},\hdots,y_ne^{-2t})||_2 \le ||(y_{k+1},\hdots,y_n)||_2 \le 1.
$$

Hence the flow line starting from $y$ stays in $V$ until time $t_y$. Now given $y \in D^k(r) \times D^{n-k}$, the flow line starting from $y$ stays in $V$ until time $t_y$, hence we have $\phi^{t_y}(y) \in \SM^{k-1} \times D^{n-k}$ and further

$$\phi^{t_y}(y) = (y_1e^{2t_y},\hdots,y_ke^{2t_y},y_{k+1}e^{-2t_y},\hdots,y_ne^{-2t_y}).$$

Moreover

$$
\begin{aligned}
||(y_{k+1}e^{-2t_y},\hdots,y_ne^{-2t_y})||_2 & = e^{-2t_y}||(y_{k+1},\hdots,y_n)||_2 \\
& \le e^{-2t_y} \\
& = ||(y_1,\hdots,y_k)||_2 \\
& \le r.
\end{aligned}
$$

Hence $D^k(r) \times D^{n-k}$ is a neighborhood of $a$ satisfying the required condition. \qedhere

\end{proof}

\textit{End of the proof of proposition \ref{proposition: stable manifolds form a stratification}.} Let $y \in \SM^{k-1} \times 0 \subset H_-$. By continuity of the flow, there exists a $r_y > 0$ and a neighborhood $W$ of $y$ in $\SM^{k-1}$ such that every flow line passing through $W \times D(r_y)$ also passes through $\Ul_{\pi(y)}$. By compactness of $\SM^{k-1}$, there exists a $r > 0$ such that every flow line passing through $\SM^{k-1} \times D^{n-k}(r)$ also passes through $\Ul_z$ for some $z \ge a$. We conclude by applying lemma \ref{structure du flot au voisinage d'un point critique}. \qedhere

\end{proof}

A first basic observation, connecting the discussion in this section with that in section \ref{section: stratified spaces, exit paths and stratified simplicial sets}, is that to every unbroken trajectory of $\xi$ there corresponds an exit path for the stratification of $X$ by stable manifolds. To make this precise, we make the following construction.

\begin{Construction}\label{construction: parametrization of trajectories by values of the function}
    Let $a, b$ by two distinct critical points of $f$ and suppose that there exists an unbroken trajectory $\gamma$ of $\xi$ connecting $a$ to $b$. Let $\widetilde{\gamma} \subset X$ be the image of $\gamma$ (in the sense of definition \ref{definition: image of a trajectory}). We endow $\widetilde{\gamma}$ with the topology induced by that on $X$. The restriction of $f$ to $\widetilde{\gamma}$ induces a homeomorphism
    
    $$
    f : \widetilde{\gamma} \overset{\simeq}{\longrightarrow} [f(b),f(a)].
    $$
    
    The inverse of this homeomorphism is denoted
    
    $$
    \fonctionsansnom{[f(b),f(a)]}{X,}{s}{\gamma(s).}
    $$

    By definition, it satisfies $f(\gamma(s))=s$ for every $s \in [f(b),f(a)]$. We refer to this map as the \emph{parametrization of} $\gamma$ \emph{by the values of} $f$. Further, by precomposing this map by the affine nonincreasing surjection $[0,1] \rightarrow [f(b),f(a)]$ \footnote{If $a \neq b$, this is a decreasing bijection.}, we obtain a path denoted
    
    $$
    \iota_{a,b}(\gamma) : [0,1] \rightarrow X,
    $$
    
    which we refer to as the \emph{path in $X$ associated with $\gamma$}. These constructions extend to the case $a=b$ in a straightforward way.
\end{Construction}

\begin{Lemmas}\label{lemma: unbroken trajectories give exit paths}
    Let $a, b$ be two critical points of $f$ and assume that there exists an unbroken trajectory $\gamma$ of $\xi$ connecting $a$ to $b$. Then the path in $X$ associated with $\gamma$ (construction \ref{construction: parametrization of trajectories by values of the function}) is an exit path with respect to the stratification of $X$ by the stable manifolds of $\xi$.
\end{Lemmas}

\begin{proof}
    In the case where $a=b$, $\iota_{a,b}(\gamma)$ is constant and is therefore an exit path. Let us assume that $a \neq b$. We have $\iota_{a,b}(\gamma)(0) = a \in X_a$ and $\iota_{a,b}(\gamma)(1) = b \in X_b$. For every $s \in ]0,1[$, $\iota_{a,b}(\gamma)(s)$ belongs to the image of a flow line of $\xi$ defined on $(-\infty,+\infty)$ whose limit at $+ \infty$ is $b$, hence $\iota_{a,b}(\gamma)(s) \in X_b$. This proves that $\iota_{a,b}(\gamma)$ is an exit path.
\end{proof}

\begin{Remark}\label{remark: map from unbroken trajectories to exit paths}
    By lemma \ref{lemma: unbroken trajectories give exit paths}, we have, for every pair of critical points $a,b$ of $f$, an inclusion

    $$
    \iota_{a,b} : \Mi(a,b) \rightarrow \Exit(a,b).
    $$

    By remark \ref{remark: different topologies on morph spaces} below, it is continuous. We will study this inclusion in section \ref{section: spaces of exit paths as spaces of pseudo gradient trajectories}, where we will prove that it is a weak homotopy equivalence.
\end{Remark}

\begin{Remark}\label{remark: choice for definition of iota a b}

    In the case when $a < b$, the map $\iota_{a,b}$ is defined by using the affine decreasing bijection $[0,1] \rightarrow [f(b),f(a)]$, but we could have defined it using any decreasing homeomorphism between $[0,1]$ and $[f(b),f(a)]$ instead. The map $\Mi(a,b) \rightarrow \Exit(a,b)$ obtained from any such homeomorphism is homotopic to $\iota_{a,b}$.

\end{Remark}

The following proposition implies that Lurie's exit path $\infty$-category construction (definition \ref{definition: stratified singular simplicial set} and theorem \ref{theorem: stratified simplicial set of conically stratified space is infty cat}) applies to the stratification of $X$ by the stable manifolds of $(f,\xi)$. We will outline a proof of it in section \ref{section: Whitney stratifications} (see corollary \ref{corollary: morse stratification is conical}).

\begin{Proposition}\label{proposition: the stratification by stable manifolds is conical}
    The stratification of $X$ by the stable manifolds of $(f,\xi)$ is conical (in the sense of definition \ref{definition: conically stratified space}).
\end{Proposition}

\subsubsection{The flow category of a Morse-Smale pair}\label{section: the flow category of a Morse-Smale pair}

Assume $X$ is closed and let $(f,\xi)$ be a Morse-Smale pair on $X$. We defined previously the notion of trajectory of $\xi$ (definition \ref{definition: space of unbroken trajectories}). The goal of this section is to a introduce a notion of composition between the trajectories of $\xi$.

\begin{Definition}\label{definition: broken trajectory}
    A \emph{broken trajectory} of $\xi$ connecting $a$ to $b$ is a sequence $\gamma_0, \hdots, \gamma_n$ of non-constant unbroken trajectories of $\xi$, such that $n \ge 1$ and there exist critical points $a_0 < a_1 < \hdots < a_{n+1}$ such that $a_0 = a$, $a_{n+1} = b$ and for every $0 \le i \le n$ one has $\gamma_i \in \Mi(a_i,a_{i+1})$. We say that the broken trajectory is \emph{broken} at $a_1, a_2, \hdots, a_{n-1}$ and $a_n$.
\end{Definition}

\begin{Example}\label{example: borken trajectory}

    See figure \ref{figure: example of broken simplex} for an example of broken trajectory.

\end{Example}

\begin{Construction}\label{construction: parametrization of all trajectories by values of the function}

Construction \ref{construction: parametrization of trajectories by values of the function} extends to broken trajectories as follows: given a broken trajectory $\gamma_0,\hdots,\gamma_n$ connecting $a$ to $b$ as in definition \ref{definition: broken trajectory} the \emph{parametrization of} $\gamma$ \emph{by the values of} $f$ is the map

$$
[f(b),f(a)] \rightarrow X
$$

sending $s$ to $\gamma_i(s)$ if $s \in [f(a_{i+1}),f(a_i)]$. In other words, the embedding

$$
\fonctionsansnom{\Mi(a,b)}{C^0([f(b),f(a)],X),}{\gamma}{(s \mapsto \gamma(s))}
$$

extends to an injection from the set of all (broken and unbroken) trajectories connecting $a$ to $b$, by concatenation on the right-hand side.
\end{Construction}

\begin{Definition}\label{definition: space of trajectories between two critical points}
    The set of broken and unbroken trajectories of $\xi$ connecting $a$ to $b$, endowed with the topology induced by the compact-open topology on $C^0([f(b),f(a)],X)$, is called the \emph{space of trajectories of} $\xi$ \emph{connecting} $a$ \emph{to} $b$, and is denoted $\Ml(a,b)$.
\end{Definition}

\begin{Remark}\label{remark: different topologies on morph spaces}

The restriction of this topology to $\Mi(a,b)$ coincides with the one previously defined; see \cite[Theorem B.4]{CalleLiuClassifyingSpace} for a proof.

\end{Remark}

\begin{Construction}\label{construction: concatenation of trajectories}
    Given three critical points $a, b$ and $c$, the trajectories from $a$ to $b$ and from $b$ to $c$ can be concatenated in a straighforward way, thus defining a continuous map
    $$
    \Ml(a,b) \times \Ml(b,c) \rightarrow \Ml(a,c).
    $$

    This concatenation operation is (strictly) associative, by definition of broken trajectories.
\end{Construction}

\begin{Definition}\label{Catégorie de flot}
    The \emph{flow category} of $(X,f,\xi)$ is the topological category $\Ml$ whose set of objects is the set of critical points of $f$, whose space of morphisms between $a$ and $b$ is $\Ml(a,b)$ and whose identity morphisms are the constant trajectories. The composition operation in $\Ml$ is defined as the concatenation of trajectories of $\xi$, as in construction \ref{construction: concatenation of trajectories}.
\end{Definition}

This terminology was introduced by Cohen-Jones-Segal in \cite{CohenJonesSegalFloer}, and is also used in the foundational paper by Abouzaid-Blumberg \cite{AbouzaidBlumbergI}.

The smooth structures on the spaces of unbroken trajectories connecting critical points (remark \ref{remark: smoothness and dimension of spaces of unbroken trajectories}) can be extended to the spaces of possibly broken trajectories as now recall. We first need a definition.

\begin{Definition}\label{definition: k stratum of manifold with corners}

Let $M$ be a smooth manifold with corners of dimension $n$. For every $m \in M$ there exists a unique integer $0 \le k \le n$ such that a neighborhood of $m$ in $M$ is diffeomorphic to $[0,+\infty)^k \times \R^{n-k}$. The point $m$ is said to belong to the $k$-\emph{stratum} of $M$.

\end{Definition}

The following result is proved in \cite[Proposition 2.11]{LatourExistence} and \cite[Theorem 3.3]{QinModuli}.

\begin{Proposition}\label{proposition: structure of smooth manifold with corners on space of broken traj}
    For every pair of critical points $(a,b)$ the space $\Ml(a,b)$ is compact and can be endowed with a structure of smooth manifold with corners such that the following conditions are satisfied.
    
    \begin{enumerate}[label=(\roman*)]
        \item The $k$-stratum is the set of trajectories that are a concatenation of $k + 1$ unbroken non-constant trajectories (in other terms, that are broken exactly $k$ times).
        
        \item The smooth structure on $\Ml(a,b)$ extends that on $\Mi(a,b)$, in other words, the interior of $\Ml(a,b)$ is the smooth manifold $\Mi(a,b)$.
        
        \item These smooth structures are compatible in the sense that for every triple of critical points $a,b,c$, the composition map
        
        $$\Ml(a,b) \times \Ml(b,c) \rightarrow \Ml(a,c)$$
        
        is a smooth embedding.
    \end{enumerate}
\end{Proposition}

As a consequence we obtain:

\begin{Corollary}\label{Corollary: space of broken trajectories equivalent to space of trajectories}
    Given two critical points $a$ and $b$, the inclusion of the subspace of unbroken trajectories between $a$ and $b$, into $\Ml(a,b)$, is a homotopy equivalence.
\end{Corollary}

\begin{proof}
    By proposition \ref{proposition: structure of smooth manifold with corners on space of broken traj}, this is the inclusion of the interior of a smooth manifold with corners, hence the result.
\end{proof}

\begin{Remark}\label{remark: trajectories as closure of unbroken trajectories}
    Proposition \ref{proposition: structure of smooth manifold with corners on space of broken traj} also implies that the image of $\Ml(a,b)$ in $C^0([f(b),f(a)],X)$ is the \emph{closure} of that of $\Mi(a,b)$.
\end{Remark}

\begin{Remark}\label{remark: smooth structure on morph spaces}

In this paper we will not need the existence of smooth structures on the morphism spaces of $\Ml$. We will only need the structures of topological manifold with boundary.

\end{Remark}

\subsection{Precise formulation of the main theorem}\label{section: precise formulation of the main theorem}

Using sections \ref{section: Left Kan extensions}, \ref{section: infinity categories}, \ref{section: stratified spaces, exit paths and stratified simplicial sets} and \ref{section: Morse-Smale pairs and flow categories} we have all the ingredients to formulate precisely the main theorem of this paper. To this end, we start by recalling the main objects involved.

Let $X$ be a smooth closed manifold, $(f,\xi)$ a Morse-Smale pair on $X$. Let $\Ml$ the associated flow category, regarded as a topological category, and $\Nl(\Ml)$ the homotopy coherent nerve of $\Ml$.

Let $A$ be the set of critical points of $f$ endowed with the Smale ordering. Endow $X$ with the $A$-stratification by the stable manifolds of $\xi$ and denote by $\Sing_A(X)$ the associated exit path $\infty$-category. Recall that $\Sing_A(X)$ is endowed with the structure of an $A$-stratified simplicial set, i.e., a morphism $\Sing_A(X) \rightarrow N(A)$.

The precise statement of our theorem will require two additional observations.

Firstly, $\Nl(\Ml)$ is also naturally equipped with a morphism of simplicial sets $\Nl(\Ml) \rightarrow N(A)$ defined as follows. Every $n$-simplex $\sigma : \Delta^n \rightarrow \Nl(\Ml)$ determines a sequence of $n+1$ critical points $\ag =(a_0,\hdots,a_n)$ and, for every $0 \le i \le n-1$, a morphism $\gamma_i \in \Ml(a_i,a_{i+1})$. In particular, $\Ml(a_i,a_{i+1})$ is nonempty and therefore the sequence $\ag$ is increasing. The morphism $\Nl(\Ml) \rightarrow N(A)$ is defined by mapping $\sigma$ to $\ag$.

Secondly, suppose given two critical points $a,b$ of $f$. By proposition \ref{proposition: morphism spaces of a topological category and its coherent nerve} we have a natural homotopy equivalence between the morphism space of $\Nl(\Ml)$ between $a$ and $b$ and the Kan complex $\Sing(\Ml(a,b))$. By proposition \ref{proposition: space of morphism is space of exit paths} we have a natural isomorphism between the morphism space of $\Sing_A(X)$ between $a$ and $b$ and the Kan complex $\Sing(\Exit(a,b))$. Recall from remark \ref{remark: map from unbroken trajectories to exit paths} the map $\iota_{a,b} : \Mi(a,b) \rightarrow \Exit(a,b)$. We will prove in section \ref{section: spaces of exit paths as spaces of pseudo gradient trajectories} that this is a weak homotopy equivalence. Since the inclusion map $\Mi(a,b) \rightarrow \Ml(a,b)$ is a homotopy equivalence (corollary \ref{Corollary: space of broken trajectories equivalent to space of trajectories}), the induced morphism of Kan complexes $\Sing(\Mi(a,b)) \rightarrow \Sing(\Ml(a,b))$ is a homotopy equivalence. Altogether, we have an isomorphism in the homotopy category of spaces between the Kan complexes of morphisms of $\Nl$ and $\Sing_A(X)$ between $a$ and $b$.

\begin{Th}[Precise formulation of theorem \ref{Main theorem}]\label{theorem: main theorem reformulated}

There exists a canonical equivalence of $\infty$-categories $\Nl(\Ml) \rightarrow \Sing_A(X)$ that satisfies the following properties.

\begin{enumerate}[label=(\roman*)]
    \item This equivalence commutes with the functors to $N(A)$ on both sides.
    
    \item Suppose given two critical points $a$ and $b$. At the level of the morphism spaces between $a$ and $b$, this equivalence induces the isomorphism constructed above.
\end{enumerate}

\end{Th}

\begin{Remark}\label{remark: properties of the equivalence}

    Suppose given a functor $G : \Nl(\Ml) \rightarrow \Sing_A(X)$ that satisfies property (i). The meaning of property (ii) for such a functor is not obvious since we do not assume that $G(a) = a$ and $G(b)=b$, so let us clarify it. The functor $G$ determines a map from the morphism space of $\Nl(\Ml)$ between $a$ and $b$ to the morphism space of $\Sing_A(X)$ between $G(a)$ and $G(b)$. Now, by condition (i), $G(a)$ belongs to the $a$-stratum. The space of morphisms of $\Sing_A(X)$ between $a$ and $G(a)$ is the space of paths in $X_a$ starting at $a$ and ending at $G(a)$. In particular, every such morphism is an isomorphism, and since $X_a$ is contractible, this space of morphisms is contractible. Therefore $a$ and $G(a)$ are \emph{canonically isomorphic} as objects of $\Sing_A(X)$. The same discussion holds with $b$ instead of $a$. Consequently, the morphism spaces of $\Sing_A(X)$ between $a$ and $b$ and between $G(a)$ and $G(b)$ are canonically homotopy equivalent. The functor $G$ thus determines a map between the morphism spaces of $\Nl(\Ml)$ and $\Sing_A(X)$ between $a$ and $b$ that only depends on a contractible space of extra choices.
 
\end{Remark}

The proof of this theorem is the purpose of sections \ref{section: spaces of exit paths as spaces of pseudo gradient trajectories} through \ref{section: turning a semi simplicial map into functor}, and will be given in section \ref{section: end of the proof}. We close this section with two remarks.

\begin{Remark}\label{remark: opposite of flow category}
    We could consider alternatively the stratification of $X$ by the \emph{unstable manifolds} of $\xi$. These are defined to be the stable manifolds of $-\xi$. Notice that $-\xi$ is a negative pseudo-gradient vector field for $-f$ and that the Smale transversality condition for $\xi$ is equivalent to the Smale transversality condition for $-\xi$. Notice also that the flow category associated with $(X,-f,-\xi)$ is canonically isomorphic to $\Ml^{\op}$. Theorem \ref{Main theorem} therefore implies that there exists an equivalence between the $\infty$-category defined by $\Ml^{\op}$ and the exit path $\infty$-category associated with the stratification of $X$ by the unstable manifolds of $\xi$, that satisfies the corresponding properties.
\end{Remark}

\begin{Remark}\label{remark: naturality of the equivalence}

    In the statement of the theorem, by "canonical" we mean that the method that we follow exhibits an equivalence of $\infty$-categories that is unique up to natural equivalence.

\end{Remark}

\newpage

\section{Spaces of exit paths as spaces of pseudo-gradient trajectories}\label{section: spaces of exit paths as spaces of pseudo gradient trajectories}

Let $X$ be a closed smooth manifold, $(f,\xi)$ a Morse-Smale pair on $X$ and $A$ the poset of critical points of $X$. Let $a, b \in A$. Recall that in section \ref{section: stratification by stable manifolds} we constructed a continuous inclusion

$$ \iota_{a,b} : \Mi(a,b) \rightarrow \Exit(a,b) $$

of the space of unbroken trajectories of $\xi$ from $a$ to $b$, into the space of exit paths from $a$ to $b$ with respect to the stratification $X \rightarrow A$ by the stable manifolds of $\xi$ (construction \ref{construction: parametrization of trajectories by values of the function} and lemma \ref{lemma: unbroken trajectories give exit paths}).

The main goal of this section is to prove the following proposition:

\begin{Proposition}\label{proposition: equivalence between space of trajectories and space of exit paths}
    The map $\iota_{a,b}$ is a weak homotopy equivalence.
\end{Proposition}

Recall that the space of morphisms of the exit path $\infty$-category $\Sing_A(X)$ between $a$ and $b$ has the homotopy type of $\Exit(a,b)$ (proposition \ref{proposition: space of morphism is space of exit paths}). Recall also that the space of morphism of $\Ml$ between $a$ and $b$ has the homotopy type of $\Mi(a,b)$ (corollary \ref{Corollary: space of broken trajectories equivalent to space of trajectories}). Together with proposition \ref{proposition: equivalence between space of trajectories and space of exit paths}, this provides evidence for our main theorem \ref{Main theorem}. This is also an essential ingredient of our proof of it.

In section \ref{section: Whitney conditions and cylindrical structures}, we review the definition of the Whitney conditions and their geometric interpretation. We deduce from the latter the fact that the Whitney condition (b) for stratifications with two strata implies the existence of a cylindrical structure, in the sense of section \ref{section: geometric description of spaces of exit paths} (proposition \ref{proposition: cylindrical structure for Whitney stratifications}). In section \ref{section: Whitney stratifications}, we give the general definition of a Whitney stratification and we state a result of Nicolaescu asserting that the Smale transversality condition for the pseudo-gradient $\xi$ is equivalent to the condition that the partition of $X$ by the stable manifolds of $\xi$ is a Whitney stratification (theorem \ref{theorem: Smale if and only if Whitney}). We also extend our discussion from the first section and outline a proof of the fact that Whitney stratifications are cylindrical along each of their strata (corollary \ref{corollary: cylindrical structure for general Whitney stratifications}), and conical (corollary \ref{corollary:  Whitney stratifications are conical}.) Finally, section \ref{section: equivalence between unbroken trajectories and exit paths} mostly consists of the proof of proposition \ref{proposition: equivalence between space of trajectories and space of exit paths}. The argument combines properties of Whitney stratifications with the description of the space of exit paths between the two strata of a cylindrically stratified space provided by proposition \ref{proposition: space of exit paths as boundary of regular neighborhood}.

\subsection{Whitney conditions and cylindrical structure}\label{section: Whitney conditions and cylindrical structures}

In this section, we give a review of the theory of Whitney stratifications, and prove that Whitney stratifications with two strata are cylindrical. We follow \cite[Section 4]{NicolaescuInvitation} as well as \cite{TopologicalStability}.

Although we will not use it, we begin by giving the classical definition of the Whitney conditions (a) and (b).

\begin{Definition}\label{definition: Whitney conditions}
    Let $n \ge 0$ and let $Y, Y'$ be two disjoint submanifolds of $\R^n$.
    
    \begin{itemize}
        \item[$-$] We say that the pair $(Y,Y')$ satisfies the \emph{Whitney regularity condition (a)} at a point $y \in Y \cap \overline{Y'}$ if, for any sequence $y_n' \in Y'$ such that 
        
        \begin{itemize}
        
        \item[$\bullet$] $y_n' \underset{n \rightarrow + \infty}{\longrightarrow} y$, and
        
        \item[$\bullet$] the sequence of tangent spaces $T_{y_n'}(Y')$ converges to some vector subspace $T_{\infty}$ of $\R^n$,
        
    \end{itemize}
        
    we have $T_y Y \subset T_{\infty}$.
    
    \item[$-$] We say that the pair $(Y,Y')$ satisfies the \emph{Whitney regularity condition (b)} at a point $y \in Y \cap \overline{Y'}$ if, for any sequence $(y_n,y_n') \in Y \times Y'$ such that
    
    \begin{itemize}
        \item[$\bullet$] $y_n, y_n' \underset{n \rightarrow + \infty}{\longrightarrow} y$,
        
        \item[$\bullet$] the sequence of lines $\R(y_n'-y_n)$ converges to some line $l_{\infty}$, and
        
        \item[$\bullet$] the sequence of tangent spaces $T_{y_n'} Y'$ converges to some vector subspace $T_{\infty}$ of $\R^n$,
    
    \end{itemize}
    we have $l_{\infty} \subset T_{\infty}$.
    
    \item[$-$] We say that the pair $(Y,Y')$ \emph{satisfies the Whitney regularity condition (a) or (b) along} $Y$ if it satisfies this condition at any point of $Y \cap \overline{Y'}$.
    \end{itemize}
\end{Definition}

\begin{Remark}\label{remark: b implies a, and independant on coordinates}
    Condition (b) implies condition (a) (\cite[Proposition 2.4]{MatherNotes}). Moreover, one can check that these two conditions are independent of the choice of coordinates in the following sense. Given $y \in Y \cap \overline{Y'}$ and a diffeomorphism $\varphi$ between an open neighborhood $U$ of $y$ in $\R^n$, and an open subset $V$ of $\R^n$, the Whitney regularity condition (a) (resp. (b)) is satisfied by the pair $(Y \cap U,Y' \cap U)$ at $y$ if and only if it is satisfied by its image by $\varphi$, at the image of $y$. This allows to extend the definition of the Whitney regularity conditions to the case of two disjoint submanifolds of any smooth manifold.
\end{Remark}

These two conditions have a geometrically insightful interpretation, which is the one that we will use in practice. To explain it, we need to introduce the following notion of tubular neighborhood.

\begin{Definition}\label{definition: tubular neighborhood}
    Let $M$ be a smooth manifold and $Y$ a smooth submanifold of $M$. A \emph{tubular neighborhood} of $Y$ in $M$ consists of the following datum:
    
    \begin{itemize}
        \item a vector bundle $E \overset{\pi}{\longrightarrow} Y$ over $Y$, endowed with a fiberwise Euclidean structure $g$, that is, $g$ is a smooth family of scalar products on the fibers of $E$,
        
        \item a smooth function $\varepsilon: Y \rightarrow (0,+\infty)$ called the \emph{width function},
        
        \item an open neighborhood $T$ of $Y$ in $M$ and a diffeomorphism
        
        $$
        \left\{ v \in E \mid g(v,v) < \varepsilon(\pi(v))^2 \right\} \overset{\varphi}{\longrightarrow} T
        $$
        
        that restricts to the identity between the $0$-section of $E$, and $Y$.
    \end{itemize}
    
    Given a tubular neighborhood of $Y$ in $M$, the function
    
    $$
    \fonctionsansnom{E}{[0,+\infty)}{v}{g(v,v)}
    $$
    
    induces a function $\rho : T \rightarrow [0,+\infty)$. The map $\pi : E \rightarrow Y$ induces a map $T \rightarrow Y$, which we also denote by $\pi$. These are respectively called the \emph{radial function} and the \emph{projection} associated with the tubular neighborhood.
\end{Definition}

\begin{Remark}\label{remark: vector bundle in tubular neighborhood is normal bundle}
    Note that, in this definition of a tubular neighborhood, the vector bundle $E$ is necessarily isomorphic to the normal bundle of $Y$ in $M$.
\end{Remark}

\begin{Notation}\label{notation: sub tubular neighborhood}
    Given a tubular neighborhood of $Y$ in $M$ as in definition \ref{definition: tubular neighborhood}, and a smooth function $\alpha : Y \rightarrow (0,+\infty)$ such that $\alpha < \varepsilon$, we set
    
    $$
    T_{\alpha} = \{ t \in T \mid \rho(t) < \alpha(\pi(t))^2 \}.
    $$
    
    Given $y \in Y$ and $0 < r < \alpha(y)$, we denote
    
    $$
    T_{\alpha}(y,r) = \{ t \in T_{\alpha} \mid \pi(t)=y \text{ and } \rho(t) \le r^2 \}
    $$
    
    and
    
    $$
    \partial T_{\alpha}(y,r) = \{ t \in T_{\alpha}(y,r) \mid \rho(t) = r^2 \}.
    $$
    
    Note that $T_{\alpha}(y,r)$ and $\partial T_{\alpha}(y,r)$ are respectively diffeomorphic to a closed disc and its boundary.
\end{Notation}

\begin{Thms}\label{theorem: existence of tubes}
Let $M$ be a smooth manifold and $Y$ a smooth submanifold of $M$. Then there exists a tubular neighborhood of $Y$ in $M$.
\end{Thms}

\begin{proof}
    See \cite{TopologicalStability}, chapter II, theorem 1.6.
\end{proof}

The following result gives a geometric interpretation to the Whitney regularity conditions (a) and (b).

\begin{Proposition}\label{proposition: geometric interpretation of Whitney}
Let $M$ be a smooth manifold and $Y,Y'$ two disjoint submanifolds of $M$ such that $Y \subset \overline{Y'}$. Then the following hold.

\begin{itemize}
    \item[$-$] If the pair $(Y,Y')$ satisfies the Whitney regularity condition (a) along $Y$, then for every tubular neighborhood of $Y$ in $M$ (in the sense of definition \ref{definition: tubular neighborhood}), there exists a smooth function $\alpha : Y \rightarrow (0,+\infty)$ satisfying $ \alpha < \varepsilon$, such that for every $y \in Y$ and every $0 < r < \alpha(y)$, $T_{\alpha}(y,r)$ intersects $Y'$ and this intersection is transverse.
    
    \item[$-$] If the pair $(Y,Y')$ satisfies the Whitney regularity condition (b) along $Y$, then for every tubular neighborhood of $Y$ in $M$ (in the sense of definition \ref{definition: tubular neighborhood}), there exists a smooth function $\alpha : Y \rightarrow (0,+\infty)$ satisfying $ \alpha < \varepsilon$, such that for every $y \in Y$ and every $0 < r < \alpha(y)$, $\partial T_{\alpha}(y,r)$ intersects $Y'$ and this intersection is transverse.
\end{itemize}
\end{Proposition}

\begin{Remark}\label{remark: transverse intersection implies submersion}
    The two conditions of proposition \ref{proposition: geometric interpretation of Whitney} are respectively equivalent to the condition that the restriction
    
    $$
    Y' \cap T_{\alpha} \overset{\pi}{\longrightarrow} Y
    $$
    
    is a submersion, and the condition that the restriction
    
    $$
    Y' \cap T_{\alpha} \xrightarrow{\pi \times \rho} Y \times (0,+\infty)
    $$
    
    is a submersion. In particular, if condition (b) is satisfied and we further assume that $Y = \overline{Y'} \backslash Y'$, then the corresponding restriction is a proper submersion and therefore a fiber bundle over its image, by Ehresmann's theorem.
\end{Remark}

\begin{proof}[Proof of proposition \ref{proposition: geometric interpretation of Whitney}]
    See \cite{TopologicalStability}, chapter II, lemma 2.3 for the case of the Whitney condition (b). The case of condition (a) is similar. \qedhere
\end{proof}

\begin{Remark}\label{remark: geometric condition equivalent to Whitney}
    Conversely, it was proven by Trotman in \cite{TrotmanGeometric} that certain transversality conditions analogous to the ones of proposition \ref{proposition: geometric interpretation of Whitney} characterize the Whitney regularity conditions.
\end{Remark}

Now let $M$ be a smooth manifold and $Y, Y'$ two submanifolds of $M$ such that $Y = \overline{Y'} \backslash{Y'}$. Let $Z = \overline{Y'}$ and let us define a stratification of $Z$ by $0 < 1$ by declaring $Z_0 = Y$ and $Z_1=Y'$. The remaining of this section will consist in proving the classical fact that if $(Y,Y')$ satisfies the Whitney regularity condition (b) along $Y$, then $Z$ is cylindrically stratified (in the sense of definition \ref{definition: cylindrically stratified space}). This will take the more precise form of proposition \ref{proposition: cylindrical structure for Whitney stratifications} below. 

In the next section, we will explain how this result is generalized to arbitrary Whitney stratifications (ie with possibly more than two strata). This generalization, combined with the description of spaces of exit paths in cylindrical stratifications developed in section \ref{section: geometric description of spaces of exit paths}, will be used in section \ref{section: equivalence between unbroken trajectories and exit paths} to prove proposition \ref{proposition: equivalence between space of trajectories and space of exit paths}.

We begin by introducing some notation.

\begin{Notation}\label{notation: cylindrical neighborhood for Whitney stratifications}
    
    By theorem \ref{theorem: existence of tubes}, there exists a tubular neighborhood of $Y$ in $M$; we fix one such from now on. Let $\alpha : Y \rightarrow (0,+\infty)$ be a smooth function satisfying the second condition of proposition \ref{proposition: geometric interpretation of Whitney}. By remark \ref{remark: transverse intersection implies submersion}, the restriction
    
    $$
    Y' \cap T_{\alpha} \xrightarrow{\pi \times \rho} Y \times (0,+\infty)
    $$
    
    is a fiber bundle over its image. Let us denote this image by $U$, which can be explicitly described as
    
    $$
    U = \{(y,r) \mid 0< r < \alpha(y)^2\}.
    $$
    
    We denote by
    
    $$
    p : Y' \cap T_{\alpha} \rightarrow U
    $$
    
    this fiber bundle.
    
    Let us introduce the following space
    
    $$
    L_{\alpha/2} = \left\{t \in Y' \cap T \mid \rho(t) = \left(\frac{\alpha(\pi(t))}{2}\right)^2 \right\}.
    $$
    
    The restriction of $\pi$ to $L_{\alpha/2}$, which we denote
    
    $$
    \pi_{\alpha/2} : L_{\alpha/2} \rightarrow Y
    $$
    
    is a fiber bundle over $Y$. Indeed, consider the inclusion
    
    $$
    \fonction{i}{Y}{U}{y}{\left(y,\left(\frac{\alpha(y)}{2}\right)^2 \right).}
    $$
    
    Then $\pi_{\alpha/2}$ can be described as the pullback
    
    $$
    \xymatrix{
    L_{\alpha/2} \ar[r] \ar[d]_{\pi_{\alpha/2}} & Y' \cap T_{\alpha} \ar[d]^{p} \\
    Y \ar[r]^{i} & U.
    }
    $$
    
    Finally let us introduce the following neighborhood of $Z_0$ in $Z$
    
   $$
    \widetilde{V}_{\alpha / 2} = \left\{t \in Z \cap T \mid 0 \le \rho(t) \le \left(\frac{\alpha(\pi(t))}{2}\right)^2 \right\}
    $$
    
    and denote by $U_{\alpha/2}$ the image of $\widetilde{V}_{\alpha / 2}$ by $p$, explicitly described as
    
    $$
    U_{\alpha / 2} = \left\{(y,r) \in U \mid 0 \le r \le \left(\frac{\alpha(y)}{2}\right)^2\right\}.
    $$
    
    We denote by $h$ the homeomorphism
             
     $$
    \begin{array}{ccc}
    Y \times [0,1] & \overset{\simeq}{\longrightarrow} & U_{\alpha/2} \\
    (y,s) & \longmapsto & \left(y,s\left(\frac{\alpha(y)}{2}\right)^2\right).
    \end{array}
    $$
    
    \end{Notation}
    
    \begin{Proposition}\label{proposition: cylindrical structure for Whitney stratifications}
        There exists a homeomorphism between the cylinder of $\pi_{\alpha/2}$ and $\widetilde{V}_{\alpha/2}$ such that the following conditions hold.
        
        \begin{itemize}
            \item It is the identity on $Y$. In particular, $Z$ is cylindrically stratified.
            
            \item It is the identity on $L_{\alpha/2}$, identified with $L_{\alpha/2} \times \{1\} \subset \Cyl(\pi_{\alpha/2})$.
            
            \item The following diagram commutes
            
            $$
            \xymatrix{
            \Cyl(\pi_{\alpha/2}) \ar[r]^-{\simeq} \ar[d]_{R} & \widetilde{V}_{\alpha/2} \ar[d]^{p} \\
            Y \times [0,1] \ar[r]_-{h} & U_{\alpha/2}.
            }
            $$
        \end{itemize}
    \end{Proposition}
    
    \begin{proof}
    
    For every fiber bundle $p : E \rightarrow Y \times (0,1]$, denoting by $F$ the restriction of $E$ to $Y \times \{1\}$, there exists an isomorphism of fiber bundles
    
    $$
    \xymatrix{
    F \times (0,1] \ar[rd]_{p \times \id} \ar[rr]^{\simeq} && E \ar[ld]^p \\
    & Y \times (0,1]
    }
    $$
    
    that restricts to the identity of $F$ over $Y \times \{1\}$.
    
    We apply this to the restriction of $p$ to $U_{\alpha/2}$. We get an isomorphism of fiber bundles
    
    $$
    \xymatrix{
    L_{\alpha/2} \times (0,1] \ar[d]_{\pi_{\alpha/2} \times \id} \ar[r]^-{\simeq}_-{\varphi} & \widetilde{V}_{\alpha/2} \backslash Y \ar[d]^p \\
    Y \times (0,1] \ar[r]_-h & U_{\alpha/2} \backslash (Y \times \{0\}).
    }
    $$
    
    We extend $\varphi$ to a map $Y \cup L_{\alpha/2} \times (0,1] \rightarrow \widetilde{V}_{\alpha/2}$ by declaring it to be the identity on $Y$. This defines a map $\Cyl(\pi_{\alpha/2}) \rightarrow \widetilde{V}_{\alpha/2}$, which one can check to be a homeomorphism. It fulfills the desired conditions. \qedhere
    
    \end{proof}

\subsection{Whitney stratifications}\label{section: Whitney stratifications}

The general definition of a Whitney stratification is the following.

\begin{Definition}\label{definition: Whitney stratified space}
    Let $M$ be a smooth manifold and $Z$ a locally closed subset of $M$. A \emph{Whitney stratification of $Z$} is a partition of $Z$ into nonempty connected smooth submanifolds of $M$ called the \emph{strata}
    
    $$
    Z = \bigsqcup_{s \in S} Z_s
    $$
    
    such that the following two conditions hold:
    
    \begin{itemize}
        \item The family $\{Z_s, s \in S\}$ is locally finite, in the sense that every $z \in Z$ has a neighborhood in $Z$ that intersects only finitely many strata,
        
        \item for every $s, t \in S$ and every $z \in Z_s \cap \overline{Z_t}$, the pair $(Z_s,Z_t)$ satisfies the Whitney regularity condition (b) at $z$.
    \end{itemize}
    
\end{Definition}

\begin{Notation}\label{notation: Whitney stratification by set of strata}
    We will refer to a Whitney stratification of $Z$ by indicating its set of strata $S$.
\end{Notation}

The notion of Whitney stratification is relevant for us because of the following result of Nicolaescu.

\begin{Thms}[Nicolaescu]\label{theorem: Smale if and only if Whitney}

Let $\xi$ be a negative pseudo-gradient adapted to a Morse function $f : X \rightarrow \R$, where $X$ is closed. The following two conditions are equivalent:

\begin{enumerate}[label=(\roman*)]
    \item $\xi$ satisfies the Smale transversality condition,
    
    \item the partition of $X$ by the stable manifolds of $\xi$ is a Whitney stratification of $X$ (seen as a subset of itself).
\end{enumerate}
    
\end{Thms}

\begin{proof}
See \cite{NicolaescuInvitation}, theorems 4.32 and 4.33.
\end{proof}

The purpose of the remaining of this section is to generalize the discussion of section \ref{section: Whitney conditions and cylindrical structures} to Whitney stratifications with possibly more than two strata. To begin with, the following result establishes a connection between the notion of Whitney stratification, and the notion of stratified topological space from section \ref{section: stratified topological spaces and exit paths}.

\begin{Proposition}\label{proposition: Whitney stratified space is stratified space}
    Let $M$ be a smooth manifold and $Z$ a locally closed subset of $M$ endowed with a Whitney stratification $S$. The relation $\le$ on $S$ defined by
     
    $$
    s \le t \Leftrightarrow X_s \cap \overline{X_t} \neq \varnothing
    $$
    
    is a partial order. Moreover, the map $Z \rightarrow S$ that sends each point to the stratum to which it belongs is a stratification of $Z$ by $(S, \le)$, in the sense of definition \ref{definition: stratified space}.
\end{Proposition}

The proof will be based on the following lemma.

\begin{Lemmas}\label{lemma: frontier condition}
    Suppose we have a Whitney stratification $S$ of $Z$ and two elements $s, t \in S$ such that $Z_s \cap \overline{Z_t} \neq \varnothing$. Then $Z_s \subset \overline{Z_t}$.
\end{Lemmas}

\begin{proof}
See \cite{TopologicalStability}, chapter II, corollary 5.7.
\end{proof}

\begin{proof}[Proof of proposition \ref{proposition: Whitney stratified space is stratified space}]
    Let us first prove that $\le$ is a partial order. Reflexivity is immediate and transitivity follows from lemma \ref{lemma: frontier condition}. For antisymmetry, given $s,t \in S$ such that $s \le t$ and $t \le s$, we wish to prove that $t = s$. By lemma \ref{lemma: frontier condition}, we have $Z_t \subset \overline{Z_s}$ and $Z_s \subset \overline{Z_t}$. Since the strata are disjoint, it suffices to prove that $Z_s \cap Z_t \neq \varnothing$. Let $z \in Z_t$. Since $Z_t$ is a smooth submanifold of $M$, there exists a neighborhood $U$ of $z$ in $M$ such that $Z_t \cap U$ is closed in $U$. Since $Z_t \subset \overline{Z_s}$, we have $U \cap Z_s \neq \varnothing$. Since $Z_s \subset \overline{Z_t}$ and $Z_t \cap U$ is closed in $U$, we have $Z_t \cap Z_s \cap U \neq \varnothing$, in particular $Z_t \cap Z_s$ is not empty.
    
    Let us now prove that the partition
    
    $$
    Z = \bigsqcup_{s \in S} Z_s
    $$
    
    is a stratification of $Z$ by the poset $(S, \le)$. We need to prove that, for every $s \in S$, $Z_{\ge s}$ is open in $Z$. Let $t \ge s$ and $z \in Z_t$. There exists a neighborhood of $z$ in $Z$ intersecting a finite number of strata. Hence there also exists a neighborhood $U$ of $z$ in $Z$ having the property that, for every $u \in S$ such that $Z_t \cap \overline{Z_u} = \varnothing$, we have $U \cap Z_u = \varnothing$. Hence for every $u \in S$ such that $U \cap Z_u \neq \varnothing$, we have $u \ge t$. Hence $U \subset Z_{\ge t}$ and therefore also $U \subset Z_{\ge s}$.
\end{proof}

It is an essential feature of Whitney stratifications that the nice properties that hold in the two strata case, which we reviewed in the previous section, extend in the presence of more than two strata. To make this precise, we introduce the following definition. Notice that this definition is given in the topological framework, although in this section we have worked in the smooth category so far.

\begin{Definition}\label{definition: stratified fiber bundle}
    Let $Y$ be a connected topological space. A \emph{stratified fiber bundle} over $Y$ is a topological fiber bundle $p : E \rightarrow Y$ with fiber $F$, such that $E$ and $F$ are both endowed with a stratification and the following holds: there exists an open covering $\Ul$ of $Y$ with local trivializations $p^{-1}(U) \simeq U \times F$ compatible with the stratifications for every $U \in \Ul$ (where $p^{-1}(U)$ is endowed with the stratification induced from that of $E$, and $U \times F$ is endowed with the product stratification).
\end{Definition}

The following is a generalization of Ehresmann's theorem to the Whitney stratified setting, known as Thom's first isotopy lemma. It first appeared in \cite{ThomEnsemblesEtMorphismes}.

\begin{Thms}\label{theorem: Thom's first isotopy lemma}
    Let $M$ be a smooth manifold and $Z$ a locally closed subset of $N$ endowed with a Whitney stratification $S$. Suppose we have a connected smooth manifold $Y$ and a smooth surjective map $g : M \rightarrow Y$ such that for every $s \in S$, $g_{|Z_s}$ is a submersion and $g_{\left| \overline{Z_s} \cap Z \right.}$ is proper (where $\overline{Z_s}$ denotes the closure of $Z_s$ in $M$). Then $g_{|Z}$ is a stratified fiber bundle.
\end{Thms}

\begin{proof}
See \cite{TopologicalStability}, chapter 2, theorem 5.2.
\end{proof}

Using this theorem, one shows the existence of cylindrical structures for Whitney stratifications with many strata in the sense of definition \ref{definition: cylindrical stratification with many strata}, generalizing theorem \ref{proposition: cylindrical structure for Whitney stratifications}. Namely, let $M$ be a smooth manifold and $Z$ a locally closed subset of $M$ endowed with a Whitney stratification $S$. Suppose given $s \in S$ and let $T$ be a tubular neighborhood of $Z_s$ in $M$. Let $\alpha : Z_s \rightarrow (0,+\infty)$ be a smooth function such that for every $z \in Z_s$, $s'>s$ and $0 < r < \alpha(z)$, the sphere $T_{\alpha}(z,r)$ intersects $Z_{s'}$ and this intersection is transverse (we refer to notation \ref{notation: sub tubular neighborhood} for explanation of the notation). Furthermore, since the partition of $Z$ into strata is locally finite, there exists a neighborhood $V$ of $Z_s$ in $Z$ such that $V = Z_{\ge s} \cap V$, and therefore it is possible to choose $\alpha$ such that $Z \cap T_{\alpha} = Z_{\ge s} \cap T_{\alpha}$, which we do. Denote by

$$
p : Z_{>s} \cap T_{\alpha} \rightarrow Z_s \times (0,+\infty)
$$

the map induced by the product of the projection $\pi$ and radius map $\rho$ of $T$. Denote by $U$ the image of $p$. Define further

$$
    L_{\alpha/2} = \left\{t \in Z_{>s} \cap T \mid \rho(t) = \left(\frac{\alpha(\pi(t))}{2}\right)^2 \right\}
    $$
    
and

 $$
    \widetilde{V}_{\alpha / 2} = \left\{t \in Z_{>s} \cap T \mid 0 \le \rho(t) \le \left(\frac{\alpha(\pi(t))}{2}\right)^2 \right\}.
    $$
    
For every $s'>s$, we denote $L_{\alpha/2,s'} = L_{\alpha/2} \cap Z_{s'}$ and $\widetilde{V}_{\alpha/2,s'} = \widetilde{V}_{\alpha/2} \cap Z_{s'}$. Finally, we denote by $\pi_{\alpha/2} : L_{\alpha/2} \rightarrow Z_s$ (resp. $\pi_{\alpha/2,s'} : L_{\alpha/2,s'} \rightarrow Z_s$) the restriction of $\pi$ to $L_{\alpha/2}$ (resp. $L_{\alpha/2,s'}$), postcomposed with the projection to $Z_s$.

\begin{Proposition}\label{proposition: stratified fiber bundles in Whitney stratifications}
    The map $p$ is a stratified fiber bundle over $U$ (where its source is endowed with the $S_{>s}$ stratification induced from that of $Z_{>s}$). Consequently, so is the map $\pi_{\alpha/2,s}$. Therefore, for every $s'>s$, the map $\pi_{\alpha/2,s'}$ is a fiber bundle over $U$.
\end{Proposition}

\begin{proof}
The assertion of the first sentence follows from Thom's first isotopy lemma (theorem \ref{theorem: Thom's first isotopy lemma}). The assertions of the other two sentences are consequences of the first one. \qedhere
\end{proof}

As in the previous section, we denote by $U_{\alpha/2}$ the image of $\widetilde{V}_{\alpha / 2}$ by $p$, explicitly described as
    
    $$
    U_{\alpha / 2} = \left\{(z,r) \in U \mid 0 \le r \le \left(\frac{\alpha(z)}{2}\right)^2\right\}.
    $$
    
    We denote by $h$ the homeomorphism
             
     $$
    \begin{array}{ccc}
    Z_s \times [0,1] & \overset{\simeq}{\longrightarrow} & U_{\alpha/2} \\
    (z,s) & \longmapsto & \left(z,s\left(\frac{\alpha(z)}{2}\right)^2\right).
    \end{array}
    $$

The proof of the following corollary is then obtained in a way analogous to that of proposition \ref{proposition: cylindrical structure for Whitney stratifications}.

\begin{Corollary}\label{corollary: cylindrical structure for general Whitney stratifications}
        There exists a stratified homeomorphism between the cylinder of $\pi_{\alpha/2}$ and $\widetilde{V}_{\alpha/2}$ such that the following conditions hold.
        
        \begin{itemize}
            \item It is the identity on $Z_s$. In particular, $Z$ is cylindrically stratified along $Z_s$.
            
            \item It is the identity on $L_{\alpha/2}$, identified with $L_{\alpha/2} \times \{1\} \subset \Cyl(\pi_{\alpha/2})$.
            
            \item The following diagram commutes
            
            $$
            \xymatrix{
            \Cyl(\pi_{\alpha/2}) \ar[r]^-{\simeq} \ar[d]_{R} & \widetilde{V}_{\alpha/2} \ar[d]^{p} \\
            Z_s \times [0,1] \ar[r]_-h & U_{\alpha/2}.
            }
            $$
        \end{itemize}
    \end{Corollary}
    
    Consequently, for every $s' >s$, a similar statement holds for the $(0<1)$-stratified space $Z_s \cup \widetilde{V}_{\alpha/2,s'}$. In particular, the $(0<1)$-stratified spaces formed by two consecutive strata in a Whitney stratification are cylindrically stratified. More precisely:
    
    \begin{Corollary}\label{corollary: cyclindrical structures for two strata in Whitney}
    For every $s' > s$ there exists a stratified homeomorphism between the cylinder of $\pi_{\alpha/2,s'}$ and $Z_s \cup \widetilde{V}_{\alpha/2,s'}$ such that the following conditions hold.
    
    \begin{itemize}
            \item It is the identity on $Z_s$. In particular, $Z_s \cup \widetilde{V}_{\alpha/2,s'}$ is cylindrically stratified.
            
            \item It is the identity on $L_{\alpha/2,s'}$, identified with $L_{\alpha/2,s'} \times \{1\} \subset \Cyl(\pi_{\alpha/2,s'})$.
            
            \item The following diagram commutes
            
            $$
            \xymatrix{
            \Cyl(\pi_{\alpha/2,s'}) \ar[r]^-{\simeq} \ar[d]_{R} & Z_s \cup \widetilde{V}_{\alpha/2,s'} \ar[d]^{p} \\
            Z_s \times [0,1] \ar[r]_-h & U_{\alpha/2}.
            }
            $$
        \end{itemize}
    \end{Corollary}

Another corollary of proposition \ref{proposition: stratified fiber bundles in Whitney stratifications} is the following.

\begin{Corollary}\label{corollary:  Whitney stratifications are conical}
    The stratified topological space $(Z,S)$ is conically stratified.
\end{Corollary}

\begin{proof}[Outline of proof]

We place ourselves in the setting of proposition \ref{proposition: stratified fiber bundles in Whitney stratifications} and given $z \in Z_s$, we prove that $Z$ is conically stratified at $z$. Denote by $l_{\alpha/2}$ the fiber of $\pi_{\alpha/2}$ over $z$, regarded as a $S_{>s}$-stratified topological space. Let $V$ be a neighborhood of $z$ in $Z_s$ that is diffeomorphic to $\R^k$ for some positive integer $k$. Since $V$ is contractible, the restriction of $\pi_{\alpha/2}$ to $V$ is a trivial stratified fiber bundle. Applying a similar argument as in the proof of corollary \ref{corollary: cylindrical structure for general Whitney stratifications}, we obtain a stratified homeomorphism between an open neighborhood of $z$ in $Z$ and $V \times C(l_{\alpha/2})$. \qedhere

\end{proof}

\begin{Corollary}\label{corollary: morse stratification is conical}
    Let $\xi$ be a Morse-Smale negative pseudo-gradient vector field for a Morse function on a smooth closed manifold $X$. The stratification of $X$ by the stable manifolds of $\xi$ is conical.
\end{Corollary}

\begin{proof}

Combine theorem \ref{theorem: Smale if and only if Whitney} with corollary \ref{corollary:  Whitney stratifications are conical}. \qedhere

\end{proof}

\begin{Remark}\label{remark: cone like description of Whitney stratification}
    The outline of proof of corollary \ref{corollary:  Whitney stratifications are conical} shows that a space $l$ providing a local cone-like description of a Whitney stratification near a point $z \in Z_s$ can be obtained as follows. Consider a neighborhood $V$ of $z$ in $Z_s$ and, given a tubular neighborhood $T$ of $V$ in the ambient manifold, with projection $\pi$ and radial function $\rho$, take $l$ to be the intersection of $Z$ with the sphere of radius $\varepsilon$ centered at $z$, for $\varepsilon > 0$ sufficiently close to $0$. Example \ref{example: example of Morse stratification is conical} presents some pictures illustrating this remark.
\end{Remark}

\begin{Remark}\label{remark: conditions on pair of strata}
    Corollary \ref{corollary:  Whitney stratifications are conical} is an example of a property of a stratification involving possibly more than two strata, that holds only assuming conditions on \emph{pairs} of strata. In this sense, it is analogous to proposition \ref{proposition: Sing_A infinity category in terms of starting point evaluation}.
\end{Remark}

\subsection{The equivalence}\label{section: equivalence between unbroken trajectories and exit paths}

This section is devoted to the proof of proposition \ref{proposition: equivalence between space of trajectories and space of exit paths}.

To begin with, note that if $\Exit(a,b) = \varnothing$, then $\Mi(a,b) = \varnothing$ as well, so the statement holds. We therefore assume from now on that $\Exit(a,b) \neq \varnothing$. In particular we have $a \in \overline{X_b}$ and so, by lemma \ref{lemma: frontier condition} combined with corollary \ref{corollary: morse stratification is conical}, we have $X_a \subset \overline{X_b}$. We first treat the case when $a = b$.

\begin{Lemmas}\label{lemma: equivalence in case a=b}
    Proposition \ref{proposition: equivalence between space of trajectories and space of exit paths} holds in the case when $a=b$.
\end{Lemmas}
    
\begin{proof}
    In that case the map $\iota_{a,b} : \Mi(a,b) \rightarrow \Exit(a,b)$ is identified with the map
    
    $$
    * \rightarrow \Omega_b (X_b)
    $$
    
    associated with the constant loop at $b$. Since $X_b$ is contractible, $\Omega_b (X_b)$ is contractible as well, so this map is a homotopy equivalence, and in particular a weak homotopy equivalence, as desired.
\end{proof}
    
Assume from now on that $a \neq b$. Since $a \in \overline{X_b}$, we have $f(a) > f(b)$. Let us briefly recall the construction of the map $\iota_{a,b} : \Mi(a,b) \rightarrow \Exit(a,b)$ (construction \ref{construction: parametrization of trajectories by values of the function}). Suppose that $\Mi(a,b)$ is nonempty \footnote{It actually follows from the proposition that $\Mi(a,b)$ is nonempty when $\Exit(a,b)$ is nonempty.} and let $\gamma \in \Mi(a,b)$. To $\gamma$ is associated a map $[f(b),f(a)] \rightarrow X$, referred to as the parametrization of $\gamma$ by the values of $f$ (construction \ref{construction: parametrization of trajectories by values of the function}). The exit path $\iota_{a,b}(\gamma)$ is the path obtained by precomposing this map with the affine decreasing bijection $[0,1] \rightarrow [f(b),f(a)]$.
    
For the purpose of the present proof, we slightly generalize this construction as follows. For every $f(b) < s < f(a)$, we consider the restriction to $[s,f(a)] \subset [f(b),f(a)]$ of the parametrization of $\gamma$ by the values of $f$, and we define a path $\iota_{a,b}^s(\gamma)$ by precomposing this restriction with the affine decreasing bijection $[0,1] \rightarrow [s,f(a)]$. This is an exit path starting at $a$ and ending in the $b$-stratum. This therefore defines a continuous map 
    
$$\iota_{a,b}^s : \Mi(a,b) \rightarrow \Exit(a,X_b).$$

Furthermore, we fix a Morse chart at $a$. We assume that the latter identifies an open neighborhood $U$ of $a$ in $X$ with an open ball of $\R^n$ centered at $0$. Denoting by $k$ the Morse index of $a$, this chart identifies the restriction of $f$ to $U$ with the function

$$
(x_1, \hdots, x_n) \mapsto f(a)-\sum_{i=1}^k x_i^2 + \sum_{i=k+1}^n x_i^2.
$$

We endow $U$ with the stratification induced from that of $X$, in particular we denote $U_a = U \cap X_a$ and $U_b = U \cap X_b$.
    
Note that there exists a neighborhood $W$ of $f(a)$ in $[f(b),f(a)]$ such that if $s \in W$, then the path $\iota_{a,b}^s(\gamma)$ is contained in $U$ for every $\gamma$. Fix such a $W$ and assume that $s \in W \backslash \{f(a)\}$. In this case the map $\iota_{a,b}^s$ factors through a map
    
    $$
    \tilde{\iota}_{a,b}^s : \Mi(a,b) \rightarrow \Exit(a,U_b).
    $$
    
    Consider the square
    
    $$
    \xymatrix{
    \Mi (a,b) \ar[r]^-{\iota} \ar[d]_-{\tilde{\iota}_{a,b}^s} & \Exit(a,b) \ar[d]^-{i} \\
    \Exit(a,U_b) \ar[r]_-{j} & \Exit(a,X_b).
    }
    $$

    Here, $U_b \cup \{a\}$ and $X_b \cup \{a\}$ are endowed with the $(0<1)$-stratifications induced by the $A$-stratification of $X$, that is, whose $0$-strata are $\{a\}$. In particular, $\Exit(a,U_b)$ (resp. $\Exit(a,X_b)$) denotes the space of paths starting at $a$ and lying in $U_b$ (resp. $X_b$) at all positive times, and the maps $i$ and $j$ are the natural inclusions.
    
    The proof of proposition \ref{proposition: equivalence between space of trajectories and space of exit paths} follows from the following four lemmas.
    
    \begin{Lemmas}\label{lemma: i weak homotopy equivalence}
    
        The map $i$ is a weak homotopy equivalence.
    
    \end{Lemmas}
    
    \begin{Lemmas}\label{lemma: j weak homotopy equivalence}
    
        The map $j$ is a homotopy equivalence.
    
    \end{Lemmas}
    
    \begin{Lemmas}\label{lemma: square commutes up to homotopy}
    
        This square is homotopy commutative.
    
    \end{Lemmas}
    
    \begin{Lemmas}\label{lemma: iota_s weak homotopy equivalence}
    
        There exists a neighborhood $W'$ of $f(a)$ in $[f(b),f(a)]$ contained in $W$ and such that for every $s \in W' \backslash \{f(a)\}$, the map $\tilde{\iota}_{a,b}^s$ is a weak homotopy equivalence.

    \end{Lemmas}
    
We complete the proof of proposition \ref{proposition: equivalence between space of trajectories and space of exit paths} by proving these lemmas.
    
\begin{proof}[Proof of lemma \ref{lemma: i weak homotopy equivalence}]

Consider the evaluation at the endpoint map

$$
\Exit(a,X_b) \rightarrow X_b.
$$

This is a fibration whose fiber over $b$ is

$$
\Exit(a,b) \overset{i}{\longrightarrow} \Exit(a,X_b).
$$

Since $X_b$ is contractible (proposition \ref{proposition: stable and unstable manifolds are R^k}), $i$ is a weak homotopy equivalence, as desired. \qedhere

\end{proof}

\begin{proof}[Proof of lemma \ref{lemma: j weak homotopy equivalence}]

Apply lemma \ref{lemma: exit paths in a neighborhood} to the space $X_b \cup \{a\}$ endowed with the  $(0 < 1)$-stratification induced by the stratification of $X$. \qedhere

\end{proof}

\begin{proof}[Proof of lemma \ref{lemma: square commutes up to homotopy}]

Note that $i \circ \iota = \iota_{f(b)}$ and $j \circ \tilde{\iota}_{a,b}^s = \iota_{a,b}^s$. The map

\begin{equation*}
\fonctionsansnom{\Mi(a,b) \times [f(b),s]}{\Exit(a,X_b)}{(\gamma,v)}{\iota_{a,b}^v(\gamma).}
\end{equation*}

is a homotopy between $i \circ \iota$ and $j \circ \tilde{\iota}_{a,b}^s$. \qedhere

\end{proof}

\begin{proof}[Proof of lemma \ref{lemma: iota_s weak homotopy equivalence}]

    We note first that it suffices to prove that \emph{there exists} some $s \in W \backslash \{f(a)\}$ such that $\tilde{\iota}_{a,b}^s$ is a weak homotopy equivalence. Indeed, for every $s \le t < f(a)$, the following map is a homotopy between $\tilde{\iota}_{a,b}^s$ and $\tilde{\iota}_{a,b}^{t}$
    
    $$
    \fonctionsansnom{\Mi(a,b) \times [s,t]}{\Exit(a,U)}{(\gamma,v)}{\tilde{\iota}_{a,b}^v(\gamma),}
    $$
    
    and therefore $\tilde{\iota}_{a,b}^{t}$ is a weak homotopy equivalence as well. Hence $W' = [s,f(a)]$ is a neighborhood of $f(a)$ in $[f(b),f(a)]$ that satisfies the desired condition.

   We now aim to prove that there exists such a $s$. Given $s \in W \backslash \{f(a)\}$, let us introduce the following subspace of $X$ 
   
   $$L_s = X_b \cap W^u(a) \cap f^{-1}(s).$$ 
   
    Recall that $W^u(a)$ denotes the \emph{unstable manifold} of $a$ (definition \ref{definition: stable and unstable manifolds}). In other words, $L_s$ is the subspace of those $x \in X$ lying on an unbroken trajectory from $a$ to $b$ and such that $f(x) = s$. This space is naturally homeomorphic to $\Mi(a,b)$ and in the sequel of the proof, we make the identification $\Mi(a,b) = L_s$.
   
   Replacing $W$ by a smaller neighborhood of $f(a)$ in $[f(b),f(a)]$ if necessary, we may assume that for every $s \in W$, we have that $U$ contains $W^u(a) \cap f^{-1}(s)$. The latter is identified via the Morse chart on $U$ with
   
   $$
   \left\{ (x_1,\hdots,x_k,0,\hdots,0) \mid \sum_{i=1}^k x_i^2 = f(a)-s \right\}.
   $$
   
   For every $s \in W \backslash \{f(a)\}$, let us introduce the following neighborhood of $a$ in the unstable manifold of $a$
   
   $$
   D_s = \left\{ (x_1,\hdots,x_k,0,\hdots,0) \mid \sum_{i=1}^k x_i^2 \le f(a)-s \right\}.
   $$
   
   Let us recall the following construction and notation from section \ref{section: geometric description of spaces of exit paths}. Suppose we have a topological space $Z$ and a closed point $z_0 \in Z$, and consider the $(0<1)$-stratification of $Z$ whose $0$-stratum is $\{z_0\}$. Suppose further that there exists a space $L$ and a homeomorphism $h$ between the closed cone of $L$, denoted $C(L)$, and a neighborhood of $z_0$ in $Z$, and this homeomorphism is compatible with the $(0 < 1)$-stratification of $C(L)$ whose $0$-stratum is reduced to the apex of the cone. The map $h$ allows to define the following map
   
   $$
   \fonction{\alpha_h}{L}{\Exit(Z_0,Z_1),}{(l,t)}{h(l,t).}
   $$

   Coming back to the case of the space $U_b \cup \{a\}$ endowed with the $(0 < 1)$-stratification whose $0$-stratum is reduced to the point $a$, we have a map $\tilde{\iota}_{a,b}^s : L_s \rightarrow \Exit(a,U_b)$ which we wish to prove is a weak homotopy equivalence. This map coincides with the map $\alpha_h$ obtained from the following homeomorphism $h$ between $C(L_s)$ and $D_s \cap (X_b \cup \{a\})$:
   
   $$\left\{
   \begin{array}{ccc}
       * & \mapsto & a  \\
       (\gamma,t) \in L_s \times (0,1] & \mapsto & \iota_{a,b}^s(\gamma)(t). 
   \end{array}
   \right.
   $$
   
   Note that this homeomorphism restricts to the identity on $L_s \times \{1\}$. By lemma \ref{lemma: link to exit paths map does not depend on homeomorphism}, for any homeomorphism $g : C(L_s) \simeq D_s \cap (X_b \cup \{a\})$ such that $g_{|L_s \times \{1\}} = \id$ and $g(*) = \{a\}$, the two maps $\iota^s_{a,b}, \alpha_g : L_s \rightarrow \Exit(a,U_b)$ are homotopic. Therefore the condition that $\iota_{a,b}^s$ is a weak homotopy equivalence is equivalent to the condition that \emph{there exists} such a homeomorphism $g$ for which the map $\alpha_g : L_s \rightarrow \Exit(a,U_b)$ is a weak homotopy equivalence.
   
   The sequel of the proof will consist in showing that there exists $s \in W \backslash \{f(a)\}$ for which such a homeomorphim exists. We will do this by combining the properties of cylindrical structures for Whitney stratifications that we have established in section \ref{section: Whitney conditions and cylindrical structures} (in particular proposition \ref{proposition: cylindrical structure for Whitney stratifications}) with the properties of the spaces of exit paths in cylindrically stratified spaces that we have established in section \ref{section: geometric description of spaces of exit paths} (in particular proposition \ref{proposition: space of exit paths as boundary of regular neighborhood}).
   
   From now on, we identify the manifolds $U_a$, $U \cap W^u(a)$ and $U_b$ with submanifolds of $\R^n$ via the Morse chart on $U$. Note that $a$ is identified with $0$, $U_a$ is identified with an open neighborhood of $0$ in $0 \times \R^{n-k} \subset \R^n$, and $U \cap W^u(a)$ is identified with an open neighborhood of $0$ in $\R^k \times 0 \subset \R^n$. Our argument will make use of the following tubular neighborhood (in the sense of definition \ref{definition: tubular neighborhood}) of $U_a$ in $\R^n$. We take the fiber bundle $E$ to be $\R^k \times U_a$ endowed with the projection to $U_a$; we regard $E$ as a subspace of $\R^n$. We endow $E$ with the metric $g$ defined to be the Euclidean metric on the fibers. We take $\varepsilon$ to be constant equal to $1$ and we consider the following neighborhood of $U_a$ in $\R^n$
   
   $$
   T = \left\{ (x_1,\hdots,x_n) \in \R^n \mid \sum_{i = 1}^k x_i^2 < 1 \text{ and } (x_{k+1},\hdots,x_n) \in U_a \right\}.
   $$
   
   Finally, the homeomorphism that we consider between $T$ and
   
   $$
   \{v \in E \mid g(v,v) < 1\}
   $$
   
   is simply the identity. The radial function and the projection associated with this tubular neighborhood are respectively the maps
   
   $$
   \fonction{\rho}{T}{[0, + \infty)}{(x_1,\hdots,x_n)}{\sum_{i=1}^k x_i^2}
   $$
   
   and
   
   $$
   \fonction{\pi}{T}{U_a}{(x_1,\hdots,x_n)}{(x_{k+1},\hdots,x_n).}
   $$
   
  By theorem \ref{theorem: Smale if and only if Whitney}, the stratification of $X$ by the stable manifolds of $\xi$ is Whitney. In particular, its restriction to $U$ is Whitney. We apply corollary \ref{corollary: cyclindrical structures for two strata in Whitney} to this setting (and we refer to section \ref{section: Whitney stratifications} for explanation of the notation).
    
    By corollary \ref{corollary: cyclindrical structures for two strata in Whitney}, there exists a homeomorphism $G : \Cyl(\pi_{\alpha/2,b}) \xrightarrow{\simeq} U_a \cup \widetilde{V}_{\alpha/2,b}$ which is the identity on $U_a$ and on $L_{\alpha/2,b} \times \{1\}$, and such that the following diagram commutes 
    
    $$
    \xymatrix{
    \Cyl(\pi_{\alpha/2,b}) \ar[r]^-{G} \ar[d] & U_a \cup \widetilde{V}_{\alpha/2,b} \ar[d] \\
    U_a \ar[r]_-{=} & U_a.
    }
    $$
    
    Let $s = f(a) - (\frac{\alpha(a)}{2})^2$. The fiber of $\pi_{\alpha/2,b} : L_{\alpha/2,b} \rightarrow U_a$ over $a$ can be described as follows:
    
    $$
    \begin{aligned}
    (\pi_{\alpha/2,b})^{-1}(a) & = \left\{t \in U_b \mid \rho(t) = \left(\frac{\alpha(\pi(t))}{2}\right)^2 \text{ and } \pi(t) = a \right\} \\
    & = \left\{ (t_1, \hdots, t_n) \in U_b \mid \sum_{i=1}^n t_i^2 = \left(\frac{\alpha(a)}{2}\right)^2 \text{ and } t_{k+1} = \hdots = t_n = 0 \right\} \\
    & = \left\{ t \in U_b \cap W^u(a)  \mid \sum_{i=1}^k t_i^2 = \left(\frac{\alpha(a)}{2}\right)^2  \right\} \\
    & = \left\{ t \in U_b \cap W^u(a)  \mid f(t) = s \right\} \\
    & = L_s.
    \end{aligned}
    $$
    
    Similarly, the fiber of the map $\Cyl(\pi_{\alpha/2,b}) \rightarrow U_a$ (resp. $U_a \cup \widetilde{V}_{\alpha/2,b}\rightarrow U_a$) over $a$ is $C(L_s)$ (resp. $D_s \cap (X_b \cup \{a\})$). In particular, the homeomorphism $G$ restricts to a homeomorphism $g : C(L_s) \xrightarrow{\simeq} D_s \cap (X_b \cup \{a\})$ such that $g_{|L_s \times \{1\}} = \id$ and $g(*) = \{a\}$. We conclude the proof by showing that this homeomorphism fulfills the desired condition, ie that the associated map
    
    $$
    \alpha_g : L_s \rightarrow \Exit(a,U_b)
    $$
    
    is a weak homotopy equivalence.
    
    Consider the following commutative diagram
    
    $$
    \xymatrix{
    L_s \ar[r]^-{\alpha_g} \ar[d] & \Exit(a,U_b) \ar[d] \\
    L_{\alpha/2,b} \ar[d]_-{\pi_{\alpha/2,b}} \ar[r]^-{\alpha_G} & \Exit(U_a, U_b) \ar[d]^-{\ev_0} \\
    U_a \ar[r]_-{=} & U_a.
    }
    $$
    
    Corollary \ref{corollary: cyclindrical structures for two strata in Whitney} and proposition \ref{proposition: space of exit paths as boundary of regular neighborhood} combined show that the map $\alpha_G$ is a weak homotopy equivalence. By proposition \ref{proposition: stratified fiber bundles in Whitney stratifications}, the left vertical column is a fiber bundle. Combining proposition \ref{proposition: the stratification by stable manifolds is conical}, theorem \ref{theorem: stratified simplicial set of conically stratified space is infty cat} and proposition \ref{proposition: Sing_A infinity category in terms of starting point evaluation}, we obtain that the right column is a Serre fibration. Recall that we assumed that the Morse chart on $U$ identifies $U$ with an open ball centered at $0$ in $\R^n$, in particular $U_a$ is contractible. Altogether, this proves that $\alpha_g$ is a weak homotopy equivalence, as desired. \qedhere

\end{proof}

The proof of proposition \ref{proposition: equivalence between space of trajectories and space of exit paths} is complete. \hfill\qed

As a consequence we have the following equivalent characterizations of the Smale ordering on the set of critical points of $f$.

\begin{Corollary}\label{corollary: equivalent characterizations of partial order on A}
Let $a_0, a_1$ be two critical points of $f$. The following conditions are equivalent:

\begin{enumerate}[label=(\roman*)]

    \item $a_0 \le a_1$ in the sense of definition \ref{definition: partial order on set of critical points}, i.e., there exists an unbroken trajectory of $\xi$ from $a_0$ to $a_1$,
    
    \item there exists an exit path from $a_0$ to $a_1$,
    
    \item the $a_0$-stratum of $X$ intersects the closure of the $a_1$-stratum of $X$.
    
     \item the $a_0$-stratum of $X$ is contained in the closure of the $a_1$-stratum of $X$.
    
\end{enumerate}

\end{Corollary}

\begin{proof}
    The equivalence (i) $\Leftrightarrow$ (ii) follows from proposition \ref{proposition: equivalence between space of trajectories and space of exit paths}. The equivalence (ii) $\Leftrightarrow$ (iii) is true for any conically stratified space, and therefore follows from corollary \ref{corollary: morse stratification is conical}. The equivalence (iii) $\Leftrightarrow$ (iv) follows from theorem \ref{theorem: Smale if and only if Whitney} together with lemma \ref{lemma: frontier condition}.
\end{proof}

\newpage

\section{The homotopy coherent nerve of a topological category}\label{section: the homotopy coherent nerve}

In this section, we define a functor $F : \sSet \rightarrow \Cat_{\Top}$ called \textit{the geometric realization of the Leitch-Cordier functor}. Our motivation is that the \textit{homotopy coherent nerve functor} $\Nl : \Cat_{\Top} \rightarrow \sSet$, introduced in section \ref{section: infinity categories}, is right adjoint to $F$. In particular for every topological category $\Cj$ we have by adjunction

$$
\begin{aligned}
    \Nl(\Cj)_n & = \Hom_{\sSet}(\Delta^n,\Nl(\Cj)) \\
    & = \Hom_{\Cat_{\Top}}(F(\Delta^n),\Cj).
\end{aligned}
$$

In section \ref{section: definition of the homotopy coherent nerve} we define the functor $F$ and discuss the adjunction between $F$ and $\Nl$. The topological categories $F(\Delta^n)$ are defined in a somehow unenlightening way, as the geometric realization of simplicial categories whose morphism simplicial sets are nerves of posets. In section \ref{section: Study of the geometric realization of the Leitch-Cordier functor}, we study $F$ in more detail, and give a more explicit description of the simplices of the homotopy coherent nerve of a topological category. One consequence of our discussion is that an $n$-simplex of $\Nl(\Cl)$ is a \emph{homotopy coherent composition} of $n$ composable morphisms of $\Cj$.

\subsection{Definition}\label{section: definition of the homotopy coherent nerve}

We begin by defining the Leitch-Cordier functor, following \cite[Section 1.1.5]{HigherTopos}. We start by defining it at the level of the category $\Delta$.

Recall from definition \ref{definition: nerve of a poset} that to every poset is associated a category, and the the nerve of the poset is the nerve of this category. Moreover, the associated functor $\Pos \rightarrow \Cat$ as well as the nerve functor $\Cat \rightarrow \sSet$ both commute with products. The nerve functor on posets $\Pos \rightarrow \sSet$, which is the composition of the two latter, therefore also commutes with products.

\begin{Definition}\label{definition: simplicially enriched categories associated with simplices}
    Let $n \ge 0$ be an integer. We define a functor $\Cr : \Delta \rightarrow \Cat_{\sSet}$ as follows.

\begin{itemize}
\item[$\bullet$] The set of objects of $\Cr (\Delta^n)$ is $\{0, \hdots, n\}$.
\item[$\bullet$] For $0 \le i,j \le n$,

$$
\Cr(\Delta^n)(i,j)=
\left\{
    \begin{array}{ll}
        \varnothing & \mbox{if } i > j, \\
        \text{nerve of } P_{n,i,j} & \mbox{if } i \le j,
    \end{array}
\right.
$$

where $P_{n,i,j}$ denotes the set $\left\{ E \subset \{0,\hdots,n\} \text{ : } i,j \in E \text{ and } \forall k \in E, \, i \le k \le j \right\}$ partially ordered by the inclusion.
        
\item[$\bullet$] The composition

$$
\Cr(\Delta^n)(i_0,i_1)\times \hdots \times \Cr(\Delta^n)(i_{k-1},i_k) \rightarrow \Cr(\Delta^n)(i_0,i_k)
$$

is defined to be the nerve of the map of posets

$$
\fonctionsansnom{P_{n,i_0,i_1}\times \hdots \times P_{n,i_{k-1},i_k}}{P_{n,i_0,i_k},}{(I_1,\hdots,I_k)}{I_1 \cup \hdots \cup I_k.}
$$

\end{itemize}

Let $f : \Delta^n \rightarrow \Delta^m$ be an arrow in $\Delta$. The functor between simplicially enriched categories $\Cr(f) : \Cr(\Delta^n) \rightarrow \Cr(\Delta^m)$ is defined as follows:

\begin{itemize}
\item[$\bullet$] It is equal to $f$ on the objects.

\item[$\bullet$] For $0 \le i \le j \le n$, the morphism of simplicial sets

$$
\Cr(f)(i,j) : \Cr(\Delta^n)(i,j) \rightarrow \Cr(\Delta^m)(f(i),f(j))
$$

is defined to be the nerve of the map of posets

$$
\fonctionsansnom{P_{n,i,j}}{P_{m,f(i),f(j)},}{I}{f(I).}
$$

\end{itemize}
\end{Definition}

\begin{Notation}\label{notation: poset enriched category associated with simplices}
    For later use, we record the following notation. For every integer $n \ge 0$, the simplicially enriched category $\Cr(\Delta^n)$ from definition \ref{definition: simplicially enriched categories associated with simplices} is defined by applying the nerve functor to the posets of morphisms of a category enriched over posets. We denote the latter by $\PR(\Delta^n)$.
\end{Notation}

\begin{Definition}\label{definition: Leitch Cordier functor}

As the category of small simplicial categories admits all small colimits (appendix \ref{appendix: colimits of simplicial categories}), the functor $\Cr$ admits a left Kan extension to $\sSet = \Pre(\Delta)$ (construction \ref{construction: left Kan extension}). We still denote it $\Cr$ and refer to it as the \textit{Leitch-Cordier functor}.

\end{Definition}

The terminology comes from the fact that the functor $\Cr$ was constructed by Leitch in $\cite{LeitchHomotopy}$ and further studied by Cordier in $\cite{CordierDiagramme}$.

\begin{Notation}\label{notation: right adjoint of Leitch Cordier functor}
     The right adjoint of the Leitch-Cordier functor (construction \ref{construction: right adjoint of left Kan extension}) is denoted $\mathscr{N}$.
\end{Notation}

Let $C$ be a simplicial category. Since the geometric realization functor commutes with products when considered with values in the category of compactly generated topological spaces, there is a topological category $|C|$ associated to $C$, defined by

$$\left\{
    \begin{array}{ll}
        \Obj (|C|) = \Obj (C), \\
        |C|(X,Y) = |C(X,Y)|.
    \end{array}
\right.$$

This defines a functor $|-| : \Cat_{\sSet} \rightarrow \Cat_{\Top}$. Moreover, since the singular simplicial set functor commutes with products when taken with source the category of compactly generated topological spaces, to every topological category $\Cj$ is associated a simplicial category $\Sing(\Cj)$ defined by

$$\left\{
    \begin{array}{ll}
        \Obj (\Sing(\Cj)) = \Obj (\Cj), \\
        \Sing(\Cj)(X,Y) = \Sing(\Cj(X,Y)).
    \end{array}
\right.$$

The singular-realization adjunction between $\Top$ and $\sSet$ then holds between $\Cat_{\Top}$ and $\Cat_{\sSet}$ as well.

\begin{Definition}\label{definition: geometric realization of the Leitch-Cordier functor}
    The composition

    $$
    |-| \circ \Cr : \sSet \rightarrow \Cat_{\Top}
    $$

    is called the \emph{geometric realization of the Leitch-Cordier functor} and is denoted by $F$.
\end{Definition}

It admits the following right adjoint:

\begin{Definition}\label{definition: homotopy coherent nerve}
 The composition

    $$
    \mathscr{N} \circ \Sing : \Cat_{\Top} \rightarrow \sSet
    $$

    is called the \emph{homotopy coherent nerve functor} and is denoted by $\Nl$.
\end{Definition}

We will study the topological categories $F(\Delta^n)$ in detail in section \ref{section: Study of the geometric realization of the Leitch-Cordier functor}. In the sequel of the paper, we will also need descriptions of the images by $F$ of certain simplicial sets different from standard simplices (namely in the proof of lemma \ref{lemma: image of an inner horn by F}). To that effect, we will use the description of the restriction of $F$ to $\Delta$ given in section \ref{section: Study of the geometric realization of the Leitch-Cordier functor} below, together with the fact that $F$ commutes with colimits.

As explained in the introduction of this section, the adjunction also implies that for every topological category $\Cj$ we have the following formula for $\Nl(\Cj)$:

$$
\fonction{\Nl(\Cj)}{\Delta^{\mathrm{op}}}{\Set,}{\Delta^n}{\Hom_{\Cat_{\Top}}(F(\Delta^n),\Cj).}
$$

An $n$-simplex of $\Nl(\Cj)$ can be thought of as a homotopy coherent composition of $n$ composable morphisms in $\Cj$. In the next section, we give properties of $F$ that make this idea precise.

\subsection{Study of the geometric realization of the Leitch-Cordier functor}\label{section: Study of the geometric realization of the Leitch-Cordier functor}

Let us denote by $I$ the interval $[0,1]$. We begin by defining some homeomorphisms that we will use throughout the paper:

\[\label{homeomorphisms}
\tag{H}
F(\Delta^n)(i,j) \simeq \left\{
    \begin{array}{ll}
        * & \mbox{if } i=j, \\
        I^{j-i-1} & \mbox{if } i < j.
    \end{array}
\right.
\]

By definitions \ref{definition: Leitch Cordier functor} and \ref{definition: geometric realization of the Leitch-Cordier functor}, $F(\Delta^n)(i,i)$ and $F(\Delta^n)(i,i+1)$ are both the geometric realization of a trivial simplicial set. Assume $i + 1 < j$ and consider the poset $P_{n,i,j}$ from definition \ref{definition: simplicially enriched categories associated with simplices}. For all $i < k < j$ consider the morphism of posets

$$
\fonction{p_k}{P_{n,i,j}}{0 < 1,}{E}{\left\{
    \begin{array}{ll}
        0 & \mbox{if } k \notin E \\
        1 & \mbox{if } k \in E.
    \end{array}
\right.}
$$

The product of the $p_k$'s is then an isomorphism of posets. Hence after passing to the nerve, and then to the geometric realization, we get a homeomorphism

$$
F(\Delta^n)(i,j) \simeq \prod_{k=i+1}^{j-1} |N(0 < 1)|.
$$

Note that this product is taken in the category of compactly generated topological spaces (remark \ref{remark: product in the category of compactly generated spaces}), but in this case there is no difference with the associated usual product in the category of topological spaces.

Now, $N(0<1)$ is canonically isomorphic to $\Delta^1$, and therefore $|N(0 < 1)|$ is canonically identified with $|\Delta^1|$, which by definition \ref{definition: geometric realization at the level of simplex category} is

$$
|\Delta^1| = \left\{ (t_0,t_1) \in \R^2 \mid t_0 + t_1 = 1 \text{ and } t_0, t_1 \ge 0 \right\}.
$$

We define a homeomorphism between $I$ and $|N(0 < 1)|$ by

$$
\fonctionsansnom{I}{|N(0 < 1)|,}{t}{(1-t,t).}
$$

This defines the homeomorphisms \eqref{homeomorphisms}.

\begin{Proposition}\label{proposition: image of the faces by F}
Let $0 \le i \le n$ and $\delta_n^i : \Delta^{n-1} \rightarrow \Delta^n$ be the corresponding face morphism of $\Delta$ (notation and terminology \ref{notation and terminology about simplicial sets}). Under the homeomorphisms \eqref{homeomorphisms}, the map $F(\delta_n^i)(0,n-1) : F(\Delta^{n-1})(0,n-1) \rightarrow F(\Delta^n)(\delta_n^i(0),\delta_n^i(n-1))$ is identified with the inclusion of the face $\{t_i = 0\}$ for $0 < i < n$, and is the identity for $i = 0$ and $i = n$.
\end{Proposition}

\begin{proof}
For $0 < i < n$, the map $F(\delta_n^i)(0,n-1) : F(\Delta^{n-1})(0,n-1) \rightarrow F(\Delta^n)(0,n)$ induces under \eqref{homeomorphisms} a map $I^{n-2} \rightarrow I^{n-1}$ and the statement is equivalent to the commutativity of the following three diagrams:

$$
\xymatrix{
I^{n-2} \ar[r] \ar[rd]_{\equiv 0} & I^{n-1} \ar[d]^{t_i} \\
& I,}
$$

$$
\text{for } l<i:
\xymatrix{
I^{n-2} \ar[r] \ar[rd]_{t_l} & I^{n-1} \ar[d]^{t_l} \\
& I,}
$$

$$
\text{for } l>i:
\xymatrix{
I^{n-2} \ar[r] \ar[rd]_{t_{l-1}} & I^{n-1} \ar[d]^{t_l} \\
& I.}
$$

This follows respectively from the commutativity of the following diagrams of posets, where we use the notation from the definition of the homeomorphisms \eqref{homeomorphisms}:

$$
\xymatrix{
P_{n-1,0,n-1} \ar[r]^{\delta_n^i} \ar[rd]_{\equiv 0} & P_{n,0,n} \ar[d]^{p_i} \\
& 0 < 1,}
$$

$$
\text{for } l<i:
\xymatrix{
P_{n-1,0,n-1} \ar[r]^{\delta_n^i} \ar[rd]_{p_l} & P_{n,0,n} \ar[d]^{p_l} \\
& 0 < 1,}
$$

$$
\text{for } l > i:
\xymatrix{
P_{n-1,0,n-1} \ar[r]^{\delta_n^i} \ar[rd]_{p_{l-1}} & P_{n,0,n} \ar[d]^{p_l} \\
& 0 < 1.}
$$

The proof in the cases $i = 0$ and $i = n$ is similar.
\end{proof}

\begin{Remark}\label{map induced by an inclusion of consecutive elements}
    The cases $i=0$ and $i = n$ in proposition \ref{proposition: image of the faces by F} imply that for every $0 \le i \le j \le n$, the morphism

    $$
    \fonction{\theta}{\Delta^{j-i}}{\Delta^n}{k}{k + i}
    $$

    induces a map $F(\theta)(0,j-i)$ which is identified by the homeomorphisms \eqref{homeomorphisms} with the identity map.
\end{Remark} 

\begin{Proposition}\label{proposition: compositions in F(n)}
    The composition in $F(\Delta^n)$ is described as follows: for $0 < k < n$, under the homeomorphisms \eqref{homeomorphisms}, the composition map $F(\Delta^n)(0,k) \times F(\Delta^n)(k,n) \rightarrow F(\Delta^n)(0,n)$ is identified with the inclusion of the face $\{t_k=1\}$ of $I^{n-1}$.

    More generally, for $0 \le i < k <j \le n$, under the homeomorphisms \eqref{homeomorphisms}, the composition map $F(\Delta^n)(i,k) \times F(\Delta^n)(k,j) \rightarrow F(\Delta^n)(i,j)$ is identified with the inclusion of the face $\{t_{k-i}=1\}$ of $I^{j-i-1}$.
\end{Proposition}

\begin{proof}

By remark \ref{map induced by an inclusion of consecutive elements}, the second part of the proposition is a consequence of the first part. Let us prove the first part. We use the notation from the definition of the homeomorphisms \eqref{homeomorphisms}.

The composition map $F(\Delta^n)(0,k) \times F(\Delta^n)(k,n) \rightarrow F(\Delta^n)(0,n)$ induces under \eqref{homeomorphisms} a map $I^{n-2} \rightarrow I^{n-1}$ and the statement is equivalent to the commutativity of the following three diagrams:

$$
\xymatrix{
I^{n-2} \ar[r] \ar[rd]_{\equiv 1} & I^{n-1} \ar[d]^{t_k} \\
& I,}
$$

$$
\text{for } l<k:
\xymatrix{
I^{n-2} \ar[r] \ar[rd]_{t_l} & I^{n-1} \ar[d]^{t_l} \\
& I,}
$$

$$
\text{for } l>k:
\xymatrix{
I^{n-2} \ar[r] \ar[rd]_{t_{l-1}} & I^{n-1} \ar[d]^{t_l} \\
& I.}
$$

The inclusion of $\{1\}$ in $I$ is induced, after passing to the nerve and then to the geometric realization, by the inclusion of posets $1 \hookrightarrow 0 < 1$. Hence, the commutativity of the first diagram follows from the commutativity of the following diagram of posets:

$$
\xymatrix{
P_{n,0,k} \times P_{n,k,n} \ar[r]^-{\cup} \ar[rd]_{\equiv 1} & P_{n,0,n} \ar[d]^{p_k} \\
& 0 < 1.}
$$

The commutativity of the remaining two diagrams follows respectively from the commutativity of the following two diagrams of posets:

$$
\text{for } l<k:
\xymatrix{
P_{n,0,k} \times P_{n,k,n} \ar[r]^-{\cup} \ar[rd]_{p_l} & P_{n,0,n} \ar[d]^{p_l} \\
& 0 < 1,}
$$

$$
\text{for } l > k:
\xymatrix{
P_{n,0,k} \times P_{n,k,n} \ar[r]^-{\cup} \ar[rd]_{p_{l}} & P_{n,0,n} \ar[d]^{p_l} \\
& 0 < 1.}
$$

\end{proof}

\begin{Definition}\label{definition: broken and unbroken morphism of F(n)}
For $0 \le i < k < j \le n$, a morphism $\gamma \in F(\Delta^n)(i,j)$ is said to be \emph{broken} at $k$ if it can be written as the composition of an element of $F(\Delta^n)(k,j)$ with an element of $F(\Delta^n)(i,k)$. It is said to be \textit{unbroken} if it cannot be written as a composition of two morphisms that are both different from the identity.
\end{Definition}

\begin{Corollary}\label{corollary: writing as composition of unbroken morphisms}
    Every morphism of $F(\Delta^n)$ can be written in a unique way as a composition of unbroken morphisms.
\end{Corollary}

\begin{proof}
    By remark \ref{map induced by an inclusion of consecutive elements}, it suffices to treat the case of a morphism $\gamma \in F(\Delta^n)(0,n)$. Using the homeomorphism $F(\Delta^n)(0,n) \simeq I^{n-1}$ from \eqref{homeomorphisms} we write $\gamma = (t_1,\hdots,t_{n-1})$. By proposition \ref{proposition: compositions in F(n)}, letting $0 < i_1 < \hdots < i_k < n$ be the $i$'s such that $t_i = 1$, $i_{k+1} := n$, and $\gamma_j \in F(\Delta^n)(i_j,i_{j+1})$ the preimage of $(t_{i_j}, \hdots, t_{i_{j+1}})$ by the homeomorphism $F(\Delta^n)(i_j,i_{j+1}) \simeq I^{i_{j+1}-i_j-1}$, we have that $\gamma_j$ is unbroken and $\gamma = \gamma_k \circ \hdots \circ \gamma_0$. Moreover, again by proposition \ref{proposition: compositions in F(n)}, any decomposition of $\gamma$ as a composition of unbroken morphisms is uniquely determined by these properties.
\end{proof}

Because of proposition \ref{proposition: compositions in F(n)}, we can think of an $n$-simplex of $\Nl(\Cj)$ as a \emph{homotopy coherent composition} of $n$ composable morphisms of $\Cj$. We now explain this in detail for $n = 0$ to $n = 3$.

$n=0$: $F(\Delta^0)$ is the topological category with one object and no morphism different from the identity, so a $0$-simplex of $\Nl(\Cj)$ is the datum of an object of $\Cj$.

$n=1$: $F(\Delta^1)$ is the topological category with two objects and one morphism between them. Hence a $1$-simplex of $\Nl(\Cj)$ is the datum of two objects of $\Cj$ and a morphism between them.

$n=2$: A $2$-simplex $\sigma \in \Nl(\Cj)_2$ is the datum of:

\begin{itemize}
\item[$\bullet$] Three objects $X$, $Y$ and $Z$, the vertices of $\sigma$, respectively equal to the images of the three objects $0$, $1$ and $2$ of $F(\Delta^2)$.

\item[$\bullet$] Three morphisms, the faces of $\sigma$, as in the following diagram.

\begin{center}

\begin{tikzpicture}

    \node (X) at (-1.2,0) {$X$};
    \node (Y) at (1.2,0) {$Y$};
    \node (Z) at (0,2) {$Z$};
    \draw[->] (X)--(Y) node[midway, below] {$\delta$};
    \draw[->] (X)--(Z);
    \draw[->] (Y)--(Z);
    \draw node at (0.9,1.2) {$\gamma$};
    \draw node at (-0.9,1.2) {$\theta$};
        
\end{tikzpicture}

\end{center}

The remaining datum can be thought of as making this diagram commutative \emph{up to coherent homotopy} in $\Cj$.

\item[$\bullet$] By proposition \ref{proposition: compositions in F(n)} and proposition \ref{proposition: image of the faces by F}, a continuous map $h : I \rightarrow \Cj(X,Z)$ such that $h(0)= \theta$ and $h(1)=\gamma \circ \delta$, that is, a homotopy between $\theta$ and $\gamma \circ \delta$.
\end{itemize}

$n=3$: A $3$-simplex $\sigma \in \Nl(\Cj)_3$ is the datum of:

\begin{itemize}

\item[$\bullet$] Four objects $X$, $Y$, $Z$ and $T$, the vertices of $\sigma$, respectively equal to the images of the three objects $0$, $1$, $2$ and $3$ of $F(\Delta^3)$.

\item[$\bullet$] Six morphisms in $\Cj$ as in the following diagram.

\begin{center}
    
\begin{tikzpicture}
        \draw node at (0,-0.2) {$\alpha$};
        \draw node at (0.3,1) {$\eta$};
        \draw node at (-0.8,1.2) {$\mu$};
        \draw node at (1.1,1.5) {$\gamma$};
        \draw node at (0.1,0.4) {$\nu$};
        \draw node at (1.7,0.1) {$\beta$};
        \node (X) at (-1.2,0) {$X$};
        \node (Y) at (1.2,0) {$Y$};
        \node (T) at (0,2) {$T$};
        \node (Z) at (1.9,0.6) {$Z$};
        
        \draw[->] (X)--(Y);
        \draw[->] (X)--(Z);
        \fill[white] (0.95,0.416) circle (1.5pt);
        \draw[->] (Y)--(Z);
        \draw[->] (X)--(T);
        \draw[->] (Y)--(T);
        \draw[->] (Z)--(T);
\end{tikzpicture}
    
\end{center}

The remaining data can be thought of as making this diagram commutative \emph{up to coherent homotopy} in $\Cj$.

\item[$\bullet$] Four $2$-simplices, the faces of $\sigma$, whose boundaries are the four triangles in the above diagram. They come with four continuous maps:

\begin{itemize}

\item a homotopy $h_1 : I \rightarrow \Cj(X,Z)$ between $\eta$ and $\beta \circ \alpha$,

\item a homotopy $h_2 : I \rightarrow \Cj(Y,T)$ between $\theta$ and $\gamma \circ \beta$,

\item a homotopy $h_3 : I \rightarrow \Cj(X,T)$ between $\mu$ and $\gamma \circ \eta$,

\item a homotopy $h_4 : I \rightarrow \Cj(X,T)$ between $\mu$ and $\theta \circ \alpha$.

\end{itemize}

\item[$\bullet$] By proposition \ref{proposition: compositions in F(n)} and proposition \ref{proposition: image of the faces by F}, a continuous map $c : I^2 \rightarrow \Cj (X,T)$ whose restriction to the boundary of $I^2$ is as in the following picture. Here, $h_{\alpha}$ and $h_{\gamma}$ denote the constant maps from $I$ to $\alpha$ and $\gamma$ respectively.

\begin{center}

\begin{tikzpicture}
    \draw (-1.5,1.5)--(1.5,1.5)--(1.5,-1.5)--(-1.5,-1.5)--cycle;
    \draw node at (-1.7,1.7) {$\gamma \circ \eta$};
    \draw node at (0,1.7) {$h_{\gamma} \circ h_1$};
    \draw node at (1.9,1.7) {$\gamma \circ \beta \circ \alpha$};
    \draw node at (2.1,0) {$h_2 \circ h_{\alpha}$};
    \draw node at (1.7,-1.7) {$\theta \circ \alpha$};
    \draw node at (0,-1.7) {$h_4$};
    \draw node at (-1.7,-1.7) {$\mu$};
    \draw node at (-1.7,0) {$h_3$};
    \end{tikzpicture}

\end{center}

\end{itemize}

Let us finally give a description of the maps induced by the degeneracy maps in $\Delta$, under the homeomorphisms \eqref{homeomorphisms}.

\begin{Proposition}\label{proposition: image of codegeneracies by F(n)}
    Let $n \ge 2$, $0 \le i \le n$ and $\sigma_n^i : \Delta^{n+1} \rightarrow \Delta^n$ be the corresponding degeneracy map in $\Delta$ (notation and terminology \ref{notation and terminology about simplicial sets}). Under the homeomorphisms \eqref{homeomorphisms} the map $F(\sigma_n^i)(0,n+1) : F(\Delta^{n+1})(0,n+1) \rightarrow F(\Delta^n)(0,n)$ is described as follows

    $$
    \fonctionsansnom{I^n}{I^{n-1}}{(t_1,\hdots,t_n)}{\left\{
    \begin{array}{ll}
         (t_2,\hdots,t_n) & \mbox{if } i=0,  \\
         (t_1,\hdots,t_{i-1},\textup{max}(t_i,t_{i+1}),t_{i+2},\hdots,t_n) & \mbox{if } 0 < i < n, \\
         (t_1,\hdots,t_{n-1}) & \text{if } i=n.
    \end{array} 
    \right.}
    $$
    
\end{Proposition}

\begin{proof}
    The proof is similar to the proofs of propositions \ref{proposition: compositions in F(n)} and \ref{proposition: image of the faces by F}, using the fact that the following morphism of posets

    $$
    \begin{array}{ccc}
(0 < 1 ) \times (0 < 1) & \longrightarrow & 0 < 1 \\
(0,0) & \longmapsto & 0 \\
(0,1) & \longmapsto & 1 \\
(1,0) & \longmapsto & 1 \\
(1,1) & \longmapsto & 1 \end{array}
    $$

    induces, after passing to the geometric realization and identifying $|\Delta^1|$ and $I$, the following map

$$
\fonctionsansnom{I^2}{I}{(t_1,t_2)}{\text{max}(t_1,t_2).}
$$
    
\end{proof}

\begin{Corollary}\label{corollary: image of an unbroken morphism}
    Let $\varphi : \Delta^n \rightarrow \Delta^m$ be a morphism in $\Delta$ and $0 \le i \le j \le n$. Let $\gamma \in F(\Delta^n)(i,j)$ be an unbroken morphism (in the sense of definition \ref{definition: broken and unbroken morphism of F(n)}). Then $F(\varphi)(i,j)(\gamma) \in F(\Delta^m)(\varphi(i),\varphi(j))$ is unbroken.
\end{Corollary}

\begin{proof}
    Since face and degeneracy maps generate the morphisms of $\Delta$, it suffices to treat the case when $\varphi$ is a face or degeneracy map, which we assume from now on. By remark \ref{map induced by an inclusion of consecutive elements}, it suffices to treat the case when $i=0$ and $j=n$, which we also assume from now on. First, if $\varphi(n) = \varphi(0)$ then the result follows from the fact that $F(\Delta^m)(\varphi(0),\varphi(n)) = \{\id\}$. Otherwise, we have $n > 0$ and $\varphi(n) > \varphi(0)$, and using the homeomorphisms \eqref{homeomorphisms} we write $\gamma = (t_1,\hdots,t_{n-1})$ and $\varphi(0,n)(\gamma) = (s_{\varphi(0)},\hdots,s_{\varphi(n)-1})$. By propositions \ref{proposition: image of the faces by F} and \ref{proposition: image of codegeneracies by F(n)}, if $t_i \neq 1$ for all $1 \le i \le n-1$ then $s_j \neq 1$ for all $\varphi(0) \le i \le \varphi(n)-1$. The desired result then follows from proposition \ref{proposition: compositions in F(n)}.
\end{proof}

\newpage

\section{\texorpdfstring{The cubical stratified geometric realization}{The cubical stratified geometric realization}}\label{section: cubical stratified geometric realization}

Let $(f,\xi)$ be a Morse-Smale pair on a smooth closed manifold $X$, let $A$ be the set of critical points of $f$ and $\Ml$ the associated flow category (see section \ref{section: Morse-Smale pairs and flow categories} for the definitions of these objects). Given $\sigma \in \Nl(\Ml)$, denote $\sigma (i) = a_i$. Identifying $F(\Delta^n)(0,n)$ and $I^{n-1}$ using the homeomorphisms \eqref{homeomorphisms}, we get a map $\sigma(0,n) : I^{n-1} \rightarrow \Ml(a_0,a_n)$. Moreover, any trajectory $\gamma \in \Ml(a_0,a_n)$ admits a natural parametrization by the interval $[f(a_n),f(a_0)]$, denoted $s \mapsto \gamma(s)$, and defined by the condition that $f(\gamma(s))=s$ (constructions \ref{construction: parametrization of trajectories by values of the function} and \ref{construction: parametrization of all trajectories by values of the function}). This defines an evaluation map

$$
\Ml(a_0,a_n) \times [f(a_n),f(a_0)] \longrightarrow X.
$$

\begin{Notation}\label{notation: map f sigma}
    Combining these two maps we get a map

    $$
    \fonction{f_{\sigma}}{ I^{n-1} \times [f(a_n),f(a_0)]}{X}{(\tu,s)}{\sigma(0,n)(\tu)(s).}
    $$
\end{Notation}

When trying to pass from the singular cube $f_{\sigma}$ in $X$ to a singular stratified $n$-simplex in $X$, we encounter two difficulties.

\begin{enumerate}[label=(\roman*)]

    \item  As explained in section \ref{section: detailed summary of the proof}, $f_{\sigma}$ might visit more intermediary strata between $a_0$ and $a_n$ than only the $X_{a_i}$'s. This is because nothing prevents the image of an unbroken morphism of $F(\Delta^n)$ by $\sigma$ to be a broken trajectory of $\xi$.

    \item The domain of $f_{\sigma}$ cannot be identified with a simplex on the nose.
    
\end{enumerate}

Recall that we endow $A$ with the partial order given by the relation of being connected by a trajectory of $\xi$, and we denote by $\Delta_A$ the simplex category of the poset $A$ (see \ref{definition: partial order on set of critical points} and \ref{definition: simplex category of a poset} for precise definitions). In this section, we resolve difficulty (i) and start addressing difficulty (ii) as follows.

\begin{itemize}
    \item In section \ref{section: non stratified functor C A}, we define an equivalence relation $\sim_{\ag}$ on $I^{n-1} \times [f(a_n),f(a_0)]$ for every $\Delta^{\ag} \in \Delta_A$. The map $f_{\sigma}$ is compatible with this equivalence relation. Moreover, we define a functor

   $$
    \fonction{C_A^{\Top}}{\Delta_A}{\Top}{\Delta^{\ag}}{I^{n-1} \times [f(a_n),f(a_0)] / \sim_{\ag}.}
    $$

    The images under $C_A^{\Top}$ of the arrows of $\Delta_A$ are constructed from the images under the functor $F$ of the arrows of $\Delta$. In section \ref{section: non stratified functor C A}, the equivalence relation $\sim_{\ag}$ arises as a way to address (ii), but also to ensure functoriality of $C_A^{\Top}$.

    \item In section \ref{section: lift to stratified spaces}, we lift $C_A^{\Top}$ to $\Top_A$. We denote the resulting functor by $C_A$ and refer to it as the \emph{cubical stratified geometric realization} \footnote{Note that the source category of this functor is the simplex category of $A$. The adjective "cubical" here refers to the fact that the values of $C_A$ are constructed as quotients of cubes.}. Asking the arrow $f_{\sigma}$ to be compatible with the stratifications resolves (i). This condition is equivalent to the condition that $\sigma$ maps unbroken morphisms to unbroken trajectories (proposition \ref{proposition: simplex unbroken iff image of unbroken is unbroken}).
    
\end{itemize}

We will later finish resolving difficulty (ii) as follows.

\begin{itemize}
    
    \item In section \ref{section: comparison between realizations}, we prove a comparison theorem between $C_A$ and the usual stratified geometric realization $|-|_A$ (theorem \ref{theorem: isomorphism between |-|'_A and C_A}). Namely, we introduce a \emph{quotient} of $|-|_A$, denoted $|-|_A'$, and prove that $|-|_A'$ and $C_A$ are isomorphic in restriction to the subcategory $\Delta_A^+ \subset \Delta_A$ generated by face morphisms.
\end{itemize}

\subsection{The functor to topological spaces}\label{section: non stratified functor C A}

In this section, we define $C_A^{\Top}$ as a functor to topological spaces (not stratified). Consider two finite increasing sequences of elements of $A$, $\ag = [a_0 \le \hdots \le a_n]$ and $\bg = [b_0 \le \hdots \le b_m]$. Suppose that there exists an arrow $\Phi : \Delta^{\ag} \rightarrow \Delta^{\bg}$ in $\Delta_A$, and denote by $\phi : \Delta^n \rightarrow \Delta^m$ the underlying arrow in $\Delta$. 

Note first that the existence of $\Phi$ implies that $b_0 \le a_0$ and $a_n \le b_m$, in particular, the interval $[f(a_n),f(a_0)]$ is included in the interval $[f(b_m),f(b_0)]$. Secondly, we have a map

$$
F(\varphi)(0,n) : F(\Delta^n)(0,n) \rightarrow F(\Delta^m)(\varphi(0),\varphi(n)).
$$

In the case when $\varphi$ is not surjective, choosing two elements $\delta \in F(\Delta^m)(0,\varphi(0))$ and $\theta \in F(\Delta^m)(\varphi(n),m)$ defines a map

$$
\fonctionsansnom{F(\Delta^n)(\varphi(0),\varphi(n))}{F(\Delta^m)(0,m)}{\gamma}{\theta \circ F(\varphi)(0,n)(\gamma) \circ \delta.}
$$

However, it is not possible to make choices of $\theta$ and $\delta$ for all $\varphi$ in such a way that $\Delta^n \mapsto F(\Delta^n)(0,n)$ becomes a \textit{functor} out of $\Delta$. Indeed, consider for example the following commutative diagram in $\Delta$

$$
\xymatrix{
& \Delta^2 \ar[rd]^{\delta_3^2} \\
\Delta^1 \ar[ru]^{\delta_2^2} \ar[rd]_{\delta_2^2} && \Delta^3. \\
& \Delta^2 \ar[ru]_{\delta_3^3}
}
$$

Since there is just one morphism in $F(\Delta^2)(1,2)$ and one morphism in $F(\Delta^3)(2,3)$, there is only one possible choice in the procedure above. Then, using the homeomorphism \eqref{homeomorphisms} and propositions \ref{proposition: compositions in F(n)} and \ref{proposition: image of the faces by F}, we see that the map $I^0 \rightarrow I^2$ obtained by passing through the top of the diagram is $0 \mapsto (1,0)$ and by passing through the bottom is $0 \mapsto (1,1)$.

However, we will now see that the arrow

$$f_{\sigma} : F(\Delta^n)(0,n) \times [f(a_n),f(a_0)] \rightarrow X$$

obtained from any $\sigma \in \Nl(\Ml)_n$ factors through the quotient of $F(\Delta^n)(0,n) \times [f(\sigma_n),f(\sigma_0)]$ by the following relation.

\begin{Definition}\label{definition: relation on cubes}
    For every object $\Delta^{\ag} \in \Delta_A$, denoting the length of $\ag$ by $n$, we define a relation $\sim_{\ag}$ on $F(\Delta^n) \times [f(a_n),f(a_0)]$ as

    $$
    (\nu \circ \gamma \circ \eta,s) \sim_{\ag} (\nu' \circ \gamma \circ \eta',s)
    $$

    for every $0 \le i \le j \le n$ and $s$ such that $f(a_j) \le s \le f(a_i)$ and $\eta, \eta' \in F(\Delta^n)(0,i)$, $\gamma \in F(\Delta^n)(i,j)$, $\nu, \nu' \in F(\Delta^n)(j,n)$.
\end{Definition}

\begin{Lemmas}\label{lemma: f sigma is compatible with equivalence relation on cube}
    Let $n$ be a nonnegative integer and $\sigma$ an $n$-simplex of $\Nl(\Ml)$. Let $\ag$ be the sequence defined as $a_i = \sigma(i)$. The map $f_{\sigma}$ from notation \ref{notation: map f sigma} is compatible with $\sim_{\ag}$.
\end{Lemmas}

\begin{proof}
Let $(\nu \circ \gamma \circ \eta,s)$ and $(\nu' \circ \gamma \circ \eta',s)$ be as in definition \ref{definition: relation on cubes}. We wish to prove that they have the same image by $f_{\sigma}$. By definition we have

$$
\begin{aligned}
f_{\sigma}(\nu \circ \gamma \circ \eta,s) & = \sigma(\nu \circ \gamma \circ \eta)(s) \\
& = (\sigma(\nu) \circ \sigma(\gamma) \circ \sigma(\eta))(s).
\end{aligned}
$$

The trajectories $\sigma(\eta)$, $\sigma(\gamma)$ and $\sigma(\nu)$ connect $a_0$ to $a_i$, $a_i$ to $a_j$ and $a_j$ to $a_n$ respectively. Since $f(a_j) \le s \le f(a_i)$ we have

$$
(\sigma(\nu) \circ \sigma(\gamma) \circ \sigma(\eta))(s) = \sigma(\gamma)(s),
$$

and similarly

$$
f_{\sigma}(\nu' \circ \gamma \circ \eta',s) = \sigma(\gamma)(s).
$$

Hence $f_{\sigma}(\nu \circ \gamma \circ \eta,s) = f_{\sigma}(\nu' \circ \gamma \circ \eta',s)$, as desired.
\end{proof}

\begin{Lemmas}\label{lemma: relation on cube is an equivalence relation}
    The relation $\sim_{\ag}$ is an equivalence relation for every finite increasing sequence $\ag$ of elements of $A$.
\end{Lemmas}

\begin{proof}
    This follows from corollary \ref{corollary: writing as composition of unbroken morphisms}.
\end{proof}

We can now introduce the main definition of this section.

\begin{Definition}\label{definition: definition of the functor C_A to Top}

The functor $C_A^{\Top} : \Delta_A \rightarrow \Top$ is defined on objects as

$$
C_A^{\Top}(\Delta^{\ag}) = F(\Delta^n)(0,n) \times [f(a_n),f(a_0)] / \sim_{\ag},
$$

and given a morphism $\Phi : \Delta^{\ag} \rightarrow \Delta^{\bg}$ in $\Delta_A$, denoting its image in $\Delta$ by $\varphi : \Delta^n \rightarrow \Delta^m$:

$$
\fonction{C_A^{\Top}(\Phi)}{C_A^{\Top}(\Delta^{\ag})}{C_A^{\Top}(\Delta^{\bg})}{(\gamma,s)}{(\theta \circ F(\varphi)(0,n)(\gamma) \circ \delta,s)}
$$

for any choice of $\theta \in F(\Delta^m)(\phi(n),m)$ and $\delta \in F(\Delta^n)(0,\phi(0))$.

\end{Definition}

\begin{Proposition}\label{proposition: C A is a functor}

$C_A^{\Top}$ is a functor.

\end{Proposition}

\begin{proof}

Observe that $C_A^{\Top}(\Delta^{\ag})$ is the coequalizer of the two maps \footnote{I owe this formulation to Lukas Waas.}

\[\label{coequalizer} \tag{C}
\bigsqcup_{0 \le i' \le i \le j \le j' \le n} F(\Delta^n)(i',i) \times F(\Delta^n)(i,j)\times [f(a_j),f(a_i)] \times F(\Delta^n)(j,j') \rightrightarrows \bigsqcup_{0 \le i \le j \le n} F(\Delta^n)(i,j) \times [f(a_j),f(a_i)],
\]

where the first arrow is the projection to the middle factor and the second arrow is induced by the inclusion $\iota : [f(a_j),f(a_i)] \subseteq [f(a_{j'}),f(a_{i'})]$ and by the composition map

$$
F(\Delta^n)(i',i) \times F(\Delta^n)(i,j) \times F(\Delta^n)(j,j') \rightarrow F(\Delta^n)(i',j').
$$

For every $0 \le i' \le i \le j \le j' \le n$, $\varphi$ induces maps

$$
\begin{aligned}
F(\varphi)(i',i) & \times F(\varphi)(i,j) \times \mathrm{id} \times F(\varphi)(j,j') : F(\Delta^n)(i',i) \times F(\Delta^n)(i,j)\times [f(a_j),f(a_i)] \times F(\Delta^n)(j,j') \\
& \longrightarrow F(\Delta^m)(\varphi(i'),\varphi(i)) \times F(\Delta^m)(\varphi(i),\varphi(j))\times [f(b_{\varphi(j)}),f(b_{\varphi(i)})] \times F(\Delta^m)(\varphi(j),\varphi(j'))
\end{aligned}
$$

and

$$
F(\varphi)(i,j) \times \mathrm{id} : F(\Delta^n)(i,j)\times [f(a_j),f(a_i)] \longrightarrow F(\Delta^m)(\varphi(i),\varphi(j))\times [f(b_{\varphi(j)}),f(b_{\varphi(i)})]
$$

which make diagram (\ref{coequalizer}), and therefore the corresponding coequalizer, functorial in $\ag$. The arrows between coequalizers obtained in this way coincide with the ones of definition \ref{definition: definition of the functor C_A to Top}, and this proves that $C_A^{\Top}$ is a functor. \qedhere

\end{proof}

\begin{Remark}\label{remark: functor defined by cohen jones segal}

    The restriction of the functor $C_A^{\Top}$ to the full subcategory of $\Delta_A$ whose objects are \emph{strictly increasing} sequences was already considered by Cohen, Jones and Segal in \cite[Section 3]{CohenJonesSegalMorse}.

\end{Remark}

\begin{Remark}\label{remark: C_A does not factor through simplex category}
    For every length $1$ sequence $\ag = [a_0 \le a_1] \in \Delta_A$, the image of $\Delta^{\ag}$ by $C_A^{\Top}$ is identified with the interval $[f(a_1),f(a_0)]$. Now suppose that there exist $a,b \in A$ such that $a < b$, or equivalently, that $\dim X \ge 1$. Then $C_A^{\Top}(\Delta^{[a=a]})$ and $C_A^{\Top}(\Delta^{[a < b]})$ are not homeomorphic. In particular, the functor $C_A^{\Top}$ does not factor through the forgetful functor $\Delta_A \rightarrow \Delta$, contrary to the functor $|-|_A$.
\end{Remark}

\begin{Notation}\label{notation: left Kan extension of unstratified cubes}
    We denote by $|-|_{C_A^{\Top}}$ the left Kan extension of the functor $C_A^{\Top} : \Delta_A \rightarrow \Top$ to $\sSet_A$.
\end{Notation}

We finally show that the functorial structure of $C_A^{\Top}$ is, in some sense, compatible with the simplicial structure of $\Nl(\Ml)$. To make this precise, we introduce the following notation.

\begin{Notation}\label{notation: unstratified singular simplicial set obtained from C_A}
    We denote by $\Sing_{C_A^{\Top}} : \Top \rightarrow \sSet_A$ the right adjoint of the left Kan extension of $C_A^{\Top} : \Delta_A \rightarrow \Top$ to $\sSet_A$ (construction \ref{construction: right adjoint of left Kan extension}).
\end{Notation}

Recall that for every topological space $Y$, $\Sing_{C_A^{\Top}}(Y)$ is an $A$-stratified simplicial set whose underlying simplicial set is given by

    $$\Sing_{C_A^{\Top}}(Y)_n = \{c : C_A^{\Top}(\Delta^{\ag}) \rightarrow X \text{ continuous, for some } \ag \text{ of length } n\},$$

with the simplicial structure induced by the functorial structure of $C_A^{\Top}$.

In particular, for every nonnegative integer $n$ there is a map of sets

$$
\fonction{\Sl_n}{\Nl(\Ml)_n}{\Sing_{C_A^{\Top}}(X)_n}{\sigma}{f_{\sigma}.}
$$

\begin{Proposition}\label{Proposition: compatibility of simplicial structures}
    These maps determine a morphism of simplicial sets

    $$
    \Sl : \Nl(\Ml) \rightarrow \Sing_{C_A^{\Top}}(X).
    $$
\end{Proposition}

\begin{proof}
    Let $\varphi : \Delta^n \rightarrow \Delta^m$ be a morphism in $\Delta$. We want to justify that the following diagram commutes

    $$\xymatrix{
    \Nl(\Ml)_m \ar[r]^-{\Sl_m} \ar[d]_{\Nl(\Ml)(\varphi)} & \Sing_{C_A^{\Top}}(X)_m \ar[d]^{\Sing_{C_A^{\Top}}(\varphi)} \\
    \Nl(\Ml)_n \ar[r]_-{\Sl_n} & \Sing_{C_A^{\Top}}(X)_n.
    }$$

    Let $\sigma \in \Nl(\Ml)_m$. Define $\Delta^{\ag} \in \Delta_A$ as $a_i = \sigma(i)$. Define another sequence $\ag \circ \varphi$ of length $n$ by the condition that the following diagram of posets commutes

    $$
    \xymatrix{
    [n] \ar[rr]^{\varphi} \ar[rd]_{\ag \circ \varphi} && [m] \ar[ld]^{\ag} \\
    & A.
    }
    $$
    
    Denote by $\Phi$ the morphism $\Delta^{\ag \circ \varphi} \rightarrow \Delta^{\ag}$ determined by this diagram. We have:

    $$
    \begin{aligned}
        \Sl_n(\Nl(\Ml)(\varphi)(\sigma)) & = \Sl_n(\sigma \circ F(\varphi)) \\
        & = f_{\sigma \circ F(\varphi)}
    \end{aligned}
    $$

    and

    $$
    \begin{aligned}
        \Sing_{C_A^{\Top}}(\varphi)(\Sl_m(\sigma)) & = \Sing_{C_A^{\Top}}(\varphi)(f_{\sigma}) \\
        & = f_{\sigma} \circ C_A^{\Top}(\Phi).
    \end{aligned}
    $$

    For every $(\gamma,s) \in C_A^{\Top}(\ag \circ \varphi)$ we have

    $$
    f_{\sigma \circ F(\varphi)}(\gamma,s) = \sigma(F(\varphi)(\gamma))(t).
    $$

    But since $f(a_{\varphi(n)}) \le s \le f(a_{\varphi(0)})$ we also have for every $\alpha \in F(\varphi(n),m)$ and $\beta \in F(0,\varphi(0))$
    
    $$
    \begin{aligned}
    f_{\sigma \circ F(\varphi)}(\gamma,s) & = \sigma(\alpha \circ F(\varphi)(\gamma) \circ \beta)(t) \\
    & = f_{\sigma}\circ C_A^{\Top}(\Phi)(\gamma,s).
    \end{aligned}
    $$

    This concludes the proof.
\end{proof}

\subsection{The lift to stratified topological spaces}\label{section: lift to stratified spaces}

In this section, we introduce a lift of the functor $C_A^{\Top}$, defined in section \ref{section: non stratified functor C A}, to stratified topological spaces. This takes the form of the following theorem.

\begin{Thms}\label{theorem: stratifications on cubes}
    There exists a unique family of $A$-stratifications
    
    $$(\pi_{\ag} : F(\Delta^{l(\ag)})(0,l(\ag)) \times [f(a_{l(\ag)}),f(a_0)] \rightarrow A)_{\Delta^{\ag} \in \Delta_A}$$
    
    such that for every $\Delta^{\ag} \in \Delta_A$, denoting the length of $\ag$ by $n$, the following conditions hold:

    \begin{enumerate}[label=(\roman*)]
        \item For $n=0$, the only stratum corresponds to $a_0$.

        \item If $s \neq f(a_0)$ and $\gamma$ is unbroken (definition \ref{definition: broken and unbroken morphism of F(n)}) then $(\gamma,s) \in F(\Delta^n)(0,n) \times [f(a_n),f(a_0)]$ is in the $a_n$-stratum.

        \item The stratification map $\pi_{\ag}$ is compatible with the equivalence relation $\sim_{\ag}$ (definition \ref{definition: relation on cubes}) and the induced stratifications on the values of $C_A^{\Top}$ yield a lift of $C_A^{\Top}$ to $\Top_A$.
    \end{enumerate}

    Moreover, $\pi_{\ag}$ is explicitly described as follows: given $(\gamma,s) \in F(\Delta^n)(0,n) \times [f(a_n),f(a_0)]$, by proposition \ref{corollary: writing as composition of unbroken morphisms} we can write $\gamma$ as a composition of unbroken morphisms $\gamma_k \circ \hdots \circ \gamma_0$, with $\gamma_j \in F(\Delta^n)(i_j,i_{j+1})$. Then:

    \begin{itemize}
        \item if $s=f(a_0)$, then $\pi_{\ag}(\gamma,s) = a_0$;

        \item otherwise, there exists a unique $l$ such that $f(a_{i_{l+1}}) \le s < f(a_{i_l})$, and then $\pi_{\ag}(\gamma,s) = a_{i_{l+1}}$.
    \end{itemize}
\end{Thms}

\begin{Example}\label{example: stratified cubes}

Examples of stratifications $\pi_{\ag}$ are given by figure \ref{figure: stratifications of the 1-cube} for $\ag$ of length 1, by figure \ref{figure: stratifications of the 2-cube} for $\ag$ of length 2, and by figure \ref{figure: some stratifications of 3-cubes} for $\ag$ of length 3. Examples of the corresponding stratifications on $C_A^{\Top}(\Delta^{\ag})$ are given by figure \ref{figure: quotients of stratified 2-cubes} for $\ag$ of length 2, and by figures \ref{figure: quotient of cube for a0<a1<a2<a3}, \ref{figure: quotient of cube for a0<a1=a2<a3} and \ref{figure: quotient of cube for a0<a1=a2=a3} for $\ag$ of length 3.

\end{Example}

\begin{proof}

The proof is done in two steps:

\textit{\underline{Step 1:}} We first prove that if there exist $A$-stratifications $(\pi_{\ag} : F(\Delta^n)(0,n) \times [f(a_n),f(a_0)] \rightarrow A)_{\Delta^{\ag} \in \Delta_A}$ satisfying conditions (i), (ii) and (iii), then they are given by the expression in the statement of the theorem.

\textit{\underline{Step 2:}} We then prove that these expressions actually define stratifications that satisfy these conditions.

\textit{\underline{Proof of step 1:}} Assume that there exists a family of $A$-stratifications that satisfy the conditions. Let $\Delta^{\ag} \in \Delta_A$, and denote the length of $\ag$ by $n$. Given $\gamma, \delta \in F(\Delta^n)(0,n)$, one can write $\gamma = \gamma \circ \id_{0}$ and $\delta = \delta \circ \id_{0}$, so $(\gamma,f(a_0)) \sim_{\ag} (\delta,f(a_0))$. Hence the image of $F(\Delta^n)(0,n) \times \{f(a_0)\}$ in $C_A^{\Top}(\Delta^{\ag})$ is a point. This point is the image of $C_A^{\Top}(\Delta^{[a_0]})$ under the image by $C_A^{\Top}$ of the arrow $\Delta^{[a_0]} \rightarrow \Delta^{\ag}$ determined by the inclusion of the first element of the sequence $\ag$. By conditions (i) and (iii) we must have $\pi_{\ag}(F(\Delta^n)(0,n) \times \{f(a_0)\}) = \{a_0\}$.

Otherwise, let $(\gamma,s) \in F(\Delta^n)(0,n) \times [f(a_n),f(a_0)]$ such that $s \neq f(a_0)$. Write $\gamma$ as a composition of unbroken morphisms and consider the integer $l$ such that $f(a_{i_{l+1}}) \le s < f(a_{i_l})$, as in the statement of the theorem. Let us define a morphism of posets

$$
\fonction{\theta}{[i_{j+1}-i_j]}{[n]}{k}{k+i_j.}
$$

Let us also define a sequence $\ag \circ \theta$ of length $i_{j+1}-i_j$, by the condition that the following diagram of posets commutes

$$
\xymatrix{
[i_{j+1}-i_j] \ar[rd]_-{\ag \circ \theta} \ar[rr]^-{\theta} && [n] \ar[ld]^-{\ag} \\
& A.
}
$$

Denote by $\Theta : \Delta^{\ag \circ \theta} \rightarrow \Delta^{\ag}$ the arrow in $\Delta_A$ determined by this diagram.

Let us abuse notation and denote by $\gamma_l$ the element of $F(\Delta^{i_{j+1}-i_j})(0,i_{j+1}-i_j)$ which is mapped to $\gamma_l$ by $F(\theta)(0,i_{l+1}-i_l)$. Then, using the same notation for elements of the quotient by $\sim_{\ag}$ and their representatives, we have, by definition of $C_A^{\Top}$:

$$
C_A^{\Top}(\Theta)(\gamma_l,s) = (\gamma,s).
$$

By conditions (ii) and (iii), we must have $\pi_{\ag}(\gamma,s) = a_{i_{l+1}}$.

\textit{\underline{Proof of step 2:}} Let us define the $\pi_{\ag}$'s as in the statement of the theorem. We first prove that these define stratifications, that is, that these are continuous. Given $0 \le i \le n$, recall that we denote by $A_{\ge a_i} \subset A$ the open subset of those elements greater or equal to $a_i$. We want to show that $\pi_{\ag}^{-1}(A_{\ge a_i})$ is open in $I^{n-1} \times [f(a_n),f(a_0)]$. If $i=0$ this holds since $\pi_{\ag}^{-1}(U_{a_0}) = I^{n-1} \times [f(a_n),f(a_0)]$. Otherwise, let $x = (t_1, \hdots, t_{n-1},s)$ such that $\pi_{\ag}(x) \ge a_i$. Let $0 < i_1 < \hdots < i_k < n$ be the $i$'s such that $t_i=1$ (if there are none, we let $k=0$). There exists a unique $l$ such that $f(a_{i_{l+1}}) \le s < f(a_{i_l})$ and we have $\pi_{\ag}(x) = a_{i_{l+1}}$. Consider the open subset $U = \{s_{i_l+1} \neq 1, s_{i_l+2} \neq 1, \hdots, s_{i_{l+1}-1} \neq 1 \}$. Then $U \times [f(a_n),f(a_{i_l})[$ is a neighborhood of $x$ in $\pi_{\ag}^{-1}(A_{\ge a_{i_{l+1}}})$, and therefore also in $\pi_{\ag}^{-1}(A_{\ge a_i})$. This proves that $\pi_{\ag}^{-1}(A_{\ge a_i})$ is open, as desired.

We now show that the $\pi_{\ag}$'s satisfy the three conditions. For conditions (i) and (ii), this is immediate. To prove that $\pi_{\ag}$ is compatible with $\sim_{\ag}$, consider two equivalent elements of $F(\Delta^n)(0,n)$,

$$
(\nu \circ \delta \circ \eta,s) \sim_{\ag} (\nu' \circ \delta \circ \eta',s),
$$

where $0 \le i \le j \le n$, $f(a_j) \le s \le f(a_i)$ and $\eta, \eta' \in F(\Delta^n)(0,i)$, $\delta \in F(\Delta^n)(i,j)$, $\nu, \nu' \in F(\Delta^n)(j,n)$. The expression of $\nu \circ \delta \circ \eta$ (resp. $\nu' \circ \delta \circ \eta'$) as a composition of unbroken morphisms is the concatenation of those of $\nu$, $\delta$ and $\eta$ (resp. $\nu'$, $\delta$ and $\eta'$). Hence, letting $\theta$ be the morphism of posets

$$
\fonction{\theta}{[j-i]}{[n]}{k}{k+i}
$$

we have

$$
\pi_{\ag}(\nu \circ \delta \circ \eta,s) = \pi_{\ag}(\nu' \circ \delta \circ \eta',s) = \pi_{\ag \circ \theta}(\delta,s).
$$

Therefore $\pi_{\ag}$ is compatible with $\sim_{\ag}$. To prove that this defines a lift of $C_A^{\Top}$ to $\Top_A$, consider an arrow $\Phi$ in $\Delta_A$, determined by a commutative diagram of posets

$$
\xymatrix{
[n] \ar[rr]^{\varphi} \ar[rd]_-{\ag} && [m] \ar[ld]^-{\bg} \\
& A.
}
$$

Let $(\gamma,s) \in F(\Delta^n)(0,n)$ and $\gamma = \gamma_k \circ \hdots \circ \gamma_0$ be the expression of $\gamma$ as a composition of unbroken morphisms. Choosing $\alpha \in F(\Delta^m)(0,\varphi(0))$, $\beta \in F(\Delta^m)(\varphi(n),m)$ and using the same notation for the element $(\gamma,s)$ and its image in the quotient by $\sim_{\ag}$, we have by definition

$$
C_A^{\Top}(\Phi)(\gamma,s)=(\beta \circ F(\varphi)(0,n)(\gamma) \circ \alpha,s).
$$

But, by corollary \ref{corollary: image of an unbroken morphism}, $F(\varphi)(0,n)(\gamma_k) \circ \hdots \circ F(\varphi)(0,n)(\gamma_0)$ is the expression of $F(\varphi)(0,n)(\gamma)$ as a composition of unbroken morphisms. Hence $\pi_{\bg}(\beta \circ F(\varphi)(0,n)(\gamma) \circ \alpha,s) = \pi_{\ag}(\gamma,s)$, in other words, $C_A^{\Top}(\Phi)$ is compatible with stratifications, as desired. \qedhere

\end{proof}

\begin{Remark}

Given an arbitrary poset $A$, the construction of the functor $C_A : \Delta_A \rightarrow \Top_A$ can be reproduced verbatim as soon as $f$ is endowed with a strictly decreasing map $f : A \rightarrow \R$. We will come back to this in section \ref{section: Comparison between stratified categories and exit path categories}.

\end{Remark}

\begin{Definition}\label{Definition: Sing_C_A}
    We denote by $C_A$ the lift of $C_A^{\Top}$ to $\Top_A$ provided by theorem \ref{theorem: stratifications on cubes}. We denote by $|-|_{C_A}$ the left Kan extension of $C_A$ to $\sSet_A$ (construction \ref{construction: left Kan extension}), and we refer to this functor as the \emph{cubical stratified geometric realization} functor.
\end{Definition}

\begin{Remark}\label{remark: Kan extensions of stratified and not stratified cubes}
    Since the forgetful functor $\Top_A \rightarrow \Top$ commutes with colimits, the corestriction of this left Kan extension to $\Top$ coincides with the functor $|-|_{C_A^{\Top}} : \sSet_A \rightarrow \Top$ from the previous section.
\end{Remark}

\begin{Remark}\label{remark: explicit description of cubical realization}
    Recall that the functor $|-|_{C_A}$ is determined (up to canonical isomorphism) by the conditions that its restriction to $\Delta_A$ is $C_A$, and it commutes with colimits. Let us describe more explicitly this left Kan extension, following remark \ref{Remark: concrete description of left Kan extension}. Let $S$ be an $A$-stratified simplicial set, regarded as a functor $\Delta_A^{\op} \rightarrow \Set$. For every $\Delta^{\ag} \in \Delta_A$, denote by $S_{\ag}$ the image of $\Delta^{\ag}$ by $S$. The stratified geometric realization of $S$ is obtained from the disjoint union of stratified topological spaces
    
    $$
    \bigsqcup_{\Delta^{\ag} \in \Delta_A} S_{\ag} \times C_A(\Delta^{\ag})
    $$
    
    by modding out by the equivalence relation generated by the relations: $(S(\Phi)(\tau),x) \sim (\tau, C_A(\Phi)(x))$ for every arrow $\Phi : \Delta^{\ag} \rightarrow \Delta^{\bg} \in \Delta_A$, every $\tau \in S_{\bg}$ and every $x \in C_A(\Delta^{\ag})$.
\end{Remark}

In construction \ref{construction: filtration on usual realization} and definition \ref{definition: filtration on usual realization}, we introduced a natural transformation $f_A$ between the corestriction of $|-|_A$ to $\Top$, and the constant functor at $\R$ on $\sSet_A$. We now introduce an analogous natural transformation, between the corestriction of $|-|_{C_A}$ to $\Top$, and the constant functor at $\R$ on $\sSet_A$.

\begin{Construction}\label{construction: filtration on C A}
    For every finite increasing sequence $\ag = [a_0 \le \hdots \le a_n]$ of elements of $A$, the projection map

    $$
    I^{n-1} \times [f(a_n),f(a_0)] \rightarrow [f(a_n),f(a_0)]
    $$

    is compatible with the equivalence relation $\sim_{\ag}$, thus defining a map

    $$C_A(\Delta^{\ag}) \rightarrow \R$$

    by passing to the quotient by $\sim_{\ag}$ and composing with the inclusion $[f(a_n),f(a_0)] \rightarrow \R$.

    Moreover, for every arrow $\Phi : \Delta^{\ag} \rightarrow \Delta^{\bg}$ in $\Delta_A$, the following diagram commutes

    $$
    \xymatrix{
    C_A(\Delta^{\ag}) \ar[rr]^{C_A(\Phi)} \ar[rd] && C_A(\Delta^{\bg}) \ar[ld] \\
    & \R.
    }
    $$
    
   We have thus constructed a natural transformation $C_A \rightarrow \underline{\R}$ of functors $\sSet_A \rightarrow \Top$. Furthermore, the latter extends to a natural transformation between the left Kan extensions of these functors to $\sSet_A$.
\end{Construction}

\begin{Definition}\label{definition: filtration on C A}
    We refer to the natural transformation $|-|_{C_A} \rightarrow \underline{\R}$ of functors from $\sSet_A$ to $\Top$, constructed in \ref{construction: filtration on C A}, as the \emph{filtration} on $|-|_{C_A}$ induced by $f$, and we denote it by $f_{C_A}$.
\end{Definition}

\begin{Remark}\label{remark: explicit description of filtration on cubical realization}
    Given an $A$-stratified simplicial set $S$, an explicit description of the map $f_{C_A}(S) : |S|_{C_A} \rightarrow \R$ can be obtained from remark \ref{remark: explicit description of filtration}, replacing the usual geometric realization by $C_A$.
\end{Remark}

\subsection{\texorpdfstring{The flow coherent nerve of $\Ml$}{The flow coherent nerve of M}}\label{section: The flow coherent nerve of M}

Recall that in section \ref{section: non stratified functor C A} we introduced a functor $C_A^{\Top} : \Delta_A \rightarrow \Top$ as well as a morphism of simplicial sets $\Nl(\Ml) \rightarrow \Sing_{C_A^{\Top}}(X)$. The latter reflects the fact that to every $n$-simplex of $\Nl(\Ml)$ is associated a singular $n$-cube in $X$. In section \ref{section: lift to stratified spaces}, we introduced a lift $C_A$ of $C_A^{\Top}$ to $\Top_A$. In this section, we introduce the simplicial subset of $\Nl(\Ml)$ formed by those simplices for which the corresponding cube in $X$ is compatible with stratifications (definition \ref{definition: flow coherent nerve}). We also give an alternative description of it, which involves stratifications of the morphism spaces of $\Ml$ (proposition \ref{proposition: unbroken iff compatible with stratifications on morph spaces}).

We start with the following definition.

\begin{Definition}\label{definition: Sing C A}
    We denote by $\Sing_{C_A} : \Top_A \rightarrow \sSet_A$ the right adjoint of the cubical stratified geometric realization. Given an $A$-stratified topological space $Y$, we refer to $\Sing_{C_A}(Y)$ as the \emph{simplicial set of stratified cubes of} $Y$.
    
\end{Definition}

Given a $A$-stratified topological space $Y$, we will denote the simplicial set associated with the $A$-stratified simplicial set $\Sing_{C_A}(Y)$ using the same notation, i.e. we will denote

$$
\Sing_{C_A}(X)_n = \bigsqcup_{\substack{\ag \in \Delta_A\\ l(\ag)=n}} \Sing_{C_A}(X)_{\ag}.
$$

Recall from construction \ref{construction: right adjoint of left Kan extension} that, for every $A$-stratified topological space $Y$, $\Sing_{C_A}(Y)$ is the $A$-stratified simplicial set given by

    $$\Sing_{C_A}(Y)_n = \{ c : C_A(\Delta^{\ag}) \rightarrow Y  \text{ stratified for some } \ag \in N(A)_n \},$$

with structural morphism $\Sing_{C_A}(X) \rightarrow N(A)$ given by $(f : C_{A}(\Delta^{\ag}) \rightarrow X) \mapsto \ag$.

\begin{Remark}\label{remark: stratified cubical simplicial set simplicial subset of unstratified one}
For every $A$-stratified topological space $Y$, the simplicial set underlying $\Sing_{C_A}(Y)$ is a simplicial subset of the simplicial set $\Sing_{C_A^{\Top}}(Y)$ defined in section \ref{section: non stratified functor C A}.
\end{Remark}

\begin{Definition}\label{definition: flow coherent nerve}

The \emph{flow coherent nerve} of $\Ml$, denoted $S_{\Ml}$, is the pullback

$$
\xymatrix{
S_{\Ml} \ar[d] \ar[r] & \Nl(\Ml) \ar[d]^-{\Sl} \\
\Sing_{C_A}(X) \ar[r] & \Sing_{C_A^{\Top}}(X).
}
$$
    
\end{Definition}

\begin{Notation}\label{notation: functors from flow coherent nerve}
    We denote by
    
    $$
    i_{\Ml} : S_{\Ml} \rightarrow \Nl(\Ml) \qquad \text{and} \qquad s_{\Ml}^{(X,A)} : S_{\Ml} \rightarrow \Sing_{C_A}(X)
    $$
    
    the two morphisms appearing in the definition of $S_{\Ml}$.
\end{Notation}

The following terminology will be motivated by proposition \ref{proposition: simplex unbroken iff image of unbroken is unbroken} below.

\begin{Definition}\label{definition: unbroken simplices and diagrams}
    A simplex of $S_{\Ml}$ is called an \emph{unbroken simplex} of $\Nl(\Ml)$. More generally, a diagram $K \rightarrow \Nl(\Ml)$, where $K \in \sSet$, is said to be \emph{unbroken} if it factors through $S_{\Ml}$.
\end{Definition}

\begin{Remark}\label{remark: reformulation of definition of flow coherent nerve}
    The two previous definitions admits the following reformulation: an unbroken simplex of $\Ml$ is a functor $\sigma : F(\Delta^n) \rightarrow \Ml$ such that the associated map $f_{\sigma} : C_A(\Delta^{\ag}) \rightarrow X$ (where the sequence $\ag = [a_0 \le \hdots \le a_n]$ is defined as $a_i = \sigma(i)$) preserves stratifications.
\end{Remark}

The condition that a simplex of $\Nl(\Ml)$ is unbroken can be alternatively formulated as follows:

\begin{Proposition}\label{proposition: simplex unbroken iff image of unbroken is unbroken}
    An $n$-simplex $\sigma : F(\Delta^n) \rightarrow \Ml$ is unbroken if and only if the image of every unbroken morphism of $F(\Delta^n)$ (in the sense of definition \ref{definition: broken and unbroken morphism of F(n)}) by $\sigma$ is an unbroken trajectory.
\end{Proposition}

\begin{proof}
    Define an increasing sequence $\ag$ by $a_i = \sigma(i)$. The map $f_{\sigma}$ is stratified if and only if for every $(\gamma,s) \in F(\Delta^n)(0,n) \times [f(a_n),f(a_0)]$, writing $\gamma = \gamma_k \circ \hdots \circ \gamma_0$ as a composition of unbroken morphisms with $\gamma_j \in F(\Delta^n)(i_j,i_{j+1})$, and letting $l$ such that $f(a_{l+1}) \le s < f(a_l)$, we have that $\sigma(\gamma)(s)$ is in the stratum of $a_{i_{l+1}}$. Since $\sigma(\gamma)(s) = \sigma(\gamma_l)(s)$, the latter is equivalent to $\sigma(\gamma_l)$ being unbroken, which proves the statement.
\end{proof}

\begin{Example}\label{example: broken morphism iff corresponding morphism is unbroken}
    By proposition \ref{proposition: simplex unbroken iff image of unbroken is unbroken}, a trajectory of $\xi$ is unbroken if and only if the corresponding morphism $\Delta^1 \rightarrow \Nl(\Ml)$ is unbroken in the sense of definition \ref{definition: unbroken simplices and diagrams}.
\end{Example}

We also have the following useful characterization of unbroken simplices.

\begin{Corollary}\label{corollary: condition of being unbroken in terms of faces}
    An $n$-simplex of $\Nl(\Ml)$ is unbroken if and only if its restrictions to all faces of $\Delta^n$ are unbroken, and it maps the interior of $F(\Delta^n)(0,n)$ into the interior of $\Ml(a_0,a_n)$.
\end{Corollary}

\begin{proof}
   Let $\gamma$ be an unbroken morphism of $F(\Delta^n)$. By proposition \ref{proposition: image of the faces by F}, either $\gamma$ is in the image of $F(\delta)$ for some face $\delta$ with target $\Delta^n$, or $\gamma$ is in the interior of $F(\Delta^n)(0,n)$, and is in particular unbroken. Moreover, the interior of $\Ml(a_0,a_n)$ consists precisely of the unbroken trajectories that connect $a_0$ to $a_n$. We conclude by applying proposition \ref{proposition: simplex unbroken iff image of unbroken is unbroken}.
\end{proof}

We finish this section with an alternative characterization of the condition of being unbroken for a simplex $\sigma \in \Nl(\Ml)_n$. We will see that it can be formulated as a condition of compatibility with stratifications for the map $\sigma(0,n) : F(\Delta^n)(0,n) \rightarrow \Ml(\sigma(0),\sigma(n))$. To do this, we begin by introducing stratifications on the morphism spaces of $F(\Delta^n)$ and $\Ml$.

\begin{Definition}\label{definition: subdivision of a poset}
    The \textit{subdivision of} $A$ is the full subcategory of $\Delta_A$ whose objects are the finite sequences of elements of $A$ which are \emph{strictly increasing}. It is denoted by $\sd A$.
\end{Definition}

This category satisfies the property that there is precisely one object in each isomorphism class, and at most one arrow between two given objects. Consequently, the relation

$$a \le b \text{ if and only if } \Hom_{\sd A}(a,b) \neq \varnothing$$

is a partial order on the set of objects of $\sd A$. Moreover, the inclusion $\sd A \rightarrow \Delta_A$ admits a left adjoint which we now introduce.

\begin{Construction}\label{construction: left adjoint of inclusion of subdivision}
    Let $\ag = [a_0 \le \hdots \le a_n]$ be a finite increasing sequence of elements of $A$. We denote by $E(\ag)$ the ordered sequence of distinct elements of $\ag$. This defines a functor $E : \Delta_A \rightarrow \sd A$, which is left adjoint to the inclusion $\sd A \rightarrow \Delta_A$.
\end{Construction}

Fix $\Delta^{\ag} \in \Delta_A$, with $\ag$ of length $n$, and define

$$
\fonction{\pi}{\Ml(a_0,a_n)}{\Obj (\sd A)}{\gamma}{\text{the ordered sequence of critical points at which } \gamma \text{ is broken along with } a_0 \text{ and } a_n.}
$$

Moreover, for a given $\gamma \in F(\Delta^n)(0,n)$, let $\gamma_k \circ \hdots \circ \gamma_0$ be the writing of $\gamma$ as a composition of unbroken morphisms, with $\gamma_j \in F(\Delta^n)(i_j,i_{j+1})$, and define

$$
\fonction{\tau}{F(\Delta^n)(0,n)}{\Obj (\sd A)}{\gamma}{E(a_0 \le a_{i_1} \le \hdots \le a_{i_k} \le a_n)}
$$

\begin{Proposition}\label{proposition: stratifications of morphism spaces}
    The maps $\pi$ and $\tau$ define stratifications of the morphism spaces $\Ml(a_0,a_n)$ and $F(\Delta^n)(0,n)$ by $(\sd A)^{\mathrm{op}}$.
\end{Proposition}

\begin{proof}

If $\ag$ is constant there is nothing to prove, so we suppose that is not, i.e. that $a_0 < a_n$.

Note that, for every $\bg \in \sd A$, we have $((\sd A)^{\op})_{\ge \bg} = (\sd A)_{\le \bg}$, and this set can be described as the set of subsequences of $\bg$. We want to prove that for every $\bg$, $\pi^{-1}((\sd A)_{\le \bg})$ and $\tau^{-1}((\sd A)_{\le \bg})$ are open subsets of $\Ml(a_0,a_n)$ and $F(\Delta^n)(0,n)$ respectively. Let $l$ be the length of $\bg$.

Firstly, $\pi^{-1}((\sd A)_{\le \bg}) = \Ml(a_0,a_n) \backslash Z$ where $Z$ is the set of trajectories that are broken at a point which is not part of $\bg$. In other words $Z$ is the union over the sequences of the form $a_0 < a_1 < \hdots < a_k = a_n$ such that $k \ge 2$ and $a_i \notin \bg$ for every $1 \le i \le k-1$, of the images of the composition maps

$$
\Ml(a_0,a_1) \times \hdots \times \Ml(a_{k-1},a_n) \rightarrow \Ml(a_0,a_n).
$$

By proposition \ref{proposition: structure of smooth manifold with corners on space of broken traj}, these images are closed in $\Ml(a_0,a_n)$. Since there is a finite number of such sequences, $Z$ is closed in $\Ml(a_0,a_n)$, hence $\pi^{-1}((\sd A)_{\le \bg})$ is open, as desired.

Secondly, let us identify $F(\Delta^n)(0,n)$ with $I^{n-1}$. Then $\tau^{-1}((\sd A)_{\le \bg}) = I^{n-1} \backslash Z'$ where $Z'$ is the set of those $(t_1,\hdots,t_{n-1})$ such that there exists $1 \le j \le n-1$ such that $a_j$ is not part of $\bg$ and $t_j = 1$. This is a closed subset of $I^{n-1}$, hence $\tau^{-1}((\sd A)_{\le \bg})$ is open, as desired.
\end{proof}

\begin{Proposition}\label{proposition: unbroken iff compatible with stratifications on morph spaces}
    Let $\sigma$ be an $n$-simplex of $\Nl(\Ml)$. Define $\Delta^{\ag} \in \Delta_A$ by $a_i = \sigma(i)$ and endow $F(\Delta^n)(0,n)$ and $\Ml(a_0,a_n)$ with the corresponding stratifications by $(\sd A)^{\mathrm{op}}$. Then $\sigma$ is unbroken if and only if $\sigma(0,n)$ is compatible with the stratifications.
\end{Proposition}

\begin{proof}
    $\sigma(0,n)$ is compatible with the stratifications if and only if, for every $\gamma \in F(\Delta^n)(0,n)$, $\sigma(\gamma)$ is broken at exactly as many points as $\gamma$. Writing $\gamma$ as a composition of unbroken morphisms $\gamma = \gamma_k \circ \hdots \circ \gamma_0$, this is equivalent to the condition that $\sigma(\gamma_j)$ is unbroken for every $0 \le j \le k$. We conclude by applying proposition \ref{proposition: simplex unbroken iff image of unbroken is unbroken}.
\end{proof}

\begin{Remark}\label{remark: unbroken simplices defined for any stratified category}
    One feature of the equivalent definition of the flow coherent nerve of $\Ml$ provided by proposition \ref{proposition: unbroken iff compatible with stratifications on morph spaces}, is that it does not rely on the stratified space $X$, but only on the topological category $\Ml$, together with extra structure on its morphism spaces. In section \ref{section: stratified categories as exit path categories} we will exploit this property to define the flow coherent nerve for a more general class of topological categories which we will introduce.
\end{Remark}

\newpage

\section{Comparison between the cubical and standard stratified geometric realizations}\label{section: comparison between realizations}

It follows from the definition of $S_{\Ml}$ (definition \ref{definition: flow coherent nerve}) that there are two morphisms of simplicial sets 

$$\Nl(\Ml) \leftarrow S_{\Ml} \rightarrow \Sing_{C_A}(X).$$

Recall that our ultimate goal is to compare the $\infty$-categories $\Nl(\Ml)$ and $\Sing_A(X)$. With this in mind, we state and prove in this section a comparison theorem between the usual stratified geometric realization $|-|_A : \Delta_A \rightarrow \Top_A$, and its cubical counterpart $C_A$ (theorem \ref{theorem: isomorphism between |-|'_A and C_A}). Recall that for every nonnegative integer $n$ we have

$$
\Sing_A(X)_n = \{\sigma : |\Delta^{\ag}|_A \rightarrow X \text{ stratified for some } \ag \text{ of length } n\}
$$

and

$$
\Sing_{C_A}(X)_n = \{\sigma : C_A(\Delta^{\ag}) \rightarrow X \text{ stratified for some } \ag \text{ of length } n\},
$$

with the simplicial structures induced by the functorial structures of $|-|_A$ and $C_A$ respectively. In particular, any comparison result between the functors $|-|_A$ and $C_A$ gives rise to a comparison result between the simplicial sets $\Sing_A(X)$ and $\Sing_{C_A}(X)$.

This section consists of three subsections. In section \ref{section: modifying geometric realization}, we observe that a mere isomorphism of functors does not exist, and we remedy this by introducing a modification of the functor $|-|_A$, as well as restricting ourselves to a subcategory of $\Delta_A$. In section \ref{section: proof of comparison between realizations}, we state and prove the comparison result between $|-|_A$ and $C_A$ (theorem \ref{theorem: isomorphism between |-|'_A and C_A}). The proof relies on another result, theorem \ref{theorem: stratified simplices and cubes as balls}, whose proof is the subject of section \ref{section: proof of the ball theorem for realizations}. To summarize, here is a roadmap of logical dependencies leading to theorem \ref{theorem: isomorphism between |-|'_A and C_A}:

$$
\text{Theorem } \ref{theorem: filtration on cube is trivial} \Rightarrow \text{Theorem } \ref{theorem: stratified simplices and cubes as balls} \Rightarrow \text{Theorem } \ref{theorem: isomorphism between |-|'_A and C_A}.
$$

\subsection{Modifying the standard stratified geometric realization}\label{section: modifying geometric realization}

Ideally, we would like the functors $|-|_A$ and $C_A$ to be \textit{isomorphic}. But observe that if $\ag$ is a constant finite sequence of elements of $A$, then $C_A(\Delta^{\ag})$ is a point, but on the other hand, if $\ag$ is of length at least $1$, then $|\Delta^{\ag}|_A$ is not a point and therefore $C_A(\Delta^{\ag})$ and $|\Delta^{\ag}|_A$ are not homeomorphic. The goal of the present section is to introduce a modification of the functor $|-|_A$ which better compares to $C_A$.

\begin{Construction}\label{construction: functor |-|_A'}

Given a finite increasing sequence $\ag$ of elements of $A$ of length $n$, we say that a constant subsequence of $\ag$ of consecutive elements $\ag_{p,q} = [a_p = a_{p+1} = \hdots = a_q]$ is \textit{of maximal length} when, if $p > 0$ then $a_{p-1} \neq a_p$ and if $q < n$ then $a_q \neq a_{q+1}$. The inclusion of such a subsequence in $\ag$ induces a map $|\Delta^{\ag_{p,q}}|_A \rightarrow |\Delta^{\ag}|_A$ in $\Top_A$. Let us declare any two elements in the image of this map to be equivalent. This defines an equivalence relation on $|\Delta^{\ag}|_A$. We denote by $|\Delta^{\ag}|_A'$ the quotient of $|\Delta^{\ag}|_A$ by this equivalence relation.

\end{Construction}

\begin{Proposition}\label{proposition: modified stratified realization is a functor}
    These equivalence relations on the values of $|-|_A$ are compatible with the stratifications by $A$ as well as with the maps induced by the morphisms of $\Delta_A$ by functoriality. We have therefore constructed a functor $|-|_A' : \Delta_A \rightarrow \Top_A$.
\end{Proposition}

\begin{proof}

    Fix a constant subsequence of $\ag$ of consecutive elements $\ag_{p,q} = [a_p = a_{p+1} = \hdots = a_q]$ of maximal length. Firstly, the image of the map $|\Delta^{\ag_{p,q}}|_A \rightarrow |\Delta^{\ag}|_A$ consists of the convex combinations of the form
    
    $$
    \sum_{k=p}^q t_ke_k.
    $$
    
    This image is therefore contained in the $a_p$-stratum, and so the equivalence relation on $|\Delta^{\ag}|_A$ is compatible with the stratification by $A$.
    
    Secondly, consider a morphism $\Phi$ in $\Delta_A$ determined by a commutative diagram of posets of the form
    
    $$
    \xymatrix{
    [n] \ar[rd]_-{\ag} \ar[rr]^-{\varphi} && [m] \ar[ld]^-{\bg} \\
    & A.
    }
    $$
    
    According to the condition that this diagram is commutative, we have $b_{\varphi(p)} = a_p$ and $b_{\varphi(q)} = a_q$, and therefore $b_{\varphi(p)} = b_{\varphi(q)}$. The subsequence $\bg_{\varphi(p),\varphi(q)} = [b_{\varphi(p)} \le b_{\varphi(p)+1} \le \hdots \le b_{\varphi(q)}] \subseteq \bg$ is thus constant, and therefore contained in a constant subsequence of consecutive elements and of maximal length. We deduce that all the points in the image of the map $|\Delta^{\bg_{\varphi(p),\varphi(q)}}|_A \rightarrow |\Delta^{\bg}|_A$ are mapped to the same point by the quotient map $|\Delta^{\bg}| \rightarrow |\Delta^{\bg}|_A'$. On the other hand, the image of the composition $|\Delta^{\ag_{p,q}}|_A \rightarrow |\Delta^{\ag}| \xrightarrow[]{|\Phi|_A} |\Delta^{\bg}|$ is contained in the image of $|\Delta^{\bg_{\varphi(p),\varphi(q)}}|_A$. We conclude that $|\Phi|_A$ is compatible with the equivalence relations defining $|\Delta^{\ag}|_A'$ and $|\Delta^{\bg}|_A'$, as desired. \qedhere

\end{proof}

\begin{Example}\label{example: quotients of simplex}

    Figures \ref{figure: quotiented stratified 2-simplices} and \ref{figure: quotiented stratified 3-simplices} in section \ref{section: detailed summary of the proof} provide examples of the image by $|-|_A'$ of sequences of length $2$ and $3$.

\end{Example}

Observe that this addresses the difficulty raised at the beginning of this section. Indeed, if $\ag$ is constant, then $|\Delta^{\ag}|_A'$ is a point.

\begin{Notation}\label{notation: left Kan extension of |-|_A' and its right adjoint}
    The left Kan extension of $|-|_A' : \Delta_A \rightarrow \Top_A$ to $\sSet_A$ (construction \ref{construction: left Kan extension}) is also denoted $|-|_A'$. Its right adjoint (construction \ref{construction: right adjoint of left Kan extension}) is denoted $\Sing_{A}'(-)$.
\end{Notation}

\begin{Remark}\label{explicit description of modified realization}
    An explicit description of this left Kan extension can be obtained from remark \ref{remark: explicit description of cubical realization} by replacing $C_A$ by $|-|_A'$. Furthermore, for every $A$-stratified topological space $Y$, $\Sing_A'(Y)$ is the $A$-stratified simplicial set explicitly described as

    $$
    \Sing_A'(Y)_n = \{\sigma : |\Delta^{\ag}|_A' \rightarrow Y \text{ stratified for some } \ag \text{ of length } n\}.
    $$
\end{Remark}

\begin{Remark}\label{remark: natural transformation |-|_A -> |-|_A'}
    The quotient maps from $|\Delta^{\ag}|_A$ to $|\Delta^{\ag}|_A'$ for every $\Delta^{\ag} \in \Delta_A$ define a natural transformation $|-|_A \rightarrow |-|_A'$ of functors from $\Delta_A$ to $\Top_A$, which induces a natural transformation of functors from $\sSet_A$ to $\Top_A$, and a natural transformation $\Sing_A'(-) \rightarrow \Sing_A(-)$. Furthermore, for every $A$-stratified simplicial set $S$, the natural transformation $|-|_A \rightarrow |-|_A'$ exhibits $|S|_A'$ as a \emph{quotient} of $|S|_A$. 
\end{Remark}

\begin{Remark}\label{remark: modification of |-|_A is harmless}
This modification of $|-|_A$ will be somehow harmless for us. We will see that if $(Y,A)$ is a stratified space such that $\Sing_A(Y)$ is an $\infty$-category then $\Sing_A'(Y)$ is one too (proposition \ref{proposition: if Sing_A infinity category then Sing_A' too}). Then we will see that the fact that the strata of $(X,A)$ are contractible implies that the functor $\Sing_A'(X) \rightarrow \Sing_A(X)$ is an equivalence when restricted to a suitable subcategory (proposition \ref{proposition: restricted version into unrestricted is equivalence}).
\end{Remark}

\begin{Proposition}\label{proposition: filtration descends to modification}
    The filtration on $|-|_A$ induced by $f$ (definition \ref{definition: filtration on usual realization}) factors through the natural transformation $|-|_A \rightarrow |-|_A'$.
\end{Proposition}

\begin{proof}
    
    Recall that this filtration is denoted by $f_A$, and is defined by left Kan extension from a natural transformation $|-|_A \rightarrow \underline{\R}$ of functors from $\Delta_A$ to $\Top$. To prove the statement, it suffices to prove that the latter natural transformation factors through the natural transformation $|-|_A \rightarrow |-|_A'$ of functors from $\Delta_A$ to $\Top$. This amounts to proving that for every $\Delta^{\ag} \in \Delta_A$ and every constant subsequence $\ag_{p,q} = [a_p = a_{p+1} = \hdots = a_q] \subset \ag$, the map $f_A(\ag) : |\Delta^{\ag}|_A \rightarrow \R$ is constant on the image of the map $|\Delta^{\ag_{p,q}}| \rightarrow |\Delta^{\ag}|_A$. But the latter image consists of the convex combinations of the form
    
    $$
    \sum_{k=p}^q t_k e_k.
    $$
    
    The image by $f_A(\ag)$ of a convex combination of this form is, by construction \ref{construction: filtration on usual realization}:
    
    $$
    \begin{aligned}
    \sum_{k=p}^q t_kf(a_k) & = \left( \sum_{k=i}^k t_k \right) f(a_p) \\
    & = f(a_p).
    \end{aligned}
    $$
    
    The restriction of $f_A(\ag)$ to this image is therefore constant to $f(a_p)$, and in particular constant, as desired.
\end{proof}

\begin{Definition}\label{definition: filtration on modified realization}
    We refer to the natural transformation $|-|_A' \rightarrow \underline{\R}$ obtained from $f_A$, as the \emph{filtration} on $|-|_A'$ induced by $f$, and denote it by $f_A'$.
\end{Definition}

\begin{Remark}\label{remark: filtration on modified realization explicitly}
    Given an $A$-stratified simplicial set $S$, an explicit description of the map $f_A'(S) : |S|_A' \rightarrow \R$ can be obtained from remark \ref{remark: explicit description of filtration} by replacing $|-|_A$ by $|-|_A'$. Alternatively, regarding $|S|_A'$ as a quotient of $|S|_A$, one can describe $f_A'(S)$ as the map induced by $f_A(S)$ on the quotient.
\end{Remark}

It turns out that there still does not exist an isomorphism between the functors $|-|_A'$ and $C_A$. A similar observation is made in \cite[Remark 2.2.2.6]{HigherTopos}.

\begin{Proposition}\label{proposition: functors not isomorphic}
   
    Assume that $\dim (X) \ge 1$. Then the functors $|-|_A', C_A : \Delta_A \rightarrow \Top_A$ are not isomorphic.
    
\end{Proposition}

\begin{proof}

    The condition that $\dim(X) \ge 1$ is equivalent to the condition that there exist two elements $a,b \in A$ such that $a < b$. Consider the two sequences $\ag_2 = [a = a < b]$ and $\ag_3 = [a=a=a < b]$ in $\Delta_A$. We have two morphisms $\Sigma_2^0$ and $\Sigma_2^1$ in $\Delta_A$ determined by the following two commutative diagrams of posets
    
    $$
\Sigma_2^0= \left(
\vcenter{
\hbox{
\xymatrix{
[3] \ar[rd]_-{\ag_3} \ar[rr]^{\sigma_2^0} && [2] \ar[ld]^-{\ag_2} \\
& A
}
}}
\right)
\qquad \text{and} \qquad
\Sigma_2^1= \left(
\vcenter{
\hbox{
\xymatrix{
[3] \ar[rd]_-{\ag_3} \ar[rr]^{\sigma_2^1} && [2] \ar[ld]^-{\ag_2} \\
& A
}
}}
\right).
$$

    The map
    
    $$
    |\sigma_2^0| \times |\sigma_2^1| : |\Delta^3| \rightarrow |\Delta^2| \times |\Delta^2|
    $$
    
    is injective, so that the map
    
    $$
    |\Sigma_2^0|_A' \times |\Sigma_2^1|_A' : |\Delta^{\ag_3}|_A' \rightarrow |\Delta^{\ag_2}|_A' \times |\Delta^{\ag_2}|_A'
    $$
    
    is injective as well. On the other hand, by definition of $C_A$ and proposition \ref{proposition: image of codegeneracies by F(n)}, the map $C_A(\Sigma_2^0)$ is obtained from the map
    
    $$
    \fonction{g}{I^2 \times [f(b),f(a)]}{I \times [f(b),f(a)]}{(t_1,t_2,s)}{(t_2,s)}
    $$
    
    by modding out $I^2 \times [f(b),f(a)]$ by $\sim_{\ag_3}$ and $I \times [f(b),f(a)]$ by $\sim_{\ag_2}$, while the map $C_A(\Sigma_2^1)$ is obtained from the map
    
    $$
    \fonction{g'}{I^2 \times [f(b),f(a)]}{I \times [f(b),f(a)]}{(t_1,t_2,s)}{(\max(t_1,t_2),s)}
    $$
    
    by modding out $I^2 \times [f(b),f(a)]$ by $\sim_{\ag_3}$ and $I \times [f(b),f(a)]$ by $\sim_{\ag_2}$. However, the product of the maps $g$ and $g'$ is not injective on the interior of $I^2 \times [f(b),f(a)]$, and the equivalence relation $\sim_{\ag_3}$ identifies distinct elements only when these both belong to the boundary of $I^2 \times [f(b),f(a)]$. In particular, the map
    
    $$
    C_A(\Sigma_2^0) \times C_A(\Sigma_2^1) : C_A(\Delta^{\ag_3}) \rightarrow C_A(\Delta^{\ag_2}) \times C_A(\Delta^{\ag_2})
    $$
    
    is not injective, and thus an isomorphism of functors between $C_A$ and $|-|_A'$ cannot exist. \qedhere

\end{proof}

We will show in the next section that, however, there exists an isomorphism between $|-|_A'$ and $C_A$ after restriction to a suitable subcategory. We finish this section by introducing this subcategory.

\begin{Notation}\label{notation: Delta_A^+}
    Let $\Delta_A^+$ denote the subcategory of $\Delta_A$ defined as:

    $$
    \Delta_A^+ = \left\{ 
    \begin{array}{ll}
        \Obj & = \Obj({\Delta_A}) \\
        \Hom(\Delta^{\ag},\Delta^{\bg}) & = \{\phi \in \Hom_{\Delta_A}(\Delta^{\ag},\Delta^{\bg}) \text{ such that the associated arrow of } \Delta \text{ is injective} \}.
    \end{array}
    \right.
    $$
\end{Notation}

The morphisms of $\Delta_A^+$ are generated by the following subset of morphisms of $\Delta_A$:

\begin{Definition}\label{definition: face of stratified simplex}
    A \textit{face morphism} of $\Delta_A$ is a morphism in $\Delta_A$ coming from a commutative diagram of posets of the form

    $$
    \xymatrix{
    [n-1] \ar[rr]^{\delta_n^i} \ar[rd] && [n] \ar[ld]^{\ag} \\
    & A
    }
    $$

    for some $n \ge 1$ and $0 \le i \le n$. This map is called the \emph{face} of $\Delta^{\ag}$ \emph{opposite to the vertex} $i$.
\end{Definition}

\begin{Definition}\label{definition: face of image of simplex by three functors}
    Let $\Delta^{\ag} \in \Delta_A$. The image by $|-|_A$ (resp. $|-|_A'$, $C_A$) of the face of $\Delta^{\ag}$ opposite to the vertex $i$ is called the \emph{face of} $|\Delta^{\ag}|$ (resp. $|\Delta^{\ag}|_A'$, $C_A(\Delta^{\ag})$) \emph{opposite to the vertex} $i$.
\end{Definition}

\subsection{The isomorphism of functors}\label{section: proof of comparison between realizations}

In the previous section, we observed that the functors $|-|_A$ and $C_A$ where not isomorphic (as long as $\dim X > 0)$. We introduced a quotient of the usual stratified geometric realization, which we denoted $|-|_A'$, as well as a subcategory $\Delta_A^+ \subset \Delta_A$. The goal of this section is to prove the following comparison theorem between $|-|_A$ and $C_A$.

\begin{Thms}\label{theorem: isomorphism between |-|'_A and C_A}

There exists an isomorphism of functors between the restrictions of $|-|'_A$ and $C_A$ to $\Delta_A^+$. Moreover, any two such isomorphisms are homotopic.
    
\end{Thms}

By a homotopy between two morphisms of functors from $|-|_A'$ to $C_A$, we mean a morphism of functors $|-|_A' \times [0,1] \rightarrow C_A$ defined on $\Delta_A^+$, where $|-|_A' \times [0,1]$ is defined by $(|-|_A' \times [0,1])(\Delta^{\ag})=|\Delta^{\ag}|_A' \times [0,1]$ with the $A$-stratification induced by the first projection.

\begin{Remark}\label{remark: relation to unproven lemma of CJS}

    In \cite[Section 3]{CohenJonesSegalMorse}, Cohen, Jones and Segal considered two functors which, in our notation, can be described as the functors
    
    $$
    \sd (\Delta_A^+) \rightarrow \Top
    $$
    
    obtained from $|-|_A'$ and $C_A$ by restricting to the full subcategory $\sd (\Delta_A^+) \subseteq \Delta_A^+$ whose objects are the strictly increasing sequences, and corestricting along the forgetful functor $\Top_A \rightarrow \Top$ (note that after this restriction, $|-|_A'$ coincides with $|-|_A$). \cite[Lemma 3.5]{CohenJonesSegalMorse} claims without proof that these two functors are isomorphic. Theorem \ref{theorem: isomorphism between |-|'_A and C_A} establishes this claim.

\end{Remark}

\begin{Remark}\label{remark: we can prove compatibility with filtrations}
    We can actually prove that there exists an isomorphism which is compatible with the respective filtrations on $|-|_A'$ and $C_A$ induced by $f$ (see remark \ref{remark: triviality theorem also holds for standard realization} for more details). However, the proof of our main result does not require this stronger statement, so we will only prove theorem \ref{theorem: isomorphism between |-|'_A and C_A} for the sake of brevity.
\end{Remark}

Our strategy to prove theorem \ref{theorem: isomorphism between |-|'_A and C_A} is to proceed by induction on the length of finite sequences of consecutive elements of $A$. More precisely, assume that we know that such an isomorphism exists in restriction to the sequences of length at most $n-1$ (with $n \ge 1$); we wish to extend it to sequences of length at most $n$. We will see that for every $\Delta^{\ag} \in \Delta_A$, the image of $\partial \Delta^{\ag}$ by $|-|_A'$ (resp. $|-|_{C_A}$) can be described as a colimit over the faces of $|\Delta^{\ag}|_A'$ (resp. $C_A(\Delta^{\ag})$) (corollary \ref{corollary: image of stratified boundary is colimit of image of faces}). Recall that the notion of face, to which we refer here, was introduced in definition \ref{definition: face of stratified simplex}. In particular, if $\ag$ is of length $n$, the induction hypothesis provides a stratified homeomorphism $|\partial \Delta^{\ag}|_A' \simeq |\partial \Delta^{\ag}|_{C_A}$. Towards the goal of performing the induction step and proving that the latter extends to a stratified homeomorphism $|\Delta^{\ag}|_A' \simeq C_A(\Delta^{\ag})$, we will state theorem \ref{theorem: stratified simplices and cubes as balls}. It essentially reduces the question of finding such an extension to that of finding an extension of a homeomorphism of $\SM^{n-1}$ to a homeomorphism of $D^{n}$, which is always possible.

This section is organised as follows. We begin by proving corollary \ref{corollary: image of stratified boundary is colimit of image of faces}. We then state theorem \ref{theorem: stratified simplices and cubes as balls} and from it we deduce the main theorem of this section, theorem \ref{theorem: isomorphism between |-|'_A and C_A}. The proof of theorem \ref{theorem: stratified simplices and cubes as balls} is deferred to the next section.

As announced, our first goal is to give a description of the stratified spaces $|\partial \Delta^{\ag}|_A'$ and $|\partial \Delta^{\ag}|_{C_A}$. We start by exhibiting the simplicial set $\partial \Delta^n$ as a colimit of standard simplices. Informally speaking, $\partial \Delta^n$ is obtained by gluing together the faces of $\Delta^n$ along their common boundary faces. More formally:

\begin{Construction}\label{construction: boundary of simplex as colimit of simplices}

Let $\mathrm{Facets}(\Delta^n)$ denote the category whose objects are the nondegenerate simplices of $\Delta^n$ of dimension $n-1$ and $n-2$, and whose non-identity morphisms correspond to the inclusions of $(n-2)$-simplices as boundary faces of $(n-1)$-simplices in $\Delta^n$. This category can be described as follows.

\begin{itemize}
    \item For every integer $0 \le i \le n$, there is an object $F_i$ associated with the face of $\Delta^n$ opposite to the vertex $i$.
    
    \item For any two integers $0 \le i < j \le n$, there is one object $F_{i,j}$ associated with the common boundary face of the face opposite to the vertex $i$, and the face opposite to the vertex $j$.
    
    \item For any two integers $0 \le i < j \le n$, there are two non-identity morphisms $\phi_{i,j}^i : F_{i,j} \rightarrow F_i$ and $\phi_{i,j}^j : F_{i,j} \rightarrow F_j$, associated with the inclusion of boundary faces.
\end{itemize}

There is a functor $\Gd^n : \mathrm{Facets}(\Delta^n) \rightarrow \sSet$ defined by

$$
\left\{
\begin{array}{l}
     F_i \mapsto \Delta^{n-1}  \\
     F_{i,j} \mapsto \Delta^{n-2} \\
     \phi_{i,j}^i \mapsto \delta_{n-1}^{j-1} \\
     \phi_{i,j}^j \mapsto \delta_{n-1}^i
\end{array}
\right.
$$

and a natural transformation $\Gd^n \rightarrow \underline{\Delta^n}$ (where $\underline{\Delta^n}$ denotes here the constant functor at $\Delta^n$ on $\mathrm{Facets}(\Delta^n)$) defined by

$$
F_i \mapsto (\delta_n^i : \Delta^{n-1} \rightarrow \Delta^n).
$$

The latter factors through $\partial \Delta^n$ and we have the following result.

\end{Construction}

\begin{Proposition}\label{proposition: boundary of standard simplex as colimit of its faces}
    The natural transformation $\Gd^n \rightarrow \underline{\partial \Delta^n}$ exhibits $\partial \Delta^n$ as a colimit of $\Gd^n$. \qed
\end{Proposition}

This statement generalizes to the case where $\Delta^n$ is endowed with an $A$-stratification, as we now explain. Let $\ag$ be a finite increasing sequence of elements of $A$ of length $n$. Recall that, equivalently, this is an $A$-stratification of $\Delta^n$, as explained in example \ref{example: stratification of Delta^n associated to an element of Delta_A}. Recall also from notation \ref{notation: stratifications of simplex, boundary and horns} that we denote by $\Delta^{\ag}$ (resp. $\partial \Delta^{\ag}$) the simplicial set $\Delta^n$ (resp. $\partial \Delta^n$) endowed with this $A$-stratification (resp. its restriction to $\partial \Delta^n$). 

We then have the following generalization of construction \ref{construction: boundary of simplex as colimit of simplices} and proposition \ref{proposition: boundary of standard simplex as colimit of its faces} to the $A$-stratified setting.

\begin{Proposition}\label{proposition: boundary of stratified simplex as colimit of faces}
    There exists a unique lift of the functor $\Gd^n$ to a functor $\Gd^{\ag} : \mathrm{Facets}(\Delta^n) \rightarrow \sSet_A$ such that the natural transformation $\Gd^n \rightarrow \underline{\Delta^n}$ lifts to a natural transformation $\Gd^{\ag} \rightarrow \underline{\Delta^{\ag}}$. Moreover, the latter factors through $\partial \Delta^{\ag}$ and the natural transformation $\Gd^{\ag} \rightarrow \underline{\partial \Delta^{\ag}}$ exhibits $\partial \Delta^{\ag}$ as a colimit of $\Gd^{\ag}$.
\end{Proposition}

\begin{proof}
The first part of the proposition follows from the construction of $\Gd^n$ and the natural transformation $\Gd^n \rightarrow \underline{\Delta^n}$ (construction \ref{construction: boundary of simplex as colimit of simplices}). The second part follows from proposition \ref{proposition: boundary of standard simplex as colimit of its faces} and the fact that the forgetful functor $\sSet_A \rightarrow \sSet$ commutes with colimits (proposition \ref{proposition: stratified simplicial sets admit small limits and colimits}).
\end{proof}

\begin{Corollary}\label{corollary: image of stratified boundary is colimit of image of faces}
   The $A$-stratified spaces $|\partial \Delta^{\ag}|_A$, $|\partial \Delta^{\ag}|_A'$ and $|\partial \Delta^{\ag}|_{C_A}$ are, respectively, naturally identified with $\colim (|-|_A \circ \Gd^{\ag})$, $\colim (|-|_A' \circ \Gd^{\ag})$ and $\colim (|-|_{C_A} \circ \Gd^{\ag})$. In other words, $|\partial \Delta^{\ag}|_A$ (resp. $|\partial \Delta^{\ag}|_A'$, $|\partial \Delta^{\ag}|_{C_A}$) is obtained by gluing together the faces of $|\Delta^{\ag}|_A$ (resp. $|\Delta^{\ag}|_A'$, $C_A(\Delta^{\ag})$), in the sense of definition \ref{definition: face of image of simplex by three functors}, along their common faces. 
\end{Corollary}

\begin{proof}
     Since $|-|_A, |-|_A',|-|_{C_A} : \sSet_A \rightarrow \Top_A$ are left Kan extensions, they commute with colimits. The desired result therefore follows from proposition \ref{proposition: boundary of stratified simplex as colimit of faces}.
\end{proof}

Combining corollary \ref{corollary: image of stratified boundary is colimit of image of faces} with the following theorem, we will be able to prove theorem \ref{theorem: isomorphism between |-|'_A and C_A}.

\begin{Thms}\label{theorem: stratified simplices and cubes as balls}
    Let $G$ denote either of the two functors $|-|_A'$ and $|-|_{C_A}$, and assume that we have a non-constant finite increasing sequence $\ag = [a_0 \le \hdots \le a_n]$ of elements of $A$. There exists a homeomorphism between $D^n$ and $G(\Delta^{\ag})$ that restricts to a homeomorphism between $\SM^{n-1}$ and $G(\partial \Delta^{\ag})$. That is, there exists a commutative square
    
    $$
    \xymatrix{
    \SM^{n-1} \ar[r] \ar[d] & D^n \ar[d] \\
    G(\partial \Delta^{\ag}) \ar[r] & G(\Delta^{\ag})
    }
    $$
    
    where the two vertical arrows are homeomorphisms, the top horizontal one is the inclusion, and the bottom horizontal one is the one induced by the inclusion $\partial \Delta^{\ag} \rightarrow \Delta^{\ag}$.
    
\end{Thms}

The proof of this theorem is quite long and is the purpose of section \ref{section: proof of the ball theorem for realizations}.

\begin{proof}[Proof of theorem \ref{theorem: isomorphism between |-|'_A and C_A}]
     The statement that there exists an isomorphism of functors means that there exists a family of stratified homeomorphisms
     
     $$(|\Delta^{\ag}|_A' \simeq C_A(\Delta^{\ag}))_{\Delta^{\ag} \in \Delta_A}$$
     
     compatible with the faces as defined in \ref{definition: face of image of simplex by three functors} (because faces generate morphisms in $\Delta_A^+$). We prove that such a family exists by induction on the length of $\ag$. For $\ag$ of length $0$, $C_A(\Delta^{\ag})$ and $|\Delta^{\ag}|_A'$ are both points. Their unique strata both correspond to the element $a_0 \in A$. Hence, a stratified homeomorphism between $C_A(\Delta^{\ag})$ and $|\Delta^{\ag}|_A'$ exists (and is unique).

    Suppose that a family of stratified homeomorphisms $(|\Delta^{\ag}|_A' \simeq C_A(\Delta^{\ag}))_{\Delta^{\ag} \in \Delta_A, l(\ag) \le n-1}$ compatible with faces exists for some $n \ge 1$. To perform the induction step, we will show that one can extend this family to every $\Delta^{\ag} \in \Delta_A$ such that $l(\ag) = n$, in a way compatible with faces.

    Assume first that $\ag$ is a constant sequence. Then, again, $C_A(\Delta^{\ag})$ and $|\Delta^{\ag}|_A'$ are both points, and their unique strata correspond to the element $a_0 \in A$. Hence a stratified homeomorphism between $C_A(\Delta^{\ag})$ and $|\Delta^{\ag}|_A'$ exists (and is unique). Assume now that $\ag$ is not constant. By the induction hypothesis, there exist stratified homeomorphisms between the faces of $|\Delta^{\ag}|_A'$ and the corresponding faces of $C_A(\Delta^{\ag})$, which coincide in restriction to their common faces. In other words, by corollary \ref{corollary: image of stratified boundary is colimit of image of faces}, there exists a stratified homeomorphism $h : |\partial \Delta^{\ag}|_A' \overset{\simeq}{\rightarrow} |\partial \Delta^{\ag}|_{C_A}$. Pick two homeomorphisms $h_1 : |\Delta^{\ag}|_A' \overset{\simeq}{\rightarrow} D^n$ and $h_2 : C_A(\Delta^{\ag}) \overset{\simeq}{\rightarrow} D^n$ as in theorem \ref{theorem: stratified simplices and cubes as balls}. The composition $(h_2) \vert_{|\partial \Delta^{\ag}|_{C_A}} \circ h \circ (h_1^{-1}) \vert_{\SM^{n-1}}$ is then a homeomorphism of $\SM^{n-1}$. The latter extends to a homeomorphism $H$ of $D^n$. By construction, the homeomorphism $h_2^{-1} \circ H \circ h_1 : |\Delta^{\ag}|_A' \overset{\simeq}{\rightarrow} C_A(\Delta^{\ag})$ is compatible with faces. It restricts to the homeomorphism $h$ on $|\partial \Delta^{\ag}|_A'$, and therefore preserves the stratifications in restriction to $|\partial \Delta^{\ag}|_A'$. Moreover it maps $|\Delta^{\ag}|_A' \backslash |\partial \Delta^{\ag}|_A'$ to $C_A(\Delta^{\ag}) \backslash |\partial \Delta^{\ag}|_{C_A}$, both of which lie in the $a_n$-strata. Hence $h_2^{-1} \circ H \circ h_1$ is a stratified homeomorphism. This completes the proof of the existence of an isomorphism of functors. The fact that any two isomorphisms are homotopic follows from a similar argument, using the fact that for every positive integer $n$, given two continuous maps from $(D^n, \SM^{n-1})$ to itself, every homotopy between their restrictions to $\SM^{n-1}$ can be extended to a homotopy between them on $D^n$. \qedhere
\end{proof}

\subsection{The proof of theorem \ref{theorem: stratified simplices and cubes as balls}}\label{section: proof of the ball theorem for realizations}

This section consists of the proof of theorem \ref{theorem: stratified simplices and cubes as balls}. It is divided into five paragraphs. The first one collects some general facts about polytopes, a notion that will play an essential role in our proof. The second one consists of a summary of the proof, while the remaining three complete the proof.

\subsubsection{Polytopes}\label{section: polytopes}

In this section, we record several definitions and results about polytopes. We use the book of Ziegler (\cite{ZieglerLectures}) as a reference.

\begin{itemize}
    \item A \emph{polytope} is a nonempty topological subspace $P$ of some Euclidean space $\R^d$ satisfying the following two equivalent conditions.

    \begin{enumerate}[label=(\roman*)]

    \item $P$ is the convex hull of a finite set of points in $\R^d$.

    \item $P$ is bounded and is the intersection of a finite family of closed affine halfspaces of $\R^d$.

    \end{enumerate}

    For a proof of the equivalence between these two conditions, see \cite[Theorem 1.1]{ZieglerLectures}.
    
    \begin{Example}\label{example: realization of simplex is polytope}
        For every nonnegative integer $n$, the topological space $|\Delta^n| \subset \R^{n+1}$ is a polytope.
    \end{Example}

    \item Let $P \subset \R^d$ be a polytope. The \emph{dimension} of $P$ is the dimension of the affine subspace of $\R^d$ spanned by $P$.
    
    \item Let $P \subset \R^d$ be a polytope. Given $c \in \R^d$, we say that the inequality $\langle c,x \rangle \le c_0$ is \emph{valid} for $P$ if it is satisfied for every point $x \in P$ (where $\langle \cdot , \cdot \rangle$ denotes the Euclidean scalar product). A \emph{face} of $P$ is a subset of $P$ of the form $\{ x \in P \mid \langle c,x \rangle = c_0\}$ for some valid inequality $\langle c,x \rangle \le c_0$.
    
    \item Taking $c = c_0 = 0$, we see that $P$ is a face of itself. A face of $P$ that is different from $P$ is called a \emph{proper} face.
    
    \item The \emph{boundary} of $P$ is the union of the proper faces of $P$. It is denoted $\partial P$. The \emph{interior} of $P$ is $P \backslash \partial P$. It is denoted $\Int(P)$. If $P \subset \R^d$ and $P$ has dimension $d$, these respectively coincide with the boundary of $P$ in $\R^d$, and the interior of $P$ in $\R^d$.
    
    \item The intersection of two faces of a polytope is a face (\cite[Proposition 2.3]{ZieglerLectures}).
    
    \item A face of a polytope is a polytope (\cite[Proposition 2.3]{ZieglerLectures}). The \emph{dimension} of a face is the dimension of the corresponding polytope.
    
    \begin{Remark}\label{remark: two notions of face}
        We warn the reader of the following conflict of terminology: for every $n \ge 1$, we have already introduced a notion of face of $\Delta^n \in \sSet$ (notation and terminology \ref{notation and terminology about simplicial sets}), as well as a corresponding notion of face of $|\Delta^n| \in \Top$ (definition \ref{definition: face of image of simplex by three functors}). The latter coincides with the notion of $(n-1)$-dimension face of $|\Delta^n|$ regarded as an $n$-dimensional polytope.
    \end{Remark}
    
    \item A \emph{vertex} of a polytope is a face of dimension $0$, that is, a face that is reduced to a point. Every polytope is the convex hull of its vertices (\cite[Proposition 2.2]{ZieglerLectures}).

\end{itemize}

The following lemma about polytopes will be essential to us.

\begin{Lemmas}\label{lemma: intersection of halfline with polytope boundary}
    Let $P \subset \R^d$ be a $d$-dimensional polytope and $p$ a point belonging to the interior of $P$. Let $q \in \R^d$ be a point different from $p$. The intersection of $\partial P$ with the half-line $p + \R_+(q-p)$ is a point.
\end{Lemmas}

\begin{proof}

    We denote this half-line by $l$. It can be identified with $\R_+$ and in particular it inherits an order relation from $\R_+$.
    
    We start by proving that $l \cap \partial P \neq \varnothing$. Since $P$ is bounded, the set $\{r \in l \mid r \notin P\}$ is nonempty, and we define $m$ to be its infimum. Since the interior of $P$ is open in $\R^d$, we have $m \notin \Int(P)$, and since $P$ is closed in $\R^d$ we have $m \in P$. Hence $m \in \partial P$.
    
    We conclude by proving that $l \cap \partial P$ is reduced to a point. Let $m$ be the infimum of $l \cap \partial P$; since $l \cap \partial P$ is closed we have $m \in l \cap \partial P$. Let us write $P$ as a finite intersection of affine closed halfspaces of $\R^d$
    
    $$
    P = \bigcap_{1 \le i \le n} \{x \in \R^d \mid \langle c_i,x\rangle \le c_{0,i} \}.
    $$
    
    Since $m \in \partial P$, there exists $1 \le j \le n$ such that $\langle c_j,m \rangle = c_{0,j}$. Since $p$ lies in the interior of $P$ we have $\langle c_j,p \rangle < c_{0,j}$, and since $p \in l$ and $p < m$ we deduce that $\langle c_j,m' \rangle > c_{0,j}$ for every $m ' \in l$ such that $m' > m$. In particular, for every $m' \in l$ such that $m' > m$, we have $m' \notin P$, as desired. \qedhere

\end{proof}

We will now prove a result that makes polytopes useful for the proof of theorem \ref{theorem: stratified simplices and cubes as balls}. For this, consider a $d$-dimensional polytope $P \subset \R^d$, so that the interior of $P$ is an open subset of $\R^d$. Denote by $\SM^{d-1}$ (resp $D^d$) the Euclidean sphere (resp. closed ball) of radius $1$ centered at $0$.

\begin{Proposition}\label{proposition: a polytope is a disc}

    Let $p$ be a point belonging to the interior of $P$. The following hold.
    
    \begin{enumerate}[label=(\alph*)]
     \item The map 
     
    $$
    \fonction{\varphi}{\partial P \times [0,1]}{D^d}{(q,s)}{s \cdot \frac{q-p}{||q-p||}}
    $$
    
    induces a homeomorphism between $(C(\partial P),\partial P \times \{1\})$ and $(D^d,\SM^{d-1})$. Here, $C(\partial P)$ denotes the closed cone of $P$, i.e., the space $(\partial P \times [0,1]) / (\partial P \times \{0\})$.
     
     \item The map
     
     $$
     \fonction{\Phi_p^P}{\partial P \times \R_+}{\R^d}{(q,s)}{p+s(q-p)}
     $$
     
     induces a homeomorphism between the open cone of $\partial P$, i.e., the space $(\partial P \times \R_+)/(\partial P \times \{0\})$, and $\R^d$. 
     
     \item The latter homeomorphism restricts to a homeomorphism between $(C(\partial P),\partial P \times \{1\})$ and $(P, \partial P)$. In particular, combining this with item (a), we obtain a homeomorphism between the pairs $(P,\partial P)$ and $(D^d,\SM^{d-1})$.
     
     \end{enumerate}

\end{Proposition}

\begin{proof}

    Let us prove that the map $C(\partial P) \rightarrow D^d$ induced by $\varphi$ is a homeomorphism. Since its source is compact and its target is Hausdorff, it suffices to prove that this is a bijection. To this end, we exhibit its inverse. Let $v \in D^d$ different from $0$. By virtue of lemma \ref{lemma: intersection of halfline with polytope boundary}, the half-line starting at $p$ and passing through $p + v$ intersects $\partial P$ at exactly one point, which we denote by $q$. The map that carries $v$ to $(q,||v||)$ and that carries $0$ to the image of $\partial P \times \{0\}$ in $C(\partial P)$ is the inverse we are looking for. This completes the proof of item (a).
    
    Let us prove that $\Phi_p^P$ is surjective. First, we have $\Phi_p^P(q,0) = p$ for every $q \in \partial P$, and therefore $p$ belongs to the image of $\Phi_p^P$. Otherwise, let $x \in \R^d$ such that $x \neq p$. By lemma \ref{lemma: intersection of halfline with polytope boundary}, the half-line $p + \R_+(x-p)$ intersects $\partial P$ at exactly one point $q$. In particular, we have $\Phi_p^P(q,s)$ for some $0 < s$.
    
    Note also that the intersection of the half-line $p + \R_+(x-p)$ with $P$ is the segment $[p,q]$, so that if $x \in P$, then $s \le 1$. In particular, the image of $\partial P \times [0,1]$ by $\Phi_p^P$ is $P$.
    
    The map $\partial P \times \R_+ \rightarrow P$ factors through $(\partial P \times \R_+)/(\partial P \times \{0\})$, we now prove that this factorization is a homeomorphism. We start by showing that it is injective. Since $p \notin \partial P$, we have $(\Phi_p^P)^{-1}(p) = \partial P \times \{0\}$. We therefore wish to prove that $\Phi_p^P$ is injective in restriction to $\partial P \times (0,+\infty)$. Now, if we have two elements $(q,s),(q',s') \in \partial P \times (0,+\infty)$ such that $\Phi_p^P(q,s)=\Phi_p^P(q',s')$, then the two half-lines $p+\R_+(q-p)$ and $p+\R_+(q'-p)$ must coincide, and therefore we must have $q=q'$, by lemma \ref{lemma: intersection of halfline with polytope boundary}. Consequently, we must also have $s=s'$.
    
    To conclude, we observe that $\Phi_p^P$ induces a continuous bijection between the one point compactifications of $(\partial P \times \R_+)/(\partial P \times \{0\})$ and $\R^k$. Since these one point compactifications are Hausdorff, this bijection is a homeomorphism, and therefore $\Phi_p^P$ as well, as desired. This completes the proof of items (b) and (c). \qedhere

\end{proof}

\subsubsection{Summary of the proof}\label{section: summary of the proof of the ball theorem}

Let us fix a nonconstant sequence $\ag = [a_0 \le \hdots \le a_n]$ of elements of $A$, of length $n$. In this section, we summarize the proof of theorem \ref{theorem: stratified simplices and cubes as balls}. We start with the case $G=|-|_A'$. For each maximally constant subsequence $[a_p=a_{p+1}=\hdots=a_q] \subset \ag$, i.e., satisfying $a_{p-1} < a_p$ if $p \ge 1$ and $a_{q+1} > a_q$ is $q \le n-1$, the convex hull of the set of vertices $\{e_p,e_{p+1},\hdots,e_q \}$ of the polytope $|\Delta^n|$ is a face of $|\Delta^n|$. By virtue of the assumption that $\ag$ is not constant, this face is proper. By construction of the functor $|-|_A'$ (construction \ref{construction: functor |-|_A'}), the space $|\Delta^{\ag}|_A'$ is obtained from $|\Delta^n|$ by collapsing each of these proper faces to a point. Note that these faces are disjoint. We will prove the following lemma.

\begin{Lemmas}\label{lemma: proper face collapsing lemma}
    Let $P$ be a polytope, $F \subset P$ a proper face and $U$ a neighborhood of $F$ in $P$. There exists a map $c : P \rightarrow P$ satisfying the following properties.
    
    \begin{enumerate}[label=(\roman*)]
        \item The map $c$ is the identity outside of $U$.
    
        \item The map $c$ is constant on $F$.

        \item The resulting map $P/F \rightarrow P$ is a homeomorphism.
        
        \item This homeomorphism restricts to a homeomorphism $\partial P / F \simeq \partial P$.
    \end{enumerate}
\end{Lemmas}

This lemma will be a consequence of an improved statement (lemma \ref{lemma: improved proper face collapsing lemma}), whose proof will be given in section \ref{section: collapsing lemma}. The proof of theorem \ref{theorem: stratified simplices and cubes as balls} in the case $G=C_A$ will follow from this lemma together with the fact that an $n$-dimensional polytope is homeomorphic to an $n$-dimensional disc (proposition \ref{proposition: a polytope is a disc}); the details can be found in section \ref{section: collapsing lemma}.

We now summarize the proof of theorem \ref{theorem: stratified simplices and cubes as balls} in the case when $G=C_A$. As in the case when $G=|-|_A'$, the space $G(\Delta^{\ag})$ is obtained from a polytope (namely the polytope $I^{n-1} \times [f(a_n),f(a_0)]$) by modding out by an equivalence relation. However, in general, passing to the quotient does not amount to collapsing disjoint faces, nor to collapsing full faces, which leads to extra difficulties compared to the case $G=|-|_A'$. In order to establish the desired result, we will make use of the filtration on $C_A(\Delta^{\ag})$ induced by $f$ (construction \ref{construction: filtration on C A}). Recall that this is the map $f_{C_A}(\Delta^{\ag}) : C_A(\Delta^{\ag}) \rightarrow \R$ obtained from the projection map

$$
p_{\ag} : I^{n-1} \times [f(a_n),f(a_0)] \rightarrow [f(a_n),f(a_0)] \subset \R
$$

by quotienting $I^{n-1} \times [f(a_n),f(a_0)]$ by the equivalence relation $\sim_{\ag}$. We will prove the following theorem.

\begin{Thms}\label{theorem: filtration on cube is trivial}
    The filtration on $C_A(\Delta^{\ag})$ induced by $f$ satisfies the following properties.
    
    \begin{enumerate}[label=(\roman*)]
        \item The preimages of $f(a_0)$ and $f(a_n)$ by $f_{C_A}(\Delta^{\ag})$ are points.
        
        \item The map $f_{C_A}(\Delta^{\ag})$ can be identified with a $(n-1)$-dimensional disc bundle over the open interval $(f(a_n),f(a_0))$, that is, there exists a pullback square
        
        $$
        \xymatrix{
        D^{n-1} \times (f(a_n),f(a_0)) \ar[d] \ar[r] & C_A(\Delta^{\ag}) \ar[d]^-{f_{C_A}(\Delta^{\ag})} \\
        (f(a_n),f(a_0)) \ar[r] & \R
        }
        $$
        
        where the left vertical arrow is the projection and the bottom horizontal arrow is the inclusion.
        
        \end{enumerate}
        
        Given a topological space $Y$, we denote by $\Sigma (Y)$ the \emph{suspension} of $Y$, defined as $\Sigma(Y) = (Y \times [0,1] \sqcup *_1 \sqcup *_2) / (Y \times \{0\} \sim *_1, Y\times \{1\} \sim *_2)$. Consider a pullback square as in item (ii). By item (i), using a homeomorphism $[f(b),f(a)] \simeq [0,1]$, the top horizontal arrow of the latter extends to a homeomorphism $\Sigma(D^{n-1}) \simeq C_A(\Delta^{\ag})$.
        
        \begin{enumerate}[label=(\roman*), start=3]
        
        \item The homeomorphism $\Sigma(D^{n-1}) \simeq C_A(\Delta^{\ag})$ restricts to a homeomorphism between $\Sigma(\SM^{n-2})$ and $|\partial \Delta^{\ag}|_{C_A}$, that is, there exists a homeomorphism $\Sigma(\SM^{n-2}) \simeq |\partial \Delta^{\ag}|_{C_A}$ such that the following diagram commutes
        
         $$
        \xymatrix{
        \Sigma(\SM^{n-2}) \ar[r] \ar[d]_-{\simeq} & \Sigma(D^{n-1}) \ar[d]^-{\simeq} \\
        |\partial \Delta^{\ag}|_{C_A} \ar[r] & C_A(\Delta^{\ag})
        }
        $$
        
        where the top horizontal arrow is induced by the inclusion $\SM^{n-2} \rightarrow D^{n-1}$, and the bottom horizontal arrow is induced by the inclusion $\partial \Delta^{\ag} \rightarrow \Delta^{\ag}$.
    \end{enumerate}
\end{Thms}

The proof of theorem \ref{theorem: stratified simplices and cubes as balls} in the case $G=C_A$ follows from item (iii) in theorem \ref{theorem: filtration on cube is trivial} together with the fact that the pair $(\Sigma(D^{n-1}),\Sigma(\SM^{n-2}))$ is homeomorphic to the pair $(D^n,\SM^{n-1})$.

\begin{Remark}\label{remark: triviality theorem also holds for standard realization}
    Theorem \ref{theorem: filtration on cube is trivial} can be regarded as a triviality theorem for the filtration on $C_A$ induced by $f$, at the level of stratified simplices. We can actually prove that this theorem also holds after replacing $C_A$ by $|-|_A'$ and $f_{C_A}$ by $f_A'$. With this stronger statement at hand about $|-|_A'$, one can adapt the proof of theorem \ref{theorem: isomorphism between |-|'_A and C_A} in order to prove the existence of an isomorphism between the restrictions of $C_A$ and $|-|_A'$ to $\Delta_A^+$ which is compatible with the filtrations. For the sake of brevity, we omit the details in this paper.
\end{Remark}

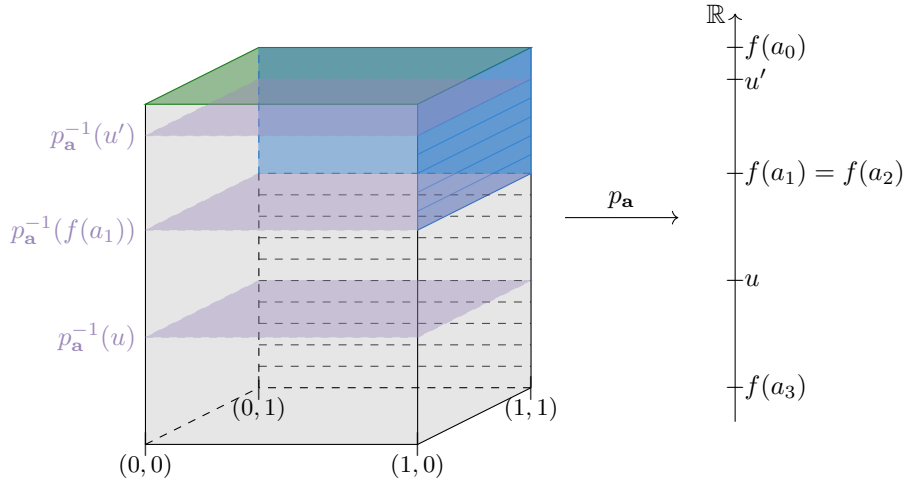
\begin{figure}[H]
    \centering
    \begin{tikzpicture}[scale=1.5]
    
        
        
            \coordinate (A) at (-1.2,-1.5);
            \coordinate (B) at (1.2,-1.5);
            \coordinate (E) at (1.2,0.4);
            \coordinate (H) at (-1.2,1.5);
            \coordinate (C) at (2.2,-1);
            \coordinate (D) at (-0.2,-1);
            \coordinate (J) at (2.2,2);
            \coordinate (I) at (1.2,1.5);
            \coordinate (K) at (-0.2,2);
            
            
            \coordinate (G) at ($(D)!0.63!(K)$);
            \coordinate (E) at ($(B)!0.63!(I)$);
            \coordinate (F) at ($(C)!0.63!(J)$);
            \coordinate (P) at ($(A)!0.63!(H)$);
            
            
            \coordinate (Q) at ($(H)!0.25!(P)$);
            \coordinate (R) at ($(I)!0.25!(E)$);
            \coordinate (S) at ($(J)!0.25!(F)$);
            \coordinate (T) at ($(K)!0.25!(G)$);
            
            
            \coordinate (L) at ($(A)!0.5!(P)$);
            \coordinate (M) at ($(B)!0.5!(E)$);
            \coordinate (N) at ($(C)!0.5!(F)$);
            \coordinate (O) at ($(D)!0.5!(G)$);
            
            
            \coordinate (Z) at (2.5,0.5);
            \coordinate (a) at (3.5,0.5);
            
            
            \coordinate (U) at (4,-1);
            \coordinate (Y) at (4,2);
            \coordinate (b) at (4,-1.3);
            \coordinate (c) at (4,2.3);
            \coordinate (W) at ($(U)!0.63!(Y)$);
            \coordinate (X) at ($(W)!0.75!(Y)$);
            \coordinate (V) at ($(U)!0.5!(W)$);

        \node[below] at (A) {\small{$(0,0)$}};
        \node at (A) {\small{$|$}};
        \node[below] at (B) {\small{$(1,0)$}};
        \node at (B) {\small{$|$}};
        \node[below] at (C) {\small{$(1,1)$}};
        \node at (C) {\small{$|$}};
        \node[below] at (D) {\small{$(0,1)$}};
        \node at (D) {\small{$|$}};
            
         \draw (H)--(A)--(B)--(E);
         \draw[dashed] (A)--(D);
         \draw[dashed] (C)--(D);
         \draw[dashed] (D)--(G);
         \draw (C)--(F);
         \draw (B)--(C);
             

            \draw[green] (H)--(I)--(J)--(K)--cycle;
            \draw[fill, green, opacity=0.4] (H)--(I)--(J)--(K)--cycle;

            \draw[blue] (I)--(E)--(F)--(J);
            \draw[fill, blue, opacity=0.4] (I)--(E)--(F)--(J)--cycle;
            \draw[dashed, blue] (F)--(G)--(K);
            \draw[fill, blue, opacity=0.4] (K)--(G)--(F)--(J)--cycle;

            \draw[fill, opacity = 0.1] (-1.2,-1.5)--(1.2,-1.5)--(2.2,-1)--(2.2,2)--(-0.2,2)--(-1.2,1.5)--cycle;
            
        
        \foreach \t in {0.1,0.2,...,0.9} {
        \draw[dashed, opacity = 0.7] 
            ($(D)!\t!(G)$) -- ($(C)!\t!(F)$);
        }
        
        \foreach \t in {0.15,0.3,...,0.9} {
        \draw[blue, opacity = 0.7] 
            ($(E)!\t!(I)$) -- ($(F)!\t!(J)$);
        }

             
            
            \draw[fill, dashed, sectioncolor, opacity = 0.4] (Q)--(R)--(S)--(T)--cycle;
            \node[left, sectioncolor] at (Q) {$p_{\ag}^{-1}(u')$};
                
                
            \draw[fill, dashed, sectioncolor, opacity = 0.4] (P)--(E)--(F)--(G)--cycle;
            \node[left, sectioncolor] at (P) {$p_{\ag}^{-1}(f(a_1))$};
                
                
            \draw[fill, dashed, sectioncolor, opacity = 0.4] (L)--(M)--(N)--(O)--cycle;
            \node[left, sectioncolor] at (L) {$p_{\ag}^{-1}(u)$};
            
        
        \draw[->] (Z)--(a) node[midway, above] {$p_{\ag}$};
            
        
            
            \draw[->] (b)--(c) node[left] {$\R$};
            
            
            \node at (Y) {$-$};
            \node at (X) {$-$};
            \node at (W) {$-$};
            \node at (V) {$-$};
            \node at (U) {$-$};
            
            
            \node[right] at (Y) {$f(a_0)$};
            \node[right] at (X) {$u'$};
            \node[right] at (W) {$f(a_1)=f(a_2)$};
            \node[right] at (V) {$u$};
            \node[right] at (U) {$f(a_3)$};

         \end{tikzpicture}
    \caption{Three level sets of $p_{\ag}$ are depicted: the one associated with the critical value $f(a_1)$, the one associated with a real number $u$ strictly comprised between $f(a_3)$ and $f(a_1)$, and the one associated with a real number $u'$ strictly comprised between $f(a_1)$ and $f(a_0)$. Recall that the stratified space $C_A(\Delta^{\ag})$ is obtained from this cube by collapsing to a point each of the two faces $I^2 \times \{f(a_0)\}$ and $I^2 \times \{f(a_3)\}$, as well as each of the horizontal segments $I \times \{1\} \times \{s\}$ for $f(a_3) \le s \le f(a_1)$, and each of the horizontal segments $\{1\} \times I \times \{s\}$ for $f(a_1) \le s \le f(a_0)$.}
    \label{figure: filtration on 3 cube}
\end{figure}

We now summarize our strategy to prove theorem \ref{theorem: filtration on cube is trivial}. Figure \ref{figure: filtration on 3 cube} illustrates the situation in the special case when $\ag$ is of the form \textcolor{green}{$a_0$} < \textcolor{blue}{$a_1$} = \textcolor{blue}{$a_2$} < $a_3$. Recall that we denote by $p_{\ag}$ the projection map $I^2 \times [f(a_3),f(a_0)] \rightarrow [f(a_3),f(a_0)]$.

Figure \ref{figure: filtration on 3 cube} illustrates the following general facts.

\begin{itemize}
    
    \item For every $f(a_n) \le s \le f(a_0)$, the equivalence relation $\sim_{\ag}$ induces an equivalence relation $\sim_{\ag,s}$ on $I^{n-1}$, by identifying $I^{n-1}$ with the subset $I^{n-1} \times \{s\} \subset I^{n-1} \times [f(a_n),f(a_0)]$. Moreover, for every $(\tu,s), (\tu',s') \in I^{n-1} \times [f(a_n),f(a_0)]$, if $(\tu,s) \sim_{\ag} (\underline{t'},s')$, then $s = s'$. In other words, the equivalence relation $\sim_{\ag}$ can be recovered from the equivalence relations $(\sim_{\ag,s})_{f(a_n) \le s \le f(a_0)}$, in the sense that the following function is a bijection
    
    $$
    \bigsqcup_{f(a_n) \le s \le f(a_0)} I^{n-1} / {\sim_{\ag,s}} \rightarrow I^{n-1} \times [f(a_n),f(a_0)] / \sim_{\ag}.
    $$
    
    \item As $s$ moves along $[f(a_n),f(a_0)]$ without encountering any critical value among the $f(a_i)$'s, the equivalence relation $\sim_{\ag,s}$ does not change. That is, for every $0 \le i < n$ such that $a_i < a_{i+1}$ and every $f(a_{i+1}) < s, s' < f(a_i)$, we have $\sim_{\ag,s} = \sim_{\ag,s'}$. Consequently, denoting this equivalence relation on $I^{n-1}$ by $\sim_{\ag,f(a_i),f(a_{i+1})}$, we have a pullback square
    
    $$
    \xymatrix{
    I^{n-1}/\sim_{\ag,f(a_i),f(a_{i+1})} \times (f(a_{i+1}),f(a_i)) \ar[d] \ar[r] & C_A(\Delta^{\ag}) \ar[d]^-{f_{C_A}(\Delta^{\ag})} \\
    (f(a_{i+1}),f(a_i)) \ar[r] & \R.
    }
    $$
    
\end{itemize}

Recall that we wish to prove that the map $f_{C_A}(\Delta^{\ag})$ can be identified with a $(n-1)$-dimensional disc bundle over $(f(a_n),f(a_0))$. In the special case depicted on figure \ref{figure: filtration on 3 cube}, we have the following.

\begin{itemize}

    \item The quotient $I^2 / \sim_{\ag,f(a_1),f(a_3)}$ can be naturally identified with a polytope as follows. Consider the convex hull of the three points $(0,0), (1,0)$ and $(0,1)$ in $\R^2$. This is a polytope, which we denote by $Q$. This polytope is also the intersection of $I^2$ with the affine halfspace $\{(x_1,x_2) \in \R^2 \mid x_1+x_2 \le 1\}$. The map depicted on the following figure is compatible with $\sim_{\ag,f(a_1),f(a_3)}$, and yields a homeomorphism between $I^2 / \sim_{\ag,f(a_1),f(a_3)}$ and $Q$.
    
    \begin{figure}[H]
        \centering
        \begin{tikzpicture}[scale = 2]
        
            
            
                \coordinate (A) at (-2,0);
                \coordinate (B) at (-1,0);
                \coordinate (C) at (-1,1);
                \coordinate (D) at (-2,1);
                \coordinate (K) at (-2,0.5);
                \coordinate (L) at (-1,0.5);
                \coordinate (M) at (-1.5,1);
                \coordinate (N) at (-1.5,0);
                
                
                \coordinate (E) at (-0.5,0.5);
                \coordinate (F) at (0.5,0.5);
                
                
                \coordinate (G) at (1,0);
                \coordinate (H) at (2,0);
                \coordinate (I) at (1,1);
                
            
            \draw (A)--(B)--(C)--(D)--cycle;
            \node[below left=-0.7mm] at (A) {\small{$(0,0)$}};
            \node[below right=-0.7mm] at (B) {\small{$(1,0)$}};
            \node[above right=-0.7mm] at (C) {\small{$(1,1)$}};
            \node[above left=-0.7mm] at (D) {\small{$(0,1)$}};
            
            
            \draw[->] (E)--(F);
            
            
            \draw (G)--(H)--(I)--cycle;
            \node[above=-0.7mm] at (I) {\small{$(0,1)$}};
            \node[below left=-0.7mm] at (G) {\small{$(0,0)$}};
            \node[below right=-0.7mm] at (H) {\small{$(1,0)$}};
            
            
            
                \foreach \t in {0.1,0.2,...,0.9} {
                 \draw[blue] 
                    ($(A)!\t!(B)$) -- ($(D)!\t!(C)$);
                }
                
                
                \foreach \t in {0.1,0.2,...,0.9} {
                 \draw[blue] 
                    ($(G)!\t!(H)$) -- (I);
                }

            
            
                \foreach \t in {0.1,0.2,...,0.9} {
                 \draw[red] 
                    ($(A)!\t!(D)$) -- ($(B)!\t!(C)$);
                }
                
                
                \foreach \t in {0.1,0.2,...,0.9} {
                 \draw[red] 
                    ($(G)!\t!(I)$) -- ($(H)!\t!(I)$);
                }

        \end{tikzpicture}
    \end{figure}
    
    This map can be precisely defined as follows. Declare $(1,1)$ to be mapped to $(0,1)$, and the other vertices of $I^2$ to be mapped to themselves. Then, write every point $(t_1,t_2) \in I^2$ as the following convex combination of the vertices of $I^2$
    
    $$
    (t_1,t_2) = (1-t_1)(1-t_2) \times (0,0) + t_1(1-t_2) \times (1,0) + (1-t_1)t_2 \times (0,1) + t_1t_2 \times (1,1)
    $$
    
    and extend to $I^2$ by linearity on such convex combinations. This way of decomposing an element $(t_1,t_2) \in I^2$ as a convex combination of the vertices of $I^2$ extends to higher dimensional cubes, yielding what we call the $C$-\emph{representation}, which is the subject of section \ref{section: C-representation}. Note that already in dimension $2$, there is not always a unique way of decomposing an element $(t_1,t_2) \in I^2$ as a convex combination of the vertices of $I^2$ (see remark \ref{remark: writing as convex combination not unique for cube}).
    
    Moreover, a homeomorphism between $I^2 / \sim_{\ag, f(a_0), f(a_1)}$ and $Q$ can be obtained in an analogous way, by declaring $(1,1)$ to be mapped to $(1,0)$ instead of $(0,1)$.

    \item Let us denote by $G$ the face of $Q$ defined as the edge connecting the vertices $(1,0)$ and $(0,1)$. The equivalence relation $\sim_{\ag, f(a_1)}$ is generated by the two equivalence relations $\sim_{\ag, f(a_0),f(a_1)}$ and $\sim_{\ag,f(a_1),f(a_3)}$. The quotient map $I^2 \rightarrow I^2 / \sim_{\ag,f(a_1)}$ therefore factors through the two quotient maps $I^2 \rightarrow I^2 / \sim_{\ag,f(a_0),f(a_1)}$ and $I^2 \rightarrow I^2 / \sim_{\ag,f(a_1),f(a_3)}$. This yields two maps $Q \rightarrow I^2 / \sim_{\ag,f(a_1)}$ that induce two homeomorphisms $Q / G \simeq I^2 / \sim_{\ag,f(a_1)}$. Note that these two homeomorphisms are \emph{different}.
     
\end{itemize}

This is actually an instance of the following general statement.

\begin{Proposition}\label{proposition: filtration of cube}
    
    For every $0 \le i < n$ such that $a_i < a_{i+1}$, consider the following $(n-1)$-dimensional polytope in $\R^{n-1}$:
    
    $$
    P \consel = I^{n-1} \cap \{ \sum_{j \le i} x_j \le 1 \} \cap \{ \sum_{j \ge i+1} x_j \le 1 \}.
    $$
    
    There exists a map $\qo \consel : I^{n-1} \rightarrow P \consel$ satisfying the following properties.
    
    \begin{enumerate}[label=(\roman*)]

    \item It is compatible with $\sim_{\ag,f(a_i),f(a_{i+1})}$.
    
    \item The induced map $I^{n-1} / \sim_{\ag,f(a_i),f(a_{i+1})} \rightarrow P \consel$ is a homeomorphism, thus yielding a pullback square
    
    $$
    \xymatrix{
    P \consel \times (f(a_{i+1}),f(a_i)) \ar[d] \ar[r] & C_A(\Delta^{\ag}) \ar[d]^-{f_{C_A}(\Delta^{\ag})} \\
    (f(a_{i+1}),f(a_i)) \ar[r] & \R.
    }
    $$
    
    \item There are faces $F_{f(a_{i+1})}^{\ag,+} \subseteq P \consel \supseteq F_{f(a_i)}^{\ag,-}$ such that, denoting by $\widetilde{P} \consel$ the space obtained from $P \consel \times [f(a_{i+1}),f(a_i)]$ by collapsing each of the subsets $F_{f(a_{i+1})}^{\ag,+} \times \{f(a_{i+1})\}$ and $F_{f(a_i)}^{\ag,-} \times \{f(a_i)\}$ to a point, the above pullback square extends to a pullback square

     $$
    \xymatrix{
    \widetilde{P} \consel \ar[d] \ar[r] & C_A(\Delta^{\ag}) \ar[d]^-{f_{C_A}(\Delta^{\ag})} \\
    [f(a_{i+1}),f(a_i)] \ar[r] & \R.
    }
    $$
    
    Moreover, if $f(a_i) \neq f(a_0)$, then $F_{f(a_i)}^{\ag,-}$ is proper. Likewise, if $f(a_{i+1}) \neq f(a_n)$, then $F_{f(a_{i+1})}^{\ag,+}$ is proper.
    
    \end{enumerate}
    
\end{Proposition}

The proof of this proposition is the purpose of section \ref{section: quotients of cubes as polytopes}.

\begin{Remark}\label{remark: combinatorial equivalence}

    Recall the map $f_A(\Delta^{\ag}) : |\Delta^{\ag}|_A \rightarrow \R$ (definition \ref{definition: filtration on usual realization}). Given a real number $f(a_{i+1}) < s < f(a_i)$, the fiber of $f_A(\Delta^{\ag})$ over $s$ is a polytope. This polytope is closely related to the polytope $P_{f(a_i),f(a_{i+1})}^{\ag}$ from proposition \ref{proposition: filtration of cube}. One can prove that these two polytopes are \emph{combinatorially equivalent}, i.e., that there exists a bijection between their sets of faces which is compatible with the inclusion relation. This combinatorial equivalence is described as follows. Regard $I^{n-1}$ as the subset of those points of $I^{n+1}$ whose $0^{\text{th}}$ and $n^{\text{th}}$ coordinates are equal to $1$. Given $0 \le j \le i$ and $i+1 \le k \le n$, we map the vertex of $P_{f(a_i),f(a_{i+1})}^{\ag}$ whose $j^{\text{th}}$ and $k^{\text{th}}$ coordinates are equal to $1$ (and the others are equal to $0$) to the intersection point between $f_A(\Delta^{\ag})^{-1}(s)$ and the segment connecting the vertices $e_j$ and $e_k$. See \cite[Section 2.2]{ZieglerLectures} for a discussion of the notion of combinatorial equivalence of polytopes.

\end{Remark}

Item (ii) of proposition \ref{proposition: filtration of cube}, combined with the fact that a $(n-1)$-dimensional polytope is homeomorphic to a $(n-1)$-dimensional disc, implies that the map $f_{C_A}(\Delta^{\ag})$ is a $(n-1)$-dimensional disc bundle in restriction to the complement of the set of critical values $\{f(a_n), f(a_{n-1}), \hdots, f(a_0)\} \subset [f(a_n),f(a_0)]$. In view of item (iii) of proposition \ref{proposition: filtration of cube}, in order to complete the proof, we need to answer the following question: given a polytope $P$ and a proper face $F \subset P$, how can we describe the quotient space $P \times [0,1] / (F \times \{1\})$? This will be done thanks to the following improvement of lemma \ref{lemma: proper face collapsing lemma}.

\begin{Lemmas}\label{lemma: improved proper face collapsing lemma}
    Let $P$ be a polytope, $F \subset P$ a proper face and $U$ a neighborhood of $F$ in $P$. There exists a map $h : P \times [0,1] \rightarrow P \times [0,1]$ satisfying the following properties.
    
    \begin{itemize}
        \item The map $h$ commutes with the projection to $[0,1]$. In other words, there is a $[0,1]$-family of maps $h_t : P \rightarrow P$ such that $h(x,t) = (h_t(x),t)$.
        
        \item The map $h$ is the identity outside of $U \times [0,1]$.
        
        \item The map $h_1 : P \rightarrow P$ is constant on $F$ and the resulting map $P \times [0,1] / (F \times \{1\}) \rightarrow P \times [0,1]$ is a homeomorphism.
        
        \item This homeomorphism restricts to a homeomorphism $\partial P \times [0,1] / (F \times \{1\}) \simeq \partial P \times [0,1]$.
    \end{itemize}
    
\end{Lemmas}

The proof of this lemma is the purpose of section \ref{section: collapsing lemma}.

\begin{Remark}\label{remark: improved collapsing lemma implies first one}
    Lemma \ref{lemma: improved proper face collapsing lemma} implies lemma \ref{lemma: proper face collapsing lemma}: indeed, given a map $h$ as in lemma \ref{lemma: improved proper face collapsing lemma}, the map $c=h_1$ fulfills the conditions of lemma \ref{lemma: proper face collapsing lemma}.
\end{Remark}

We deduce that the images by $f$ of the element of the sequence $\ag$ form a finite sequence of real numbers $f(a_n) =\lambda_k < \lambda_{k-1} < \hdots < \lambda_1 < \lambda_0 = f(a_0)$ such that $f_{C_A}(\Delta^{\ag})$ is a $(n-1)$-dimensional disc bundle over $(\lambda_k,\lambda_{k-1}]$, over $[\lambda_1,\lambda_0)$ as well as over $[\lambda_{j+1},\lambda_j]$ for every $1 \le j \le k-2$. Now, we have the following lemma.

\begin{Lemmas}\label{lemma: gluing trivializations over interval}
    Let $k \ge 2$ be an integer, consider two real numbers $\lambda_0, \lambda_k$ such that $\lambda_k < \lambda_0$, and suppose we are given a topological space $E$ and a map $p : E \rightarrow (\lambda_k, \lambda_0)$. Suppose further that there exists a sequence of real numbers $\lambda_k < \lambda_{k-1} < \hdots < \lambda_1 < \lambda_0$ such that $p$ is a $(n-1)$-dimensional disc bundle over $(\lambda_k,\lambda_{k-1}]$, over $[\lambda_1,\lambda_0)$ as well as over $[\lambda_{j+1},\lambda_j]$ for every $1 \le j \le k-2$. Then, $p$ is a $(n-1)$-dimensional disc bundle over $(\lambda_k,\lambda_0)$.
\end{Lemmas}

In section \ref{section: triviality theorem for filtration on cube}, we prove this lemma and collect all the ingredients to conclude the proof of theorem \ref{theorem: filtration on cube is trivial}.

\subsubsection{A collapsing lemma for polytopes}\label{section: collapsing lemma}

The goal of this section is to prove lemma \ref{lemma: improved proper face collapsing lemma}, and use this lemma to complete the proof of theorem \ref{theorem: stratified simplices and cubes as balls} in the case $G=|-|_A'$.

To begin with, the proof of lemma \ref{lemma: improved proper face collapsing lemma} will be done in three steps.

\emph{\underline{Step 1:}} Denoting by $k$ the dimension of $P$, we will prove that every neighborhood of $F$ in $\partial P$ contains a neighborhood of $F$ in $\partial P$ that is homeomorphic to $\R^{k-1}$ through a homeomorphism that identifies $F$ with a polytope in $\R^{k-1}$.

\emph{\underline{Step 2:}} We will prove an analogue of lemma \ref{lemma: improved proper face collapsing lemma} with $(P,F)$ replaced by $(\R^{k-1},Q)$, where $Q \subset \R^{k-1}$ is a polytope. Combining this and step 1, we will deduce an analogue of lemma \ref{lemma: improved proper face collapsing lemma} with $(P,F)$ replaced by $(\partial P,F)$.

\emph{\underline{Step 3:}} We will deduce the desired result from step 2.

\emph{\underline{Step 1:}} The goal of this step is to prove the following lemma

\begin{Lemmas}\label{lemma: neighborhood of face in boundary}
    Let $P \subset \R^d$ be a $k$-dimensional polytope, $F \subset P$ a proper face and $U$ a neighborhood of $F$ in $\partial P$. There exists a neighborhood $V$ of $P$ in $U$ and a homeomorphism $V \simeq \R^{k-1}$ that identifies $F$ with a polytope in $\R^{k-1}$.
\end{Lemmas}

\begin{proof}

    Up to restricting ourselves to the affine hull of $P$ in $\R^d$, we can assume that $k=d$, which we do from now on. Let $\langle c,x \rangle \le M$ be a valid inequality for $P$ such that $F = \{p \in P \mid \langle c,p \rangle = M \}$ (see section \ref{section: polytopes} for the definition of a valid inequality). We denote by $\lambda_c$ the linear form on $\R^d$ defined by $x \mapsto \langle c,x \rangle$. Given a subset $E \subseteq \R^d$ and a real number $s$, we denote by $E_{s}$ the set of those $e \in E$ such that $\lambda_c(e)=s$, so that in particular, we have $F = P_M$. Similarly, given a real number $s'$, we also use the notation $E_{> s}$, $E_{> s}^{\le s'}$.
    
    Since $F$ is a proper face of $P$, we have $P \neq P_M$. Consequently, the subset $\lambda_c(P) \subset \R$ is a closed bounded interval with nonempty interior, whose maximum is $M$. Moreover the set $\lambda_c(^c U)$ is a closed subset of $\lambda_c(P)$ which does not contain $M$. In particular, there exists a point $s \in \lambda_c(P)$ such that $s \notin \lambda_c(^c U)$ and $s < M$. Denote $V = (\partial P)_{>s}$; this is an open neighborhood of $F$ in $\partial P$ which is contained in $U$.
    
    We proceed by constructing a homeomorphism between $V$ and $\R^{d-1}$ that satisfies the desired condition. To this end, we start by justifying that there exists a point $p \in (\Int(P))_s$. For this, denote $\lambda_c(P)$ by $[m,M]$. Since $s \in (m,M)$, it suffices to show that $\lambda_c(\Int(P))=(m,M)$. Now, we observe two things, which follow from the description of $P$ as a finite intersection of affine closed halfspaces of $\R^d$.
    
    \begin{itemize}
        \item The interior of $P$ is connected. Since $\lambda_c$ is an open map, we deduce that $\lambda_c(\Int(P))$ is an open interval.
        
        \item $P$ is the closure of $\Int (P)$. Therefore, $\lambda_c(P)$ is the closure of $\lambda_c(\Int(P))$. This forces $\lambda_c(\Int(P))=(m,M)$. 
    \end{itemize}
    
    Fix $p \in (\Int(P))_s$ from now on and consider the subspace $(\R^d)_{>s}^{\le M} \subset \R^d$. There is a map $(\R^d)_{>s}^{\le M} \rightarrow (\R^d)_M$ that carries a point $x$ to the intersection of the half-line $p + \R_+(x-p)$ with $(\R^d)_M$. We denote by $\varphi$ its restriction to $V$; we regard $\varphi$ as a map $V \rightarrow \R^{d-1}$ by using an affine identification between $(\R^d)_M$ and $\R^{d-1}$. We will prove that $\varphi$ is a homeomorphism satisfying the desired property.
    
    First note that since $F \subset (\R^d)_M$, the map $\varphi$ is the identity on $F$, and therefore the image of $F$ by $\varphi$ is a polytope in $\R^{d-1}$. Second, $\varphi$ is bijective; indeed, for every point $x \in (\R^d)_M$, the half-line $p + \R_+(x-p)$ intersects $\partial P$ at a unique point by lemma \ref{lemma: intersection of halfline with polytope boundary}, and this point has to belong to $V$.
    
    To conclude that $\varphi$ is a homeomorphism, we observe that $\varphi$ extends to a bijective map between the one-point compactifications of $V$ and $\R^{d-1}$. Since these one-point compactifications are Hausdorff, this extension is a homeomorphism, and therefore $\varphi$ as well. \qedhere

\end{proof}

\emph{\underline{Step 2:}} The goal of this step is to prove two analogues of lemma \ref{lemma: improved proper face collapsing lemma}.

\begin{Lemmas}\label{lemma: Euclidean polytope collapsing lemma}
 Let $d \ge 1$ be an integer, $Q \subset \R^d$ a polytope and $U$ a neighborhood of $Q$ in $\R^d$. There exists a map $h : \R^d \times [0,1] \rightarrow \R^d \times [0,1]$ satisfying the following properties.
    
    \begin{itemize}
        \item The map $h$ commutes with the projection to $[0,1]$. In other words, there is a $[0,1]$-family of maps $h_t : \R^d \rightarrow \R^d$ such that $h(x,t)=(h_t(x),t)$.
        
        \item The map $h$ is the identity outside of $U \times [0,1]$.
        
        \item The map $h_1 : \R^d \rightarrow \R^d$ is constant on $Q$ and the resulting map $\R^d \times [0,1] / (Q \times \{1\}) \rightarrow \R^d \times [0,1]$ is a homeomorphism.
    \end{itemize}
    
\end{Lemmas}

\begin{proof}

    We will prove the statement by induction on the codimension of $Q$ in $\R^d$. We start by considering the case of codimension $0$, that is, the case when the dimension of $Q$ is $d$. In that case, choose a point $q$ belonging to the interior of $Q$ and recall the map $\Phi_q^Q : \partial Q \times \R_+ \rightarrow \R^d$ from proposition \ref{proposition: a polytope is a disc}. Fix a positive real number $\varepsilon$ such that $\Phi_q^Q(\partial Q \times [0,1+\varepsilon]) \subset U$. Using proposition \ref{proposition: a polytope is a disc}, we see that it suffices to prove the following statement: there exists a map $h : (\partial Q \times \R_+) \times [0,1] \rightarrow (\partial Q \times \R_+) \times [0,1]$ satisfying the following properties:
    
    \begin{itemize}
        \item The map $h$ commutes with the projection to $[0,1]$.
        
        \item The map $h$ is the identity outside of $\partial Q \times [0,1+\varepsilon] \times [0,1]$.
        
        \item The restriction of the map $h_1$ to $\partial Q \times [0,1]$ is the projection to $\partial Q \times \{0\}$ and, denoting by $C^{\circ}(\partial Q)$ the space $(\partial Q \times \R_+)/(\partial Q \times \{0\})$ and by $C(\partial Q)$ the subspace $(\partial Q \times [0,1])/(\partial Q \times \{0\})$ of $C^{\circ}(\partial Q)$, the resulting map $h : C^{\circ}(\partial Q) \times [0,1] / (C(\partial Q) \times \{1\}) \rightarrow C^{\circ}(\partial Q) \times [0,1]$ is a homeomorphism.
    \end{itemize}
    
    The following map satisfies the desired properties.
    
    $$
    \fonction{h}{\partial Q \times \R_+ \times [0,1]}{\partial Q \times \R_+ \times [0,1]}{(q,s,t)}{
    \left\{
    \begin{array}{ll}
        (q,s,t) & \text{if } s \ge 1+\varepsilon  \\
        (q,(1-t)s,t) & \text{if } 0 \le s \le 1 \\
        (q,\frac{1}{\varepsilon}((1+\varepsilon-s)(1-t)+(s-1)(1+\varepsilon)),t) & \text{if } 1 \le s \le 1+ \varepsilon.
    \end{array}
    \right.
    }
    $$
    
    The initial step of the induction is thus complete. Now suppose that $Q$ has codimension $k \ge 1$ in $\R^d$, and that the result holds for polytopes of codimension at most $k-1$. Recall that the dimension of $Q$ is defined to be the dimension of the affine hull of $Q$ in $\R^d$. We deduce that $Q$ is contained in an affine hyperplane of $\R^d$. There exists an affine isomorphism $\R^{d-1} \times \R \rightarrow \R^d$ that restricts to an affine isomorphism between $\R^{d-1} \times \{0\}$ and an affine hyperplane containing $Q$; we are thus reduced to the case when $Q$ is a polytope contained in $\R^{d-1} \times \{0\} \subset \R^{d-1} \times \R$. In that case, there exists a neighborhood $U'$ of $Q$ in $\R^{d-1} \times \{0\}$ as well as a real number $\varepsilon > 0$ such that $U' \times [-\varepsilon,\varepsilon] \subset U$. Up to rescaling, we may assume that $\varepsilon = 1$, which we do from now on, in order to simplify the notation. By the induction hypothesis, there exists a map $\tilde{h} : \R^{d-1} \times [0,1] \rightarrow \R^{d-1} \times [0,1]$ satisfying the condition of the lemma with respect to the polytope $Q \subset \R^{d-1}$ and the neighborhood $U'$. The following map satisfies the conditions of the lemma with respect to the polytope $Q \subset \R^{d-1} \times \R$ and the neighborhood $U' \times [-1,1]$ (and, therefore, with respect to the neighborhood $U$ as well):
    
     $$
    \fonction{h}{\R^{d-1} \times \R \times [0,1]}{\R^{d-1} \times \R \times [0,1]}{(x,s,t)}{
    \left\{
    \begin{array}{ll}
        (\tilde{h}(x,t-|s|),s,t) & \text{if } |s| \le t  \\
        (\tilde{h}(x,0),s,t) & \text{if } |s| \ge t.
    \end{array}
    \right.
    }
    $$ \qedhere

\end{proof}

\begin{Corollary}\label{corollary: collapsing lemma for face in boundary}

    Let $P$ be a polytope, $F \subset P$ a proper face and $U$ a neighborhood of $F$ in $\partial P$. There exists a map $h : \partial P \times [0,1] \rightarrow \partial P \times [0,1]$ satisfying the following properties.
    
    \begin{itemize}
        \item The map $h$ commutes with the projection to $[0,1]$. In other words, there is a $[0,1]$-family of maps $h_t : \partial P \rightarrow \partial P$ such that $h(x,t)=(h_t(x),t)$.
        
        \item The map $h$ is the identity outside of $U \times [0,1]$.
        
        \item The map $h_1 : \partial P \rightarrow \partial P$ is constant on $F$ and the resulting map $\partial P \times [0,1] / (F \times \{1\}) \rightarrow \partial P \times [0,1]$ is a homeomorphism.
    \end{itemize}

\end{Corollary}

\begin{proof}

    Denote by $k$ the dimension of $P$. By lemma \ref{lemma: neighborhood of face in boundary}, there exists a neighborhood $V$ of $P$ in $U$ and a homeomorphism $V \simeq \R^{k-1}$ that identifies $F$ with a polytope in $\R^{k-1}$. Let $K$ be a compact neighborhood of $F$ in $V$. By lemma \ref{lemma: Euclidean polytope collapsing lemma}, there exists a map $\tilde{h} : V \times [0,1] \rightarrow V \times [0,1]$ that commutes with the projection to $[0,1]$, is the identity outside of $K \times [0,1]$, is constant on $F \times \{1\}$, and such that the resulting map $(V \times [0,1])/(F \times \{1\}) \rightarrow V \times [0,1]$ is a homeomorphism. Let us define
    
    $$
    \fonction{h}{\partial P \times [0,1]}{\partial P \times [0,1]}{(q,t)}{
    \left\{
    \begin{array}{ll}
        \tilde{h}(q,t) & \text{if } q \in V  \\
        (q,t) & \text{if } q \notin V.
    \end{array}
    \right.
    }
    $$
    
    This is a continuous map; indeed, for every open subset $W \subset \partial P \times [0,1]$ we have
    
    $$
    h^{-1}(W)=(W \cap (^c K \times [0,1])) \cup \tilde{h}^{-1}(W \cap (V \times [0,1])).
    $$
    
    This map satisfies the first two the desired properties by construction. Also by construction, the map $h_1$ is constant on $F$, and the resulting map $\partial P \times [0,1] / (F \times \{1\}) \rightarrow \partial P \times [0,1]$ is a bijection. Since its source is compact and its target is Hausdorff, this is a homeomorphism. This completes the proof. \qedhere

\end{proof} 

\emph{\underline{Step 3:}} We are now ready to collect the ingredients and prove lemma \ref{lemma: improved proper face collapsing lemma}.

\begin{proof}[Proof of lemma \ref{lemma: improved proper face collapsing lemma}]

    We denote by $C(\partial P)$ the closed cone of $\partial P$, that is, the space $(\partial P \times [0,1])/(\partial P \times \{0\})$. Identifying $P$ with $C(\partial P)$ using proposition \ref{proposition: a polytope is a disc}, we are reduced to proving that for every neighborhood $U$ of $F \times \{1\}$ in $\partial P \times [0,1]$, there exists a map $h : (\partial P \times [0,1])\times [0,1] \rightarrow (\partial P \times [0,1]) \times [0,1]$ satisfying the following properties.
    
    \begin{itemize}
        \item The map $h$ commutes with the second projection to $[0,1]$.
        
        \item The map $h$ is the identity outside of $U \times [0,1]$ and $h_t((\partial P \times \{0\}) \times [0,1]) \subseteq (\partial P \times \{0\}) \times [0,1]$.
        
        \item The map $h_1$ is constant on $F \times \{1\} \times \{1\}$ and the resulting map $(C(\partial P) \times [0,1])/(F \times \{1\} \times \{1\}) \rightarrow C(\partial P) \times [0,1]$ is a homeomorphism.
        
        \item This homeomorphism restricts to a homeomorphism $((\partial P \times \{1\}) \times [0,1]) / (F \times \{1\} \times \{1\}) \simeq (\partial P \times \{1\}) \times [0,1]$.
    \end{itemize}
    
    Fix a neighborhood $V$ of $F$ in $\partial P$ and a real number $0 < \varepsilon < 1$ such that $V \times (\varepsilon,1] \subset U$. Let $\tilde{h} : \partial P \times [0,1] \rightarrow \partial P \times [0,1]$ be a map satisfying the conditions of corollary \ref{corollary: collapsing lemma for face in boundary} with respect to the face $F$ and the neighborhood $V$. The following map satisfies the desired properties.
    
    $$
    \fonction{h}{(\partial P \times [0,1]) \times [0,1]}{(\partial P \times [0,1]) \times [0,1]}{(p,s,t)}{
    \left\{
    \begin{array}{ll}
        (\tilde{h}(p,0),s,t) & \text{if } s \le \varepsilon  \\
        (\tilde{h}(p,\frac{s-\varepsilon}{1-\varepsilon}t),s,t) & \text{if }  s \ge \varepsilon.
    \end{array}
    \right.
    }
    $$ \qedhere

\end{proof}

For later use, we record the following corollary.

\begin{Corollary}\label{corollary: proper face collapsing homeomorphic to disc}

    Let $P$ be a $d$-dimensional polytope and $F \subset P$ a proper face. The pair $(P/F, \partial P / F)$ is homeomorphic to $(D^d,\SM^{d-1})$.

\end{Corollary}

\begin{proof}

    This follows from lemma \ref{lemma: proper face collapsing lemma} together with proposition \ref{proposition: a polytope is a disc}. \qedhere

\end{proof}

We conclude this section by proving theorem \ref{theorem: stratified simplices and cubes as balls} in the case $G=|-|_A'$.

\begin{proof}[Proof of theorem \ref{theorem: stratified simplices and cubes as balls} in the case $G=|-|_A'$]

    Recall that there are proper disjoint faces $F_1,\hdots,F_k \subset |\Delta^n|$ such that the space $|\Delta^{\ag}|_A'$ is obtained from $|\Delta^n|$ by collapsing each of these to a point. Take a family of disjoint open subsets $V_1,\hdots,V_k \subseteq |\Delta^n|$ such that for every $1 \le i \le k$, $V_i$ contains $F_i$. For every $1 \le i \le k$, take a compact neighborhood $L_i$ of $F_i$ in $V_i$. By lemma \ref{lemma: proper face collapsing lemma}, for every $1 \le i \le k$, take a map $c_i : |\Delta^n| \rightarrow |\Delta^n|$ such that $c_i$ is the identity outside of $L_i$, $c_i$ is constant on $F_i$ and the resulting map $|\Delta^n| / F_i \rightarrow |\Delta^n|$ is a homeomorphism that restricts to a homeomorphism $\partial |\Delta^n| / F_i \simeq \partial |\Delta^n|$. Then define a map
    
    $$
    \fonction{c}{|\Delta^n|}{|\Delta^n|}{x}{
    \left\{
    \begin{array}{ll}
        c_i(x) & \text{if } x \in L_i \text{ for some } 1 \le i \le k, \\
        x & \text{otherwise.} 
    \end{array}
    \right.
    }
    $$
    
    This is a continuous map; indeed, for every open subset $U \subset |\Delta^n|$ one has
    
    $$
    c^{-1}(U) = (U \cap \, ^c(L_1 \cup \hdots \cup L_k)) \cup c_1^{-1}(U \cap V_1) \cup \hdots \cup c_k^{-1}(U \cap V_k).
    $$
    
    By construction, the map $c$ is constant on $F_1 \cup \hdots \cup F_k$ and the resulting map $|\Delta^{\ag}|_A' \rightarrow |\Delta^n|$ is a bijection. Since its source is compact and its target is Hausdorff, this is a homeomorphism. The latter restricts to a homeomorphism between the image of $\partial |\Delta^n|$ in $|\Delta^{\ag}|_A'$, and $\partial |\Delta^n|$. We denote by $\partial |\Delta^{\ag}|_A'$ the image of $\partial |\Delta^n|$ in $|\Delta^{\ag}|_A'$. Since the pair $(|\Delta^n|, \partial |\Delta^n|)$ is homeomorphic to the pair $(D^n,\SM^{n-1})$ (proposition \ref{proposition: a polytope is a disc}), we deduce that the pair $(|\Delta^{\ag}|_A', \partial |\Delta^{\ag}|_A')$ is homeomorphic to $(D^n,\SM^{n-1})$.
    
    To complete the proof, it remains to show that the natural map $|\partial \Delta^{\ag}|_A' \rightarrow |\Delta^{\ag}|_A'$ induces a homeomorphism between $|\partial \Delta^{\ag}|_A'$ and $\partial |\Delta^{\ag}|_A'$. To this end, we will use the description of $|\partial \Delta^{\ag}|_A'$ as a colimit, provided by proposition \ref{proposition: boundary of stratified simplex as colimit of faces}. We know from that result that the image of $|\partial \Delta^{\ag}|_A'$ in $|\Delta^{\ag}|_A'$ is contained in $\partial |\Delta^{\ag}|_A'$. Now, the map $|\partial \Delta^{\ag}|_A' \rightarrow \partial |\Delta^{\ag}|_A'$ admits an inverse that one can define as follows. For every integer $0 \le i \le n$, we have the quotient map $|\Delta^{n-1}| \rightarrow |\Delta^{\ag \circ \delta_n^i}|_A'$. These maps assemble into a map $|\partial \Delta^n| \rightarrow |\partial \Delta^{\ag}|_A'$ which factors through the quotient map $\partial |\Delta^n| = |\partial \Delta^n| \rightarrow \partial |\Delta^{\ag}|_A'$. This factorization is the desired inverse. \qedhere

\end{proof}

\subsubsection{\texorpdfstring{The cubical representation of a Euclidean vector}{The C-representation of a Euclidean vector}}\label{section: C-representation}

The goal of this section and the next one is to prove proposition \ref{proposition: filtration of cube}. In this section, we introduce the notion of cubical representation of a point $(t_1,\hdots,t_N) \in \R^N$, called $C$-representation, which will be involved in the definition of the maps $\qo \consel$ from proposition \ref{proposition: filtration of cube}.

Let us fix an integer $N \ge 1$.

\begin{Notation}\label{notation: set of vertices of cube}
We denote by $V_{N}$ the set of vertices of the polytope $I^{N} \subset \R^N$. That is, $V_{N} = \{0,1\}^{N}$.
\end{Notation}
    
Every point $\tu=(t_1,\hdots,t_{N}) \in \R^{N}$ can be written as a linear combination of the elements of $V_{N}$ as follows:
    
    $$
    \tu = \sum_{v \in V_{N}} ( \prod_{j \mid v_j =0} (1-t_j) \prod_{j \mid v_j = 1} t_j  ) v.
    $$

\begin{Definition}\label{definition: C-representation}    
We call this expression the $C$-\emph{representation} of $\tu$.
\end{Definition}

\begin{Remark}\label{remark: C representation as convex combination}
The sum of the coefficients of the $C$-representation of $\tu$ is equal to $1$, that is, we have
    
    $$
    \sum_{v \in V_{N}} ( \prod_{j \mid v_j =0} (1-t_j) \prod_{j \mid v_j = 1} t_j ) = 1.
    $$
    
In particular, for every $\tu \in I^{N}$, the $C$-representation of $\tu$ is a way to write $\tu$ as a convex combination of the vertices of $I^{N}$.
\end{Remark}
    
\begin{Remark}\label{remark: writing as convex combination not unique for cube}
    The polytope $I^N$ is the convex hull of the set $\{0,1\}^N$, that is, every map $\Delta^{2^N-1} \rightarrow I^N$ that is bijective on vertices and linear, is surjective. The $C$-representation can be regarded as a section of such a map. However, as long as $N \ge 2$, such a section is \emph{not} unique. For instance, a different section can be obtained by choosing a triangulation of $I^N$ whose set of vertices is $\{0,1\}^N$, and writing every element of $I^N$ as a convex combination of the vertices of a simplex of the triangulation in which it lies.
\end{Remark}
    
\begin{Notation}\label{notation: vertices for which coefficient is nonzero}
For every $\tu \in \R^{N}$, we denote by $V_{N}(\tu)$ the subset of $V_{N}$ formed by those vertices for which the corresponding coefficient in the $C$-representation of $\tu$ is nonzero. In other words, for every $v \in V_N$, we have $v \in V_N(\tu)$ except if $v_j=0$ and $t_j=1$ or if $v_j=1$ and $t_j=0$.
\end{Notation}

\begin{Construction}\label{construction: extension by linearity of C-representations}
    
For every $d \ge 0$ and every map $g : V_{N} \rightarrow \R^d$, we denote by $\overline{g}$ the map $\R^{N} \rightarrow \R^d$ defined to be $g$ on $V_{N}$, and to be "linear on $C$-representations" in the sense that for every $\tu \in \R^{N}$, $\overline{g}(\tu)$ is defined as
    
    $$
    \overline{g} (\tu) = \sum_{v \in V_{N}} ( \prod_{j \mid v_j =0} (1-t_j) \prod_{j \mid v_j = 1} t_j  ) g(v).
    $$
    
    Denoting by $Q$ the convex hull of $g(V_{N})$, the map $\overline{g}$ restricts to a map $\overline{g} : I^N \rightarrow Q$.

\end{Construction}

\begin{Lemmas}\label{lemma: functoriality of extension to C representations}
    The assignment $g \mapsto \overline{g}$ satisfies the following functoriality properties.
    
    \begin{enumerate}[label=(\roman*)]
        \item If $\lambda : \R^d \rightarrow \R^m$ is linear, then $\overline{\lambda \circ g} = \lambda \circ \overline{g}$.
    \end{enumerate}
    
    Fix a subset $E \subseteq \{1,...,N\}$, denote by $^c\!E$ the complement of $E$ in $\{1,\hdots,N\}$, and by $V_E$ (resp. $V_{\,^c\!E}$) the set of vertices of $I^E$ (resp. $I^{\,^c\!E}$). Denote by $\pi_E$ and $\pi_{\,^c\!E}$ the projection maps $I^N \rightarrow I^E$ and $I^N \rightarrow I^{\,^c\!E}$. Then for every positive integer $d$ and every map $g : V_E \rightarrow \R^d$ we have: 
    
    \begin{enumerate}[label=(\roman*), start=2]
        
        \item $\overline{g \circ (\pi_E)_{\mid V_N}} = \overline{g} \circ \pi_E$.
        
    \end{enumerate}
    
    Finally, for every $u \in V_{\,^c\!E}$, denote by $\kappa_u$ the map $I^E \rightarrow I^N$ defined by the condition that $\pi_E \circ \kappa_u = \id$ and $\pi_{\,^c\!E} \circ \kappa_u \equiv u$. Then for every positive integer $d$ and every map $g : V_N \rightarrow \R^d$ we have:
    
    \begin{enumerate}[label=(\roman*), start=3]
        \item $\overline{g \circ (\kappa_u)_{\mid V_E}} = \overline{g} \circ \kappa_u$.
    \end{enumerate}
\end{Lemmas}

\begin{proof}

    Item (i) is straightforward. As for item (ii), we have for every $\tu \in I^N$
    
    $$
    \begin{aligned}
    \overline{g \circ (\pi_E)_{\mid V_N}} (\tu) & = \sum_{v \in V_{N}} ( \prod_{j \mid v_j = 0} (1-t_j) \prod_{j \mid v_j=1} t_j ) g \circ \pi_E(v)\\
    & = \sum_{w \in V_E} \sum_{v \in V_{N} \mid \pi(v) = w} ( \prod_{j \mid v_j = 0} (1-t_j) \prod_{j \mid v_j=1} t_j ) g(w) \\
    & = \sum_{w \in V_E} ( \prod_{j \mid w_j = 0} (1-\pi_E(\tu)_j) \prod_{j \mid w_j=1} \pi_E(\tu)_j ) (\sum_{u \in V_{\,^c\!E}} ( \prod_{j \mid u_j = 0} (1-(\pi_{\,^c\!E}(\tu))_j) \prod_{j \mid u_j=1} (\pi_{\,^c\!E}(\tu))_j )) g(w) \\
    & = \sum_{w \in V_E} ( \prod_{j \mid w_j = 0} (1-\pi_E(\tu)_j) \prod_{j \mid w_j=1} \pi_E(\tu)_j ) g(w) \\
    & \,\,\,\,\, \text{ since } \,\,\, \sum_{u \in V_{\,^c\!E}} ( \prod_{j \mid u_j = 0} (1-(\pi_{\,^c\!E}(\tu))_j) \prod_{j \mid u_j=1} (\pi_{\,^c\!E}(\tu))_j ) = 1 \\
    & = \overline{g}(\pi_E(\tu)).
    \end{aligned}
    $$

    As for item (iii), we have for every $\tu \in I^E$:
    
    $$
    \overline{g} \circ \kappa_u (\tu) = \sum_{v \in V_N}(\prod_{i \mid v_i = 0}(1-\kappa_u(\tu)_i)\prod_{i \mid v_i = 1} \kappa_u(\tu)_i)g(v).
    $$
    
    But note that for every $v \in V_N$, either $\pi_{\,^c\!E}(v) = u$ and then
    
    $$
    \prod_{i \in ^c\!E \mid v_i = 0}(1-\kappa_u(\tu)_i)\prod_{i \in ^c\!E \mid v_i = 1} \kappa_u(\tu)_i = 1,
    $$
    
    or $\pi_{^c\!E}(v) \neq u$ and then
    
    $$
    \prod_{i \in ^c\!E \mid v_i = 0}(1-\kappa_u(\tu)_i)\prod_{i \in ^c\!E \mid v_i = 1} \kappa_u(\tu)_i = 0.
    $$
    
    Therefore
    
    $$
    \begin{aligned}
    \overline{g} \circ \kappa_u (\tu) & = \sum_{\substack{v \in V_N \\ \pi_{^c\!E}(v) = u}} (\prod_{i \in E \mid v_i=0}(1- \kappa_u(\tu)_i) \prod_{i \in E \mid v_i=1} \kappa_u(\tu)_i)g(v) \\
    & = \sum_{w \in V_E} (\prod_{i \mid w_i=0} (1-t_i) \prod_{i \mid w_i=1} t_i)g(\kappa_u(w)) \\
    & = \overline{g \circ (\kappa_u)_{\mid V_E}} (\tu).
    \end{aligned}
    $$
    \qedhere

\end{proof}

\subsubsection{Quotients of cubes as polytopes}\label{section: quotients of cubes as polytopes}

The goal of this section is to proceed with the proof of proposition \ref{proposition: filtration of cube} and complete it. We start with the case $n=1$.

\begin{Lemmas}\label{lemma: proposition in case n=1}
    Proposition \ref{proposition: filtration of cube} holds in the case $n=1$.
\end{Lemmas}

\begin{proof}

    In this case $\ag$ is of the form $[a_0 < a_1]$ and $I^{n-1}$, $P_{a_0,a_1}^{\ag}$ are points. The items (i) and (ii) are therefore straightforward, and the faces $F_{a_1}^{\ag,+} = P_{a_0,a_1}^{\ag} = F_{a_0}^{\ag,-}$ fulfill the conditions of item (iii). \qedhere

\end{proof}

Throughout the rest of this section, we fix a nonconstant sequence $\ag = [a_0 \le \hdots \le a_n]$ of elements of $A$ of length $n \ge 2$. Recall that for every $f(a_n) \le s \le f(a_0)$ we denote by $\sim_{\ag,s}$ the restriction of the equivalence relation $\sim_{\ag}$ to $I^{n-1} \simeq I^{n-1} \times \{s\} \subset I^{n-1} \times [f(a_n),f(a_0)]$.

Our first goal is to prove items (i) and (ii) in proposition \ref{proposition: filtration of cube}. We begin with a lemma describing explicitly the equivalence relations $\sim_{\ag,s}$ when $s \notin \{f(a_n),f(a_{n-1}),\hdots,f(a_0)\}$.

\begin{Lemmas}\label{lemma: restriction of equivalence relation to noncritical level}

    Consider a real number $f(a_n) < s < f(a_0)$ such that $s \notin \{f(a_n),f(a_{n-1}),\hdots,f(a_0)\}$ and denote by $i$ the integer satisfying $f(a_{i+1}) < s < f(a_i)$. Regard $I^{n-1}$ as the set of maps $\tu : \{0,\hdots,n\} \rightarrow I$, $\tu = (t_0,\hdots,t_n)$, such that $t_0=t_n=1$. For every $\tu, \tu' \in I^{n-1}$, we have $\tu \sim_{\ag,s} \tu'$ if and only if there exist $0 \le m \le i$ and $i+1 \le l \le n$ such that $t_m=t_m'=1$, $t_l=t_l'=1$ and $(t_m,\hdots,t_l)=(t_m',\hdots,t_l')$.
    
    Consequently, $\sim_{\ag,s}$ is the smallest equivalence relation on $I^{n-1}$ that identifies the elements that have the same image by at least one map among the projection maps $I^{k-1} \times \{1\} \times I^{n-k-1} \rightarrow I^{n-k-1}$ for $0 < k \le i$ and the projections maps $I^{k-1} \times \{1\} \times I^{n-k-1} \rightarrow I^{k-1}$ for $i+1 \le k < n$.
    
\end{Lemmas}

\begin{proof}

    Consider two elements $\tu$ and $\tu'$ of the morphism space $F(\Delta^n)(0,n)$. By definition, we have $\tu \sim_{\ag,s} \tu'$ if and only if $(\tu,s) \sim_{\ag} (\tu',s)$, if and only if $\tu$ and $\tu'$ are respectively of the form $\nu \circ \gamma \circ \eta$ and $\nu' \circ \gamma \circ \eta'$, for some $\nu, \nu' \in F(\Delta^n)(l,n)$, $\eta, \eta' \in F(\Delta^n)(0,m)$ and $\gamma \in F(\Delta^n)(m,l)$, with $0 \le m \le i$ and $i+1 \le l \le n$. The result therefore follows from proposition \ref{proposition: compositions in F(n)}. \qedhere

\end{proof}

In particular, for every $0 \le i < n$ such that $a_i < a_{i+1}$, the equivalence relations $(\sim_{\ag,s})_{f(a_{i+1}) < s < f(a_i)}$ are indeed all the same, as stated in section \ref{section: summary of the proof of the ball theorem}. Recall that we denote by $\sim_{\ag,f(a_i),f(a_{i+1})}$ this equivalence relation.

We now fix an integer $0 \le i < n$ such that $a_{i+1} < a_i$, and construct the maps $\qo \consel$ from the statement of proposition \ref{proposition: filtration of cube}.

\begin{Construction}\label{construction: maps out of cube}

    We define $\qo \consel : I^{n-1} \rightarrow P \consel$ to be the map obtained, through construction \ref{construction: extension by linearity of C-representations}, from the following map $q \consel : V_{n-1} \rightarrow P \consel$:
    
    \begin{itemize}
        \item If the set $E_1=\{ 1 \le j \le i \mid v_j = 1 \}$ is nonempty, we let $k_1 = \max E_1$ and we declare $q \consel (v)_j = 0$ for every $1 \le j < k_1$, and $q \consel (v)_j=v_j$
        for every $k_1 \le j \le i$. Otherwise, we declare $q \consel (v)_j=v_j=0$ for every $1 \le j \le i$.
        
        \item If the set $E_2=\{ i+1 \le j \le n-1 \mid v_j = 1 \}$ is nonempty, we let $k_2 = \min E_2$ and we declare $q \consel (v)_j = 0$ for every $k_2 < j \le n-1$, and $q \consel (v)_j=v_j$
        for every $i+1 \le j \le k_2$. Otherwise, we declare $q \consel (v)_j=v_j=0$ for every $i+1 \le j \le n-1$.
        
    \end{itemize}

\end{Construction}

Note that, by construction, for every $v \in V_{}$, at most one of the coefficients $(q \consel (v)_j)_{1 \le j \le i}$ is equal to $1$ and the others are equal to $0$, and likewise for the coefficients $(q \consel (v)_j)_{i+1 \le j \le n-1}$. In particular, we have indeed $q \consel (v) \in P \consel$.

Figures \ref{figure: collapsed cube as polytope 1} and \ref{figure: collapsed cube as polytope 2} below illustrate the equivalence relation $\sim_{\ag,f(a_i),f(a_{i+1})}$, the polytope $P \consel$ and the map $\qo \consel$ in the case when $\ag$ is of the form $a_0 < a_1 < a_2 \le a_3 \le a_4$ and $i=0$ and $i=1$ respectively.

\begin{figure}[H]
    \centering
    \begin{minipage}{0.3\textwidth}
    \centering
    \begin{tikzpicture}[scale=0.3]
    
        
        
            \coordinate (A) at (-0,0);
            \coordinate (B) at (8,0);
            \coordinate (C) at (11,2);
            \coordinate (D) at (3,2);
            \coordinate (E) at (0,8);
            \coordinate (F) at (8,8);
            \coordinate (G) at (11,10);
            \coordinate (H) at (3,10);
            
            
            \node[below] at (A) {\footnotesize{$(0,0,0)$}};
            \node[left] at (E) {\footnotesize{$(0,0,1)$}};
            \node[below] at (B) {\footnotesize{$(1,0,0)$}};
            \node[below] at (D) {\footnotesize{$(0,1,0)$}};

            
            \draw (A)--(B)--(C)--(G)--(H)--(E)--cycle;
            \draw (E)--(F)--(B);
            \draw (F)--(G);
            \draw[dashed] (A)--(D)--(C);
            \draw[dashed] (D)--(H);
            
            
            \foreach \t in {0.1,0.2,...,0.9} {
            \draw[dashed] 
            ($(D)!\t!(C)$) -- ($(H)!\t!(G)$);
            }
            
            \draw[fill, opacity=0.4] (F)--(B)--(C)--(G)--cycle;
             
    \end{tikzpicture}
    \end{minipage}
    \begin{minipage}{0.1\textwidth}
    \centering
    \begin{tikzpicture}[scale=0.3]
    
        
        \draw[->] (-4,0)--(4,0) node[midway, above] {$\qo ^{\ag}_{f(a_0),f(a_1)}$};
             
    \end{tikzpicture}
    \end{minipage}
    \begin{minipage}{0.3\textwidth}
    \centering
    \begin{tikzpicture}[scale=0.3]
    
        
        
            \coordinate (A) at (-0,0);
            \coordinate (B) at (8,0);
            \coordinate (D) at (3,2);
            \coordinate (E) at (0,8);
            
            
            \node[below] at (A) {\footnotesize{$(0,0,0)$}};
            \node[left] at (E) {\footnotesize{$(0,0,1)$}};
            \node[below] at (B) {\footnotesize{$(1,0,0)$}};
            \node[below] at (D) {\footnotesize{$(0,1,0)$}};

            
            \draw (A)--(B)--(E)--cycle;
            \draw[dashed] (A)--(D)--(B);
            \draw[dashed] (D)--(E);
            
    \end{tikzpicture}
    \end{minipage}
    \caption{By lemma \ref{lemma: restriction of equivalence relation to noncritical level}, the equivalence relation $\sim_{\ag,a_0,a_1}$ on $I^3$ identifies the elements that have the same image by the projection $\{1\} \times I^2 \rightarrow *$ (these correspond to the grey part) and by the first projection $I \times \{1\} \times I \rightarrow I$ (these correspond to the vertical dashed segments). The map $\qo ^{\ag}_{f(a_0),f(a_1)}$ from construction \ref{construction: maps out of cube} carries all the vertices of the face $\{1\} \times I^2$ to $(1,0,0)$, the point $(0,1,1)$ to $(0,1,0)$ and the other vertices of $I^3$ to themselves.}
     \label{figure: collapsed cube as polytope 1}
\end{figure}
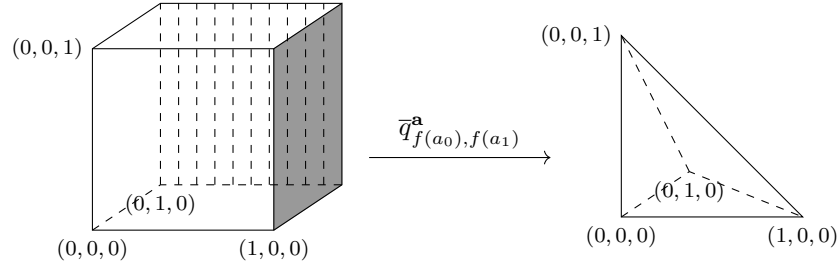

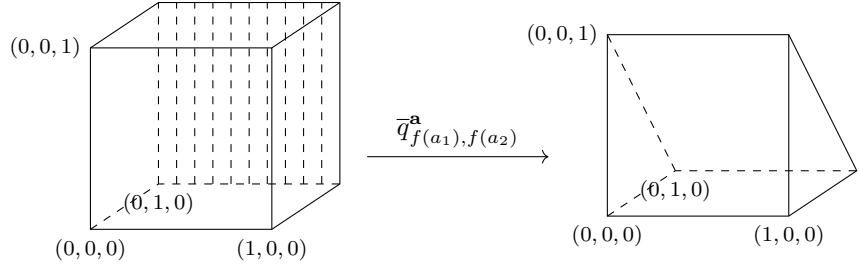
\begin{figure}[H]
    \centering
    \begin{minipage}{0.3\textwidth}
    \centering
    \begin{tikzpicture}[scale=0.3]
    
        
        
            \coordinate (A) at (-0,0);
            \coordinate (B) at (8,0);
            \coordinate (C) at (11,2);
            \coordinate (D) at (3,2);
            \coordinate (E) at (0,8);
            \coordinate (F) at (8,8);
            \coordinate (G) at (11,10);
            \coordinate (H) at (3,10);
            
            
            \node[below] at (A) {\footnotesize{$(0,0,0)$}};
            \node[left] at (E) {\footnotesize{$(0,0,1)$}};
            \node[below] at (B) {\footnotesize{$(1,0,0)$}};
            \node[below] at (D) {\footnotesize{$(0,1,0)$}};

            
            \draw (A)--(B)--(C)--(G)--(H)--(E)--cycle;
            \draw (E)--(F)--(B);
            \draw (F)--(G);
            \draw[dashed] (A)--(D)--(C);
            \draw[dashed] (D)--(H);
            
            
            \foreach \t in {0.1,0.2,...,0.9} {
            \draw[dashed] 
            ($(D)!\t!(C)$) -- ($(H)!\t!(G)$);
            }
             
    \end{tikzpicture}
    \end{minipage}
    \begin{minipage}{0.1\textwidth}
    \centering
    \begin{tikzpicture}[scale=0.3]
    
        
        \draw[->] (-4,0)--(4,0) node[midway, above] {$\qo ^{\ag}_{f(a_1),f(a_2)}$};
             
    \end{tikzpicture}
    \end{minipage}
    \begin{minipage}{0.3\textwidth}
    \centering
    \begin{tikzpicture}[scale=0.3]
    
        
        
            \coordinate (A) at (-0,0);
            \coordinate (B) at (8,0);
            \coordinate (D) at (3,2);
            \coordinate (E) at (0,8);
            
            
            \node[below] at (A) {\footnotesize{$(0,0,0)$}};
            \node[left] at (E) {\footnotesize{$(0,0,1)$}};
            \node[below] at (B) {\footnotesize{$(1,0,0)$}};
            \node[below] at (D) {\footnotesize{$(0,1,0)$}};

            
            \draw (A)--(B)--(C)--(F)--(E)--cycle;
            \draw (B)--(F);
            \draw[dashed] (A)--(D)--(C);
            \draw[dashed] (D)--(E);
            
    \end{tikzpicture}
    \end{minipage}
    \caption{By lemma \ref{lemma: restriction of equivalence relation to noncritical level}, the equivalence relation $\sim_{\ag,a_1,a_2}$ on $I^3$ identifies the elements that have the same image by the first projection $I \times \{1\} \times I \rightarrow I$ (these correspond to the vertical dashed segments). The map $\qo ^{\ag}_{f(a_1),f(a_2)}$ from construction \ref{construction: maps out of cube} carries the point $(0,1,1)$ to $(0,1,0)$, the point $(1,1,1)$ to $(1,1,0)$ and the other vertices of $I^3$ to themselves.}
     \label{figure: collapsed cube as polytope 2}
\end{figure}

The following lemma establishes item (i) in proposition \ref{proposition: filtration of cube}. This lemma is illustrated by figures \ref{figure: collapsed cube as polytope 1} and \ref{figure: collapsed cube as polytope 2}.

\begin{Lemmas}\label{lemma: cube map compatible with equivalence relation}
    The map $\qo \consel$ is compatible with $\sim_{\ag,f(a_i),f(a_{i+1})}$.
\end{Lemmas}

\begin{proof}

    According to lemma \ref{lemma: restriction of equivalence relation to noncritical level}, we wish to prove that $\qo \consel$ is compatible with the projection maps
    $I^{k-1} \times \{1\} \times I^{n-k-1} \rightarrow I^{n-k-1}$ for $1 \le k \le i$ and $I^{k-1} \times \{1\} \times I^{n-k-1} \rightarrow I^{k-1}$ for $i+1 \le k \le n-1$.
    
    For every $1 \le k \le n-1$, let us denote by $\pi_{\le k-1}$ the projection map $I^{n-1} = I^{k-1} \times I^{n-k} \rightarrow I^{k-1}$, and by $\pi_{\ge n-k-1}$ the projection map $I^{n-1} = I^k \times I^{n-k-1} \rightarrow I^{n-k-1}$. We seek to prove that when $0 \le k \le i$, the map $\qo \consel$ factors through $\pi_{\ge n-k-1}$ in restriction to $I^{k-1} \times \{1\} \times I^{n-k-1}$, and when $i+1 \le k \le n-1$, the map $\qo \consel$ factors through $\pi_{\le k-1}$ in restriction to $I^{k-1} \times \{1\} \times I^{n-k-1}$. By items (ii) and (iii) in lemma \ref{lemma: functoriality of extension to C representations}, it is enough to prove that this holds in restriction to vertices, which is indeed the case by construction of $\qo \consel$. \qedhere

\end{proof}

Our next goal is to prove item (ii) in proposition \ref{proposition: filtration of cube}. This will be achieved in several steps.

\emph{\underline{Step 1:}} We will prove that $\qo \consel$ is a local diffeomorphism in restriction to the interior of $I^{n-1}$.

\emph{\underline{Step 2:}} We will prove several lemmas concerning the images of the faces of $I^{n-1}$, and the preimages of the faces of $P \consel$, by $\qo \consel$. As a consequence, we will obtain in particular that the boundary of $I^{n-1}$ is mapped by $\qo \consel$ to the boundary of $P \consel$.

\emph{\underline{Step 3:}} We will deduce from the first two steps that the map $\qo \consel : I^{n-1} \rightarrow P \consel$ restricts to a diffeomorphism between the interiors and is surjective.

\emph{\underline{Step 4:}} We will prove that if two distinct elements of $I^{n-1}$ have the same image by $\qo \consel$, then they both belong to some codimension $1$ face of $I^{n-1}$.

\emph{\underline{Step 5:}} We will deduce item (ii) of proposition \ref{proposition: filtration of cube} from the previous steps, by induction on the length of $\ag$.

\emph{\underline{Step 1:}} The goal of this step is to prove the following lemma.

\begin{Lemmas}\label{lemma: local diffeomorphism}
    The map $\qo \consel : \R^{n-1} \rightarrow \R^{n-1}$ restricts to a local diffeomorphism on the interior of $I^{n-1}$.
\end{Lemmas}

The proof will make use of the following result.

\begin{Proposition}\label{proposition: determinant of dominant diagonal matrices}
    Let $M = (m_{i,j})$ be a real square matrix of size $N \ge 1$ such that for every $1 \le k \le N$, the following inequality holds
    
    $$
    m_{k,k} \ge \sum_{\substack{1 \le i \le N \\ i \neq k }} |m_{i,k}|.
    $$
    
    Then $\Det(M) \ge 0$. Moreover, if the strict inequality holds for every $1 \le k \le N$, then $\Det(M) > 0$.
\end{Proposition}

\begin{proof}

    This is a consequence of Ger\v{s}gorin's theorem (\cite[Theorem 1.1]{Gersgorin}) which asserts that for any complex square matrix $A=(a_{i,j})$ of size $N$ and every eigenvalue $\lambda \in \C$ of $A$, there exists an integer $1 \le k \le N$ such that the following inequality holds
    
    $$
    | \lambda - a_{k,k}| \le \sum_{\substack{1 \le i \le N \\ i \neq k }} |a_{i,k}|.
    $$
    
    If, in the assumption on $M$, one asks for the large inequality only, this implies that the real part of every eigenvalue of $M$ is nonnegative. If one asks for the strict inequality for every $1 \le k \le N$, this implies that the real part of every eigenvalue of $M$ is strictly positive. Since the determinant of $M$ is the product of its eigenvalues, and the non-real ones occur in conjugate pairs, the result follows. \qedhere

\end{proof}

\begin{proof}[Proof of lemma \ref{lemma: local diffeomorphism}]

We will prove that for every element $\tu$ in the interior of $I^{n-1}$, the Jacobian matrix of $\qo \consel$ at $\tu$ in the canonical basis of $\R ^{n-1}$, has strictly positive determinant. Let us denote by $J_{\tu}$ this matrix, by $D_{\tu}$ the differential of $\qo \consel$ at $\tu$, and by $(e_1,\hdots,e_{n-1})$ the canonical basis of $\R^{n-1}$.

For every $1 \le j \le n-1$, the $j^{\mathrm{th}}$ row of $J_{\tu}$ consists of the coordinates of the vector $D_{\tu}(e_j)$ in the canonical basis. We denote by $\pi_{\hat{j}} : I^{n-1} \rightarrow I^{n-2}$ the map that forgets the $j^{\mathrm{th}}$ coordinate. For every $v \in V_{n-2}$, the preimage of $v$ by the restriction of $\pi_{\hat{j}}$ to $V_{n-1}$ contains two elements: the one whose $j^{\mathrm{th}}$ coordinate is $0$, which we denote by $v_{j,0}$, and the other one, whose $j^{\mathrm{th}}$ coordinate is $1$, which we denote by $v_{j,1}$. We have:

$$
D_{\tu} = \sum_{v \in V_{n-2}} (\prod_{\substack{k \neq j \\ (v_{j,0})_k=0}} (1-t_k) \prod_{\substack{k \neq j \\ (v_{j,0})_k=1}} t_k ) (q \consel (v_{j,1}) - q \consel (v_{j,0})) \, \, \, \footnote{For every $v \in V_{n-2}$, the set of those $k\neq j$ such that $(v_{j,0})_k=0$ (resp. $(v_{j,0})_k=1$) is the same as the set of those $k \neq j$ such that $(v_{j,1})_k=0$ (resp. $(v_{j,1})_k=1$)}.
$$

By virtue of the assumption that $\tu$ lies in the interior of $I^{n-1}$, all the coefficients of this linear combination are strictly positive. For every $v \in V_{n-2}$, we denote by $R_j(v)$ the vector $q \consel (v_{j,1}) - q \consel (v_{j,0})$. The determinant of $J_{\tu}$ is a linear combination with strictly positive coefficients of all the elements of the following set

$$
S = \{\Det(R_1(v^1),\hdots, R_{n-1}(v^{n-1})) \mid (v^1,\hdots,v^{n-1}) \in (V_{n-2})^{n-1} \}.
$$

It is therefore enough to prove that all the elements of $S$ are nonnegative, and that they are not all equal to zero. To this end, we make two observations, which follow from the definition of $q \consel$.

\begin{itemize}
    \item If $j \le i$ and the vector $v \in V_{n-2}$ has a coefficient $1$ in some position $j < k \le i$, then $R_j(v) = 0$. Otherwise, define an integer $l$ to be the maximum of the set $\{1 \le k < j \mid (v_{j,0})_k = 1\}$ if the latter is not empty, and $l = 0$ otherwise. We have $R_j(v) = e_j - e_l$ if $l \neq 0$, and $R_j(v) = e_j$ otherwise.
    
    \item If $j \ge i+1$ and the vector $v \in V_{n-2}$ has a coefficient $1$ in some position $i+1 \le k < j$, then $R_j(v) = 0$. Otherwise, define an integer $l$ to be the minimum of the set $\{j < k \le n-1 \mid (v_{j,0})_k = 1\}$ if the latter is not empty, and $l = 0$ otherwise. We have $R_j(v) = e_j - e_l$ if $l \neq 0$, and $R_j(v) = e_j$ otherwise.
\end{itemize}

Denote by $S'$ the set of the determinants of those real square matrices $M$ of size $n-1$ satisfying the following properties:

\begin{itemize}
    \item The diagonal coefficients of $M$ are equal to $1$.
    \item Every column of $M$ of index $1 \le j \le i$ has at most one nonzero coefficient outside of the diagonal. Moreover, if such a coefficient exists, it lies above the diagonal and is equal to $-1$.
    \item Every column of $M$ of index $i+1 \le j \le n-1$ has at most one nonzero coefficient outside of the diagonal. Moreover, if such a coefficient exists, it lies below the diagonal and is equal to $-1$.
\end{itemize}

We have shown that $S = S'$. The desired result then follows from proposition \ref{proposition: determinant of dominant diagonal matrices}. \qedhere

\end{proof}

\underline{\emph{Step 2:}} The goal of this step is to prove a series of results concerning the behavior of $\qo \consel$ with respect to the proper faces of $I^{n-1}$ and $P \consel$. For illustrations of these results, we refer the reader to figures \ref{figure: collapsed cube as polytope 1} and \ref{figure: collapsed cube as polytope 2}.

\begin{Lemmas}\label{lemma: coordinate is zero then also in image}
    Let $\tu \in I^{n-1}$. For every integer $1 \le j \le n-1$, if $t_j = 0$ then $\qo \consel(\tu)_j = 0$. Conversely, if $\qo \consel (\tu)_j = 0$ and none of the coordinates of $\tu$ is equal to $1$, then $t_j = 0$.
\end{Lemmas}

\begin{proof}

    For the first implication, note that if $t_j = 0$, then every $v \in V_{n-1}(\tu)$ satisfies $v_j=0$, so $q \consel (v)_j=0$ as well by construction, and therefore $\overline{q} \consel (\tu)_j=0$. Conversely, suppose that $t_j \neq 0$ and that none of the coordinates of $\tu$ is equal to $1$. Let us prove that $\qo \consel (\tu)_j \neq 0$. For this, it suffices to exhibit a vector $v \in V_{n-1}(\underline{t})$ such that $q \consel (v)_j = 1$. Let us define $v$ as $v_j = 1$ and $v_k = 0$ for $k \neq j$. By virtue of the assumption that none of the coordinates of $\tu$ is equal to $1$, we indeed have $v \in V_{n-1}(\tu)$. We also have $q \consel(v)_j = 1$. \qedhere

\end{proof}

For the next lemma, we will adopt the following notation. We denote by $H_1$ the affine hyperplane of $\R^{n-1}$ defined as $\{\sum_{k \le i} x_k = 1\}$ and by $E_1$ the closed affine halfspace of $\R^{n-1}$ defined as $\{\sum_{k \le i} x_k \le 1\}$. Likewise, we denote by $H_2$ the affine hyperplane $\{\sum_{k \ge i+1} x_k = 1\}$ and by $E_2$ the closed affine halfspace $\{\sum_{k \ge i+1} x_k \le 1\}$.

\begin{Lemmas}\label{lemma: image of face t_j=1}
    Let $\tu \in I^{n-1}$.
    
    \begin{itemize}
        \item If $t_j = 1$ for some $1 \le j \le i$, then $\qo \consel (\tu)_j \in H_1$. Conversely, if $\qo \consel(\tu) \in H_1$, then there exists $1 \le j \le i$ such that $t_j = 1$.
        
        \item If $t_j = 1$ for some $i+1 \le j \le n-1$, then $\qo \consel (\tu)_j \in H_2$. Conversely, if $\qo \consel(\tu) \in H_2$, then there exists $i+1 \le j \le n-1$ such that $t_j = 1$.
    \end{itemize}
\end{Lemmas}

\begin{proof}

    The two cases can be treated similarly so we just do the first one.
    
    For the first implication, note that if $t_j=0$ for some $1 \le j \le n-1$, then every $v \in V_{n-1}(\tu)$ satisfies $v_j=1$. Consequently, if $1 \le j \le i$, by construction one has $q \consel (v)_j \in H_1$ for every $v \in V_{n-1}(\tu)$, and therefore $\qo \consel (\tu)_j \in H_1$ as well.
    
    Conversely, assume that $t_j \neq 1$ for every $1 \le j \le i$. We wish to prove that $\qo \consel (\tu) \notin H_1$. For this, it suffices to exhibit a vector $v \in V_{n-1}(\tu)$ such that $q \consel (v) \notin H_1$. These conditions are fulfilled by the vector $v$ defined as follows:
    
    \begin{itemize}
        \item $v_j = 0$ for every $1 \le j \le i$,
        
        \item for every $j \ge i+1$, $v_j$ is defined as follows. If $t_j=1$ then $v_j = 1$. If $t_j=0$ then $v_j = 0$. Otherwise, if neither $t_j=1$ nor $t_j = 0$, then $v_j = 1$ or $v_j = 0$ indifferently. \qedhere
    \end{itemize}

\end{proof}

\begin{Corollary}\label{corollary: boundary mapped to boundary}

    The boundary of $I^{n-1}$ is mapped by $\qo \consel$ to the boundary of $P \consel$.

\end{Corollary}

\begin{proof}

    This follows from lemmas \ref{lemma: coordinate is zero then also in image} and \ref{lemma: image of face t_j=1}. \qedhere
    
\end{proof}

\begin{Lemmas}\label{lemma: first coefficients are zero}
    Let $\tu \in I^{n-1}$.
    
    \begin{itemize}
        \item If there exists and integer $1 \le j \le i$ such that $t_j = 1$, then $\qo \consel (\tu)_k = 0$ for every $k < j$.
        
        \item If there exists and integer $i+1 \le j \le n-1$ such that $t_j = 1$, then $\qo \consel (\tu)_k = 0$ for every $k > j$.
    \end{itemize}
\end{Lemmas}

\begin{proof}

    If $t_j = 1$, then every $v \in V_{n-1}(\tu)$ satisfies $v_j = 1$. If we further assume that $j \le i$ (resp. $j \ge i+1$), we deduce that every $v \in V_{n-1}(\tu)$ satisfies $q \consel (v)_k = 0$ for every $k < j$ (resp. $k > j$), and therefore $\qo \consel (\tu)$ as well, as desired. \qedhere

\end{proof}

\begin{Lemmas}\label{lemma: max coefficient is not zero}
    Let $\tu \in I^{n-1}$.
    
    \begin{itemize}
        \item Assume that there exists an integer $1 \le j \le i$ such that $t_j=1$ and set $k = \max \{ 1 \le j \le i \mid t_j = 1 \}$. Then $\qo \consel (\tu)_k \neq 0$.
        
        \item Assume that there exists an integer $i+1 \le j \le n-1$ such that $t_j=1$ and set $k = \min \{ i+1 \le j \le n-1 \mid t_j = 1 \}$. Then $\qo \consel (\tu)_k \neq 0$.
    \end{itemize}
\end{Lemmas}

\begin{proof}

    Again, the two cases can be treated similarly, and we just do the first one.
    
    The condition that $\qo \consel (\tu) \neq 0$ is equivalent to the condition that there exists $v \in V_{n-1}(\tu)$ such that $q \consel(v)_k =1$. Consider the vector $v$ defined as follows: for every $1 \le j \le n-1$, if $t_j=1$, then $v_j=1$, otherwise, $v_j=0$. We indeed have $v \in V_{n-1}(\tu)$ by construction. Also by construction, we have $k = \max \{1 \le j \le i \mid v_j =1\}$, and therefore $q \consel(v)_k =1$, as desired. \qedhere

\end{proof}

\underline{\emph{Step 3:}} The goal of this step is to prove the following lemma.

\begin{Lemmas}\label{lemma: map q surjective and homeomorphism between interiors}

    The map $\qo \consel : I^{n-1} \rightarrow P \consel$ restricts to a homeomorphism between the interiors and is surjective.

\end{Lemmas}

\begin{proof}

    By lemma \ref{lemma: local diffeomorphism}, this map carries the interior of $I^{n-1}$ to the interior of $P \consel$. By Corollary \ref{corollary: boundary mapped to boundary}, it carries the boundary of $I^{n-1}$ to the boundary of $P \consel$. Since it is a closed map, we deduce that it is also a closed map in restriction to the interior of $I^{n-1}$. Since it is a local diffeomorphism on the interior of $I^{n-1}$, it is open on the interior of $I^{n-1}$. Since the interior of $P \consel$ is connected, we deduce that $\qo \consel (\Int(I^{n-1})) = \Int(P \consel)$. Furthermore, since $\qo \consel$ is proper, the restriction $\Int(I^{n-1}) \rightarrow \Int(P \consel)$ is proper as well. This is a surjective proper local diffeomorphism, and therefore a covering map. Since the target is contractible and the source is connected, this is a homeomorphism.
    
    Finally, since $\qo \consel (I^{n-1})$ is closed in $\R^{n-1}$, contained in $P \consel$ and contains $\Int(P \consel)$, it is equal to $P \consel$. This completes the proof. \qedhere

\end{proof}

\underline{\emph{Step 4:}} The goal of this step is to prove the following lemma. Again we refer to figures \ref{figure: collapsed cube as polytope 1} and \ref{figure: collapsed cube as polytope 2} for an illustration.

\begin{Lemmas}\label{lemma: same image then on the same face}
    Suppose we have two distinct elements $\tu, \tu' \in I^{n-1}$ such that $\qo \consel(\tu) = \qo \consel(\tu')$. Then there exists a codimension $1$ face of $I^{n-1}$ that contains both $\tu$ and $\tu'$. That is, there exists $1 \le j \le n-1$ such that $t_j=t_j'=0$ or $t_j = t_j'=1$.
\end{Lemmas}

\begin{proof}

    By corollary \ref{corollary: boundary mapped to boundary} combined with lemma \ref{lemma: map q surjective and homeomorphism between interiors}, $\tu$ and $\tu'$ must both lie on the boundary of $I^{n-1}$. Assume first that none of the coordinates of $\tu$ and $\tu'$ is equal to $1$. Then there exists $1 \le j \le n-1$ such that $t_j = 0$. By lemma \ref{lemma: coordinate is zero then also in image} we have $\qo \consel (\tu)_j = 0$ and therefore also $\qo \consel (\tu')_j = 0$. Applying lemma \ref{lemma: coordinate is zero then also in image} again, we obtain $t'_j=0$.
    
    It remains to treat the case when one of the coordinates of $\tu$ or $\tu'$ is equal to $1$. We may assume without lost of generality that there exists $1 \le j \le n-1$ such that $t_j = 1$. We only treat the case $j \le i$, since the case $j \ge i+1$ can be treated similarly. By lemma \ref{lemma: image of face t_j=1} we have $\qo \consel (\tu) \in E_1$, and therefore $\qo \consel (\tu ') \in E_1$ as well. By lemma \ref{lemma: image of face t_j=1}, there exists $1 \le j' \le i$ such that $t_{j'}' = 1$. We take $j$ and $j'$ to be such that $t_k \neq 1$ for every $j < k \le i$ and $t_k' \neq 1$ for every $j' < k \le i$.
    
    By lemma \ref{lemma: first coefficients are zero} we have $\qo \consel (\tu')_k = 0$ for every $k < j'$, and by lemma \ref{lemma: max coefficient is not zero} we have $\qo \consel(\tu)_j \neq 0$. Since $\qo \consel(\tu) = \qo \consel(\tu')$, we deduce that $j \ge j'$. Similarly we obtain $j \le j'$, and therefore $j=j'$. This completes the proof. \qedhere

\end{proof}

\underline{\emph{Step 5:}} The goal of this step is to complete the proof of item (ii) in proposition \ref{proposition: filtration of cube}. We wish to prove that the map $I^{n-1} / \sim_{\ag,f(a_i),f(a_{i+1})} \rightarrow P \consel$ induced from $\qo \consel$, is a homeomorphism. We will do this by induction on the length of $\ag$. The case of length $1$ is handled by lemma \ref{lemma: proposition in case n=1}. We thus assume that for every $m < n$, every nonconstant increasing sequence $\bg$ of elements of $A$ of length $m$, and every $0 \le k < m$ such that $b_k < b_{k+1}$, the map $I^{m-1} / \sim_{\bg,b_k,b_{k+1}} \rightarrow Q_{b_k,b_{k+1}}^{\bg}$ induced by $\qo_{b_k,b_{k+1}}^{\bg}$ is a homeomorphism.

Combining lemmas \ref{lemma: map q surjective and homeomorphism between interiors} and \ref{lemma: same image then on the same face}, we see that it remains to prove that if two elements $\tu,\tu' \in I^{n-1}$ are such that $t_j=t_j'=0$ or $t_j = t_j' = 1$ for some $1 \le j \le n-1$, and satisfy $\qo \consel (\tu) = \qo \consel (\tu')$, then $\tu \sim_{\ag,f(a_i),f(a_{i+1})} \tu'$.

It will be convenient to use the following notation.

\begin{Notation}\label{notation: polytope and map at given level}
    For every $\Delta^{\bg} \in \Delta_A$ with $\bg$ of length $m$, every $0 \le k < m$ such that $b_k < b_{k+1}$, and every $f(b_{k+1}) < s < f(b_k)$, we let $P^{\bg}_s = P_{b_k,b_{k+1}}^{\bg}$ and $\qo^{\bg}_s = \qo_{b_k,b_{k+1}}^{\bg}$.
\end{Notation}

We begin by addressing the case $t_j = t_j'=0$ for some $1 \le j \le n-1$. Consider the morphism $\Phi$ in $\Delta_A$ determined by the following commutative diagram of posets

$$
\xymatrix{
[n-1] \ar[rd]_-{\ag \circ \delta_n^j} \ar[rr]^{\delta_n^j} && [n] \ar[ld]^-{\ag} \\
& A.
}
$$

Recall that, by definition, the map $C_A(\Phi) : C_A(\Delta^{\ag \circ \delta_n^j}) \rightarrow C_A(\Delta^{\ag})$ is induced from the map

$$
F(\Delta^{n-1})(0,n-1) \times [f(a_n),f(a_0)] \xrightarrow{F(\delta_n^j)(0,n-1) \times \id} F(\Delta^n)(0,n) \times [f(a_n),f(a_0)].
$$

In particular, the latter is compatible with $\sim_{\ag \circ \delta_n^j}$ and $\sim_{\ag}$. It therefore restricts to a map $I^{n-2} \times \{s\} \rightarrow I^{n-1} \times \{s\}$ compatible with $\sim_{\ag \circ \delta_n^j,s}$ and $\sim_{\ag,s}$ for every real number $f(a_{i+1}) < s < f(a_i)$. The latter map corresponds to the inclusion of the face $\{x_j = 0 \}$ of $I^{n-1}$, by proposition \ref{proposition: image of the faces by F}. Consider the following diagram

$$
\xymatrix{
I^{n-2} \ar[d]_-{\qo _s^{\ag \circ \delta_n^j}} \ar[r] & I^{n-1} \ar[d]^-{\qo ^{\ag}_s} \\
\R^{n-2} \ar[r] & \R^{n-1},
}
$$

where the two horizontal arrows are given by the inclusion of the component $\{x_j = 0\}$. This diagram is commutative; indeed, by item (iii) in lemma \ref{lemma: functoriality of extension to C representations}, this follows from the commutativity of the corresponding diagram restricted to the vertices of $I^{n-2}$ and $I^{n-1}$, which holds by construction. By the induction hypothesis, the left vertical arrow becomes injective after passing to the quotient, which yields the desired property about the right vertical arrow.

We now treat the case $t_j = t_j' = 1$ for some $1 \le j \le n-1$. It is necessary to distinguish the two cases $j \le i$ and $j \ge i+1$. They can be handled similarly, so we will only treat the first one. We therefore assume from now on that $j \le i$. Let us denote by $\ag_{\ge j}$ the sequence $[a_j \le a_{j+1} \le \hdots \le a_n]$ and by $\Psi$ the inclusion $\ag_{\ge j} \rightarrow \ag$. Recall that, by definition, the map $C_A(\Psi) : C_A(\Delta^{\ag_{\ge j}}) \rightarrow C_A(\Delta^{\ag})$ is induced from the map

$$
\fonctionsansnom{F(\Delta^n)(j,n) \times [f(a_n),f(a_j)]}{F(\Delta^n) \times [f(a_n),f(a_0)]}{(\gamma,s)}{(\gamma \circ \delta,s)}
$$

for any $\delta \in F(\Delta^n)(0,j)$. For the sake of this proof, we will take $\delta$ to correspond to the element $(0,\hdots,0) \in I^{j-1}$. For every real number $f(a_{i+1}) < s < f(a_i)$, the map $I^{n-j-1} \rightarrow I^{n-1}$ corresponding to the inclusion of the face $\{x_1 = \hdots = x_{j-1} = 0, x_j=1 \}$ is therefore compatible with the equivalence relations $\sim_{\ag_{\ge j},s}$ and $\sim_{\ag,s}$. Using the notation of lemma \ref{lemma: functoriality of extension to C representations}, this is the map $\kappa_u$ for $E = \{j+1,\hdots,n-1\}$ and $u=(0,\hdots,0,1) \in I^j$. Note also that the map  $I^{n-j-1} / \sim_{\ag_{\ge j},s} \rightarrow I^{n-1} / \sim_{\ag,s}$ induced from $\kappa_u$ is surjective onto the image of the face $\{x_j=1\}$ in the quotient. It therefore suffices to prove that the map $\qo_s^{\ag} \circ \kappa_u$ induces an injective map on $I^{n-j-1} / \sim_{\ag_{\ge j},s}$. Denote by $\pi_{\ge j+1}$ the projection map $\R^j \times \R^{n-j-1} \rightarrow \R^{n-j-1}$. We have:

$$
\begin{aligned}
\pi_{\ge j+1} \circ \qo_s^{\ag} \circ \kappa_u & = \pi_{\ge j+1} \circ \overline{q_s^{\ag} \circ (\kappa_u)_{\mid V_{n-j-1}}} \, \,\, \text{by item (iii) in lemma \ref{lemma: functoriality of extension to C representations}} \\
& = \overline{ \pi_{\ge j+1} \circ q_s^{\ag} \circ (\kappa_u)_{\mid V_{n-j-1}}} \, \,\, \text{by item (i) in lemma \ref{lemma: functoriality of extension to C representations}} \\
& = \qo _s^{\ag_{\ge j+1}} \, \, \, \text{by construction}.
\end{aligned}
$$

Furthermore, by the induction hypothesis, the map $\qo_s^{\ag_{\ge j+1}}$ induces an injective map on $I^{n-j-1} / \sim_{\ag_{\ge j},s}$. We conclude that $\qo_s^{\ag} \circ \kappa_u$ as well, as desired.

Step 5 as well as the proof of item (ii) in proposition \ref{proposition: filtration of cube} are now complete.

We next address the proof of item (iii) in proposition \ref{proposition: filtration of cube}. We start with a lemma describing explicitly the equivalence relations $\sim_{\ag,s}$ when $s$ is a critical value.

\begin{Lemmas}\label{lemma: restriction of equivalence relation to critical level}
    Let $s \in \{f(a_n),f(a_{n-1}), \hdots,f(a_0)\}$ and let $[a_p = a_{p+1} = \hdots = a_q ] \subset \ag$ be the constant subsequence of $\ag$ formed by those critical points whose image by $f$ is $s$. Regard $I^{n-1}$ as the set of maps $\tu : \{0,\hdots,n\} \rightarrow I$ such that $t_0=t_n=1$. For every $\tu, \tu' \in I^{n-1}$, we have $\tu \sim_{\ag,s} \tu'$ if and only if there exists $0 \le m \le q$ and $p \le l \le n$ satisfying $m \le l$ and such that $t_m=t_m'=1$, $t_l=t_l'=1$ and $(t_m,\hdots,t_l)=(t_m',\hdots,t_l')$.
    
    Consequently, $\sim_{\ag,s}$ is the smallest equivalence relation on $I^{n-1}$ that identifies the elements that have the same image by at least one map among the projection maps $I^{k-1} \times \{1\} \times I^{n-k-1} \rightarrow I^{n-k-1}$  for $0 < k \le p$, the projections maps $I^{k-1} \times \{1\} \times I^{n-k-1} \rightarrow I^{k-1}$ for $q \le k < n$, and the maps $I^{k-1} \times \{1\} \times I^{n-k-1} \rightarrow *$ for $p \le k \le q$.
\end{Lemmas}

\begin{proof}

    Consider two elements $\tu$ and $\tu'$ of the morphism space $F(\Delta^n)(0,n)$. By definition, we have $\tu \sim_{\ag,s} \tu'$ if and only if $(\tu,s) \sim_{\ag} (\tu',s)$, if and only if $\tu$ and $\tu'$ are respectively of the form $\nu \circ \gamma \circ \eta$ and $\nu' \circ \gamma \circ \eta'$, for some $\nu, \nu' \in F(\Delta^n)(l,n)$, $\eta, \eta' \in F(\Delta^n)(0,m)$ and $\gamma \in F(\Delta^n)(m,l)$, with $0 \le m \le q$ and $p \le l \le n$ such that $m \le l$. The result therefore follows from proposition \ref{proposition: compositions in F(n)}. \qedhere

\end{proof}

By item (ii) of proposition \ref{proposition: filtration of cube}, we have for every real number $f(a_{i+1}) < s < f(a_i)$ a map $\qo \consel : I^{n-1} \rightarrow P \consel$, which induces a homeomorphism $P \consel \xrightarrow[]{\simeq} I^{n-1} / \sim_{\ag,s}$. Let us denote the latter by $h_s$. Moreover, by lemmas \ref{lemma: restriction of equivalence relation to noncritical level} and \ref{lemma: restriction of equivalence relation to critical level}, the two projection maps $I^{n-1} \rightarrow I^{n-1} / \sim_{\ag,f(a_{i+1})}$ and $I^{n-1} \rightarrow I^{n-1} / \sim_{\ag,f(a_i)}$ both factor through the projection map $I^{n-1} \rightarrow I^{n-1} / \sim_{\ag,f(a_i),f(a_{i+1})}$. Point (ii) of proposition \ref{proposition: filtration of cube} thus provides two maps $h_{f(a_{i+1})} : P \consel \rightarrow I^{n-1} / \sim_{\ag,f(a_{i+1})}$ and $h_{f(a_i)} : P \consel \rightarrow I^{n-1} / \sim_{\ag,f(a_i)}$. These maps assemble into a continuous map

$$
\fonctionsansnom{P \consel \times [f(a_{i+1}),f(a_i)]}{C_A(\Delta_A)}{(q,s)}{h_s(q).}
$$

Proving item (iii) amounts to proving that there are faces $F_{f(a_{i+1})}^{\ag,+} \subseteq P \consel \supseteq F_{f(a_i)}^{\ag,-}$ such that $h_{f(a_{i+1})}$ and $h_{f(a_i)}$ induce homeomorphisms $P \consel / F_{f(a_{i+1})}^{\ag,+} \simeq I^{n-1} / \sim_{\ag,f(a_{i+1})}$ and $P \consel / F_{f(a_i)}^{\ag,-} \simeq I^{n-1} / \sim_{\ag,f(a_i)}$ respectively.

The two cases of $a_i$ and $a_{i+1}$ can be treated similarly, so we only do that of $a_{i+1}$. We need to understand what identifications the equivalence relation $\sim_{\ag,f(a_{i+1})}$ makes in addition to $\sim_{\ag,f(a_i),f(a_{i+1})}$. For this, we denote by $q$ the integer such that $a_{i+1} = a_{i+2} = \hdots = a_q$, and $a_q < a_{q+1}$ in the case when $q < n$. We have the following lemma.

\begin{Lemmas}\label{lemma: further identifications at critical level}
    The map
    
    $$
    I^{n-1} / \sim_{\ag,f(a_i),f(a_{i+1})} \rightarrow I^{n-1} / \sim_{\ag,f(a_{i+1})}
    $$
    
    induces a homeomorphism between $I^{n-1} / \sim_{\ag,f(a_{i+1})}$ and the quotient of $I^{n-1} / \sim_{\ag,f(a_i),f(a_{i+1})}$ by the image of
    
    $$
    \bigcup_{i+1 \le k \le q} I^{k-1} \times \{1\} \times I^{n-k-1} \subset I^{n-1}.
    $$
    
\end{Lemmas}

\begin{proof}

    For the sake of this proof, we denote this map by $\varphi$ and this union by $\El$. We first note that $\varphi$ is constant on the image of $\El$ in $I^{n-1} / \sim_{\ag,a_1,a_{i+1}}$; indeed, by lemma \ref{lemma: restriction of equivalence relation to critical level}, the equivalence relation $\sim_{\ag,f(a_{i+1})}$ identifies all the elements of $\El$. We will complete the proof by showing that if $\tu, \tu'$ are two elements of $I^{n-1}$ such that $\tu \nsim_{\ag,a_i,a_{i+1}} \tu'$ and $\tu \sim_{\ag,f(a_{i+1})} \tu'$, then $\tu, \tu' \in \El$. The assumption that $\tu \sim_{\ag,f(a_{i+1})} \tu'$ implies, by lemma \ref{lemma: restriction of equivalence relation to critical level}, that there exist $0 \le m \le q$ and $i+1 \le l \le n$ satisfying $m \le l$ and such that $t_m = t_m' = 1$, $t_l=t_l'=1$ and $(t_m,\hdots,t_l)=(t_m',\hdots,t_l')$. But by virtue of the assumption that $\tu' \nsim_{\ag,a_i,a_{i+1}} \tu'$, we must have $m \ge i+1$, according to lemma \ref{lemma: restriction of equivalence relation to noncritical level}. In particular, $\tu,\tu' \in \El$, as desired. \qedhere

\end{proof}

The following corollary establishes item (iii) in proposition \ref{proposition: filtration of cube}. Recall that $H_2$ denotes the affine hyperplane of $\R^{n-1}$ defined as $\{\sum_{k \ge i+1} x_i =1 \}$.

\begin{Corollary}\label{corollary: }

    The map $h_{f(a_{i+1})}$ induces a homeomorphism between $I^{n-1}/\sim_{\ag,f(a_{i+1})}$ and the quotient of $P \consel$ by the face
    
    $$
    F_{f(a_{i+1})}^{\ag,+} = P \consel \cap H_2 \cap \{x_{q+1}=x_{q+2}=\hdots=x_{n-1}=0\}.
    $$

\end{Corollary}

\begin{proof}

    According to lemma \ref{lemma: further identifications at critical level}, this is equivalent to the condition that
    
    $$
    \qo \consel (\bigcup_{i+1 \le k \le q} I^{k-1} \times \{1\} \times I^{n-k-1} ) = F_{f(a_{i+1})}^{\ag,+}.
    $$
    
    The inclusion $\subseteq$ follows from lemmas \ref{lemma: image of face t_j=1} and \ref{lemma: first coefficients are zero}. The inclusion $\supseteq$ follows from lemmas \ref{lemma: image of face t_j=1} and \ref{lemma: max coefficient is not zero}. \qedhere

\end{proof}

The proof of proposition \ref{proposition: filtration of cube} is now complete.

\subsubsection{A triviality theorem for the filtration on the cubical realization}\label{section: triviality theorem for filtration on cube}

The goal of this section is to complete the proof of theorem \ref{theorem: filtration on cube is trivial} (and therefore that of theorem \ref{theorem: stratified simplices and cubes as balls} in the case $G=C_A$). We start by proving lemma \ref{lemma: gluing trivializations over interval}.

\begin{proof}[Proof of lemma \ref{lemma: gluing trivializations over interval}]

    The result also holds when considering either $[\lambda_k,\lambda_0]$, $(\lambda_k,\lambda_0]$ or $[\lambda_k,\lambda_0)$ instead of $(\lambda_k,\lambda_0)$. We will prove it in this more general version. Working by induction over $k$, we see that it suffices to treat the case $k=2$. We will do it in the case of the open interval $(\lambda_2,\lambda_0)$, but the same proof works in the other cases, by closing the intervals at the relevant places. Pick two trivializations
    
    $$
    \xymatrix{
    D^{n-1} \times (\lambda_2,\lambda_1] \ar[rr]^-{\psi_0} \ar[rd] && E \ar[ld]^-p \\
    & (\lambda_2,\lambda_1]
    }
    \,\,\, \text{ and } \,\,\,
    \xymatrix{
    D^{n-1} \times [\lambda_1,\lambda_0) \ar[rr]^-{\psi_1} \ar[rd] && E \ar[ld]^-p \\
    & [\lambda_1,\lambda_0).
    }
    $$
    
    We denote by $\varphi$ the homeomorphism of $D^{n-1}$ defined as $\psi_1(\cdot, \lambda_1)^{-1} \circ \psi_0(\cdot, \lambda_1)$. In order to construct a trivialization of $p$ over $(\lambda_2,\lambda_0)$, we use the description of the space $D^{n-1} \times (\lambda_2,\lambda_0)$ as the pushout
    
    $$D^{n-1} \times (\lambda_2,\lambda_1] \sqcup_{D^{n-1} \times \{\lambda_1\}} D^{n-1} \times [\lambda_1,\lambda_0).
    $$ 
    
    The two homeomorphisms $\psi_0$ and
    
    $$
    \fonctionsansnom{D^{n-1} \times [\lambda_1,\lambda_0)}{E}{(s,x)}{\psi_1(s,\varphi(x))}
    $$
    
    yield the desired trivialization. \qedhere

\end{proof}  

We are now ready to conclude.

\begin{proof}[Proof of theorem \ref{theorem: filtration on cube is trivial}]

    We have essentially already all the ingredients to prove items (i) and (ii). Namely, by lemma \ref{lemma: restriction of equivalence relation to critical level}, the equivalence relations $\sim_{\ag,f(a_0)}$ and $\sim_{\ag,f(a_n)}$ both identify all the elements of $I^{n-1}$, which proves item (i). As for item (ii), let us denote by $[a_0 = a_{j_0} < a_{j_1} < \hdots < a_{j_{k-1}} < a_{j_k} = a_n]$ the ordered sequence of distinct elements of $\ag$. We distinguish two cases. First, if $k=1$ then item (ii) follows from item (iii) in proposition \ref{proposition: filtration of cube} combined with the fact that a $(n-1)$-polytope is homeomorphic to a $(n-1)$-dimensional disc (proposition \ref{proposition: a polytope is a disc}). Second, assume that $k \ge 2$.  By item (iii) of proposition \ref{proposition: filtration of cube}, combined with lemma \ref{lemma: improved proper face collapsing lemma} and the fact that a $(n-1)$-dimensional polytope is homeomorphic to a $(n-1)$-dimensional disc, we obtain that the map $f_{C_A}(\Delta^{\ag})$ is a $(n-1)$-dimensional disc bundle over $(f(a_n),f(a_{j_{k-1}})]$, over $[f(a_{j_1}),f(a_0))$ as well as over $[f(a_{j_r}),f(a_{j_{r-1}})]$ for every $2 \le r \le k-1$. In that case, item (ii) of theorem \ref{theorem: stratified simplices and cubes as balls} follows from proposition \ref{lemma: gluing trivializations over interval}.
    
    It remains to prove item (iii). We start by justifying that the restriction of the homeomorphism $\Sigma (D^{n-1}) \simeq C_A(\Delta^{\ag})$ to $\Sigma(\SM^{n-2})$, yields a homeomorphism between $\Sigma (\SM^{n-2})$ and the image of the boundary of $I^{n-1} \times [f(a_n),f(a_0)]$ inside $C_A(\Delta^{\ag})$. For every $f(a_n) < s < f(a_0)$, the homeomorphism $\Sigma (D^{n-1}) \simeq C_A(\Delta^{\ag})$ restricts to a homeomorphism between $D^{n-1}$ and $I^{n-1} / \sim_{\ag,s}$. By proposition \ref{proposition: filtration of cube} combined with lemma \ref{lemma: map q surjective and homeomorphism between interiors}, there exists a polytope $P$ and a proper face $F \subset P$ such that the pair $(I^{n-1}/\sim_{\ag,s},(\partial I^{n-1})/\sim_{\ag,s})$ is homeomorphic to $(P/F, \partial P/F)$. By corollary \ref{corollary: proper face collapsing homeomorphic to disc}, we deduce that the pair $(I^{n-1}/\sim_{\ag,s},(\partial I^{n-1})/\sim_{\ag,s})$ is homeomorphic to $(D^{n-1},\SM^{n-2})$. Now, every homeomorphism of $D^{n-1}$ restricts to a homeomorphism of $\SM^{n-2}$.

    We denote by $\partial C_A(\Delta^{\ag})$ the image of $\partial (I^{n-1} \times [f(a_n),f(a_0)])$ in $C_A(\Delta^{\ag})$. What we wish to prove is therefore equivalent to the condition that the natural map $|\partial \Delta^{\ag}|_{C_A} \rightarrow C_A(\Delta^{\ag})$ induces a homeomorphism between $|\partial \Delta^{\ag}|_{C_A}$ and $\partial C_A(\Delta^{\ag})$. In order to prove this, we will use the description of $|\partial \Delta^{\ag}|_{C_A}$ as a colimit provided by proposition \ref{proposition: stratified horn as colimit of faces}.
    
    We first note that this description implies that $|\partial \Delta^{\ag}|_{C_A}$ is compact; it therefore suffices to prove the weaker statement obtained by replacing "homeomorphism" by "continuous bijection".
    
    Second, we note that every element of $|\partial \Delta^{\ag}|_{C_A}$ is mapped to an element of $\partial C_A(\Delta^{\ag})$. Indeed, this is equivalent to the condition that the faces of $C_A(\Delta^{\ag})$, in the sense of definition \ref{definition: face of image of simplex by three functors}, have image contained in $\partial C_A(\Delta^{\ag})$. Now, for every $1 \le i \le n-1$, every element of the image of the face of $C_A(\Delta^{\ag})$ opposite to the vertex $i$ is the equivalence class of an element $(\tu,s) \in I^{n-1} \times [f(a_n),f(a_0)]$ such that $t_i = 0$, by construction of the functor $C_A$ together proposition \ref{proposition: image of the faces by F}. Furthermore, for $i=0$ and $i=n$, every element of the face of $C_A(\Delta^{\ag})$ opposite to the vertex $i$ is the equivalence class of an element $(\tu,s) \in I^{n-1} \times [f(a_n),f(a_0)]$ such that $t_i = 1$, by construction of the functor $C_A$ together with proposition \ref{proposition: compositions in F(n)}.
    
    We complete the proof by constructing an inverse to the map $|\partial \Delta^{\ag}|_{C_A} \rightarrow \partial C_A(\Delta^{\ag})$. We do this by constructing a map $\varphi : \partial (I^{n-1} \times [f(a_n),f(a_0)]) \rightarrow |\partial \Delta^{\ag}|_{C_A}$ compatible with $\sim_{\ag}$. Let $(\tu,s) \in \partial (I^{n-1} \times [f(a_n),f(a_0)])$.
    
    \begin{itemize}
        \item Suppose first that $s=f(a_n)$ (resp. $s=f(a_0)$). The inclusion of the length $0$ sequence $[a_n] \rightarrow \ag$ (resp. $[a_0] \rightarrow \ag$) determines a morphism $\Delta^0 \rightarrow \Delta^{\ag}$ that factors through $\partial \Delta^{\ag}$; we declare $\varphi(\tu,s)$ to be the element of $|\partial \Delta^{\ag}|_{C_A}$ obtained by applying the functor $C_A$ to the morphism $\Delta^0 \rightarrow \partial \Delta^{\ag}$.
        
        \item Suppose otherwise that $f(a_n) < s < f(a_0)$.
        
        \begin{itemize}
            \item Suppose that there exists an integer $1 \le j \le n-1$ such that $t_j = 0$. For every $1 \le j \le n-1$ such that $t_j=0$, the element $(\pi_{\hat{j}}(\tu),s) \in I^{n-2} \times [f(a_n),f(a_0)]$ determines an element of $C_A(\Delta^{\ag \circ \delta_n^j})$ (where $\pi_{\hat{j}} : I^{n-1} \rightarrow I^{n-2}$ denotes the map that forgets the $j^{\text{th}}$ coordinate). Moreover, all these elements are mapped to the same element of $|\partial \Delta^{\ag}|_{C_A}$, which we declare to be $\varphi(\tu,s)$. 
            
            \item Suppose otherwise that $t_j \neq 0$ for every integer $1 \le j \le n-1$. In particular, there exists an integer $1 \le j \le n-1$ such that $t_j = 1$.
            
            \begin{itemize}
                \item Suppose first that there exists an integer $0 \le i < n$ such that $f(a_{i+1}) < s < f(a_i)$. Define an integer $m$ to be the maximum of the set $\{1 \le j \le i \mid t_j=1\}$ if the latter is nonempty, and $m=0$ otherwise. Define another integer $l$ to be the minimum of the set $\{i+1 \le j \le n-1 \mid t_j=1\}$ if the latter is nonempty, and $l=n$ otherwise. The inclusion of the subsequence $\ag_{\ge m}^{\le l} := [a_m \le a_{m+1} \le \hdots \le a_l ] \rightarrow \ag$ determines a morphism $\Delta^{\ag_{\ge m}^{\le l}} \rightarrow \Delta^{\ag}$. By virtue of the assumption that we made on $\tu$, we have either $m \ge 1$ or $l \le n-1$, hence this is a proper subsequence, and therefore this inclusion factors through $\partial \Delta^{\ag}$. We declare $\varphi(\tu,s)$ to be the image, by the map obtained by applying $C_A$ to this morphism, of the element of $C_A(\Delta^{\ag_{\ge m}^{\le l}})$ determined by $((t_{m+1},\hdots,t_{l-1}),s) \in I^{l-m-1} \times [f(a_l),f(a_m)]$.
                
                \item Suppose otherwise that there exists an integer $0 < i < n$ such that $s=f(a_i)$. Let $\bg = [a_p = \hdots = a_q]$ be the constant subsequence of those elements of $\ag$ whose image by $f$ is $s$. Note that $C_A(\Delta^{\bg})$ is a point. If there exists an integer $p \le j \le q$ such that $t_j=1$, we declare $\varphi(\tu,s)$ to be the element determined by applying $C_A$ to the inclusion $\Delta^{\bg} \rightarrow \Delta^{\ag}$. Otherwise, we define an integer $m$ to be the maximum of the set $\{1 \le j < p \mid t_j=1\}$ if the latter is nonempty, and $m=0$ otherwise. We define another integer $l$ to be the minimum of the set $\{q < j \le n-1 \mid t_j=1\}$ if the latter is nonempty, and $l=n$ otherwise. Then, we define $\varphi(\tu,s)$ as above.
            \end{itemize}
        \end{itemize}
    \end{itemize}
    
    One can check that the map $\varphi$ is compatible with the restriction of the equivalence relation $\sim_{\ag}$ to $\partial (I^{n-1} \times[f(a_n),f(a_n)])$, and that the induced map $\partial C_A(\Delta^{\ag}) \rightarrow |\partial \Delta^{\ag}|_{C_A}$ is the desired inverse. \qedhere

\end{proof}

\newpage

\section{\texorpdfstring{Unbroken simplices of $\Nl(\mathcal{M})$}{Unbroken simplices of N(M)}}\label{section: unbroken simplices}

In this section, we consider the flow category $\Ml$ of a Morse-Smale pair $(f,\xi)$ on a smooth closed manifold $X$. We study the flow coherent nerve of $\Ml$ (definition \ref{definition: flow coherent nerve}) denoted $S_{\Ml}$. We begin by proving in section \ref{section: flow coherent nerve is an infinity category} that this is an $\infty$-category (proposition \ref{proposition: S_M is an infinity category}). Then, in section \ref{section: equivalence between flow and homotopy coherent nerve}, we prove that the inclusion $S_{\Ml} \rightarrow \Nl(\Ml)$ is an equivalence of $\infty$-categories (proposition \ref{proposition: equivalence between flow and homotopy coherent nerve}).

In the next two sections, we carry out the comparison between $S_{\Ml}$ and $\Sing_A(X)$ initiated in section \ref{section: comparison between realizations}. In section \ref{section: filling inner horns in Sing C A (X)}, we prove that the simplicial set $\Sing_{C_A}(X)$ is an $\infty$-category (proposition \ref{proposition: Sing_C_A is an infinity category}). In section $\ref{section: comparison between S M and Sing C A (X)}$, we introduce an $\infty$-subcategory $\Sing_{C_A,r}(X) \subseteq \Sing_{C_A}(X)$ (definition \ref{definition: restricted versions}) such that the morphism $S_{\Ml} \rightarrow \Sing_{C_A}(X)$ corestricts to a morphism $S_{\Ml} \rightarrow \Sing_{C_A,r}(X)$, and we prove that the latter is an equivalence of $\infty$-categories (proposition \ref{proposition: S_M -> Sing_C_A,r is an equivalence}).

\subsection{\texorpdfstring{Filling inner horns in the flow coherent nerve of $\Ml$}{Filling inner horns in the flow coherent nerve of M}}\label{section: flow coherent nerve is an infinity category}

The goal of this section is to prove the following proposition.

\begin{Proposition}\label{proposition: S_M is an infinity category}
    The flow coherent nerve $S_{\Ml}$ of $\Ml$ is an $\infty$-category.
\end{Proposition}

The proof of proposition \ref{proposition: S_M is an infinity category} will require some preliminaries. We use the characterization of $\infty$-categories in terms of the inner horn filling property (definition \ref{definition: infinity category}). We thus wish to prove that for every $n \ge 0$ and $0 < k < n$, every lifting problem of the form

$$
\xymatrix{
\Lambda_k^n \ar[r] \ar[d] & S_{\Ml} \\
\Delta^n \ar@{-->}[ru]
}
$$

has a solution. By adjunction, to such a lifting problem is associated a functor between topological categories $F(\Lambda_k^n) \rightarrow \Ml$, and a solution is the same as an extension of this functor to $F(\Delta^n)$ such that the associated $n$-simplex of $\Ml$ is unbroken (in the sense of definition \ref{definition: unbroken simplices and diagrams}).

Before attacking the proof, let us briefly outline our strategy. In the above diagram, we can regard the morphism $\Lambda_k^n \rightarrow S_{\Ml}$ as a morphism $\Lambda_k^n \rightarrow \Nl(\Ml)$, and it follows from the fact that $\Nl(\Ml)$ is an $\infty$-category that there exists an extension of the latter to a morphism $\Delta^n \rightarrow \Nl(\Ml)$. Our goal if to justify that there exists an extension which is unbroken. In order to give an idea of what we want to prove, let us start by considering an arbitrary topological category $\Cj$ and let us outline how one can prove that $\Nl(\Cj)$ is an $\infty$-category \footnote{More generally, the homotopy coherent nerve of every simplicially enriched category such that for every objects $x,y$, the simplicial set of morphisms between $x$ and $y$ is a Kan complex, is an $\infty$-category. This is due to Cordier and Porter (\cite{CordierPorter}) and a proof if presented in \cite[Section 2.4.5]{kerodon}. This implies the result for topological categories because for every topological category $\Cj$, the simplicially enriched category $\Sing(\Cj)$, defined by $\Sing(\Cj)(x,y)=\Sing(\Cj(x,y))$, satisfies the condition.}. Consider an inner horn $\Lambda_k^n \rightarrow \Nl(\Cj)$. By adjunction, it corresponds to a functor of topological categories $F(\Lambda_k^n) \rightarrow \Cj$. As we will prove in this section, the topological category $F(\Lambda_k^n)$ can be described as follows: the canonical functor of topological categories $F(\Lambda_k^n) \rightarrow F(\Delta^n)$ induces a bijection at the level of objects, and a homeomorphism between every two morphism spaces except the ones between $0$ and $n$, where it corresponds through the homeomorphisms \eqref{homeomorphisms} to the inclusion of the \emph{hollow} $(n-1)$-\emph{cube} $\sqcap_k^{n-1} \subset I^{n-1}$, which is the subspace of $I^{n-1}$ obtained by removing the interior together with the face $\{t_k=0\}$ (lemma \ref{lemma: image of an inner horn by F}). Therefore,
filling an inner horn in $\Nl(\Cj)$ is the same as filling a hollow cube in some morphism space of $\Cj$, which is always possible since $I^{n-1}$ retracts onto every hollow cube. In order to prove proposition \ref{proposition: S_M is an infinity category}, we will follow the same approach, which will lead us to consider the problem of filling a hollow cube in some morphism space of $\Ml$ in such a way that the resulting simplex of $\Nl(\Ml)$ is unbroken. We will prove that a filling satisfying the desired property exists by using the structure of smooth manifold with corners on the morphism spaces of $\Ml$ provided by proposition \ref{proposition: structure of smooth manifold with corners on space of broken traj}.

We begin by studying the topological category $F(\Lambda_k^n)$ for every $n \ge 0$ and $0 < k < n$. We will give and prove the precise description of $F(\Lambda_k^n)$ in terms of hollow cubes that we have sketched above.

Since we studied the restriction of $F$ to $\Delta$ in section \ref{section: Study of the geometric realization of the Leitch-Cordier functor} and we know that $F$ preserves colimits, we start by writing $\Lambda^n_k$ as a colimit of standard simplices (a similar construction was given in \ref{construction: boundary of simplex as colimit of simplices} in the case of $\partial \Delta^n$). Informally speaking, $\Lambda_k^n$ is obtained by gluing together the faces of $\Delta^n$ that are different from the one opposite to the vertex $k$, along their common boundary faces. More formally:

\begin{Construction}\label{construction: horn as colimit of simplices}
    Let $n \ge 0$ and $0 \le k \le n$. We let $\mathrm{Facets}(\Lambda_k^n)$ be the category whose objects are the nondegenerate $(n-1)$-simplices of $\Delta^n$ different from the face opposite to the vertex $k$, together with the nondegenerate $(n-2)$-simplices of $\Delta^n$, and whose non-identity morphisms correspond to the inclusions of $(n-2)$-simplices as boundary faces of $(n-1)$-simplices. Recall from construction \ref{construction: boundary of simplex as colimit of simplices} the category $\mathrm{Facets}(\Delta^n)$ and the notation $F_i$, $F_{i,j}$, $\phi_{i,j}^i$, $\phi_{i,j}^j$. The category $\mathrm{Facets}(\Lambda_k^n)$ can be described as obtained from the category $\mathrm{Facets}(\Delta^n)$ by deleting the object $F_k$ as well as the morphisms with target $F_k$.

    There is a functor $\Ld_k^n : \mathrm{Facets}(\Lambda_k^n) \rightarrow \sSet$ defined by 

$$
\left\{
    \begin{array}{l}
        F_i \mapsto \Delta^{n-1} \\
        F_{i,j} \mapsto \Delta^{n-2} \\
        \phi_{i,j}^i \mapsto \delta_{n-1}^{j-1} \\
        \phi_{i,j}^j \mapsto \delta_{n-1}^i
    \end{array}
\right.
$$

and a natural transformation $\Ld_k^n \rightarrow \underline{\Delta^n}$ (where $\underline{\Delta^n}$ denotes here the constant functor at $\Delta^n$ on $\mathrm{Facets}(\Lambda_k^n)$) defined by

$$
F_i \mapsto (\delta_n^i : \Delta^{n-1} \rightarrow \Delta^n).
$$

The latter factors through $\Lambda_k^n \rightarrow \Delta^n$ and we have the following result.
\end{Construction}

\begin{Proposition}\label{proposition: horn as colimit of faces}
    For every $n \ge 0$ and $0 \le k \le n$, the natural transformation $\Ld_k^n \rightarrow \underline{\Lambda_k^n}$ exhibits $\Lambda_k^n$ as a colimit of $\Ld_k^n$. \qed
\end{Proposition}

Let us now fix an integer $n \ge 2$, and an integer $0 < k < n$. We introduce a topological subcategory of $F(\Delta^n)$ that we will prove below to be equivalent to $F(\Lambda_k^n)$ (lemma \ref{lemma: image of an inner horn by F}).

\begin{Definition}\label{definition: image of horn by F}
    Let $\Lj_k^n$ be the topological subcategory of $F(\Delta^n)$ with the same objects and with morphisms generated by the morphisms that are in the image of at least one of the functors $F(\delta_n^i)$ for $0 \le i \le n$ and $i \neq k$. The topology on the morphism spaces is the topology induced by the topology on the morphism spaces of $F(\Delta^n)$.
\end{Definition}

Before describing $\Lj_k^n$ through the homeomorphisms \eqref{homeomorphisms}, we introduce some notation.

\begin{Notation}\label{notation: hollow cube}
    For every $m \ge 0$ and every $1 \le l \le m$, we denote by $\sqcap_l^m$ (resp. $\sqcup_l^m$) the subspace of $I^m$ obtained by removing the interior of $I^m$ together with the interior of the face $I^{l-1} \times \{0\} \times I^{m-l}$ (resp. $I^{l-1} \times \{1\} \times I^{m-l}$). These are called \textit{hollow} $m$\textit{-cubes} \footnote{Note that our notation for the hollow cubes is the opposite of that used in \cite[Section 2.4.5]{kerodon}. This comes from the fact that the simplicial path categories $\mathrm{Path}[m]_{\bullet}$ of \cite{kerodon} are defined in the same way as the simplicial categories $\Cr([m])$, except that the order on the sets $P_{m,i,j}$ from definition \ref{definition: simplicially enriched categories associated with simplices} is the opposite to the one of \cite[Notation 2.4.3.1]{kerodon}.}.
\end{Notation}

Using this notation we have the following lemma.

\begin{Lemmas}\label{lemma: morphism spaces of image of horn by F}
The morphism spaces of $\Lj_k^n$ are described as follows:
    $$
\Lj_k^n(i,j) = \left\{
    \begin{array}{ll}
        F(\Delta^n)(i,j) & \mbox{if } (i,j) \neq (0,n) \\
        \sqcap_k^{n-1} & \mbox{if } (i,j)=(0,n).
    \end{array}
\right.
$$
\end{Lemmas}

\begin{proof}
The case $(i,j) \neq (0,n)$ follows from the fact that $F(\delta_n^0)$ and $F(\delta_n^n)$ induce bijections between morphism spaces, and the case $(i,j) = (0,n)$ follows from propositions \ref{proposition: image of the faces by F} and $\ref{proposition: compositions in F(n)}$.
\end{proof}

By construction, the natural transformation $F \circ \Ld_k^n \rightarrow \underline{\Delta^n}$ factors through $\underline{\Lj_k^n} \rightarrow \underline{\Delta^n}$. We obtain a natural transformation $T : F \circ \Ld_k^n \rightarrow \underline{\Lj_k^n}$ and we have

\begin{Lemmas}\label{lemma: image of an inner horn by F}
    The natural transformation $T$ exhibits $\Lj_k^n$ as a colimit of $F \circ \Ld_k^n$. Hence, since $F$ preserves colimits, we have $F(\Lambda_k^n) = \Lj_k^n$.
\end{Lemmas}

\begin{proof}
    Let $\Dj$ be a topological category together with a natural transformation $T' : F \circ \Ld_k^n \rightarrow \underline{\Dj}$. We wish to prove that there exists a unique functor between topological categories $\Gj : \Lj_k^n \rightarrow \Dj$ such that $\underline{\Gj} \circ T = T'$. 
    
    Uniqueness follows from the fact that every object of $\Lj_k^n$ is in the image of $T$, and every morphism of $\Lj_k^n$ is a composition of morphisms in the image of $T$. Let us prove the existence.

    We define the objects $X_i = \Gj(i)$ as follows: $X_0, \hdots, X_{n-1}$ are the images of the objects $0, \hdots, n-1$ of $F \circ \Ld_k^n(F_n)$ by $T'(F_n)$, and $X_n$ is the image of the object $n-1$ of $F \circ \Ld_k^n(F_0)$ by $T'(F_0)$. Note that since $T'$ is a natural transformation, the following diagram commutes
    
    $$\xymatrix@C=1.7cm@R=0.5cm{
    & F \circ \Ld_k^n (F_0) \ar[rd]^-{T'(F_0)} \\
    F \circ \Ld_k^n(F_{0,n}) \ar[ru]^-{F \circ \Ld_k^n(\phi_{0,n}^0)} \ar[rd]_-{F \circ \Ld_k^n(\phi_{0,n}^n)} && \Dj \\
    & F \circ \Ld_k^n (F_n) \ar[ru]_-{T'(F_n)}
    }$$
    
    and therefore we also have that $X_i = T'(F_0)(i-1)$ for $i > 0$.
    
    Let $i \le j$ with $(i,j) \neq (0,n)$. If $j \neq n$ we define $\Gj(i,j) : F(\Delta^n)(i,j) \rightarrow \Dj(X_i,X_j)$ to be $T'(F_0)(i,j) : F \circ \Ld_k^n(F_0)(i,j) \rightarrow \Dj(X_i,X_j)$, and if $i \neq 0$ we define it to be $T'(F_n)(i-1,j-1) : F \circ \Ld_k^n(F_n)(i-1,j-1) \rightarrow \Dj(X_i,X_j)$. The commutativity of the diagram above shows that the two definitions coincide in the cases $i \neq 0$ and $j \neq n$. Moreover, $\Gj$ is continuous on the morphism spaces, and functorial. It remains to define $\Gj(0,n)$. Let $\gamma \in \Lj_k^n(0,n)$ and $\gamma = \gamma_k \circ \hdots \gamma_0$ be the writing of $\gamma$ as a composition of unbroken morphisms. If $k \ge 1$, we define $\Gj(\gamma) = \Gj(\gamma_k) \circ \hdots \circ \Gj(\gamma_0)$. If $k = 0$ then there exists $0 < i < n$ and $\gamma' \in F(\Delta^{n-1})(0,n-1)$ such that $\gamma = F(\delta_n^i)(\gamma')$. We define $\Gj(\gamma) = T'(F_i)(\gamma')$. To see that this definition of $\Gj(\gamma)$ does not depend on $i$, assume that there exists $0 < j < n$ with $i \neq j$ and $\gamma'' \in F(\Delta^{n-1})(0,n-1)$ such that $\gamma = F(\delta_n^j)(\gamma'')$. We may assume that $i < j$. Then, by proposition \ref{proposition: image of the faces by F}, $\gamma'$ and $\gamma''$ are the images of an element $\gamma''' \in F(\Delta^{n-2})(0,n-2)$ by $F \circ \Ld_k^n(\phi_{i,j}^i)$ and $F \circ \Ld_k^n(\phi_{i,j}^j)$ respectively. Since $T'$ is a natural transformation, the following diagram commutes

    $$\xymatrix@C=1.7cm@R=0.5cm{
    & F \circ \Ld_k^n (F_i) \ar[rd]^-{T'(F_i)} \\
    F \circ \Ld_k^n(F_{i,j}) \ar[ru]^-{F \circ \Ld_k^n(\phi_{i,j}^i)} \ar[rd]_-{F \circ \Ld_k^n(\phi_{i,j}^j)} && \Dj, \\
    & F \circ \Ld_k^n (F_j) \ar[ru]_-{T'(F_j)}
    }$$

    which allows to conclude that this definition of $\Gj(\gamma)$ does not depend on $i$.
    
    Then, by proposition \ref{proposition: image of the faces by F} and $\ref{proposition: compositions in F(n)}$, $\Gj$ is continuous on each of the subspaces $(\{t_i = 0\})_{0 \le i \le n, \, i \neq k}$ and $(\{t_i = 1\})_{0 \le i \le n}$ of $\Lj_k^n(0,n)$, and therefore continuous on $\Lj_k^n(0,n)$. It is functorial over $\Lj_k^n$ and we have $\underline{\Gj} \circ T = T'$.
\end{proof}

\begin{Remark}\label{remark: description of inner horns Kerodon}
    Lemma \ref{lemma: image of an inner horn by F} can also be proved using the description of the simplicial category $\Cr(\Lambda_k^n)$ given in \cite[Proposition 2.4.5.8]{kerodon}.
\end{Remark}

We now proceed with the proof of proposition \ref{proposition: S_M is an infinity category}. Consider a lifting problem of the form

$$
\xymatrix{
\Lambda_k^n \ar[r] \ar[d] & S_{\Ml}. \\
\Delta^n \ar@{-->}[ru]
}
$$

We denote by $\varphi$ the associated functor $F(\Lambda_k^n) \rightarrow \Ml$ (obtained using the adjunction between $F$ and $\Nl$). We also denote $a_i = \varphi(i)$; this defines a sequence $[a_0 \le \hdots \le a_n ] \in \Delta_A$. By the description of $F(\Lambda_k^n)$ supplied by lemma \ref{lemma: image of an inner horn by F}, there is an associated map $\varphi(0,n) : \sqcap_k^{n-1} \rightarrow \Ml(a_0,a_n)$. By corollary \ref{corollary: condition of being unbroken in terms of faces} and lemma \ref{lemma: image of an inner horn by F}, a solution to the lifting problem is the same as a solution to the following lifting problem in $\Top$:

$$
\xymatrix{
\sqcap_k^{n-1} \ar[r]^-{\varphi(0,n)} \ar[d] & \Ml(a_0,a_n) \\
I^{n-1} \ar@{-->}[ru]
}
$$

such that $I^{n-1} \backslash \sqcap_k^{n-1}$ is mapped to the interior of $\Ml(a_0,a_n)$.

This latter condition is conveniently formulated by introducing a stratification of $\Ml(a_0,a_n)$ by the poset $0 < 1$. It is defined by declaring $\Ml(a_0,a_n)_0$ to be the subspace of broken trajectories \footnote{Notice that it can be empty.}. This is indeed a stratification because, by proposition \ref{proposition: structure of smooth manifold with corners on space of broken traj}, $\Ml(a_0,a_n)$ admits a structure of smooth manifold with corners, with boundary $\Ml(a_0,a_n)_0$. We also endow $\sqcap_k^{n-1}$ with the stratification by $0 < 1$ induced by that on $\Ml(a_0,a_n)$ and the map $\varphi(0,n)$. Finally, we endow $I^{n-1}$ with the stratification by $0 < 1$ defined by declaring $(I^{n-1})_0 = (\sqcap_k^{n-1})_0$. The solution we are looking for is then the same as a solution to the following lifting problem in $\Top_{0 < 1}$

$$
\xymatrix{
\sqcap_k^{n-1} \ar[r]^-{\varphi(0,n)} \ar[d] & \Ml(a_0,a_n). \\
I^{n-1} \ar@{-->}[ru]
}
$$

Let $p = \max \{0 \le i \le n \mid a_i = a_0\}$ and $q = \min \{0 \le i \le n \mid a_i = a_n\}$. The following lemma provides a description of the $0$-stratum of $\sqcap_k^{n-1}$.

\begin{Lemmas}\label{lemma: 0-stratum of the hollow cube}
    The $0$-stratum of $\sqcap_k^{n-1}$ is the union

    $$
    (\sqcap_k^{n-1})_0 = \bigcup_{p < i < q} I^{i-1} \times \{1\} \times I^{n-1-i}.
    $$
\end{Lemmas}

\begin{proof}

By definition, the $0$-stratum of $\sqcap_k^{n-1}$ is the subspace of those morphisms that are mapped to a broken trajectory by $\varphi$. By definition of $p$ and $q$, given $p < i < q$, all morphisms in $F(\Lambda_k^n)(0,i)$ and $F(\Lambda_k^n)(i,n)$ are mapped to a nonconstant trajectory by $\varphi$. Therefore, by proposition \ref{proposition: compositions in F(n)}, every $\gamma \in I^{i-1} \times \{1\} \times I^{n-1-i}$ is mapped to a broken trajectory by $\varphi$. Hence

$$
\bigcup_{p < i < q} I^{i-1} \times \{1\} \times I^{n-1-i} \subseteq (\sqcap_k^{n-1})_0.
$$

Conversely, let $\gamma \in \sqcap_k^{n-1}$ be a morphism lying in the interior of the face  $I^{i-1} \times \{0\} \times I^{n-1-i}$ for some $1 \le i \le n-1$ and $i \neq k$. We wish to prove that $\gamma \notin (\sqcap_k^{n-1})_0$. By proposition \ref{proposition: image of the faces by F}, $\gamma$ is the image of an unbroken morphism of $F(\Delta^{n-1})$ by $F(\delta_n^i)$. Since the $i$-th face of $\Lambda_k^n \rightarrow S_{\Ml}$ is unbroken, $\gamma$ is mapped to an unbroken trajectory by $\varphi$, by proposition \ref{proposition: simplex unbroken iff image of unbroken is unbroken}.

Let $1 \le i \le n-1$ be such that $i \le p$ or $i \ge q$, and let $\gamma$ be a morphism lying in the interior of the face $I^{i-1} \times \{1\} \times I^{n-1-i}$. We wish to prove that $\gamma \notin (\sqcap_k^{n-1})_0$. By proposition \ref{proposition: compositions in F(n)}, $\gamma$ is the concatenation of an unbroken morphism in $F(\Lambda_k^n)(0,i)$ with an unbroken morphism in $F(\Lambda_k^n)(i,n)$. Since the $0$-th and $n$-th faces of $\Lambda_k^n \rightarrow S_{\Ml}$ are unbroken, by proposition \ref{proposition: simplex unbroken iff image of unbroken is unbroken} we find that $\gamma$ is mapped by $\varphi$ to the concatenation of an unbroken trajectory with a constant trajectory (or, if $p = q$, of two constant trajectories). In particular, it is unbroken. We finally conclude that

\begin{equation*}
(\sqcap_k^{n-1})_0 = \bigcup_{p < i < q} I^{i-1} \times \{1\} \times I^{n-1-i}. \qedhere
\end{equation*}

\end{proof}

The key properties that we are going to use in order to extend $\varphi(0,n)$ to $I^{n-1}$ in $\Top_{0 < 1}$ are the following two lemmas, whose proofs we postpone to the end of this section.

\begin{Lemmas}\label{lemma: regular neighborhood of the 0-stratum of the hollow cube}
    Let $m \ge 1$, $1 \le i \le j \le m$ and

    $$
    E = \bigcup_{i \le l \le j} I^{l-1} \times \{1\} \times I^{m-l} \subset I^m.
    $$

    Denote by $\partial E$ the boundary of $E$ in $\partial I^m$. There exists a neighborhood $U$ of $E$ in $I^m$ and a homeomorphism between $(U,E)$ and $(E \times [0,1],E \times \{0\})$ that identifies $\partial E \times (0,1]$ with the intersection of $U$ with the complement of $E$ in $\partial I^m$.
\end{Lemmas}

\begin{Lemmas}\label{lemma: evaluation at 0 is a trivial fibration}
    The evaluation $\ev_0 : \Exit(\Ml(a_0,a_n)_0,\Ml(a_0,a_n)_1) \rightarrow \Ml(a_0,a_n)_0$ at the starting point is a Serre fibration and a weak homotopy equivalence. 
\end{Lemmas}

We are now ready to complete the proof of proposition \ref{proposition: S_M is an infinity category}.

\begin{proof}[Proof of proposition \ref{proposition: S_M is an infinity category}]

Suppose first that $\Ml(a_0,a_n)_0$ is empty, i.e., according to proposition \ref{lemma: 0-stratum of the hollow cube}, that $q \le p+1$. Then, a solution as we are looking for is the same as a solution to the following lifting problem in $\Top$

$$
\xymatrix{
\sqcap_k^{n-1} \ar[r]^-{\varphi(0,n)} \ar[d] & \Ml(a_0,a_n). \\
I^{n-1} \ar@{-->}[ru]
}
$$

The existence of the latter follows from the fact that $I^{n-1}$ retracts onto $\sqcap_k^{n-1}$.

We assume from now on that $\Ml(a_0,a_n)_0$ is nonempty, i.e., that $q > p+1$. We prove in three steps that there exists an extension of $\varphi(0,n)$ to $I^{n-1}$ in $\Top_{0 < 1}$.

\textit{First step:} we prove that there exists an extension to a neighborhood of $(\sqcap_k^{n-1})_0$ in $I^{k-1} \times \{0\} \times I^{n-1-k}$. Let us for a moment denote $I^{k-1} \times \{0\} \times I^{n-1-k}$ by $I^{n-2}$. We endow $I^{n-2}$ with the stratification by $0 < 1$ induced by that on $I^{n-1}$. By lemma \ref{lemma: regular neighborhood of the 0-stratum of the hollow cube}, $(I^{n-2})_0$ is included in the boundary of $I^{n-2}$. We denote by $\partial(I^{n-2})_0$ the boundary of $(I^{n-2})_0$ in $\partial I^{n-2}$. Using lemmas \ref{lemma: 0-stratum of the hollow cube} and \ref{lemma: regular neighborhood of the 0-stratum of the hollow cube}, we reduce the existence of the desired extension to the existence of a solution to a lifting problem in $\Top$ of the form

$$
\xymatrix{
\partial(I^{n-2})_0 \ar[r] \ar[d] & \Exit(\Ml(a_0,a_n)_0,\Ml(a_0,a_n)_1) \ar[d]^{\ev_0} \\
(I^{n-2})_0 \ar@{-->}[ru] \ar[r] & \Ml(a_0,a_n)_0.
}
$$

By lemma \ref{lemma: evaluation at 0 is a trivial fibration}, the right vertical arrow is a trivial fibration. The left vertical arrow is a relative cell complex, in particular a cofibration. Since trivial fibrations have the right lifting property with respect to cofibrations, we conclude that such a solution exists.

\textit{Second step:} we prove that there exists an extension to a neighborhood of $(\sqcap_k^{n-1})_0$ in $I^{n-1}$. We denote by $\partial(\sqcap_k^{n-1})_0$ the boundary of $(\sqcap_k^{n-1})_0$ in $\partial I^{n-1}$. Using lemmas \ref{lemma: 0-stratum of the hollow cube} and \ref{lemma: regular neighborhood of the 0-stratum of the hollow cube} combined with the first step, we reduce the existence of such an extension to the existence of a solution to a lifting problem in $\Top$ of the form

$$
\xymatrix{
\partial(\sqcap_k^{n-1})_0 \ar[r] \ar[d] & \Exit(\Ml(a_0,a_n)_0,\Ml(a_0,a_n)_1) \ar[d]^{\ev_0} \\
(\sqcap_k^{n-1})_0 \ar@{-->}[ru] \ar[r] & \Ml(a_0,a_n)_0.
}
$$

Again, by lemma \ref{lemma: evaluation at 0 is a trivial fibration}, the right vertical arrow is a trivial fibration. The left vertical arrow is a relative cell complex, in particular a cofibration. Since trivial fibrations have the right lifting property with respect to cofibrations, we conclude that such a solution exists.

\textit{Third step:} we prove that there exists an extension to the whole of $I^{n-1}$. By the second step, there exists $0 < \varepsilon < 1$ and an extension to

$$
\bigcup_{p \le l \le q} I^{l-1} \times [\varepsilon,1] \times I^{n-l}. 
$$

Thus, identifying $I^p \times [0,\varepsilon]^{q-p-1} \times I^{n-q}$ with $I^{n-1}$ using the identification

$$
\fonctionsansnom{[0,\varepsilon]}{[0,1]}{t}{\frac{t}{\varepsilon}}
$$

we reduce the third step to the existence of a solution to a lifting problem in $\Top$ of the form

$$
\xymatrix{
\sqcap_k^{n-1} \ar[r] \ar[d] & \Ml(a_0,a_n)_1. \\
I^{n-1} \ar@{-->}[ru]
}
$$

The existence of the latter follows from the fact that $I^{n-1}$ retracts onto $\sqcap_k^{n-1}$.
\end{proof}

We close this section with the proofs of lemmas \ref{lemma: regular neighborhood of the 0-stratum of the hollow cube} and \ref{lemma: evaluation at 0 is a trivial fibration}.

\begin{proof}[Proof of lemma \ref{lemma: regular neighborhood of the 0-stratum of the hollow cube}]
    We begin by reducing to the case where $i=1$ and $j=m$. Assume this is true in that case for every $m$, and consider the general case as in the statement of the lemma. Let $m' = j-i+1$ and

    $$E' = \bigcup_{1 \le l \le m'} I^{l-1} \times \{1\} \times I^{m'-l} \subset I^{m'}.$$

    By assumption, there exists a neighborhood $U'$ of $E'$ in $I^{m'}$ and a homeomorphism between $(U',E')$ and $(E' \times [0,1],E'\times \{0\})$ satisfying the desired condition. Since $E = I^{i-1} \times E' \times I^{m-j}$, the subspace $U = I^{i-1} \times U' \times I^{m-j}$ is a neighborhood of $E$ in $I^m$ and the product of this homeomorphism with the identity maps of $I^{i-1}$ and $I^{m-j}$ defines a homeomorphism between $(U,E)$ and $(E \times [0,1],E \times \{0\})$. Moreover the boundary of $E$ in $\partial I^m$ is the intersection of $E$ with the subspaces $I^{l-1} \times \{1\} \times I^{m-l}$ for $1 \le l \le m$ and $l \notin \{i,\hdots,j\}$, and $I^{l-1} \times \{0\} \times I^{m-l}$ for $1 \le l \le m$. Since, in addition, the first homeomorphism satisfies the desired condition, we conclude that the second does as well.

    It remains to treat the case where $i = 1$ and $j = m$. We introduce the open neighborhood $V$ of $E$ in $I^m$ defined as

    $$
    V = \bigcup_{1 \le l \le m} I^{l-1} \times ]1/2,1] \times I^{m-l}.
    $$

    The closure of $V$ in $I^m$ is

    $$
    \overline{V} = \bigcup_{1 \le l \le m} I^{l-1} \times [1/2,1] \times I^{m-l}.
    $$

    It suffices to prove that there exists a homeomorphism between $(V,E)$ and $(E \times [0,1[,E \times \{0\})$ satisfying the desired condition. Indeed, restricting such a homeomorphism to the preimage of $E \times [0,1/2]$ yields a neighborhood of $E$ and a homeomorphism as desired.

    Such a homeomorphism between $(V,E)$ and $(E \times [0,1[,E \times \{0\})$ can be constructed as follows. For every $\tu = (t_1,\hdots,t_m) \ \in E$, we let $J = \{1 \le i \le m \mid t_i \ge 1/2 \}$ (note that $J \neq \varnothing$ for all $\tu$). We consider the segment connecting $\tu$ to $\underline{t'} = (t'_1,\hdots,t'_m)$ defined as $t'_i = 1/2$ for $ i \in J$ and $t'_i = t_i$ for $i \notin J$ (see figure \ref{figure: neighborhood of E} below for an illustration in the case $m=2$). This segment admits a unique affine parametrization

    $$
    \Gamma(\tu,.) : [0,1] \rightarrow \overline{V}
    $$

    such that $\Gamma(\tu,0) = \tu$ and $\Gamma(\tu,1) = \underline{t'}$.
    
    We now argue that the map $\Gamma = E \times [0,1] \rightarrow \overline{V}$ restricts to a homeomorphism between $E \times [0,1[$ and $V$. To justify this, we prove that for every $\underline{e} \in V$, there exists a unique $\tu \in E$ such that $\underline{e}$ lies on the segment $\Gamma(\tu,.)$ (see figure \ref{figure: neighborhood of E} for an illustration). Assume such a $\tu$ exists, then let $\underline{t'} = \Gamma(\tu,1)$ and $J = \{1 \le i \le m \mid e_i \ge 1/2 \}$. We must have $t'_i = 1/2$ for $i \in J$ and $t'_i = e_i$ for $i \notin J$. Conversely, define $\underline{t'}$ like this. Since $\underline{e} \in V$, we have $\underline{e} \neq \underline{t'}$.

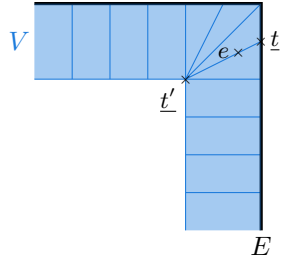
\begin{figure}[H]
    \centering
    \begin{tikzpicture}
    
    \coordinate (A) at (-1,2);
    \coordinate (B) at (2,2);
    \coordinate (C) at (2,-1);
    \coordinate (D) at (1,-1);
    \coordinate (E) at (1,1);
    \coordinate (F) at (-1,1);
    
    \coordinate (B1) at (-0.5,1);
    \coordinate (H1) at (-0.5,2);
    \coordinate (B2) at (0,1);
    \coordinate (H2) at (0,2);
    \coordinate (B3) at (0.5,1);
    \coordinate (H3) at (0.5,2);
    \coordinate (H4) at (1,2);
    \coordinate (H5) at (1.5,2);
    \coordinate (H6) at (2,1.5);
    \coordinate (H7) at (2,1);
    \coordinate (H8) at (2,0.5);
    \coordinate (H9) at (2,0);
    \coordinate (H10) at (2,-0.5);
    \coordinate (B8) at (1,0.5);
    \coordinate (B9) at (1,0);
    \coordinate (B10) at (1,-0.5);
    
    
    \draw[blue] (B1)--(H1);
    \draw[blue] (B2)--(H2);
    \draw[blue] (B3)--(H3);
    \draw[blue] (E)--(H4);
    \draw[blue] (E)--(H5);
    \draw[blue] (E)--(B);
    \draw[blue] (E)--(H6);
    \draw[blue] (E)--(H7);
    \draw[blue] (B8)--(H8);
    \draw[blue] (B9)--(H9);
    \draw[blue] (B10)--(H10);
    
    
    \draw[blue] (D)--(E)--(F);
    \draw[line width=1.2pt] (A)--(B)--(C);
    \fill[blue, opacity=0.4] (A)--(B)--(C)--(D)--(E)--(F)--cycle;
    
    \node at (2,-1.2) {$E$};
    \node[blue] at (-1.2,1.5) {$V$};
    
    
    \path (E)--(H6) coordinate[pos=0.7] (e);
    \node at (e) {\tiny{$\times$}};
    \node[anchor=east] at (e) {\small{$e$}};
    \node at (E) {\tiny{$\times$}};
    \node[anchor=north east] at (E) {\small{$\underline{t'}$}};
    \node at (H6) {\tiny{$\times$}};
    \node[anchor=west] at (H6) {\small{$\tu$}};
    
    \end{tikzpicture}
    \caption{The case $m=2$.}
    \label{figure: neighborhood of E}
\end{figure}
    
    Moreover, the intersection points of the line connecting $\underline{t'}$ to $\underline{e}$ with $E$ correspond to the $\lambda \in \R$ such that $\underline{t'} + \lambda(\underline{e}-\underline{t'}) \in E$. We have $\underline{t'} + \lambda (\underline{e}-\underline{t'}) \in E$ if and only if there exists $1 \le i \le m$ such that $(\underline{t'} + \lambda (\underline{e}-\underline{t'}))_i = 1$ and for every $1 \le j \le m$ we have $(\underline{t'} + \lambda (\underline{e}-\underline{t'}))_j \le 1$. Since for every $i \notin J$ we have $(\underline{t'} + \lambda (\underline{e}-\underline{t'}))_i = e_i < 1/2$, we must have $i \in J$. Given $i \in J$, there exists a unique $\lambda$ such that $(\underline{t'} + \lambda (\underline{e}-\underline{t'}))_i = 1$, which we denote by $\lambda_i$ and which is given by
    
    $$\lambda_i = \frac{1}{2e_i-1}.$$
    
    Together with the condition that for every $1 \le j \le m$ we have $(\underline{t'} + \lambda (\underline{e}-\underline{t'}))_j \le 1$, this implies that there exists a unique $\lambda$ such that  $\underline{t'} + \lambda(\underline{e}-\underline{t'}) \in E$, given by
    
    $$
    \lambda = \underset{i \in J}{\min} \, \lambda_i = \frac{1}{2(\underset{i \in J}{\max} \, e_i)  -1}.
    $$

    The associated intersection point is the desired $\tu$. This concludes the proof. \qedhere
\end{proof}

\begin{proof}[Proof of Lemma \ref{lemma: evaluation at 0 is a trivial fibration}]

    We first argue that it suffices to prove the statement by replacing $\Ml(a_0,a_n)$ with a neighborhood of $\Ml(a_0,a_n)_0$ in $\Ml(a_0,a_n)$. Indeed, assuming that there exists such a neighborhood, we wish to prove that for every $m \ge 0$, every lifting problem in $\Top$ of the form

    $$
    \xymatrix{
    \partial I^m \ar[r]^-b \ar[d] & \Exit(\Ml(a_0,a_n)_0,\Ml(a_0,a_n)_1) \ar[d]^{\ev_0} \\
    I^m \ar@{-->}[ru] \ar[r]_-B & \Ml(a_0,a_n)_0
    }
    $$

    has a solution. To this commutative square corresponds a map $\sqcup_{m+1}^{m + 1} \rightarrow \Ml(a_0,a_n)$ defined to be $B$ in restriction to $I^m \times \{0\}$, and defined as

    $$
    (t_1,\hdots,t_{m+1}) \mapsto b(t_1,\hdots,t_m)(t_{m+1})
    $$

    in restriction to $\partial I^m \times I$.

    Let us define stratifications of $I^{m+1}$ and $\sqcup_{m+1}^{m+1}$ by $0 < 1$, by declaring $(I^{m+1})_0$ and $(\sqcup_{m+1}^{m+1})_0$ to be equal to $I^m \times \{0\}$. Then the map $\sqcup_{m+1}^{m + 1} \rightarrow \Ml(a_0,a_n)$ is compatible with stratifications, and a solution to the above lifting problem is the same as a solution to the following lifting problem in $\Top_{0 < 1}$

    $$
    \xymatrix{
    \sqcup_{m+1}^{m+1} \ar[r] \ar[d] & \Ml(a_0,a_n). \\
    I^{m+1} \ar@{-->}[ru]
    }
    $$

    Now, by virtue of our assumption that there exists a neighborhood of $\Ml(a_0,a_n)$ on which the evaluation at $0$ from the $|0 < 1|_{0 < 1}$-homotopy link to $\Ml(a_0,a_n)$ is a trivial fibration, such a dotted arrow exists in restriction to $I^m \times [0, \varepsilon]$ for some $0 < \varepsilon \le 1$. Thus, identifying $I^m \times [\varepsilon,1]$ with $I^{m+1}$ using the identification
    
    $$
    \fonctionsansnom{[\varepsilon,1]}{[0,1]}{t}{\frac{t-\varepsilon}{1-\varepsilon}}
    $$

    we reduce the question of the existence of such a solution to that of a solution to the following lifting problem in $\Top$

    $$
    \xymatrix{
    \sqcup_{m+1}^{m+1} \ar[r] \ar[d] & \Ml(a_0,a_n)_1. \\
    I^{m+1} \ar@{-->}[ru]
    }
    $$
    
    The existence of the latter follows from the fact that $I^{m+1}$ retracts onto $\sqcup_{m+1}^{m+1}$.
    
    It remains to prove the statement by replacing $\Ml(a_0,a_n)$ with a neighborhood of $\Ml(a_0,a_n)_0$ in $\Ml(a_0,a_n)$. By proposition \ref{proposition: structure of smooth manifold with corners on space of broken traj}, $\Ml(a_0,a_n)$ admits a structure of compact smooth manifold with corners, with boundary $\Ml(a_0,a_n)_0$. Hence, $(\Ml(a_0,a_n),\Ml(a_0,a_n)_0)$ is homeomorphic to a smooth manifold with boundary. It therefore suffices to treat the case of a compact smooth manifold with boundary $M$ endowed with the stratification by $0 < 1$ defined by $M_0 = \partial M$.

    There exists a neighborhood $U$ of $\partial M$ in $M$ such that $(U,\partial M)$ is homeomorphic to $(\partial M \times [0,1], \partial M \times \{0\})$. It therefore suffices to treat the case of $Y = \partial M \times [0,1]$ endowed with the stratification by $0 < 1$ defined by $Y_0 = \partial M \times \{0\}$. To simplify the notation, we denote the space of exit paths of $Y$ from $Y_0$ to $Y_1$ by $\Hl$. Given a path $\gamma$ in $Y$, we will denote by $\gamma_1$ the projection of $\gamma$ to $\partial M$ and by $\gamma_2$ the projection of $\gamma$ to $[0,1]$.

    We first show that $\ev_0 : \Hl \rightarrow \partial M$ is a homotopy equivalence. There is a homotopy inverse of $\ev_0$ given by

    $$
    \fonction{\varphi}{\partial M}{\Hl}{m}{(t \mapsto (m,t)).}
    $$

    We clearly have $\ev_0 \circ \varphi = \mathrm{Id}$. There is a homotopy between $\varphi \circ \ev_0$ and $\mathrm{Id}$ given by

    $$
    \fonction{H}{\Hl \times [0,1]}{\Hl}{(\gamma,t)}{(s \mapsto (\gamma_1(st),(1-t)s + t\gamma_2(s))).}
    $$

    To prove that $\ev_0$ is a Serre fibration one could argue as follows: since $Y$ is conically stratified by construction, this is a consequence of theorem \ref{theorem: stratified simplicial set of conically stratified space is infty cat} and proposition \ref{proposition: Sing_A infinity category in terms of starting point evaluation} combined. We however give a complete proof in this special case. We want to show that every lifting problem of the form

    $$
    \xymatrix{
    I^{m-1} \ar[r]^h \ar[d]_i & \Hl \ar[d]^{\ev_0} \\
    I^m \ar@{-->}[ru]^l \ar[r]_g & \partial M
    }
    $$

    has a solution $l$ (where $i$ is the inclusion $I^{m-1} \simeq I^{m-1} \times \{0\} \hookrightarrow I^m$). We denote the elements of $I^{m-1}$ by $\tu$, and the elements of $I^m$ by $(\tu,t_m)$. Let us define

    $$
    \fonction{l}{I^m}{\Hl}{(\tu,t_m)}{\left( s \mapsto \left\{
    \begin{array}{cc}
        \left( g(\tu,t_m-s),\frac{s}{t_m}h(\tu)_2(s) \right) & \text{for } s < t_m \\
         (h(\tu)_1(s-t_m),h(\tu)_2(s)) & \text{for } s \ge t_m. 
    \end{array}
    \right. \right)}
    $$

    For $t_m > 0$ and $s = t_m$, the first expression is well-defined and gives $(g(\tu,0),h(\tu)_2(t_m))$, and the second expression gives $(h(\tu)_1(0),h(\tu)_2(t_m)))$. By commutativity of the square above, the two expressions agree in that case, hence $l$ is continuous. Moreover we have $l \circ i = h$ and $\ev_0 \circ l = g$, which concludes the proof.
\end{proof}

\subsection{\texorpdfstring{The equivalence between the homotopy-and flow coherent nerves of $\Ml$}{The equivalence between the flow and homotopy coherent nerves of M}}\label{section: equivalence between flow and homotopy coherent nerve}

The goal of this section is to prove the following proposition.
    
\begin{Proposition}\label{proposition: equivalence between flow and homotopy coherent nerve}
    The functor induced by the inclusion $i_{\Ml} : S_{\Ml} \rightarrow \Nl(\Ml)$ is an equivalence of $\infty$-categories.
\end{Proposition}

Our strategy is to apply the characterization of equivalences of $\infty$-categories in terms of essential surjectivity and full faithfulness (theorem \ref{theorem: full faithfulness and essential surjectivity}), and reduce the full faithfulness condition to the result that for every pair of critical points, the inclusion of the space of unbroken trajectories connecting these, into the space of broken ones, is a homotopy equivalence (corollary \ref{Corollary: space of broken trajectories equivalent to space of trajectories}).

We will use a model for the morphism spaces of an $\infty$-category which is different from that of definition \ref{definition: morphism space}. We start by introducing it, following \cite[Section 1.2.2]{HigherTopos}.

\begin{Construction}\label{construction: functor K}

    For every arrow $\varphi : \Delta^{n} \rightarrow \Delta^{m}$ in $\Delta$, there is a corresponding arrow $\Delta^{n+1} \rightarrow \Delta^{m+1}$ defined to be $\varphi$ in restriction to $\{0,\hdots,n\}$, and carrying $n+1$ to $m+1$. Consider the simplicial sets $\Delta^{n+1} / \Delta^n$ and $\Delta^{m+1} / \Delta^m$ obtained by collapsing the faces opposite to the vertices $n+1$ and $m+1$ respectively (an explicit description of these simplicial sets follows from the general description of colimits of simplicial sets given in section \ref{section: Left Kan extensions}). This arrow passes to the quotient and determines an arrow $\Delta^{n+1} / \Delta^n \rightarrow \Delta^{m+1} / \Delta^m$. This turns the assignment $\Delta^n \mapsto \Delta^{n+1} / \Delta^n$ into a functor $\Delta \rightarrow \sSet$. We denote this functor by $K$. Furthermore, for every $\Delta^n \in \Delta$, the simplicial set $K(\Delta^n)$ has two vertices: one corresponding to the collapsed face of $\Delta^{n+1}$, and one corresponding to the vertex $n+1$ of $\Delta^{n+1}$. We denote these by $x_n$ and $y_n$ respectively. Finally, we denote by $\Sing_K(-)$ the right adjoint of the left Kan extension of $K$ to $\sSet$. By definition, for every simplicial set $T$ we have
    
    $$
    \Sing_K(T)_n = \Hom(K(\Delta^n),T).
    $$
    
\end{Construction}

\begin{Definition}\label{definition: right morphism space}
    Let $T \in \sSet$ be an $\infty$-category and $x, y$ two objects of $T$. The \textit{right morphism space} of $T$ between $x$ and $y$, denoted $\Hom_{T}^R(x,y)$, is the simplicial subset of $\Sing_K(T)$ formed by those morphisms $K(\Delta^n) \rightarrow T$ sending $x_n$ to $x$ and $y_n$ to $y$. In particular, for every nonnegative integer $n$, $\Hom_{T}^R(x,y)_n$ is the set of those morphisms $\Delta^{n+1} \rightarrow T$ that are constant equal to $x$ on the face opposite to the vertex $n+1$, and that carry the vertex $n+1$ to $y$.
\end{Definition}

\begin{Proposition}\label{proposition: right morphism space is Kan complex}
    For every $x, y \in T$, $\Hom_T^R(x,y)$ is a Kan complex.
\end{Proposition}

\begin{proof}
    See \cite[Proposition 1.2.2.3]{HigherTopos}.
\end{proof}

Our next step is to state a comparison result between the notion of right morphism space, and the notion of morphism space from definition \ref{definition: morphism space}. Recall that the $n$-simplices of $\Hom_T(x,y)$ correspond to certain maps $\Delta^n \times \Delta^1 \rightarrow T$, while the $n$-simplices of $\Hom_T^R(x,y)$ correspond to certain maps $\Delta^{n+1} \rightarrow T$. The simplicial set $\Delta^n \times \Delta^1$ is naturally identified with the nerve of the poset $[n] \times [1]$, while the simplicial set $\Delta^{n+1}$ is naturally identified with the nerve of the poset $[n+1]$. There is a morphism of simplicial sets $\Delta^n \times \Delta^1 \rightarrow \Delta^{n+1}$ defined to be the nerve of the morphism of posets $[n] \times [1] \rightarrow [n+1]$ defined as $(i,0) \mapsto i$ and $(i,1) \mapsto n+1$ for every $0 \le i \le n$. These morphisms are natural in $n$ and thus give rise to a morphism of simplicial sets $\Hom_T^R(x,y) \rightarrow \Hom_T(x,y)$.

\begin{Proposition} \label{proposition: right morphism spaces are model for morphism spaces}
    This morphism of simplicial sets is natural with respect to $T$. Moreover, if $T$ is an $\infty$-category, this is a homotopy equivalence.
\end{Proposition}

\begin{proof}
    Naturality follows from the construction. For the rest, see \cite[Corollary 4.2.1.8]{HigherTopos}.
\end{proof}

\begin{Corollary}\label{corollary: full faithfulness and right morphism spaces}

    Let $F : \Cl \rightarrow \Dl$ be a functor between $\infty$-categories. The functor $F$ is fully faithful if and only if for every pair of objects $x,y \in \Cl$, the induced morphism of Kan complexes
    
    $$
    \Hom_{\Cl}^R(x,y) \rightarrow \Hom_{\Dl}^R(F(x),F(y))
    $$
    
    is a homotopy equivalence. \qed

\end{Corollary}

Our next goal is to better understand the right morphism spaces of the $\infty$-category $\Nl(\Ml)$. We closely follow \cite[Section 2.2.2]{HigherTopos}.

The description of these right morphism spaces will involve the following construction.

\begin{Construction}\label{construction: functor Q}

    For every integer $n \ge 0$, we denote by $\Ql(\Delta^n)$ the topological space $F(K(\Delta^n))(x_n,y_n)$ (recall that $F$ is the Leitch-Cordier functor and that $K(\Delta^n)$, $x_n$ and $y_n$ where introduced in construction \ref{construction: functor K}). Note that for every arrow $\varphi : \Delta^n \rightarrow \Delta^m$ in $\Delta$ one has $F(K(\Delta^n))(\varphi)(x_n) = x_m$ and $F(K(\Delta^n))(\varphi)(y_n) = y_m$. We therefore have a map $F(K(\Delta^n))(\varphi) : \Ql(\Delta^n) \rightarrow \Ql(\Delta^m)$. This turns the assignment $\Delta^n \mapsto \Ql(\Delta^n)$ into a functor $\Ql : \Delta \mapsto \Top$. We denote by $\Sing_{\Ql}$ the right adjoint of the left Kan extension of this functor to $\sSet$. Recall that for every topological space $Y$, the set of $n$-simplices of $\Sing_{\Ql}(Y)$ is $\Hom(\Ql(\Delta^n),Y)$.

\end{Construction}

Fix a nonnegative integer $n$ and two critical points $a,b \in \Ml$. By definition, an $n$-simplex of $\Hom^R_{\Nl(\Ml)}(a,b)$ is a morphism of simplicial sets

$$
K(\Delta^n) \rightarrow \Nl(\Ml)
$$

sending $x_n$ to $a$ and $y_n$ to $b$. By adjunction, this is the same as a functor of topological categories

$$
F(K(\Delta^n)) \rightarrow \Ml
$$

sending $x_n$ to $a$ and $y_n$ to $b$. In order to proceed with our analysis, we will need the following description of the topological category $F(K(\Delta^n))$.

\begin{Proposition}\label{proposition: description of F(K(n))}
    The topological category $F(K(\Delta^n))$ is described as follows.
    
    \begin{itemize}
    
        \item It has two objects $x_n$ and $y_n$.
        
        \item The endomorphism spaces of $x_n$ and $y_n$ are reduced to the corresponding identity morphisms.
        
        \item The morphism space between $x_n$ and $y_n$ is the quotient of $F(\Delta^{n+1})(0,n+1)$ by the following equivalence relation: $\gamma \circ \delta \sim \gamma \circ \delta'$ for every $0 \le i \le n+1$, $\delta, \delta' \in F(\Delta^{n+1})(0,i)$ and $\gamma \in F(\Delta^{n+1})(i,n+1)$. Equivalently, identifying $F(\Delta^{n+1})(0,n+1)$ with $I^n$, this is the quotient of $I^n$ by the following equivalence relation: $\tu \sim \tu$ if and only if there exists $0 \le i \le n$ such that $t_i = t_i' = 1$ and $(t_i,\hdots,t_n) = (t'_i,\hdots,t'_n)$.
        
        \item The morphism space between $y_n$ and $x_n$ is empty.
    \end{itemize}
    
\end{Proposition}

\begin{proof}

    Since the functor $F$ commutes with colimits, $F(K(\Delta^n))$ is canonically isomorphic to the quotient of topological categories $F(\Delta^{n+1}) / F(\Delta^n)$ induced by the inclusion $\Delta^n \rightarrow \Delta^{n+1}$ of the face opposite to the vertex $n+1$.
    
    Denote by $\Cj$ the topological category described in the statement of the proposition. There is a functor $F(\Delta^{n+1}) \rightarrow \Cj$ defined as follows.
    
    \begin{itemize}
        
        \item Each of the objects $i$ between $0$ and $n$ is carried to $x_n$, and the object $n+1$ is carried to $y_n$.
        
        \item For every $0 \le i < n+1$, an element $\gamma \in F(\Delta^{n+1})(i,n+1)$ is carried to the equivalence class of the element $\gamma \circ \delta$ in $\Cj(x_n,y_n)$, for any choice of element $\delta \in F(\Delta^{n+1})(0,i)$.
        
    \end{itemize}
    
    We leave it to the reader to check that this functor exhibits $\Cj$ as the quotient $F(K(\Delta^n))$. \qedhere

\end{proof}

\begin{Remark}\label{remark: resemblance between C_A and F(K)}
    Note the resemblance between proposition \ref{proposition: description of F(K(n))} and lemma \ref{lemma: restriction of equivalence relation to critical level}. We will study the connection between the functor $C_A$ and the functor $\Ql$ in section \ref{section: comparison between S M and Sing C A (X)}.
\end{Remark}

We now complete our analysis of the $n$-simplices of $\Hom_{\Nl(\Ml)}^R(a,b)$. As a consequence of proposition \ref{proposition: description of F(K(n))}, we obtain that such an $n$-simplex is the same as a continuous map $\Ql(\Delta^n) \rightarrow \Ml(a,b)$. These identifications are functorial in $n$; we therefore have proven:

\begin{Proposition}\label{proposition: right morphism spaces of coherent nerve}
    There is a natural isomorphism of simplicial sets
    
    $$
    \Hom_{\Nl(\Ml)}^R(a,b) \simeq \Sing_{\Ql}(\Ml(a,b)).
    $$
\end{Proposition}

\begin{Remark}\label{remark: homotopy equivalence of morphism spaces btw cat and nerve}

Note that so far, the discussion is valid after replacing $\Ml$ by any topological category. There is natural transformation of functors $\Ql \rightarrow |-|$ that induces, for every topological space $Y$, a homotopy equivalence of Kan complexes $\Sing(Y) \rightarrow \Sing_{\Ql}(Y)$; this is discussed in \cite[Section 2.2.2]{HigherTopos}. In particular, for every topological category $\Cj$ and every pair of objects $(X,Y)$ of $\Cj$ we get a natural homotopy equivalence between $\Sing(\Hom_{\Cj}(X,Y))$ and $\Hom_{\Nl(\Cj)}(X,Y)$. This is the natural homotopy equivalence referred to in proposition \ref{proposition: morphism spaces of a topological category and its coherent nerve}.

\end{Remark}

In the special case of $\Ml$, we additionally have the following result.

\begin{Proposition}\label{proposition: right morphism spaces of flow coherent nerve}
    The isomorphism of simplicial sets $\Hom_{\Nl(\Ml)}^R(a,b) \simeq \Sing_{\Ql}(\Ml(a,b))$ from proposition \ref{proposition: right morphism spaces of coherent nerve} restricts to an isomorphism of simplicial sets
    
    $$
    \Hom_{S_{\Ml}}^R(a,b) \simeq \Sing_{\Ql}(\Mi(a,b)).
    $$
\end{Proposition}

\begin{proof}

    Let $n$ be a nonnegative integer and $\sigma$ an $n$-simplex of the right morphism space of $\Nl(\Ml)$ between $a$ and $b$. Let us regard $\sigma$ as a $(n+1)$-simplex of $\Nl(\Ml)$, and equivalently, as a functor of topological categories $F(\Delta^{n+1})(0,n+1) \rightarrow \Ml$. The condition that $\sigma$ is an $n$-simplex of the right morphism space of $S_{\Ml}$ between $a$ and $b$, is equivalent to the condition that $\sigma$ is unbroken. By proposition \ref{proposition: unbroken iff compatible with stratifications on morph spaces}, this is also equivalent to the condition that $\sigma (0,n+1)$ has image contained in $\Mi(a,b)$, in other words, that the map $\Ql(\Delta^n) \rightarrow \Ml(a,b)$ associated with $\sigma$ has image contained in $\Mi(a,b)$. The latter condition is the condition that the element of $\Sing_{\Ql}(\Ml(a,b))$ associated with $\sigma$ belongs to the simplicial subset $\Sing_{\Ql}(\Mi(a,b))$, as desired. \qedhere

\end{proof}

We can summarize our discussion of this section so far as follows: the map of right morphism spaces between $a$ and $b$ induced by the functor $S_{\Ml} \rightarrow \Nl(\Ml)$ can be identified with the inclusion between Kan complexes

$$
\Sing_{\Ql}(\Mi(a,b)) \rightarrow \Sing_{\Ql}(\Ml(a,b)).
$$

In order to complete the proof of proposition \ref{proposition: equivalence between flow and homotopy coherent nerve}, we will need the following result.

\begin{Proposition}\label{proposition: Sing_Q detects weak equivalences}
    For every topological space $Y$, the simplicial set $\Sing_{\Ql}(Y)$ is a Kan complex. Furthermore, for every continuous map of topological spaces $g : Y \rightarrow Z$, the induced map
    
    $$
    \Sing_{\Ql}(Y) \rightarrow \Sing_{\Ql}(Z)
    $$
    
    is a homotopy equivalence if and only if $g$ is a weak homotopy equivalence.
\end{Proposition}

\begin{proof}

    We denote by $|-|_{\Ql}$ the left Kan extension of $\Ql$ to $\sSet$, which is also the left adjoint of $\Sing_{\Ql}$.
    
    Recall from section \ref{section: definition of the homotopy coherent nerve} that the Leitch-Cordier functor is obtained by applying the geometric realization functor, at the level of simplicial categories, to a functor denoted $\Cr : \sSet \rightarrow \Cat_{\sSet}$. In particular, the functor $\Ql$ is obtained by applying the geometric realization functor, at the level of simplicial sets, to the functor $Q : \Delta \rightarrow \sSet$ defined by $Q(\Delta^n) = \Cr(K(\Delta^n))(x_n,y_n)$. Denote by $|-|_Q$ the left Kan extension of $Q$ to $\sSet$ and by $\Sing_{Q}$ its right adjoint. By \cite[Proposition 2.2.2.9]{HigherTopos}, the adjunction
    
    $$
    \sSet
    \underset{\Sing_{Q}}{\overset{|-|_Q}{\rightleftarrows}}
    \sSet
    $$
    
    is a Quillen adjunction (where $\sSet$ is endowed with the Kan model structure). But observe that we have $|-|_{\Ql} = |-| \circ |-|_Q$ and $\Sing_{\Ql} = \Sing_Q \circ \Sing$. The adjunction between $|-|$ and $\Sing$ is also a Quillen adjunction (where $\Top$ is endowed with the classical model structure). Since a composition of Quillen adjunctions is a Quillen adjunction, we obtain that the adjunction between $|-|_{\Ql}$ and $\Sing_{\Ql}$ is Quillen. The first part of the proposition then follows from the fact that the right adjoint of a Quillen adjunction preserves fibrant objects. The second part follows from the fact that the right adjoint of a Quillen adjunction preserves weak equivalences between fibrant objects. \qedhere

\end{proof}

We are now ready to conclude this section.

\begin{proof}[Proof of proposition \ref{proposition: equivalence between flow and homotopy coherent nerve}]

    We wish to prove that this functor is essentially surjective and fully faithful. Essential surjectivity follows from the fact that this functor induces the identity on the set of objects. By corollary \ref{corollary: full faithfulness and right morphism spaces} combined with proposition \ref{proposition: right morphism spaces of flow coherent nerve}, the full faithfulness condition is equivalent to the condition that for every pair of objects $a,b \in \Ml$, the 
    morphism of Kan complexes
    
    $$
    \Sing_{\Ql}(\Mi(a,b)) \rightarrow \Sing_{\Ql}(\Ml(a,b))
    $$
    
    is a homotopy equivalence. This follows from corollary \ref{Corollary: space of broken trajectories equivalent to space of trajectories} combined with proposition \ref{proposition: Sing_Q detects weak equivalences}. \qedhere

\end{proof}

\subsection{\texorpdfstring{Filling inner horns in the simplicial set of stratified cubes of $X$}{Filling inner horns in the simplicial set of stratified cubes of X}}\label{section: filling inner horns in Sing C A (X)}

In the previous two sections, we proved that the morphism of simplicial sets $S_{\Ml} \rightarrow \Nl(\Ml)$ is an equivalence of $\infty$-categories. Our goal in the next two sections is to obtain an analogous result concerning the morphism of simplicial sets $s_{\Ml}^{(X,A)} : S_{\Ml} \rightarrow \Sing_{C_A}(X)$.

Recall that in section \ref{section: comparison between realizations}, we introduced a modification $|-|_A'$ of the functor $|-|_A$, and denoted its right adjoint by $\Sing_A'(-)$ (construction \ref{construction: functor |-|_A'} and notation \ref{notation: left Kan extension of |-|_A' and its right adjoint}). In this section, we prove for any $A$-stratified topological space $Y$, the simplicial set $\Sing_A'(X)$ is an $\infty$-category as soon as $\Sing_A(Y)$ is one (proposition \ref{proposition: if Sing_A infinity category then Sing_A' too}). Combining this with the comparison theorem between $|-|_A'$ and $C_A$ proven in section \ref{section: comparison between realizations} (theorem \ref{theorem: isomorphism between |-|'_A and C_A}), we deduce that for any $A$-stratified topological space $Y$, $\Sing_{C_A}(Y)$ is an $\infty$-category as soon as $\Sing_A(Y)$ is one (proposition \ref{proposition: Sing_C_A is an infinity category}). 

To prove these statements, we will check condition (i) in definition \ref{definition: infinity category}, and therefore consider morphisms from inner horns into $\Sing_A'(Y)$ and $\Sing_{C_A}(Y)$,

$$
\Lambda_k^n \rightarrow \Sing_A'(Y) \qquad \text{and} \qquad \Lambda_k^n \rightarrow \Sing_{C_A}(Y).
$$

We will see that, by adjunction, these correspond to stratified maps

$$
|\Lambda_k^{\ag}|_A' \rightarrow Y \qquad \text{and} \qquad |\Lambda_k^{\ag}|_{C_A} \rightarrow Y
$$

for some sequence $\ag$ of length $n$ (recall from notation \ref{notation: stratifications of simplex, boundary and horns} that we denote by $\Lambda_k^{\ag}$ the restriction to $\Lambda_k^n$ of the stratification on $\Delta^n$ induced by the sequence $\ag$). We will therefore need to better understand the stratified spaces $|\Lambda_k^{\ag}|_A'$ and $|\Lambda_k^{\ag}|_{C_A}$. Since the functors $|-|_A'$ and $|-|_{C_A}$ commute with colimits, we start by giving a description of stratified horns as colimits of their faces, generalizing construction \ref{construction: horn as colimit of simplices} and proposition \ref{proposition: horn as colimit of faces} to the stratified setting. This is analogous to the generalization of construction \ref{construction: boundary of simplex as colimit of simplices} and proposition \ref{proposition: boundary of standard simplex as colimit of its faces} to the stratified setting provided by proposition \ref{proposition: boundary of stratified simplex as colimit of faces}. Let $n$ be a nonnegative integer, $0 \le k \le n$ and $\ag$ a finite increasing sequence of elements of $A$ of length $n$.

\begin{Proposition}\label{proposition: stratified horn as colimit of faces}
    There exists a unique lift of the functor $\Ld_k^n$ to a functor $\Ld_k^{\ag} : \mathrm{Facets}(\Delta^n) \rightarrow \sSet_A$ such that the natural transformation $\Ld_k^n \rightarrow \underline{\Delta^n}$ lifts to a natural transformation $\Ld_k^{\ag} \rightarrow \underline{\Delta^{\ag}}$. Moreover, the latter factors through $\Lambda_k^{\ag}$ and the natural transformation $\Ld_k^{\ag} \rightarrow \underline{\Lambda_k^{\ag}}$ exhibits $\Lambda_k^{\ag}$ as a colimit of $\Ld_k^{\ag}$.
\end{Proposition}

\begin{proof}
The first part of the proposition follows from the construction of $\Ld_k^n$ and the natural transformation $\Ld_k^n \rightarrow \underline{\Delta^n}$. The second part follows from proposition \ref{proposition: horn as colimit of faces} and the fact that the forgetful functor $\sSet_A \rightarrow \sSet$ commutes with colimits (proposition \ref{proposition: stratified simplicial sets admit small limits and colimits}).
\end{proof}

\begin{Remark}\label{remark: image of stratified horn is colimit of image of faces}
    Since $|-|_A, |-|_A',|-|_{C_A} : \sSet_A \rightarrow \Top_A$ commute with colimits, it follows from proposition \ref{proposition: stratified horn as colimit of faces} that $|\Lambda_k^{\ag}|_A$ (resp. $|\Lambda_k^{\ag}|_A'$, $|\Lambda_k^{\ag}|_{C_A}$) is equal to $\colim (|-|_A \circ \Ld_k^{\ag})$ (resp. $\colim (|-|_A' \circ \Ld_k^{\ag})$, $\colim (|-|_{C_A} \circ \Ld_k^{\ag})$). In other words, $|\Lambda_k^{\ag}|_A$ (resp. $|\Lambda_k^{\ag}|_A'$, $|\Lambda_k^{\ag}|_{C_A}$) is obtained by gluing together the faces of $|\Delta^{\ag}|_A$ (resp. $|\Delta^{\ag}|_A'$, $C_A(\Delta^{\ag})$) different from the one opposite to the vertex $k$, along their common faces. 
\end{Remark}

We are now ready to prove that the simplicial sets $\Sing_A'(X)$ and $\Sing_{C_A}(X)$ are $\infty$-categories.

\begin{Proposition}\label{proposition: if Sing_A infinity category then Sing_A' too}
    Let $Y$ be an $A$-stratified topological space such that $\Sing_A(Y)$ is an $\infty$-category. Then the simplicial set $\Sing_A'(Y)$ introduced in notation \ref{notation: left Kan extension of |-|_A' and its right adjoint} is an $\infty$-category. In particular, $\Sing_A'(X)$ is an $\infty$-category.
\end{Proposition}

\begin{proof}
    We prove that $\Sing_A'(Y)$ satisfies the inner horn filling condition. Let $k,n$ be two integers such that $0 < k <n$ and consider a lifting problem

    $$
    \xymatrix{
    \Lambda_k^n \ar[r] \ar[d] & \Sing_A'(Y). \\
    \Delta^n \ar@{-->}[ru]
    }
    $$

    Since $N(A)$ is the nerve of a category, and therefore an $\infty$-category, the stratification on $\Lambda_k^n$ induced by that on $\Sing_A'(Y)$ extends to a stratification of $\Delta^n$. The latter corresponds to some sequence $\ag$ of length $n$. The lifting problem is then equivalent to its adjoint lifting problem in $\Top_A$:

    $$
    \xymatrix{
    |\Lambda_k^{\ag}|'_A \ar[r] \ar[d] & Y. \\
    |\Delta^{\ag}|_A' \ar@{-->}[ru]
    }
    $$
    
    Assume first that $\Lambda_k^{\ag}$ consists of one stratum only. This is equivalent to the sequence $\ag$ being constant, hence also to $\Delta^{\ag}$ consisting of one stratum only. Therefore $|\Delta^{\ag}|_A'$ is a point in that case. On the other hand, the induced stratification on every $(n-1)$-simplex of $\Lambda_k^n$ also consists of one stratum only. Consequently, by remark \ref{remark: image of stratified horn is colimit of image of faces}, $|\Lambda_k^{\ag}|_A'$ is the colimit of a diagram whose objects are points, and therefore is also a point. It follows that the lifting problem has a (unique) solution.
    
    Assume now that $\Lambda_k^{\ag}$ consists of at least two strata. The natural transformation $|-|_A \rightarrow |-|'_A$ gives rise to a lifting problem in $\Top_A$

    $$
    \xymatrix{
    |\Lambda_k^{\ag}|_A \ar[r] \ar[d] & Y \\
    |\Delta^{\ag}|_A \ar@{-->}_l[ru]
    }
    $$

    which has a solution by virtue of the assumption that $\Sing_A(Y)$ is an $\infty$-category. Let us denote by $l$ such a solution. To prove that the lifting problem we started with has a solution, it suffices to prove that $l$ is compatible with the equivalence relation on $|\Delta^{\ag}|_A$ that defines $|\Delta^{\ag}|_A'$ (construction \ref{construction: functor |-|_A'}). But observe that, since the stratification on $\Delta^{\ag}$ consists of at least two strata, two equivalent elements of $|\Delta^{\ag}|_A$ either lie on the face opposite to the vertex $0$ or on the face opposite to the vertex $n$ (in the sense of definition \ref{definition: face of image of simplex by three functors}). Since $0 < k < n$ and since, in restriction to $|\Lambda_k^{\ag}|_A$, $l$ factors through $|\Lambda_k^{\ag}|_A'$, the images of two such elements by $l$ are equal by remark \ref{remark: image of stratified horn is colimit of image of faces}. This concludes the proof.
\end{proof}

\begin{Proposition}\label{proposition: Sing_C_A is an infinity category}
    Let $Y$ be an $A$-stratified topological space such that $\Sing_A(Y)$ is an $\infty$-category. Then the simplicial set $\Sing_{C_A}(Y)$ is an $\infty$-category. In particular, $\Sing_{C_A}(X)$ is an $\infty$-category.
\end{Proposition}

\begin{proof}
    We prove that $\Sing_{C_A}(Y)$ satisfies the inner horn filling condition. Let $n \ge 0$, $0 < k < n$ and consider a lifting problem

    $$
    \xymatrix{
    \Lambda_k^n \ar[d] \ar[r] & \Sing_{C_A}(Y). \\
    \Delta^n \ar@{-->}[ru]
    }
    $$

    Since $N(A)$ is the nerve of a category, and therefore an $\infty$-category, the stratification on $\Lambda_k^n$ induced by that on $\Sing_{C_A}(Y)$ extends to a stratification of $\Delta^n$. The latter corresponds to some sequence $\ag$ of length $n$. The lifting problem is then equivalent to its adjoint lifting problem in $\Top_A$:

    $$
    \xymatrix{
    |\Lambda_k^{\ag}|_{C_A} \ar[d] \ar[r] & Y. \\
    C_A(\Delta^{\ag}) \ar@{-->}[ru]
    }
    $$

    By theorem \ref{theorem: isomorphism between |-|'_A and C_A} and remark \ref{remark: image of stratified horn is colimit of image of faces}, there exists a commutative diagram in $\Top_A$
    
    $$
    \xymatrix{
    |\Lambda_k^{\ag}|_A' \ar[d] \ar[r]^-{\simeq} & |\Lambda_k^{\ag}|_{C_A} \ar[d] \\
    |\Delta^{\ag}|_A' \ar[r]^-{\simeq} & C_A(\Delta^{\ag})
    }
    $$
    
    where the horizontal arrows are stratified homeomorphisms and the vertical ones are induced by the inclusion $\Lambda_k^{\ag} \rightarrow \Delta^{\ag}$ in $\sSet_A$. In particular, the lifting problem is equivalent to a lifting problem in $\Top_A$ of the form

    $$
    \xymatrix{
    |\Lambda_k^{\ag}|_A' \ar[d] \ar[r] & Y, \\
    |\Delta^{\ag}|_A' \ar@{-->}[ru]
    }
    $$

    and therefore, by adjunction, to a lifting problem in $\sSet_A$ of the form

    $$
    \xymatrix{
    \Lambda_k^n \ar[d] \ar[r] & \Sing_A'(Y). \\
    \Delta^n \ar@{-->}[ru]
    }
    $$

    Combining proposition \ref{proposition: if Sing_A infinity category then Sing_A' too} with the assumption that $\Sing_A(Y)$ is an $\infty$-category, we obtain that the latter has a solution. This finishes the proof.
\end{proof}

\subsection{\texorpdfstring{The equivalence between the flow coherent nerve of $\Ml$ and the restricted simplicial set of stratified cubes of $X$}{Comparison between the flow coherent nerve of M and the simplicial set of stratified cubes in X}}\label{section: comparison between S M and Sing C A (X)}

In the previous section, we proved that the simplicial sets $\Sing_{C_A}(X)$ and $\Sing_A'(X)$ are $\infty$-categories. On the other hand, theorem \ref{theorem: isomorphism between |-|'_A and C_A} immediately implies the following:

\begin{Proposition}\label{proposition: Sing_C_A and Sing_A are isomorphic as semi-simplicial sets}
    There exists an isomorphism between the semi-simplicial sets underlying the simplicial sets $\Sing_{C_A}(X)$ and $\Sing_A'(X)$, which is the identity on the set of objects (identified with $X$). \qed
\end{Proposition}

In view of proposition \ref{proposition: Sing_C_A and Sing_A are isomorphic as semi-simplicial sets}, and the fact that there are morphisms of simplicial sets

$$
s_{\Ml}^{(X,A)} : S_{\Ml} \rightarrow \Sing_{C_A}(X) \qquad \text{and} \qquad \Sing_A'(X) \rightarrow \Sing_A(X),
$$

our plan is to use $\Sing_{C_A}(X)$ and $\Sing_A'(X)$ as intermediaries to compare $S_{\Ml}$ and $\Sing_A(X)$. In this section, we start by studying the functor $s_{\Ml}^{(X,A)} : S_{\Ml} \rightarrow \Sing_{C_A}(X)$. The main result is that it becomes an equivalence of $\infty$-categories after corestriction to some simplicial subset of the target (proposition \ref{proposition: S_M -> Sing_C_A,r is an equivalence}). We start by justifying the necessity of this corestriction.

Let $x,y$ be two points of $X$ that are isomorphic as objects of the $\infty$-category $\Sing_{C_A}(X)$.
Applying the functor $\Sing_{C_A}(X) \rightarrow A$, we see that the points $x$ and $y$ lie in the same stratum, and further, denoting this stratum by $X_a$, we see that $x$ and $y$ are isomorphic as objects of the $\infty$-category $\Sing_{C_A}(X)_a$. In addition, recall that the images by $C_A$ of elements of $\Delta_A$ that correspond to constant sequences, are points. Consequently the $n$-simplices of the $a$-stratum of $\Sing_{C_A}(X)$ correspond to the points of $X_a$. More precisely, there is an isomorphism of simplicial sets

$$
\Sing_{C_A}(X)_a \simeq \bigsqcup_{p \in X_a} \Delta^0 .
$$

Therefore we must have $x=y$. In other words, there are as many isomorphism classes of objects in $\Sing_{C_A}(X)$, as there are points in $X$. On the other hand, the isomorphism classes of objects of $S_{\Ml}$ correspond to the critical points of $f$. The $\infty$-categories $S_{\Ml}$ and $\Sing_{C_A}(X)$ are therefore not equivalent as soon as $\dim X \ge 1$. To remedy this, we corestrict ourselves to the following simplicial subset of $\Sing_{C_A}(X)$.

\begin{Definition}\label{definition: restricted version for C_A}
    We define the \textit{restricted version} of $\Sing_{C_A}(X)$ to be the simplicial subset of $\Sing_{C_A}(X)$ formed by those simplices whose vertices are critical points of $f$. We denote it by $\Sing_{C_A,r}(X)$.
\end{Definition}

The restricted version of $\Sing_{C_A}(X)$ satisfies the following two properties.

\begin{Proposition}
    The simplicial set $\Sing_{C_A,r}(X)$ is an $\infty$-category.
\end{Proposition}

\begin{proof}
    Since the property that the vertices are critical points of $f$ is preserved by inner horn filling in $\Sing_{C_A}(X)$, this follows from the fact that $\Sing_{C_A}(X)$ is an $\infty$-category (proposition \ref{proposition: Sing_C_A is an infinity category}).
\end{proof}

\begin{Proposition}\label{proposition: flow coherent nerve factors through restricted cubes}
    The morphism of simplicial sets $s_{\Ml}^{(X,A)} : S_{\Ml} \rightarrow \Sing_{C_A}(X)$ factors through the inclusion $\Sing_{C_A,r}(X) \rightarrow \Sing_{C_A}(X)$.
\end{Proposition}

\begin{proof}
This statement is equivalent to saying that $(S_{\Ml})_0 \rightarrow \Sing_{C_A}(X)_0$ has image contained in the set of critical points of $f$, which is indeed the case.
\end{proof}

The main result of this section is the following.

\begin{Proposition}\label{proposition: S_M -> Sing_C_A,r is an equivalence}
    The functor $s_{\Ml}^{(X,A)} : S_{\Ml} \rightarrow \Sing_{C_A,r}(X)$ is an equivalence of $\infty$-categories.
\end{Proposition}

As for the proof of proposition \ref{proposition: equivalence between flow and homotopy coherent nerve}, our strategy is to prove that this functor is essentially surjective and fully faithful. Again, we will use the right morphisms spaces, introduced in section \ref{section: equivalence between flow and homotopy coherent nerve}, as models for the morphism spaces of these two $\infty$-categories. This time, our plan is to reduce the full faithfulness condition to the condition that for every pair of critical points $(a,b)$, the inclusion, denoted $\iota_{a,b}$, of the space of unbroken trajectories between $a$ and $b$ into the space of exit paths between $a$ and $b$ (construction \ref{construction: parametrization of trajectories by values of the function}), is a weak homotopy equivalence (proposition \ref{proposition: equivalence between space of trajectories and space of exit paths}).

Fix two critical points $a,b \in \Ml$ such that $a < b$. Recall from definition \ref{definition: right morphism space} that the right morphism space of $\Sing_{C_A,r}(X)$ between $a$ and $b$ is a simplicial subset of $\Sing_K(\Sing_{C_A,r}(X))$. In this description, we forget the stratification on $\Sing_{C_A,r}(X)$ and regard it as a mere simplicial set.

We can actually be more precise in this situation. To this end, we introduce a lift of the functor $K$ from construction \ref{construction: functor K}, to $\sSet_A$.

\begin{Construction}\label{construction: lift of functor K}

    Let $\varphi : \Delta^n \rightarrow \Delta^m$ be a morphism in $\Delta$, and consider the sequence $[a = \hdots = a < b]$ of length $n+1$, as well as the sequence $[a = \hdots = a < b]$ of length $m+1$. These define $A$-stratifications of the simplicial sets $\Delta^{n+1}$ and $\Delta^{m+1}$ respectively. These stratifications pass to the quotients $K(\Delta^n)$ and $K(\Delta^m)$, and are compatible with the map $K(\varphi)$. We denote the resulting functor $\Delta \rightarrow \sSet_A$ by $K_{a,b}$. We denote by $\Sing_{K_{a,b}}$ the right adjoint of the left Kan extension of this functor to $\sSet$. Recall that for every $A$-stratified simplicial set $Y$, $\Sing_{K_{a,b}}(Y)$ is the simplicial set whose set of $n$-simplices is $\Hom(K_{a,b}(\Delta^n),Y)$.

\end{Construction}

The right morphism space of $\Sing_{C_A,r}(X)$ between $a$ and $b$ can be described as follows.

\begin{Proposition}\label{proposition: right morphism space of restricted version}
    
    The simplicial set $\Hom_{\Sing_{C_A,r}(X)}^R(a,b)$ is identified with the following simplicial subset
    
    $$
    \Sing_{K_{a,b}}(\Sing_{C_A,r}(X)) \subset \Sing_K(\Sing_{C_A,r}(X)).
    $$
    
\end{Proposition}

\begin{proof}

    By definition, an $n$-simplex $\sigma : K(\Delta^n) \rightarrow \Sing_{C_A,r}(X)$ of the right-hand side of this inclusion, belongs to this right morphism space if and only if it sends $x_n$ to $a$ and $y_n$ to $b$. At the level of $0$-simplices, the $a$ and $b$-strata of the restricted version of $\Sing_{C_A}(X)$ are reduced to $a$ and $b$ respectively. The latter condition is therefore equivalent to the condition that $\sigma$ sends $x_n$ to the $a$-stratum, and $y_n$ to the $b$-stratum. By remark \ref{remark: stratification determined by vertices}, this is equivalent to the condition that $\sigma$ is compatible with the $A$-stratification on $K_{a,b}(\Delta^n)$, as desired. \qedhere

\end{proof}

Our next goal is to use the functor $K_{a,b}$ to formulate the precise relationship between the functor $\Ql$ (constructed in \ref{construction: functor Q}) and the functor $C_A$. For this, consider the length $n+1$ sequence of the form $\ag = [a = \hdots = a < b]$. Recall from section \ref{section: lift to stratified spaces} that the image of $\Delta^{\ag}$ by $C_A$ can be regarded as the quotient of $I^n \times [f(b),f(a)]$ by the equivalence relation $\sim_{\ag}$, endowed with the $(a < b)$-stratification whose $a$-stratum is the image of $I^n \times \{f(a)\}$ in this quotient. Recall from section \ref{section: proof of the ball theorem for realizations} that for any given real number $f(b) < s < f(a)$, the equivalence relation $\sim_{\ag}$ restricts to an equivalence relation $\sim_{\ag,s}$ on $I^n$, and propositions \ref{lemma: restriction of equivalence relation to noncritical level} and \ref{proposition: description of F(K(n))} yield an identification $I^n / \sim_{\ag,s} = \Ql(\Delta^n)$. Moreover, making this identification, $C_A(\Delta^{\ag})$ is identified with the quotient of the topological space $\Ql(\Delta^n) \times [f(b),f(a)]$ obtained by collapsing each of the two subsets $\Ql(\Delta^n) \times \{f(a)\}$ and $\Ql(\Delta^n) \times \{f(b)\}$ to distinct points, endowed with the $(a < b )$-stratification whose $a$-stratum is the image of $\Ql(\Delta^n) \times \{f(a)\}$ in this quotient. In the sequel, we will denote this space by

$$
(\Ql(\Delta^n) \times [f(b),f(a)] \sqcup * \sqcup *') / (\Ql(\Delta^n) \times \{f(a)\} \sim *, \Ql(\Delta^n) \times \{f(b)\} \sim *').
$$

Observe also that the image of $\Delta^{\ag}$ by $C_A$ coincides with the image of $K_{a,b}(\Delta^n)$ by $C_A$. Indeed, $K_{a,b}(\Delta^n)$ is the quotient of $\Delta^{\ag}$ by the face opposite to the vertex $n+1$, which corresponds to the constant subsequence $[a = \hdots = a] \subset \ag$ of length $n$. Since the image of the latter by $C_A$ is reduced to a point, and $C_A$ commutes with colimits, the natural map $C_A(\Delta^{\ag}) \rightarrow C_A(K_{a,b}(\Delta^n))$ is an isomorphism.

Combining these observations, we obtain an isomorphism of functors

$$
C_A \circ K_{a,b} \simeq (\Ql \times [f(b),f(a)] \sqcup * \sqcup *') / (\Ql \times \{f(a)\} \sim *, \Ql \times \{f(b)\} \sim *')
$$

where the right-hand side is endowed with the $(a < b)$-stratification whose $a$-stratum is the image of $\Ql \times \{f(a)\}$ in the quotient, and with the functorial structure induced from that of $\Ql$ together with the identity on $[f(b),f(a)]$. This discussion leads to the following proposition.

\begin{Proposition}\label{proposition: right morphism space as space of exit paths}
    
    There exists a (strictly) commutative square of Kan complexes
    
    $$
    \xymatrix{
    \Hom_{\Sing_{C_A,r}(X)}^R(a,b) \ar[r]^-{\simeq} & \Sing_{\Ql}(\Exit(a,b)) \\
    \Hom_{S_{\Ml}}^R(a,b) \ar[u] \ar[r]_-{\simeq} & \Sing_{\Ql}(\Mi(a,b)) \ar[u]_-{\Sing_{\Ql}(\iota_{a,b})}
    }
    $$
    
    where the left vertical arrow is the one induced by the inclusion of simplicial sets $S_{\Ml} \rightarrow \Sing_{C_A,r}(X)$, the bottom horizontal arrow is the one provided by proposition \ref{proposition: right morphism spaces of coherent nerve}, and the symbol $\simeq$ denotes an isomorphism of simplicial sets.
    
\end{Proposition}

\begin{proof}

    Let us exhibit the top horizontal arrow. To begin with, proposition \ref{proposition: right morphism space of restricted version} provides an inclusion of simplicial sets
    
    $$
    \Hom_{\Sing_{C_A,r}(X)}^R(a,b) \subset \Sing_{C_A \circ K_{a,b}}(X),
    $$
    
    where $\Sing_{C_A \circ K_{a,b}}$ is the functor $\Top_A \rightarrow \sSet_A$ obtained as the right adjoint of the left Kan extension of $C_A \circ K_{a,b}$ to $\sSet_A$. This inclusion identifies the left-hand side with the simplicial subset of the right-hand side formed by those simplices whose vertices are $a$ and $b$. We denote it by $\Sing_{C_A \circ K_{a,b},r}(X)$.
    
    Furthermore, combining the affine decreasing bijection $[0,1] \rightarrow [f(b),f(a)]$ with the isomorphism of functors previously constructed, we obtain an isomorphism of functors
    
    $$
    C_A \circ K_{a,b} \simeq (\Ql \times [0,1] \sqcup * \sqcup *') / (\Ql \times \{0\} \sim *, \Ql \times \{1\} \sim *').
    $$
    
    Let $n$ be a nonnegative integer and endow $\Ql(\Delta^n) \times [0,1]$ with the $(a < b)$-stratification whose $a$-stratum is $\Ql(\Delta^n) \times \{0\}$. Combining the above observations, we obtain an identification
    
    $$
    \begin{aligned}
    \Hom_{\Sing_{C_A,r}(X)}^R(a,b)_n
    & = \{ c \in \Hom_{\Top_A}(\Ql(\Delta^n) \times [0,1], X) \mid c \text{ is constant equal to } a \text{ and } b \\
    & \qquad \text{on } \Ql(\Delta^n) \times \{0\} \text{ and } \Ql(\Delta^n) \times \{1\} \text{ respectively} \} \\
    & = \Hom_{\Top}(\Ql(\Delta^n),\Exit(a,b)) \\
    &= \Sing_{\Ql}(\Exit(a,b))_n.
    \end{aligned}
    $$

    These identifications are functorial in $n$, thus yielding the top horizontal isomorphism in the statement of the proposition. Since this isomorphism is constructed using the same homeomorphism $[f(b),f(a)] \simeq [0,1]$ as in the construction of the map $\iota_{a,b} : \Mi(a,b) \rightarrow \Exit(a,b)$, the resulting square diagram is strictly commutative. \qedhere
    
\end{proof}

We are now ready to collect all the ingredients to conclude this section.

\begin{proof}[Proof of proposition \ref{proposition: S_M -> Sing_C_A,r is an equivalence}]

    We prove that this functor is essentially surjective and fully faithful. Essential surjectivity follows from the fact that this functor induces the identity on the set of objects. For full faithfulness, consider two critical points $a,b \in A$. We wish to prove that the morphism of Kan complexes

    $$
    \Hom_{S_{\Ml}}^R(a,b) \rightarrow \Hom_{\Sing_{C_A,r}(X)}^R(a,b)
    $$
    
    is a homotopy equivalence.
    
    If $a \nleq b$, this follows from the fact that both Kan complexes are empty. If $a=b$, this follows from the fact that these are both isomorphic to $\Delta^0$. If $a < b$, then by proposition \ref{proposition: right morphism space as space of exit paths}, this is equivalent to the condition that the morphism of Kan complexes
    
    $$
    \Sing_{\Ql}(\iota_{a,b}) : \Sing_{\Ql}(\Mi(a,b)) \rightarrow \Sing_{\Ql}(\Exit(a,b))
    $$
    
    is a homotopy equivalence. This follows from proposition \ref{proposition: equivalence between space of trajectories and space of exit paths} combined with proposition \ref{proposition: Sing_Q detects weak equivalences}. \qedhere

\end{proof}

\newpage

\section{\texorpdfstring{Turning a semi-simplicial map into a functor of $\infty$-categories}{Turning a semi-simplicial map into a functor of infinity categories}}\label{section: turning a semi simplicial map into functor}

Let us summarize the current state of the proof of theorem \ref{Main theorem}. Recall that our ultimate goal is to obtain a zigzag of equivalences of $\infty$-categories of the following form:

\[
\Nl(\Ml) \leftarrow S_{\Ml} \rightarrow \Sing_{C_A,r}(X) \rightarrow \Sing'_{A,r}(X) \rightarrow \Sing_A(X).
\]

We have introduced in section \ref{section: lift to stratified spaces} a functor $C_A : \Delta_A \rightarrow \Top_A$, and, in section \ref{section: The flow coherent nerve of M}, a simplicial subset $S_{\Ml} \subset \Nl(\Ml)$ endowed with a morphism $S_{\Ml} \rightarrow \Sing_{C_A}(X)$. We have proven, in sections \ref{section: flow coherent nerve is an infinity category} and \ref{section: equivalence between flow and homotopy coherent nerve}, that $S_{\Ml}$ is an $\infty$-category and the inclusion $S_{\Ml} \rightarrow \Nl(\Ml)$ is an equivalence. In section \ref{section: comparison between realizations}, we have introduced a quotient $|-|_A'$ of the standard stratified geometric realization $|-|_A$, yielding a simplicial subset $\Sing_A'(X) \subset \Sing_A(X)$, and we have proven that the semi-simplicial sets underlying $\Sing_{C_A}(X)$ and $\Sing_A'(X)$ are isomorphic. In sections \ref{section: filling inner horns in Sing C A (X)} and \ref{section: comparison between S M and Sing C A (X)}, we have shown that $\Sing_{C_A}(X)$ is an $\infty$-category, we have introduced an $\infty$-subcategory $\Sing_{C_A,r}(X) \subset \Sing_{C_A}(X)$, called the restricted version, such that there is a factorization $S_{\Ml} \rightarrow \Sing_{C_A,r}(X)$, and we have proven the latter to be an equivalence. It therefore remains to compare the $\infty$-categories $\Sing_A(X)$ and $\Sing_{C_A,r}(X)$.

\subsection{\texorpdfstring{The restricted versions of $\Sing_A(X)$ and $\Sing_A'(X)$}{The restricted versions of Sing A(X) Sing A'(X)}}\label{section: restricted versions}

We start be introducing analogous restricted versions of the simplicial sets $\Sing_A'(X)$ and $\Sing_A(X)$.

\begin{Definition}\label{definition: restricted versions}
    The \textit{restricted version} of $\Sing_A(X)$ (resp. $\Sing_A'(X)$) is the simplicial subset of $\Sing_A(X)$ (resp. $\Sing_A'(X)$) formed by those simplices whose vertices are critical points of $f$. We denote it by $\Sing_{A,r}(X)$ (resp. $\Sing_{A,r}'(X)$).
\end{Definition}

\begin{Proposition}\label{proposition: restricted versions are infinity categories}
    The restricted versions of $\Sing_A(X)$ and $\Sing_A'(X)$ are $\infty$-categories.
\end{Proposition}

\begin{proof}

    Recall that $\Sing_A(X)$ is an $\infty$-category, and we proved in proposition \ref{proposition: if Sing_A infinity category then Sing_A' too} that $\Sing_A'(X)$ is one too. The desired result follows from the characterization of $\infty$-categories in terms of the inner horn filling property and the fact that the property of having critical points as vertices is preserved by inner horn filling. \qedhere

\end{proof}

The morphism of simplicial sets $\Sing_A'(X) \rightarrow \Sing_A(X)$ restricts to a morphism of simplicial sets $\Sing_{A,r}'(X) \rightarrow \Sing_{A,r}(X)$. Postcomposing the latter with the inclusion $\Sing_{A,r}(X) \rightarrow \Sing_A(X)$, we obtain a morphism of simplicial sets $\Sing_{A,r}'(X) \rightarrow \Sing_A(X)$.

\begin{Proposition}\label{proposition: restricted version into unrestricted is equivalence}
    
    The morphism of simplicial sets $\Sing_{A,r}'(X) \rightarrow \Sing_A(X)$ is an equivalence of $\infty$-categories.
    
\end{Proposition}

\begin{proof}

    We prove this by proving that each of the two morphisms of simplicial sets $\Sing_{A,r}'(X) \rightarrow \Sing_{A,r}(X)$ and $\Sing_{A,r}(X) \rightarrow \Sing_A(X)$ are equivalences of $\infty$-categories, by proving that these are essentially surjective and fully faithful.
    
    We start with the first one. In that case, essential surjectivity follows from the fact that this morphism is the identity at the level of objects (identified with the set of critical points of $f$). For full faithfulness, suppose given two critical points $a,b$. First, if $a \nleq b$, then the morphism spaces of both $\infty$-categories are empty. Second, if $a=b$, we have on the one hand $\Hom_{\Sing_{A,r}'(X)}^R(a,a)=\Delta^0$, and on the other hand
    
    $$
    \begin{aligned}
    \Hom_{\Sing_{A,r}(X)}^R(a,a) & = \Hom_{\Sing_A(X)}^R(a,a) \\
    & = \Hom_{\Sing(X_a)}^R(a,a)
    \end{aligned}
    $$
    
    and the latter is homotopy equivalent to $\Delta^0$, since $X_a$ is contractible. Lastly, if $a < b$, the inclusion $\Sing_{A,r}'(X) \rightarrow \Sing_{A,r}(X)$ induces an equality between the right morphism spaces between $a$ and $b$.
    
    We now treat the case of the second morphism of simplicial sets. In that case, full faithfulness follows from the fact that this morphism induces an equality between morphism spaces (in the sense of definition \ref{definition: morphism space}). As for essential surjectivity, this functor corresponds, at the level of objects, to the inclusion of the set of critical points of $f$ into the set of points of $X$. We therefore need to prove that every point $x \in X$ is isomorphic, as an object of the $\infty$-category $\Sing_A(X)$, to a critical point. Denote by $a$ the critical point such that $x \in X_a$. Since $X_a$ is connected, there exists a path in $X_a$ starting at $x$ and ending at $a$. The corresponding morphism of $\Sing_A(X)$ is an isomorphism between $x$ and $a$. This completes the proof. \qedhere

\end{proof}

\subsection{End of the proof}\label{section: end of the proof}

The goal of this section is to complete the proof of our main theorem as well as that of the corollaries from section \ref{section: Morse homology and constructible sheaves}.

In view of proposition \ref{proposition: restricted version into unrestricted is equivalence}, all that remains to be done is to compare the $\infty$-categories $\Sing_{C_A,r}(X)$ and $\Sing_{A,r}'(X)$. Now, theorem \ref{theorem: isomorphism between |-|'_A and C_A} guarantees the existence of an isomorphism of functors between $|-|_A'$ and $C_A$ in restriction to $\Delta_A^+$, and every such isomorphism determines an isomorphism of semi-simplicial sets between $\Sing_{C_A,r}(X)$ and $\Sing_{A,r}'(X)$.

Suppose given two $\infty$-categories $\Cl$ and $\Dl$ and a morphism of semi-simplicial sets $G : \Cl \rightarrow \Dl$. We can think of $G$ as assigning a morphism of $\Dl$ to every morphism of $\Cl$ in a way compatible with composition but not necessarily unital, in the sense that $G$ does not necessarily carry identity morphisms to identity morphisms. In order to complete our zigzag of equivalences of $\infty$-categories, we are led to the following question: when can this assignment be made unital?

The answer that we will use is that given in \cite{HiroFunctors}. To explain it, recall the category $\Delta^+ \subset \Delta$ and the category of semi-simplicial sets $\sSet^+ = \Pre(\Delta^+)$. The composition

$$
j : \Delta^+ \rightarrow \Delta \rightarrow \sSet
$$

induces a functor $j_! : \sSet^+ \rightarrow \sSet$ by left Kan extension. The right adjoint of $j_!$ is the functor $j^* : \sSet \rightarrow \sSet^+$ that associates, to a simplicial set, its underlying semi-simplicial set. Given an $\infty$-category $\Cl$ we denote by $\Cl_+$ the simplicial set $j_!(j^*(\Cl))$. By adjunction, the morphism $G$ determines a morphism of \emph{simplicial sets} $G_+ : \Cl_+ \rightarrow \Dl$. The unit of the adjunction also determines a morphism of simplicial sets $\Cl_+ \rightarrow \Cl$, and we have the following theorem (\cite[Theorem 1.4]{HiroFunctors}):

\begin{Thms}\label{theorem: Hiro}

    Suppose that for every $x \in \Cl$, $G(s_0^0(x))$ is an equivalence. Then there exists a functor $G' : \Cl \rightarrow \Dl$ such that the following diagram is commutative in the homotopy category of $\infty$-categories
    
    $$
    \xymatrix{
    \Cl_+ \ar[r]^-{G_+} \ar[d] & \Dl. \\
    \Cl \ar[ru]_-{G'}
    }
    $$
    
    Moreover, such a functor is unique up to natural equivalence.

\end{Thms}

\begin{Remark}\label{remark: C + not always infty category}

The simplicial set $\Cl_+$ is not always an $\infty$-category. In order to associate to it an $\infty$-category, one takes a fibrant replacement for the Joyal model structure. We refer the reader to \cite{HiroFunctors} for details.

\end{Remark}

Let us fix inverse isomorphisms of semi-simplicial sets $G : \Sing_{C_A,r}(X) \leftrightarrows \Sing_{A,r}'(X) : T$ as provided by theorem \ref{theorem: isomorphism between |-|'_A and C_A}. The hypotheses of theorem \ref{theorem: Hiro} are satisfied by these two morphisms, so this theorem provides two inverse equivalences of $\infty$-categories

$$
G' : \Sing_{C_A,r}(X) \leftrightarrows \Sing_{A,r}'(X) : T'.
$$

We now have a complete zigzag of equivalences of $\infty$-categories between $\Nl(\Ml)$ and $\Sing_A(X)$.

\begin{Remark}\label{remark: homotopic isomorphisms of functors give homotopic equivalences}

    By theorem \ref{theorem: isomorphism between |-|'_A and C_A}, any two isomorphisms from $|-|_A'$ to $C_A$ in restriction to $\Delta_A^+$ are \emph{homotopic}. Consequently, by \cite[Corollary 1.6]{HiroFunctors}, the two equivalences of $\infty$-categories between $\Sing_{C_A,r}(X)$ and $\Sing_{A,r}'(X)$ obtained from these two isomorphisms are naturally equivalent.

\end{Remark}

We now complete the proof of our main theorem in its precise form given in section \ref{section: precise formulation of the main theorem}.

\begin{proof}[Proof of the main theorem]

By Theorem \ref{theorem: Hiro}, for every critical point $a$, $G'(a)$ and $T'(a)$ are isomorphic to $a$, i.e., belong to the $a$-stratum. Up to replacing these functors by naturally equivalent ones, we can assume that $G'(a)=a$ and $T'(a)=a$ for every critical point $a$, which we do. Fix two critical points $a$ and $b$. At the level of the right morphism spaces between $a$ and $b$, our zigzag of equivalences induces a zigzag of Kan complexes

$$
\begin{aligned}
\Sing_{\Ql}(\Ml(a,b)) \xleftarrow{\Sing_{\Ql}(i)} \Sing_{\Ql}(\Mi(a,b)) \xrightarrow{\Sing_{\Ql}(\iota_{a,b})} \Sing_{\Ql}(\Exit(a,b)) \leftrightarrows &  \, \Hom^R_{\Sing_{A,r}'(X)}(a,b) \\
& = \Hom^R_{\Sing_A(X)}(a,b)
\end{aligned}
$$

where $i : \Mi(a,b) \rightarrow \Ml(a,b)$ is the inclusion. On the one hand, there is a natural transformation $\Ql \rightarrow |-|$ of functors $\Delta \rightarrow \Top$, that induces a homotopy equivalence of Kan complexes $\Sing(\Ml(a,b)) \rightarrow \Sing_{\Ql}(\Ml(a,b))$ (see remark \ref{remark: homotopy equivalence of morphism spaces btw cat and nerve}). On the other hand, there is a homotopy equivalence of Kan complexes

$$
\Hom^R_{\Sing_A(X)}(a,b) \rightarrow \Hom_{\Sing_A(X)}(a,b) = \Sing(\Exit(a,b)).
$$

We obtain from these homotopy equivalences and this zigzag a homotopy equivalence

$$
\Sing(\Ml(a,b)) \simeq \Sing(\Exit(a,b)).
$$

We want to show that the latter coincides up to homotopy with the homotopy equivalence induced by the zigzag

$$
\Sing(\Ml(a,b) \xleftarrow{\Sing(i)} \Sing(\Mi(a,b)) \xrightarrow{\Sing(\iota_{a,b})} \Sing(\Exit(a,b)).
$$

To this end, recall from construction \ref{construction: functor K} the functor $K : \Delta \rightarrow \sSet$ and for every positive integer $n$, the two vertices $x_n$ and $y_n$ of $K(\Delta^n)$. We denote by $\Kl : \Delta \rightarrow \Top$ the functor obtained by composing $K$ with the geometric realization functor. We note that there is an isomorphism of functors

$$
\Kl \simeq (|-| \times [0,1] \sqcup *_1 \sqcup *_2)/(|-| \times \{0\} \sim *_1, |-| \times \{1\} \sim *_2)
$$

that carries $x_n$ to $0$ and $y_n$ to $1$ for every positive integer $n$. This isomorphism of functors determines an isomorphism of simplicial sets

$$
\Hom^R_{\Sing_A(X)}(a,b) \simeq \Sing(\Exit(a,b)) = \Hom_{\Sing_A(X)}(a,b)
$$

that coincides up to homotopy with the homotopy equivalence above. All that remains to be done is therefore to show that the morphism of Kan complexes

$$
\Phi : \Hom^R_{\Sing_A(X)}(a,b) \simeq \Sing(\Exit(a,b)) \rightarrow \Sing_{\Ql}(\Exit(a,b))
$$

induced by the functor $T'$ coincides up to homotopy with the morphism induced by the natural transformation $\Ql \rightarrow |-|$. We denote this morphism by $U$. Now, there exists an isomorphism between the restrictions of the functors $\Ql$ and $|-|$ to $\Delta^+$. Such an isomorphism determines an isomorphism of semi-simplicial sets 

$$
V : \Sing(\Exit(a,b)) \rightarrow \Sing_{\Ql}(\Exit(a,b))
$$

and by construction, $\Phi$ coincides up to homotopy with the morphism of Kan complexes $V'$ provided by theorem \ref{theorem: Hiro}. Now, there exists a homotopy between the isomorphism of functors $\Ql \simeq |-|$ and the restriction of the natural transformation $\Ql \rightarrow |-|$ to the subcategory $\Delta^+ \subset \Delta$. This homotopy induces a morphism of semi-simplicial sets

$$
H : \Sing(\Exit(a,b)) \times \Delta^1 \rightarrow \Sing_{\Ql}(\Exit(a,b))
$$

such that $H_0=U$ and $H_1 = V$. By \cite[Corollary 1.6]{HiroFunctors}, $U'$ and $V'$ coincide up to homotopy. Finally, since $U$ is already a morphism of simplicial sets, $U$ and $U'$ coincide up to homotopy. We conclude that $U$ and $\Phi$ coincide up to homotopy, as desired. \qedhere

\end{proof}

We finish by giving the proofs of the corollaries of our main result stated in section \ref{section: Morse homology and constructible sheaves}.

\begin{proof}[Proof of corollary \ref{corollary: infinity category of constructible sheaves}]

    This result is proven in \cite[Theorem A.9.3]{HigherAlgebra} in the case when the coefficients are in the $\infty$-category of spaces. This particular case is used in \cite[Theorem B.9]{OrsnesStratified} to deduce the case of coefficients in an arbitrary compactly generated $\infty$-category. Another approach to this result is presented in \cite{PortaTeyssierExodromy} (see in particular Theorem 5.17). Roughly speaking, the exodromy equivalence is obtained as follows. Given a constructible sheaf $\Fl$, one associates to it the functor that carries every point $x$ to the stalk $\Fl_x$, and every exit path $\gamma$ to the monodromy of $\Fl$ along $\gamma$. For a neat explanation of how to carry out this perspective rigorously, we refer to the introduction of \cite{PortaTeyssierExodromy}. \qedhere

\end{proof}

\begin{proof}[Proof of corollary \ref{corollary: localization of infinity category associated with M}]

    It is proven in \cite[Corollary A.9.4]{HigherAlgebra} that the inclusion of simplicial sets
    
    $$
    \Sing_A(X) \rightarrow \Sing(X)
    $$
    
    is a weak homotopy equivalence. This is equivalent to the condition that for every Kan complex $K$, precomposition with this inclusion induces a homotopy equivalence
    
    $$\underline{\Hom}(\Sing(X),K) \rightarrow \underline{\Hom}(\Sing_A(X), K).$$
    
    Since both sides of this map are Kan complexes, this map is an equivalence of $\infty$-categories, which yields the desired result. We also refer to \cite[Corollary 1.3]{PortaTeyssierExodromy} for another approach to this result. \qedhere

\end{proof}

\begin{proof}[Proof of corollary \ref{corollary: fundamental groupoid as localization of flow category}]

    Suppose we are given a functor of $\infty$-categories $\Cl \rightarrow \Dl$ exhibiting $\Dl$ as the localization of $\Cl$ with respect to the class of all morphisms. This formally implies that the same property holds for the functor induced at the level of homotopy categories. Recall also that there is a natural isomorphism of categories between $\mathrm{h} (\Nl(\Ml))$ and $\pi_0(\Ml)$ (remark \ref{remark: homotopy category of topological category versus quasicategory}). Moreover, there is a natural isomorphism of categories between $\mathrm{h} \Sing(X)$ and $\Pi_1(X)$. In particular, the functor $\Nl(\Ml) \rightarrow \Sing(X)$ from corollary \ref{corollary: localization of infinity category associated with M} induces a functor $\pi_0(\Ml) \rightarrow \Pi_1(X)$ that exhibits the right-hand side as the localization of the left-hand side with respect to the class of all morphisms. It remains to check that this functor carries every homotopy class of (possibly broken) gradient trajectory to the homotopy class of the corresponding path in $X$. This follows by inspecting our zigzag of equivalences between $\Nl(\Ml)$ and $\Sing_A(X)$. \qedhere

\end{proof}

\newpage

\section{Stratified categories as exit path categories}\label{section: stratified categories as exit path categories}

Let $\Ml$ be a topological category (i.e., for us, a topologically enriched category) and $Y$ a stratified topological space whose stratified singular simplicial set is an $\infty$-category (when $Y$ is metrizable, an equivalent formulation of the latter condition is given by proposition \ref{proposition: Sing_A infinity category in terms of starting point evaluation}).

The goal of this section is to explain how the techniques developed in this paper can be reproduced to prove that, under certain further assumptions, the $\infty$-category defined by $\Ml$ is equivalent to the $\infty$-category defined by the stratified singular simplicial set of $Y$. The main theorem of this section (theorem \ref{theorem: main theorem extended}) does not refer to Morse theory, but instead, to an extra structure on $\Ml$ inspired from Morse flow categories, which we refer to as a \emph{structure of stratified category}. In the particular case of Morse flow categories, this theorem recovers our main theorem \ref{Main theorem}. This section is intended to be readable as a way to understand the proof of theorem \ref{Main theorem}. The proofs of most of the results of this section are straightforward adaptations of proofs given in the previous sections, in which case these are simply referred to, and not reproduced.

We will study other applications of the main theorem of this section in future work.

\subsection{Stratified categories}\label{section: stratified categories}

Our goal in this section is to introduce an extra structure on $\Ml$ that allows to speak of unbroken morphisms of $\Ml$ and, more generally, to define the flow coherent nerve of $\Ml$. To this end, we will follow the equivalent definition of the flow coherent nerve of a Morse flow category provided by proposition \ref{proposition: unbroken iff compatible with stratifications on morph spaces}, which makes use of \emph{stratifications} on the spaces of morphisms. We will introduce the notion of \emph{stratified category} to be a topological category whose topological spaces of morphisms are coherently endowed with certain stratifications. In the next section, we will define and study the flow coherent nerve of a stratified category.

In order to formulate the precise definition of a stratified category, we start by explaining a construction that associates, to every poset, a category enriched over posets. This construction can also been found in \cite[Construction 2.3.5.1]{kerodon}.

\begin{Construction}\label{construction: poset enriched category from poset}
    Let $A$ be a poset. We construct from $A$ a poset-enriched category $P_A$, called the \emph{path category} of $A$, as follows. The set of objects of $P_A$ is $A$. Given $a,b \in A$, the poset of morphisms of $P_A$ between $a$ and $b$ is the full subposet of $(\sd A)^{\op}$ spanned by those sequences starting with $a$ and finishing with $b$. In other words, this is the set of strictly increasing sequences of elements of $A$ of the form
    
    $$
    a=a_0 < a_1 < \hdots < a_{n-1} < a_n=b
    $$
    
    endowed with the partial order defined by $\ag \le \ag'$ if and only if $\ag'$ is a subsequence of $\ag$. To define the composition operation of $P_A$, consider two sequences
    
    $$
    \ag = [a=a_0 < \hdots < a_n =b ] \in P_A(a,b) \qquad \text{and} \qquad \bg = [b=b_0 < \hdots < b_m =c] \in P_A(b,c).
    $$
    
    The composition of $\ag$ and $\bg$ is defined as the concatenation of $\ag$ and $\bg$, i.e., as the sequence
    
    $$
    \bg \circ \ag = [a_0 < \hdots < a_n < b_1 < \hdots < b_m].
    $$
\end{Construction}

The construction of the path category of a poset allows us to make the following definition. This definition was partly inspired by the definition of abstract flow categories by Abouzaid-Blumberg in \cite[Sections 2 and 3]{AbouzaidBlumbergI}.

\begin{Definition}\label{definition: stratified category}
   Let $\Ml$ be a topological category and $A$ the set of objects of $\Ml$. A \emph{structure of stratified category} on $\Ml$ consists of the datum of a partial order on $A$ together with a functor of topological categories $\Strat : \Ml \rightarrow P_A$, where $P_A$ is regarded as a topological category by endowing each of the posets $P_A(a,b)$ with the Alexandrov topology (definition \ref{definition: Alexandrov topology}). Moreover, $\Strat$ is assumed to induce the identity on objects. The triple consisting of the topological category $\Ml$, the poset $A$ and the functor $\Strat$ is called a \emph{stratified category} and denoted $(\Ml,A,\Strat)$.
\end{Definition}

\begin{Remark}\label{remark: partial order on objects of stratified category}
Suppose we are given a stratified category $(\Ml,A,\Strat)$ and two objects $a,b$ of $\Ml$. We have $P_A(a,b) \neq \varnothing$ if and only if $a \le b$. In particular, since there is a map $\Strat(a,b) : \Ml(a,b) \rightarrow P_A(a,b)$, if $\Ml(a,b) \neq \varnothing$ then we must have $a \le b$.
\end{Remark}

\begin{Notation}\label{notation: top stratum of morph space of strat cat}
Suppose we are given a stratified category $(\Ml,A,\Strat)$ and two objects $a,b \in A$. If the poset $P_A(a,b)$ is nonempty, i.e., if $a \le b$, then it has a maximal element. If $a=b$, this is the sequence $[a]$, and if $a < b$, this is the sequence $[a < b]$. We denote by $\Mi(a,b)$ the stratum of $\Ml(a,b)$ associated with this maximal element, with respect to the stratification defined by the map $\Strat(a,b) : \Ml(a,b) \rightarrow P_A(a,b)$.
\end{Notation}

\begin{Definition}\label{definition: unbroken morphisms of a stratified category}
    Suppose we are given a stratified category $(\Ml,A,\Strat)$ and a morphism $\gamma$ of $\Ml$ from an object $a$ to an object $b$. The morphism $\gamma$ is said to be \emph{unbroken} if $\gamma \in \Mi(a,b)$.
\end{Definition}

We have encountered two sources of examples of stratified categories in this paper.

\begin{Example}\label{example: stratified category associated with a poset}

    Let $A$ be a poset. Applying the nerve functor to the morphism posets of $P_A$, one associates to $P_A$ a simplicially enriched category. Applying the geometric realization functor to the simplicial sets of morphisms of the latter, one associates to $P_A$ a topological category, which we denote by $|P_A|$. We can endow $|P_A|$ with the structure of a stratified category as follows. For every $a,b \in A$, the nerve of the poset $P_A(a,b)$ is endowed with a tautological $P_A(a,b)$-stratification (defined by the identity morphism $N(P_A(a,b)) \rightarrow N_A(P_A(a,b))$), which in turn endows the geometric realization $|N(P_A(a,b))|$ with the structure of a $P_A(a,b)$-stratified topological space. We let $\Strat_A(a,b)$ denote the stratification map $|N(P_A(a,b))| \rightarrow P_A(a,b)$. These stratifications are compatible with compositions. We have thus a stratified category $(|P_A|,A,\Strat_A)$.

\end{Example}

\begin{Remark}\label{remark: stratified category structure on F(n)}

    In section \ref{section: the homotopy coherent nerve} we studied the left adjoint of the homotopy coherent nerve functor, which we denoted $F : \sSet \rightarrow \Cat_{\Top}$. For every nonnegative integer $n$, the topological category $F(\Delta^n)$ is defined by applying the nerve functor, and then the geometric realization functor, to the posets of morphisms of a poset-enriched category $\PR(\Delta^n$) (notation \ref{notation: poset enriched category associated with simplices}). On the other hand, consider the stratified category obtained from example \ref{example: stratified category associated with a poset} in the special case when $A$ is the ordered set $[n]=(0<1< \hdots < n)$. The poset-enriched categories $P_{[n]}$ and $\PR(\Delta^n)$ are almost the same except that we made the choice to define the partial orders on the morphism sets of $\PR(\Delta^n)$ to be the \emph{opposite} of those on the morphism sets of $P_{[n]}$. Example \ref{example: stratified category associated with a poset} endows $P_{[n]}$ with a structure of stratified category; we can also endow $F(\Delta^n)$ with the structure of a stratified category as follows. For every $0 \le i,j \le n$, $\PR(i,j)$ is the set $\{E \subset \{0,\hdots,n\} : i,j \in E \text{ and } \forall k \in E, i \le k \le j \}$ partially ordered by the inclusion relation, while $P_{[n]}(i,j)$ is (naturally isomorphic to) the same set partially ordered by the \emph{reverse inclusion} relation. There is a natural isomorphism between these two posets defined in both directions by mapping $E$ to $\{i,j\} \cup \{i < k < j \mid k \notin E\}$. This defines an isomorphism of poset-enriched categories $\PR(\Delta^n) \simeq P_{[n]}$. The structure of stratified category $(|P_{[n]}|,[n],\Strat_{[n]})$ thus induces a structure of stratified category on $F(\Delta^n)$, which we denote by $(F(\Delta^n),[n],\Strat_n)$. One can check that the notion of unbroken morphism of $F(\Delta^n)$ from definition \ref{definition: unbroken morphisms of a stratified category} coincides with that of definition \ref{definition: broken and unbroken morphism of F(n)}, i.e., a morphism of $F(\Delta^n)$ is unbroken if and only if it cannot be written as a composition of two morphisms that are both different from the identity.

\end{Remark}

\begin{Example}\label{example: stratified category associated with a flow category}

    Let $(f,\xi)$ be a Morse-Smale pair on a smooth closed manifold. Denote by $\Ml$ the associated flow category and by $A$ the associated poset of critical points. In section \ref{section: The flow coherent nerve of M} we defined for every $a, b \in A$ a stratification map $\Ml(a,b) \rightarrow P_A(a,b)$, which associates to a trajectory between $a$ and $b$ the ordered sequence of critical points at which it is broken, including $a$ and $b$ (see proposition \ref{proposition: stratifications of morphism spaces}). Denoting this stratification map by $\Strat_{\Ml}(a,b)$, we have a stratified category $(\Ml,A,\Strat_{\Ml})$. The notion of unbroken morphism of $\Ml$ from definition \ref{definition: unbroken morphisms of a stratified category} coincides with the notion of unbroken trajectory of $\xi$ connecting two critical points.

\end{Example}

We now introduce the notion of functor between stratified categories. We start by a construction which associates, to a morphism of posets, a functor between the corresponding path categories.

\begin{Construction}\label{construction: functor induces on path categories}

    Let $\varphi : A \rightarrow B$ be a morphism of posets. We define a functor between poset-enriched categories $P_{\varphi} : P_A \rightarrow P_B$ as follows.
    
    \begin{itemize}
        \item $P_{\varphi}$ is equal to $\varphi$ at the level of objects.
        
        \item Given $a,b \in A$ and a sequence $\ag = [a = a_0 < a_1 < \hdots < a_n = b] \in P_A(a,b)$, we define $P_{\varphi}(a,b)(\ag)$ as the ordered sequence of distinct elements of the sequence $[\varphi(a_0) \le \varphi(a_1) \le \hdots \le \varphi(a_n)]$.
    \end{itemize}

\end{Construction}

\begin{Definition}\label{definition: functor of stratified categories}

Let $(\Ml,A,\Strat)$ and $(\Nl,B,\mathrm{Strat'})$ be two stratified categories. A \emph{functor of stratified categories} between these is a functor of topological categories $G : \Ml \rightarrow \Nl$ satisfying the following conditions.

\begin{itemize}
    \item The map $A \rightarrow B$ induced by $G$ at the level of objects is order-preserving. We denote the corresponding morphism of posets by $G_{\mathrm{obj}}$.
    
    \item The following diagram of topological categories commutes
    
    $$
    \xymatrix@C=3em{
    \Ml \ar[r]^-{G} \ar[d]_-{\Strat} & \Nl \ar[d]^-{\mathrm{Strat'}} \\
    P_A \ar[r]_-{P_{G_{\mathrm{obj}}}} & P_B.
    }
    $$
\end{itemize}

\end{Definition}

\begin{Remark}\label{remark: condition that functor is order-preserving}

    Recall from remark \ref{remark: partial order on objects of stratified category} the implication $\Ml(a,b) \neq \varnothing \Rightarrow a \le b$. If, in definition \ref{definition: functor of stratified categories}, we assume that this implication is actually an equivalence in the case of $\Ml$, then the first condition of the definition is automatically satisfied. 

\end{Remark}

\begin{Remark}\label{remark: commutativity condition for functor of stratified categories}

    The square diagram in the second condition of definition \ref{definition: functor of stratified categories} is always commutative at the level of objects. The condition that this diagram is commutative is therefore equivalent to the condition that for every pair of objects $(a,b)$ of $\Ml$, the following diagram of topological spaces commutes
    
    $$
    \xymatrix@C=5em{
    \Ml(a,b) \ar[r]^-{G(a,b)} \ar[d]_-{\Strat(a,b)} & \Nl(G(a),G(b)) \ar[d]^-{\mathrm{Strat'}(G(a),G(b))} \\
    P_A(a,b) \ar[r]_-{P_{G_{\mathrm{obj}}}(a,b)} & P_B(G(a),G(b)).
    }
    $$

\end{Remark}

\begin{Lemmas}\label{lemma: composition of functors of stratified categories}

    The composition of two functors of stratified categories is a functor of stratified categories.

\end{Lemmas}

\begin{proof}

    Let $(\Ml,A,\Strat)$, $(\Nl,B,\mathrm{Strat'})$ and $(\Pl, C, \Strat '')$ be three stratified categories, $G$ a functor of stratified categories between $(\Ml,A,\Strat)$ and $(\Nl,B,\mathrm{Strat'})$ and $H$ a functor of stratified categories between $(\Nl,B,\mathrm{Strat'})$ and $(\Pl, C, \Strat '')$. Let $(a,b)$ be a pair of objects of $\Ml$. The left and right squares of the following diagram are commutative
    
    $$
    \xymatrix@C=7em{
    \Ml(a,b) \ar[r]^-{G(a,b)} \ar[d]_-{\Strat(a,b)} & \Nl(G(a),G(b)) \ar[r]^-{H(G(a),G(b))} \ar[d]^-{\mathrm{Strat'}(G(a),G(b))} & \Pl(H(G(a)),H(G(b))) \ar[d]^-{\Strat''(H(G(a)),H(G(b)))} \\
    P_A(a,b) \ar[r]_-{P_{G_{\mathrm{obj}}}(a,b)} & P_B(G(a),G(b)) \ar[r]_-{P_{H_{\mathrm{obj}}}(G(a),G(b))} & P_C(H(G(a)),H(G(b)))
    }
    $$
    
    and therefore the outer square is commutative as well. The statement follows from the fact that the composition of $P_{G_{\mathrm{obj}}}(a,b)$ and $P_{H_{\mathrm{obj}}}(G(a),G(b))$ is equal to $P_{(H \circ G)_{\mathrm{obj}}}(a,b)$. \qedhere

\end{proof}

\begin{Example}\label{example: functor of stratified categories}

    Let $(\Ml,A,\Strat_{\Ml})$ be the stratified category associated with a Morse-Smale pair on a smooth closed manifold, as in example \ref{example: stratified category associated with a flow category}. Let $\sigma$ be an $n$-simplex of the homotopy coherent nerve of $\Ml$, regarded as a functor of topological categories $F(\Delta^n) \rightarrow \Ml$. Endow $F(\Delta^n)$ with the structure of stratified category $(F(\Delta^n),[n],\Strat_n)$ as in remark \ref{remark: stratified category structure on F(n)}. The sequence $\sigma(0),\sigma(1),\hdots,\sigma(n)$ is increasing in $A$, in other words, $\sigma$ is compatible with the partial orders at the level of objects. We denote this sequence by $\ag = [a_0 \le \hdots \le a_n]$. Furthermore, according to proposition \ref{proposition: unbroken iff compatible with stratifications on morph spaces}, the simplex $\sigma$ is unbroken if and only if the following diagram commutes
    
    $$
    \xymatrix@C=4em{
    F(\Delta^n)(0,n) \ar[d]_-{\Strat_n(0,n)} \ar[r]^-{\sigma(0,n)} & \Ml(a_0,a_n) \ar[d]^-{\Strat_{\Ml}(a_0,a_n)} \\
    P_{[n]}(0,n) \ar[r]_-{P_{\ag}(0,n)} & P_A(a_0,a_n).
    }
    $$
    
    In other words, the condition that $\sigma$ is unbroken is equivalent to the fact that the condition of remark \ref{remark: commutativity condition for functor of stratified categories} is satisfied by $\sigma$ for the pair of objects $(0,n)$ only. One can check that, actually, $\sigma$ is unbroken if and only if $\sigma$ is a functor of stratified categories (see proposition \ref{proposition: weaker characterization of unbroken simplices} for a more general statement).

\end{Example}

\subsection{The flow coherent nerve of a stratified category}\label{section: flow coherent nerve of stratified category}

The goal of this section is to introduce the flow coherent nerve of a stratified category, and to give conditions under which it is an $\infty$-category, and under which it is an $\infty$-category equivalent to the homotopy coherent nerve of the underlying topological category.

Inspired by example \ref{example: functor of stratified categories}, we introduce the following definition. Recall that we denote by $F : \sSet \rightarrow \Cat_{\Top}$ the left adjoint of the homotopy coherent nerve functor. The functor $F$ is studied in detail in section \ref{section: the homotopy coherent nerve}.

\begin{Definition}\label{definition: flow coherent nerve of stratified category}

    Let $(\Ml,A,\Strat)$ be a stratified category. The \emph{flow coherent nerve} of $(\Ml,A,\Strat)$, denoted $S_{\Ml}$, is the simplicial subset of the homotopy coherent nerve of $\Ml$ whose $n$-simplices are the functors of topological categories $F(\Delta^n) \rightarrow \Ml$ that are functors of stratified categories. Here, $F(\Delta^n)$ is endowed with the structure of stratified category described in remark \ref{remark: stratified category structure on F(n)}. A simplex of $\Nl(\Ml)$ that belongs to $S_{\Ml}$ is called an \emph{unbroken simplex} of $\Nl(\Ml)$.

\end{Definition}

The fact that $S_{\Ml}$ is indeed a simplicial set, i.e., stable by face and degeneracy operators in $\Nl(\Ml)$, can be argued as follows. The structures of stratified categories on the topological categories $F(\Delta^n)$ are defined in such a way that for every morphism $\varphi : \Delta^n \rightarrow \Delta^m$ in $\Delta$, the functor of topological categories $F(\varphi) : F(\Delta^n) \rightarrow F(\Delta^m)$ is a functor of stratified categories. Now, if we have an unbroken $m$-simplex of $\Nl(\Ml)$ given by a functor of stratified categories $F(\Delta^m) \rightarrow \Ml$, its image by the operator $\Nl(\Ml)_m \rightarrow \Nl(\Ml)_n$ induced by $\varphi$ is obtained by precomposing with $F(\varphi)$. The statement follows because a composition of functors of stratified categories is a functor of stratified categories (lemma \ref{lemma: composition of functors of stratified categories}).

\begin{Example}\label{example: flow coherent nerve}

    By example \ref{example: functor of stratified categories}, the definition of the flow coherent nerve of a stratified category recovers the definition of the flow coherent nerve of the flow category of a Morse-Smale pair.

\end{Example}

We have the following characterization of unbroken simplices of a stratified category.

\begin{Proposition}\label{proposition: flow coherent nerve of stratified category and unbroken morphisms}
    Let $(\Ml,A,\Strat)$ be a stratified category. An $n$-simplex $\sigma \in \Nl(\Ml)_n$ is unbroken if and only if the image of every unbroken morphism of $F(\Delta^n)$ by the functor $\sigma : F(\Delta^n) \rightarrow \Ml$ is an unbroken morphism of $\Ml$.
\end{Proposition}

\begin{proof}

    In order to prove that the condition is necessary, assume that $\sigma$ is unbroken, consider two integers $0 \le i \le j \le n$ and an unbroken morphism $\gamma \in F(\Delta^n)(i,j)$. We wish to prove that $\sigma(i,j)(\gamma)$ is unbroken. This is straightforward if $i=j$; we therefore assume that $i < j$. Since $\sigma$ is unbroken the following diagram commutes
    
    $$
    \xymatrix@C=4em{
    F(\Delta^n)(i,j) \ar[d]_-{\Strat_n(i,j)} \ar[r]^-{\sigma(i,j)} & \Ml(a_i,a_j) \ar[d]^-{\Strat_{\Ml}(a_i,a_j)} \\
    P_{[n]}(i,j) \ar[r]_-{P_{\ag}(i,j)} & P_A(a_i,a_j).
    }
    $$
    
    The condition that $\gamma$ is unbroken means that $\gamma \in F(\Delta^n)(i,j)_{[i < j]}$ and therefore $\sigma(i,j)(\gamma)$ belongs to $\Ml(a_i,a_j)_{P_{\ag}(i,j)([i < j])}$. But $P_{\ag}(i,j)([i < j])$ is the maximal element of $P_A(a_i,a_j)$, and thus the image of $\gamma$ by $\sigma$ is unbroken, as desired.
    
    Conversely, assume that the image of every unbroken morphism of $F(\Delta^n)$ by $\sigma$ is unbroken. We wish to prove that for any two integers $0 \le i \le j \le n$ the above diagram commutes. Consider an arbitrary morphism $\gamma \in F(\Delta^n)(i,j)$. Let $\mathbf{i} = [i = i_0 < i_1 < \hdots < i_k =j]$ be the image of $\gamma$ by $\Strat_{[n]}(i,j)$. By corollary \ref{corollary: writing as composition of unbroken morphisms} there exists a family of unbroken morphisms $(\gamma_l \in F(\Delta^n)(i_l,i_{l+1}))_{0 \le l < k}$ whose composition is $\gamma$. By assumption the image of $\gamma_l$ by $\sigma$ is unbroken for every $0 \le l < k$, in other words, $\sigma(i_l,i_{l+1})(\gamma_l)$ belongs to the top stratum of $\Ml(a_{i_l},a_{i_{l+1}})$ for every $0 \le l < k$. By functoriality of $\sigma$ and $\Strat_{\Ml}$, this is equivalent to the condition that $\Strat_{\Ml}(a_i,a_j)(\sigma(i,j)(\gamma)) = P_{\ag}(i,j)(\mathbf{i})$. This proves the desired commutativity condition. \qedhere

\end{proof}

\begin{Remark}\label{remark: unbroken simplex can map borken to unbroken}

    It can happen that an unbroken $n$-simplex $\sigma \in \Nl(\Ml)_n$ sends certain broken morphisms to broken ones and certain broken morphisms to unbroken ones. Consider for example the flow category $\Ml$ of a Morse-Smale pair such that there exists a strictly increasing sequence of four critical points $a_0 < a_1 < a_2 < a_3$. Consider an unbroken $3$-simplex $\sigma \in \Nl(\Ml)_3$ such that $\sigma(0) = a_0$, $\sigma(1)=a_1$ and $\sigma(2)=\sigma(3)=a_3$ \footnote{Such a $3$-simplex exists because the flow coherent nerve of $\Ml$ is an $\infty$-category, by proposition \ref{proposition: S_M is an infinity category}.}. Then every morphism in $F(\Delta^n)(0,3)$ which is broken at $1$ is mapped by $\sigma$ to a broken trajectory, while every morphism in $F(\Delta^n)(0,3)$ which is broken only at $2$ is mapped by $\sigma$ to an unbroken trajectory.

\end{Remark}

Unbroken simplices of $\Nl(\Ml)$ are in fact characterized by an a priori weaker condition than that of the definition.

\begin{Proposition}\label{proposition: weaker characterization of unbroken simplices}

    Suppose we have an $n$-simplex $\sigma \in \Nl(\Ml)_n$. Define an increasing sequence $\ag$ of elements of $A$ by $a_i = \sigma(i)$. The simplex $\sigma$ is unbroken if and only if the following diagram commutes
    
    $$
    \xymatrix@C=4em{
    F(\Delta^n)(0,n) \ar[d]_-{\Strat_n(0,n)} \ar[r]^-{\sigma(0,n)} & \Ml(a_0,a_n) \ar[d]^-{\Strat_{\Ml}(a_0,a_n)} \\
    P_{[n]}(0,n) \ar[r]_-{P_{\ag}(0,n)} & P_A(a_0,a_n).
    }
    $$
    
    In other words, the condition that the functor $\sigma : F(\Delta^n) \rightarrow \Ml$ is a functor of stratified categories is equivalent to the condition of remark \ref{remark: commutativity condition for functor of stratified categories} for the pair of objects $(0,n)$ only.

\end{Proposition}

\begin{proof}

    We apply the characterization of unbroken simplices provided by proposition \ref{proposition: flow coherent nerve of stratified category and unbroken morphisms}. Assuming the commutativity of this diagram, we wish to show that the image by $\sigma$ of every unbroken morphism is unbroken. Let $0 \le i \le j \le n$ be two integers and $\gamma \in F(\Delta^n)(i,j)$ an unbroken morphism. Take another two unbroken morphisms $\delta \in F(\Delta^n)(0,i)$ and $\theta \in F(\Delta^n)(j,n)$. We get a morphism $\theta \circ \gamma \circ \delta \in F(\Delta^n)(0,n)$. The assumption that the diagram is commutative implies that the image by $\sigma$ of each of the three morphisms $\theta$, $\gamma$, $\delta$ is unbroken. In particular, $\gamma$ is unbroken, as desired. \qedhere

\end{proof}

The proof of proposition \ref{proposition: S_M is an infinity category} adapts to prove the following criterion for the flow coherent nerve of a stratified category to be an $\infty$-category.

\begin{Proposition}\label{proposition: when is the flow coherent nerve an infinity category}
    Let $(\Ml,A,\Strat)$ be a stratified category. For every pair of objects $(a,b)$ of $\Ml$ such that $a \le b$, endow the morphism space $\Ml(a,b)$ with the stratification by $0 < 1$ obtained by composing the stratification map $\Strat(a,b) : \Ml(a,b) \rightarrow P_A(a,b)$ with the morphism of posets $P_A(a,b) \rightarrow (0<1)$ that maps the maximal element to $1$ and every other sequence to $0$. In other words, this is the $(0<1)$-stratification whose $1$-stratum is $\Mi(a,b)$. \\[0pt]
    Suppose that the starting point evaluation map $\Exit(\Ml(a,b)_0,\Ml(a,b)_1) \rightarrow \Ml(a,b)_0$ is a Serre fibration and a weak homotopy equivalence. Then the flow coherent nerve of $\Ml$ is an $\infty$-category. \qed
\end{Proposition}

Furthermore, the discussion of section \ref{section: equivalence between flow and homotopy coherent nerve} adapts to yield the following description of the right morphism spaces of the flow coherent nerve of a stratified category. There is a functor $\Ql : \Delta \rightarrow \Top$ and for every topological category $\Ml$ and every two objects $a,b$ of $\Ml$ a natural isomorphism of simplicial sets

$$
\Hom_{\Nl(\Ml)}^R(a,b) \simeq \Sing_{\Ql}(\Ml(a,b)).
$$

Further, if $\Ml$ is endowed with a structure of stratified category $(\Ml,A,\Strat)$ this isomorphism restricts to an isomorphism of simplicial sets

$$
\Hom_{S_{\Ml}}^R(a,b) \simeq \Sing_{\Ql}(\Mi(a,b)).
$$

In particular, since a map of topological spaces is a weak homotopy equivalence if and only if its image by $\Sing_{\Ql}$ is a Kan equivalence (proposition \ref{proposition: Sing_Q detects weak equivalences}), we have the following.

\begin{Proposition}\label{proposition: when is flow coherent nerve equivalent to coherent nerve}
    Let $(\Ml,A,\Strat)$ be a stratified category whose flow coherent nerve $S_{\Ml}$ is an $\infty$-category. The following two conditions are equivalent.
    
    \begin{enumerate}[label=(\roman*)]
        \item The inclusion 
        morphism of simplicial sets $S_{\Ml} \rightarrow \Nl(\Ml)$ is an equivalence of $\infty$-categories.
        
        \item For every pair of objects $(a,b)$ of $\Ml$, the inclusion map $\Mi(a,b) \rightarrow \Ml(a,b)$ is a weak homotopy equivalence.
    \end{enumerate} \qed
\end{Proposition}
 
\subsection{Comparison between stratified categories and exit path categories}\label{section: Comparison between stratified categories and exit path categories}

In this section, we build upon the previous two sections to present a framework in which the methods of this paper apply to prove that a stratified category is equivalent, as an $\infty$-category, to the exit path $\infty$-category of some stratified space.

Suppose we are given a topological space $Y$, a poset $A$ and a stratification $\pi : Y \rightarrow A$. We will set up a framework in which we can compare the flow coherent nerve of some stratified category, to the stratified singular simplicial set of $(Y,A)$. To begin with, we will need two further assumptions.

\begin{itemize}
    \item We suppose that we have a strictly decreasing map $f : A \rightarrow \R$ \footnote{We could just as well require $f$ to be strictly increasing; our choice is made in order to be consistent with our convention of considering \emph{negative} pseudo-gradients in Morse theory.} (not every poset admits such a map, see remark \ref{remark: choice of decreasing map} for a discussion).
    
    \item Up to restricting ourselves to the full subposet of $A$ spanned by the image of $\pi$, we may assume that $\pi$ is surjective, which we do. We also assume that a set-theoretic section of $\pi$ is given (such a section always exists). Such a section is the same as a choice, for every $a \in A$, of a point in $Y_a$. We identify $A$ with a subset of $Y$ using this section, i.e., we have $a \in Y_a$.
\end{itemize}

\begin{Remark}\label{remark: choice of decreasing map}

    For cardinality reasons, not every poset admits a strictly decreasing map to $\R$ (or, equivalently, a strictly increasing one). An example of a totally ordered set that has the same cardinality as $\R$ and does not admit a strictly decreasing map to $\R$ is given in \cite[Section 2.1]{FoundationsofMeasurement}. On the other hand, every countable poset admits a strictly decreasing map to $\R$; this is proven in \cite[Section 2.1, Theorem 1]{FoundationsofMeasurement} for totally ordered sets and the proof can be adapted to the case of partially ordered sets.

\end{Remark}

Inspired by Morse theory, we start by introducing a triple $(\Ml_f(Y),A,\Strat)$ which is a stratified category under mild assumptions on $A$. The set of objects of the topological category $\Ml_f(Y)$ is $A$ and its spaces of morphisms are the spaces of "piecewise exit paths" in the $A$-stratified space $Y$ that are parameterized using $f$. More precisely, the topological category $\Ml_f(Y)$ is constructed as follows.

\begin{Construction}\label{construction: topological category built from stratified space}

    The set of objects of $\Ml_f(Y)$ is $A$ and for every $a,b \in A$, the space of morphisms of $\Ml_f(Y)$ from $a$ to $b$ is the topological subspace of $C^0([f(b),f(a)],Y)$, endowed with the compact-open topology, formed by those paths $\gamma$ satisfying the following condition: $\gamma(f(b)) = b$, $\gamma(f(a)) = a$ and there exists a strictly increasing sequence $a=a_0 < a_1 < \hdots < a_n = b$ of elements of $A$ such that for every integer $0 \le i <n$ we have $\gamma([f(a_{i+1}),f(a_i)) \subseteq Y_{a_{i+1}}$ and $\gamma(f(a_{i+1}))=a_{i+1}$. If we are given $c \in A$ and $\gamma \in \Ml_f(Y)(a,b)$, $\gamma' \in \Ml_f(Y)(b,c)$, the composition $\gamma' \circ \gamma \in \Ml_f(Y)(a,c)$ is defined as follows: given a real number $f(c) \le s \le f(a)$, the point $(\gamma' \circ \gamma) (s)$ is defined as $\gamma(s)$ if $s \ge f(b)$ and $\gamma'(s)$ if $s \le f(b)$.

\end{Construction}

\begin{Remark}\label{remark: flow category is subcategory of stratified category}

    Consider the case of a Morse-Smale pair $(f,\xi)$ on a smooth closed manifold $X$, whose flow category we denote by $\Ml$. Denote the poset of critical points of $f$ by $A$ and endow $X$ with the $A$-stratification by the stable manifolds of $\xi$. Denote the restriction of $f$ to $A$ by $f$ as well. Recall from construction \ref{construction: parametrization of all trajectories by values of the function} the map $\Ml(a,b) \rightarrow C^0([f(b),f(a)],X)$ for every $a,b \in A$. These maps yield a functor of topological categories $\Ml \rightarrow \Ml_f(X)$.

\end{Remark}

\begin{Notation}\label{notation: sequence of breaking points is unique}

    Given $a,b \in A$ and $\gamma \in \Ml_f(Y)$, there exists a \emph{unique} sequence $a=a_0 < a_1 < \hdots < a_n = b$ associated with $\gamma$ satisfying the conditions stated in construction \ref{construction: topological category built from stratified space}. Carrying $\gamma$ to this sequence yields a map $\Strat_f(Y)(a,b) : \Ml_f(Y)(a,b) \rightarrow P_A(a,b)$.

\end{Notation}

The following lemma gives conditions under which the triple $(\Ml_f(Y),A,\Strat)$ is a stratified category.

\begin{Lemmas}\label{lemma: stratified category from stratified space}

    For every $a,b \in A$ we denote by $[a,b]$ the set of those elements $c \in A$ such that $a \le c \le b$. Suppose that:
    
    \begin{enumerate}[label=(\roman*)]
        \item for every $a,b \in A$, the set $[a,b]$ is finite and
    
        \item for every $a \in A$, $\{a\}$ is a closed subset of $Y$.
    \end{enumerate}
    
    Then the partial order on $A$ and the maps $\Strat_f(Y)(a,b)$ from notation \ref{notation: sequence of breaking points is unique} endow $\Ml_f$ with the structure of a stratified category.

\end{Lemmas}

\begin{proof}

    We need to verify that $\Strat_f(Y) : \Ml_f(Y) \rightarrow P_A$ is a functor of topological categories. Functoriality holds without assumptions (i) and (ii) and is straightforward. It remains to show that for every $a,b \in A$, the map $\Strat_f(Y)(a,b) : \Ml_f(Y)(a,b) \rightarrow P_A(a,b)$ is continuous (where the right-hand side is endowed with the Alexandrov topology). In other words, we wish to show that for every sequence of the form $\ag = [a=a_0 < a_1 < \hdots < a_n =b]$, the subset $\Ml_f(Y)(a,b)_{\ge \ag} \subseteq \Ml_f(Y)(a,b)$ is open. Let $E = [a,b] \backslash \{a_0,a_1,\hdots,a_n\} \subseteq A$, $K = f(E) \subseteq [f(b),f(a)]$ and $U = X \backslash \{E\} \subseteq Y$. We have:
    
    \[
    \begin{aligned}
    \Ml_f(Y)(a,b)_{\ge \ag} & = \{ \gamma \in \Ml_f(Y)(a,b) \mid \forall c \in E, \gamma(f(c)) \neq c  \} \\
    & = \{ \gamma \in \Ml_f(Y)(a,b) \mid \forall c \in E, \gamma(f(c)) \notin E  \} \\
    & = \{ \gamma \in \Ml_f(Y)(a,b) \mid \gamma(K) \subseteq U \}. \\
    \end{aligned}
    \]
    
    By the assumptions (i) and (ii), $K$ is compact and $U$ is open in $Y$. By definition of the compact-open topology, $\Ml_f(Y)(a,b)_{\ge \ag}$ is open in $\Ml_f(Y)(a,b)$, as desired. \qedhere

\end{proof}

\begin{Remark}\label{remark: top stratum for stratified category of piecewise exit paths}

    For every $a,b \in \Ml$, the subspace of unbroken morphisms $\Mi_f(Y)(a,b) \subseteq \Ml_f(Y)(a,b)$ (definition \ref{definition: unbroken morphisms of a stratified category}) is the subspace of those $\gamma : [f(b),f(a)] \rightarrow Y$ such that $\gamma(f(b))=b$, $\gamma(f(a))=a$ and $\gamma([f(b),f(a)) \subseteq X_b$. If $a < b$, precomposing by the affine decreasing bijection $[0,1] \rightarrow [f(b),f(a)]$ yields a homeomorphism between $\Mi_f(Y)(a,b)$ and $\Exit(a,b)$.

\end{Remark}

\begin{Example}\label{example: functor is stratified in Morse case}

    The functor of topological categories from remark \ref{remark: flow category is subcategory of stratified category} is a functor of stratified categories, where the flow category is endowed with the structure of stratified category described in example \ref{example: stratified category associated with a flow category}.

\end{Example}

Now suppose that the conditions of lemma \ref{lemma: stratified category from stratified space} are satisfied, so that we have a stratified category

$$(\Ml_f(Y),A,\Strat_f(Y))$$

associated with the $A$-stratified space $Y$ and the map $f : A \rightarrow \R$. Suppose that we are given a stratified category $(\Ml,A, \Strat)$ with the same underlying poset $A$. In order to compare its homotopy coherent nerve to $\Sing_A(Y)$, we will assume that we have a functor of stratified categories $G : (\Ml,A,\Strat) \rightarrow (\Ml_f(Y),A,\Strat_f(Y))$ that induces the identity at the level of objects (in the Morse case, this functor is given by example \ref{example: functor is stratified in Morse case}). This yields a zigzag of simplicial sets

$$
\Nl(\Ml) \leftarrow S_{\Ml} \rightarrow S_{\Ml_f(Y)}
$$

where the left arrow is the inclusion and the right one is induced by $G$.

Next, we can associate to the simplices of $S_{\Ml_f(Y)}$ some kind of "stratified cubes" of $Y$, as follows. Consider more generally a simplex $\sigma \in \Nl(\Ml_f(Y))_n$ regarded as a functor of topological categories $\sigma : F(\Delta^n) \rightarrow \Ml_f(Y)$ and denote $\sigma(i) = a_i$ for every $0 \le i \le n$. Since, for every $a,b \in A$, the condition that $\Ml_f(Y)(a,b)$ is nonempty implies that $a \le b$, the sequence $\ag = [a_0,a_1,\hdots,a_n]$ is increasing.

There is a map $\sigma(0,n) : I^{n-1} \rightarrow \Ml_f(Y)(a_0,a_n)$ associated with $\sigma$. Furthermore, since $\Ml_f(Y)(a_0,a_n)$ is a subspace of $C^0([f(a_n),f(a_0)],Y)$, there is an evaluation map

$$
\Ml_f(Y)(a_0,a_n) \times [f(a_n),f(a_0)] \rightarrow Y.
$$

Combining the two latter maps we get a singular $n$-cube in $Y$ which we denote by

$$
f_{\sigma} : I^{n-1} \times [f(a_n),f(a_0)] \rightarrow Y.
$$

Since our ultimate goal is to compare $\Nl(\Ml)$ and $\Sing_A(Y)$ we have two requirements at this point: we would like the domain of $f_{\sigma}$ to be comparable to an $n$-simplex, and we would like $f_{\sigma}$ to be compatible with some stratification. To this end we define, using the morphism of posets $f : A \rightarrow \R$, a functor $C_A^{\Top} : \Delta_A \rightarrow \Top$ as in section \ref{section: non stratified functor C A} \footnote{Notice the abuse of notation here, since this functor depends on $A$ but also on $f$.}. Briefly, this is done by introducing, for every $\ag \in \Delta_A$, an equivalence relation $\sim_{\ag}$ on $I^{n-1} \times [f(a_n),f(a_0)]$ which is compatible with the maps $f_{\sigma}$. The topological space $C_A^{\Top}(\Delta^{\ag})$ is then defined as the quotient of $I^{n-1} \times [f(a_n),f(a_0)]$ by $\sim_{\ag}$. The left Kan extension of $C_A^{\Top}$ to $\sSet_A$ admits a right adjoint denoted $\Sing_{C_A^{\Top}}(-) : \Top \rightarrow \sSet_A$ and there is a morphism of simplicial sets $\Nl(\Ml_f(Y)) \rightarrow \Sing_{C_A^{\Top}}(Y)$ defined by assigning, to a simplex $\sigma$, the map induced by $f_{\sigma}$ on the quotient.

Furthermore, the discussion of section \ref{section: lift to stratified spaces} holds for $C_A^{\Top}$ as well (with the same proofs). In other words, $C_A^{\Top}$ lifts to $\Top_A$ in a way which is unique under the natural conditions of theorem \ref{theorem: stratifications on cubes}, yielding a functor $C_A : \Delta_A \rightarrow \Top_A$ and, taking the right adjoint of the left Kan extension, a functor $\Sing_{C_A}(-) : \Top_A \rightarrow \sSet_A$. Furthermore, for every $\sigma \in \Nl(\Ml_f(Y))_n$, the condition that $\sigma$ is unbroken, i.e., belongs to the flow coherent nerve $S_{\Ml_f(Y)}$, is equivalent to the condition that the associated $n$-simplex of $\Sing_{C_A^{\Top}}(Y)$ is compatible with stratifications, i.e., lies in the simplicial subset $\Sing_{C_A}(Y)$. Said differently, the following is a pullback square of simplicial sets

$$
\xymatrix{
S_{\Ml_f(Y)} \ar[r] \ar[d] & \Nl(\Ml_f(Y)) \ar[d] \\
\Sing_{C_A}(Y) \ar[r] & \Sing_{C_A^{\Top}}(Y).
}
$$

We have thus extended our zigzag to a zigzag of simplicial sets

$$
\Nl(\Ml) \leftarrow S_{\Ml} \rightarrow S_{\Ml_f(Y)} \rightarrow \Sing_{C_A}(Y).
$$

In order to compare the stratified simplicial sets $\Sing_{C_A}(Y)$ and $\Sing_A(Y)$, we compare the functors $C_A$ and $|-|_A$ in the same way as in section \ref{section: comparison between realizations}. We introduce the quotient $|-|_A'$ of the standard stratified geometric realization $|-|_A$ as in section \ref{section: modifying geometric realization}, and denote by $\Sing_A'(-) : \Top_A \rightarrow \sSet_A$ the corresponding right adjoint of left Kan extension. Theorem \ref{theorem: isomorphism between |-|'_A and C_A}, stating that the functors $|-|_A'$ and $C_A$ are isomorphic in restriction to the subcategory $\Delta_A^+ \subseteq \Delta_A$ generated by face maps, holds with the same proof. We have thus an isomorphism of semi-simplicial sets between $\Sing_{C_A}(Y)$ and $\Sing_A'(Y)$.

The main theorem of this section is the following.

\begin{Thms}\label{theorem: main theorem extended}

    Suppose that we have a stratified category $(\Ml,A,\Strat)$ together with a functor of topological categories $G : \Ml \rightarrow \Ml_f(Y)$ satisfying the following conditions:
    
    \begin{enumerate}[label=(\roman*)]
        \item The conditions of lemma \ref{lemma: stratified category from stratified space} are satisfied and G is a functor of stratified categories inducing the identity at the level of objects.
        
        \item For every $a,b \in A$, the map $\Mi(a,b) \rightarrow \Mi_f(Y)(a,b)$ induced by $G(a,b)$ is a weak homotopy equivalence.
        
        \item $Y$ has weakly contractible strata.
    \end{enumerate}
    
    Suppose further the following:
    
    \begin{enumerate}[label=(\roman*), start=4]
        \item The stratified simplicial set $\Sing_A(Y)$ is an $\infty$-category.
        
        \item The flow coherent nerve of $(\Ml,A,\Strat)$ is an $\infty$-category and the inclusion $S_{\Ml} \rightarrow \Nl(\Ml)$ is an equivalence of $\infty$-categories (see section \ref{section: flow coherent nerve of stratified category} for a discussion of these conditions).
    \end{enumerate}
    
    Then there exists an equivalence of $\infty$-categories between $\Nl(\Ml)$ and $\Sing_A(Y)$ that satisfies the following properties.
    
    \begin{itemize}
        
        \item It is compatible with the functors to $N(A)$ on both sides.
        
        \item For every $a,b \in A$, at the level of the morphism spaces between $a$ and $b$, this equivalence induces the isomorphism of homotopy types obtained from the weak homotopy equivalences $\Mi(a,b) \rightarrow \Ml(a,b)$ and $\Mi(a,b) \rightarrow \Mi_f(Y)(a,b) = \Exit(a,b)$.

    \end{itemize}

\end{Thms}

\begin{Remark}\label{remark: contractibility assumption}

In order to remove the contractibility assumption (iii) from this theorem, it would be interesting to consider categories \emph{internal to the category of topological spaces} instead of topological categories. The difference is that the set of objects of a category internal to $\Top$ is also endowed with a topology. We expect that the category $\Ml_f$ should then be replaced by a construction analogous to \emph{stratification categories} introduced by Miller in \cite[Section 8]{MillerStrongly} \footnote{I am grateful to David Chataur for pointing out this reference.}.

\end{Remark}

\begin{Example}\label{example: extended main theorem recovers main theorem}

    In the Morse case, the functor of topological categories from remark \ref{remark: flow category is subcategory of stratified category} satisfies the conditions of theorem \ref{theorem: main theorem extended}. The latter therefore recovers our main theorem (formulated in section \ref{theorem: main theorem reformulated}).

\end{Example}

\begin{proof}

By condition (i) and the preceding discussion, we have a zigzag of simplicial sets

$$
\Nl(\Ml) \leftarrow S_{\Ml} \rightarrow \Sing_{C_A}(Y)
$$

as well as an isomorphism of semi-simplicial sets between $\Sing_{C_A}(Y)$ and $\Sing_A'(Y)$.

Condition (v) says that the morphism $S_{\Ml} \rightarrow \Nl(\Ml)$ is an equivalence of $\infty$-categories. On the other hand, condition (iv) implies that $\Sing_{C_A}(Y)$ is an $\infty$-category, with the same proof as in section \ref{section: filling inner horns in Sing C A (X)}.

However, the morphism $S_{\Ml} \rightarrow \Sing_{C_A}(Y)$ has no chance to be essentially surjective in general, for the same reason as that given at the beginning of section \ref{section: comparison between S M and Sing C A (X)}. Namely, by construction of $C_A$ the simplices of $\Sing_{C_A}(Y)$ that are entirely contained in a given stratum are just points, and so there are as many isomorphism classes of objects of the $\infty$-category $\Sing_{C_A}(Y)$ as there are points in $Y$. However, the set of objects of $S_{\Ml}$ is $A$ only. We remedy this by introducing the \emph{restricted version} $\Sing_{C_A,r}(Y) \subseteq \Sing_{C_A}(Y)$ defined to consist of those simplices whose vertices belong to $A$; it is stable by inner horn filling and is therefore an $\infty$-category. Now the functor $S_{\Ml} \rightarrow \Sing_{C_A,r}(Y)$ is essentially surjective, since it is a bijection at the level of objects. It is also fully faithful under the conditions of the theorem; this can be proven by considering the maps that it induces between right morphism spaces, as follows.

Firstly, suppose we have two objects $a,b \in A$ such that $a < b$. Remark \ref{remark: top stratum for stratified category of piecewise exit paths} provides a homeomorphism between $\Mi_f(Y)(a,b)$ and $\Exit(a,b)$. Combining it with the map $\Mi(a,b) \rightarrow \Mi_f(Y)(a,b)$ induced by $G(a,b)$, we get a map denoted $\iota_{a,b} : \Mi(a,b) \rightarrow \Exit(a,b)$ which, by condition (ii), is a weak homotopy equivalence. Proposition \ref{proposition: right morphism space as space of exit paths} then holds in this setting with the same proof, namely there exists a top horizontal isomorphism of simplicial sets rendering the following commutative square of Kan complexes (strictly) commutative

$$
    \xymatrix{
    \Hom_{\Sing_{C_A,r}(Y)}^R(a,b) \ar[r]^-{\simeq} & \Sing_{\Ql}(\Exit(a,b)) \\
    \Hom_{S_{\Ml}}^R(a,b) \ar[u] \ar[r]_-{\simeq} & \Sing_{\Ql}(\Mi(a,b)), \ar[u]_-{\Sing_{\Ql}(\iota_{a,b})}
    }
$$

thus proving that the left vertical arrow of this diagram is a homotopy equivalence.

Secondly, suppose given an object $a \in A$. By condition (ii), the map $\Mi(a,a) \rightarrow \Mi_f(Y)(a,a)$ is a weak homotopy equivalence. In this case, we have $\Mi_f(Y)(a,a) = \Ml_f(Y)(a,a)$ and this space is the subspace of $C^0(\{f(a)\},Y)$ formed by those maps carrying $f(a)$ to $a$; in other words, this space is a point, and therefore $\Mi(a,a)$ is weakly contractible. Furthermore we have $\Hom_{S_{\Ml}}^R(a,a) \simeq \Sing_{\Ql}(\Mi(a,a))$ and $ \Hom_{\Sing_{C_A,r}(Y)}^R(a,a) \simeq \Delta^0$. The map $\Hom_{S_{\Ml}}^R(a,a) \rightarrow \Hom_{\Sing_{C_A,r}(Y)}^R(a,a)$ is thus a homotopy equivalence.

Lastly, suppose given two objects $a,b \in A$ such that $a \nleq b$. Then both of the right morphism spaces of $S_{\Ml}$ and $\Sing_{C_A,r}(Y)$ between $a$ and $b$ are empty, in particular there exists a unique map between them, which is a weak homotopy equivalence. This completes the proof of the fact that $S_{\Ml} \rightarrow \Sing_{C_A,r}(Y)$ is fully faithful, and an equivalence of $\infty$-categories.

In order to complete the construction of a zigzag of equivalences of $\infty$-categories between $\Nl(\Ml)$ and $\Sing_A(Y)$, we introduce the restricted versions

$$
\Sing_{A,r}'(Y) \subseteq \Sing_A'(Y) \qquad \text{and} \qquad \Sing_{A,r}(Y) \subseteq \Sing_A(Y),
$$

also defined as those simplicial subsets consisting of those simplices whose vertices belong to $A$. Condition (iv) implies that the simplicial set $\Sing_A'(Y)$ is an $\infty$-category; this can be proven in the same way as proposition \ref{proposition: if Sing_A infinity category then Sing_A' too}. Since these restricted versions are stable by inner horn filling, these are $\infty$-categories. There are two morphisms of simplicial sets

$$
\Sing_{A,r}'(Y) \rightarrow \Sing_{A,r}(Y) \rightarrow \Sing_A(Y),
$$

both of which are equivalences of $\infty$-categories; this can be proven in the same way as proposition \ref{proposition: restricted version into unrestricted is equivalence} using the fact that the strata of $Y$ are weakly contractible (condition (iii)). Finally, the isomorphism of semi-simplicial sets between $\Sing_{C_A}(Y)$ and $\Sing_A'(Y)$ restricts to an isomorphism of semi-simplicial sets between $\Sing_{C_A,r}(Y)$ and $\Sing_{A,r}'(Y)$. Consequently, theorem \ref{theorem: Hiro} provides an equivalence of $\infty$-categories between $\Sing_{C_A,r}(Y)$ and $\Sing_{A,r}'(Y)$. We have therefore shown that there exists a zigzag of equivalences of $\infty$-categories between $\Nl(\Ml)$ and $\Sing_A(Y)$. The fact that the induced equivalence of $\infty$-categories between $\Nl(\Ml)$ and $\Sing_A(Y)$ satisfies the stated properties is proved as in the Morse case in section \ref{section: end of the proof}. \qedhere

\end{proof}

\newpage

\appendix

\section{Colimits of simplicially enriched categories}\label{appendix: colimits of simplicial categories}

The goal of this appendix is to give a proof of the following statement.

\begin{Proposition}\label{proposition: simplicial categories admit small colimits}
    The category of small simplicial categories admits all small colimits.
\end{Proposition}

\begin{proof}
Let us denote by $\Cat$ the category of small categories. Recall that we denote by $\Cat_{\sSet}$ the category of small simplicial categories, which, for us, is the category of small categories enriched over simplicial sets. Let us denote by $s \Cat$ the functor category $\Fun(\Delta^{\op}, \Cat)$. Consider the functor

$$
\fonctionsansnom{\Cat_{\sSet}}{s \Cat}{\Cl}{
[n] \mapsto
\left\{
\begin{array}{ll}
    \Obj & = \Obj(\Cl) \\
    \Hom_{\Cl_n}(X,Y) & = \Hom_{\Cl}(X,Y)_n. 
\end{array}
\right.
}
$$

It induces an equivalence between $\Cat_{\sSet}$ and the full subcategory of $s \Cat$ spanned by those functors $\Delta^{\op} \rightarrow \Cat$ having the property that the face and degeneracy maps yield bijections at the level of objects in $\Cat$. It is therefore enough to prove that this property is preserved by colimits in $s \Cat$, and that $s \Cat$ admits all small colimits. By remark \ref{remark: presheaf categories admit all limits and colimits}, it is enough to prove that $\Cat$ admits all small colimits, and that the forgetful functor

$$
\fonctionsansnom{\Cat}{\Set}{C}{\Obj \, C}
$$

commutes with colimits. This follows from corollary 4.5.16 and its proof in \cite{RiehlCategory}. \qedhere
\end{proof}

\newpage

\bibliographystyle{alpha}
\bibliography{bibliographie}

@book{BDHOMorse,
  author={Barraud, Jean-Fran{\c{c}}ois and Damian, Mihai and Humili{\`e}re, Vincent and Oancea, Alexandru},
  title={Morse Homology with Differential Graded Coefficients},
  publisher={Birkhäuser},
  year={2025},
  series={Progress in Mathematics},
  volume={360}
}

@article{MaderWaasSamples,
  title={From samples to persistent stratified homotopy types},
  author={M{\"a}der, Tim and Waas, Lukas},
  journal={Journal of Applied and Computational Topology},
  volume={8},
  number={3},
  pages={761--838},
  year={2024},
  publisher={Springer}
}

@article{AbouzaidBlumbergI,
  title={Foundation of {F}loer {h}omotopy {t}heory {I}: Flow categories},
  author={Abouzaid, Mohammed and Blumberg, Andrew J.},
  journal={preprint},
  year={2024},
  note={Available at: \url{https://arxiv.org/abs/2404.03193}}
}

@article{CohenJonesSegalMorse,
  title={Morse theory and classifying spaces},
  author={Cohen, Ralph L. and Jones, John D.S. and Segal, Graeme B.},
  journal={preprint},
  year={1995},
  note={Available at: \url{https://www.math.toronto.edu/mgualt/Morse%20Theory/CohenJonesSegal.pdf}}
}

@article{CohenJonesSegalFloer,
  title={Floer’s infinite dimensional {M}orse theory and homotopy theory},
  author={Cohen, Ralph L and Jones, John DS and Segal, Graeme B},
  journal={The Floer memorial volume},
  pages={297--325},
  year={1995},
  publisher={Springer}
}

@book{HigherTopos,
  title={Higher topos theory},
  author={Lurie, Jacob},
  year={2009},
  publisher={Princeton University Press}
}

@book{AudinDamianMorseFloer,
  title={Morse theory and Floer homology},
  author={Audin, Mich{\`e}le and Damian, Mihai},
  volume={2},
  year={2014},
  publisher={Springer}
}

@book{NicolaescuInvitation,
  title={An invitation to Morse theory},
  author={Nicolaescu, Liviu I.},
  year={2007},
  publisher={Springer}
}

@unpublished{HigherAlgebra,
  title={\emph{Higher algebra}},
  author={Lurie, Jacob},
  year={2017},
  note={Available at: \url{https://www.math.ias.edu/~lurie/papers/HA.pdf}}
}

@article{CordierDiagramme,
  title={Sur la notion de diagramme homotopiquement coh{\'e}rent},
  author={Cordier, Jean-Marc},
  journal={Cahiers de topologie et g{\'e}om{\'e}trie diff{\'e}rentielle},
  volume={23},
  number={1},
  pages={93--112},
  year={1982}
}

@article{LeitchHomotopy,
  title={The homotopy commutative cube},
  author={Leitch, R.D.},
  journal={Journal of the London Mathematical Society},
  volume={2},
  number={1},
  pages={23--29},
  year={1974},
  publisher={Oxford University Press}
}

@article{WaasPresenting,
  title={Presenting the topological stratified homotopy hypothesis},
  author={Waas, Lukas},
  journal={preprint},
  year={2024},
  note={Available at: \url{https://arxiv.org/abs/2403.07686}}
}

@misc{kerodon,
  author       = {Jacob Lurie},
  title        = {Kerodon},
  howpublished = {\url{https://kerodon.net}}
}

@article{QuinnHomotopically,
  title={Homotopically stratified sets},
  author={Quinn, Frank},
  journal={Journal of the American Mathematical Society},
  volume={1},
  number={2},
  pages={441--499},
  year={1988},
  publisher={JSTOR}
}

@article{PortaTeyssierExodromy,
  title={Topological exodromy with coefficients},
  author={Porta, Mauro and Teyssier, Jean-Baptiste},
  journal={preprint},
  year={2022},
  note={Available at: \url{https://arxiv.org/abs/2211.05004}}
}

@article{WaasWoolfYokuraStratifications,
  title={On stratifications and poset-stratified spaces},
  author={Waas, Lukas and Woolf, Jon and Yokura, Shoji},
  journal={preprint},
  year={2024},
  note={Available at: \url{https://arxiv.org/abs/2407.17690}}
}

@article{DuggerSheaves,
  title={Sheaves and homotopy theory},
  author={Dugger, Daniel},
  journal={preprint},
  note={Available at: \url{https://ncatlab.org/nlab/files/cech.pdf}},
  year={1999}
}

@article{QinModuli,
  title={On moduli spaces and {CW} structures arising from {M}orse theory on {H}ilbert manifolds},
  author={Qin, Lizhen},
  journal={Journal of Topology and Analysis},
  volume={2},
  number={04},
  pages={469--526},
  year={2010},
  publisher={World Scientific}
}

@article{LatourExistence,
  title={Existence de l-formes ferm{\'e}es non singulieres dans une classe de cohomologie de de {R}ham},
  author={Latour, Fran{\c{c}}ois},
  journal={Publications Math{\'e}matiques de l'IH{\'E}S},
  volume={80},
  pages={135--194},
  year={1994}
}

@phdthesis{DouteauThese,
  title={{\'E}tude homotopique des espaces stratifi{\'e}s},
  author={Douteau, Sylvain},
  year={2019},
  school={Universit{\'e} de Picardie Jules Verne},
  note = {Available at: \url{https://arxiv.org/abs/1908.01366}}
}

@book{RiehlCategory,
  title={Category theory in context},
  author={Riehl, Emily},
  year={2017},
  publisher={Courier Dover Publications}
}

@book{TopologicalStability,
  title={Topological stability of smooth mappings},
  author={Gibson, Christopher G. and Wirthm{\"u}ller, Klaus and Du Plessis, Andrew A. and Looijenga, Eduard J.N.},
  volume={552},
  year={2006},
  publisher={Springer}
}

@book{ZieglerLectures,
  title={Lectures on polytopes},
  author={Ziegler, G{\"u}nter M.},
  volume={152},
  year={2012},
  series={Graduate {T}exts in {M}athematics},
  publisher={Springer}
}

@article{TrotmanGeometric,
  title={Geometric versions of {W}hitney regularity for smooth stratifications},
  author={Trotman, David J.A.},
  journal={Annales scientifiques de l'{\'E}cole Normale Sup{\'e}rieure},
  volume={12},
  number={4},
  pages={453--463},
  year={1979}
}

@article{ThomEnsemblesEtMorphismes,
  title={Ensembles et morphismes stratifi{\'e}s},
  author={Ren{\'e} Thom},
  journal={Bulletin of the American Mathematical Society},
  year={1969},
  volume={75},
  pages={240-284},
  url={https://api.semanticscholar.org/CorpusID:123199335}
}

@article{MatherNotes,
  title={Notes on Topological Stability},
  author={John N. Mather},
  journal={Bulletin of the American Mathematical Society},
  year={2012},
  volume={49},
  pages={475-506},
  url={https://api.semanticscholar.org/CorpusID:22650297}
}

@article{CalleLiuClassifyingSpace,
  title={On the classifying space of a {M}orse flow category},
  author={Maxine E. Calle and Fangji Liu},
  year={2026},
  journal={preprint},
  note={Available at: \url{https://arxiv.org/abs/2603.23695}}
}

@book{Gersgorin,
  title={Ger{\v{s}}gorin and his circles},
  author={Varga, Richard S.},
  volume={36},
  year={2011},
  publisher={Springer {S}eries in {C}omputational {M}athematics}
}

@article{CordierPorter,
    title={Vogt's theorem on categories of homotopy coherent diagrams},
    author={Cordier, Jean-Marc and Porter, Timothy},
    journal={Math. Proc. Cambridge Philos. Soc.},
    volume={100},
    year={1986},
    number={1},
    pages={65--90}
}

@book{FoundationsofMeasurement,
    title={Foundations of {M}easurement, {V}olume 1: {A}dditive and {P}olynomials {R}epresentations},
    author={H.Krantz, David and Luce, R.Duncan and Suppes, Patrick and Tversky, Amos},
    publisher={Academic Press},
    year={1971}
}

@misc{ClassifyingSpaceVSGroupCompletion,
  title={Classifying space versus group completion for {M}orse flow categories},
  author={Fourel, Colin},
  journal={},
  year={},
  note={In preparation}
}

@book{QuillenHomotopicalAlebra,
  author    = {Daniel Quillen},
  title     = {Homotopical Algebra},
  series    = {Lecture Notes in Mathematics},
  volume    = {43},
  publisher = {Springer},
  year      = {1967}
}

@article {MilnorGeometricRealization,
    author = {Milnor, John},
     title = {The geometric realization of a semi-simplicial complex},
   journal = {Ann. of Math. (2)},
    volume = {65},
      year = {1957},
     pages = {357--362}
}

@article{JoyalEquivalence,
  title={Quasi-categories vs simplicial categories},
  author={Joyal, Andr{\'e}},
  journal={preprint},
  year={2007},
  note={Available at: \url{https://www.math.uchicago.edu/~may/IMA/Incoming/Joyal/QvsDJan9(2007).pdf}}
}

@inproceedings{HolsteinMorita,
  title={Morita cohomology},
  author={Holstein, Julian VS},
  booktitle={Mathematical Proceedings of the Cambridge Philosophical Society},
  volume={158},
  pages={1--26},
  year={2015},
  organization={Cambridge University Press}
}

@article{HiroFunctors,
  title={Functors (between infinity categories) that aren’t strictly unital},
  author={Tanaka, Hiro Lee},
  journal={Journal of Homotopy and Related Structures},
  volume={13},
  number={2},
  pages={273--286},
  year={2018},
  publisher={Springer}
}

@article{TrygveAlice,
  title={Structured flow categories and twisted presheaves},
  author={Hedenlund, {A}lice and {O}ldervoll, {T}rygve {P}oppe},
  journal={preprint},
  year={2026},
  note={Available at: \url{https://arxiv.org/abs/2603.29576}}
}

@article{OrsnesStratified,
  title={The stratified homotopy type of the reductive {B}orel--{S}erre compactification},
  author={{\O}rsnes Jansen, Mikala},
  journal={International Mathematics Research Notices},
  volume={2023},
  number={19},
  pages={16394--16452},
  year={2023},
  publisher={Oxford University Press}
}

@article{MillerStrongly,
  title={Strongly stratified homotopy theory},
  author={Miller, David},
  journal={Transactions of the American Mathematical Society},
  volume={365},
  number={9},
  pages={4933--4962},
  year={2013}
}

@inproceedings{CorneaHomotopicalII,
  title={Homotopical {D}ynamics {I}{I}: Hopf invariants, smoothings and the {M}orse complex},
  author={Cornea, Octavian},
  booktitle={Annales scientifiques de l'Ecole normale sup{\'e}rieure},
  volume={35},
  pages={549--573},
  year={2002}
}

@incollection{BarraudCorneaHomotopical,
  title={Homotopical dynamics in symplectic topology},
  author={Barraud, Jean-Fran{\c{c}}ois and Cornea, Octav},
  booktitle={Morse theoretic methods in nonlinear analysis and in symplectic topology},
  pages={109--148},
  year={2006},
  publisher={Springer}
}

@book{HoveyModelCategories,
  title={Model categories},
  author={Hovey, Mark},
  volume={63},
  year={1999},
  publisher={American Mathematical Society}
}

@article{BarraudCorneaSerre,
  title={Lagrangian intersections and the {S}erre spectral sequence},
  author={Barraud, Jean-Fran{\c{c}}ois and Cornea, Octav},
  journal={Annals of mathematics},
  pages={657--722},
  volume={166},
  year={2007}
}

\end{document}